\documentclass[preprint,12pt]{elsarticle}

\usepackage{newtxtext,newtxmath}
\usepackage{amsmath}
\usepackage{xcolor}
\usepackage{subcaption}
\usepackage{algpseudocode}

\journal{Journal of Computational Physics}

\newcommand{\RE}{\mathbb{R}}


\begin{document}

\begin{frontmatter}

\title{Data-driven reduced order modeling for parametric PDE eigenvalue problems using Gaussian process regression}

\author[label1]{Fleurianne Bertrand}
\author[label2,label3]{Daniele Boffi}
\author[label2,label4]{Abdul Halim}
\address[label1]{organization={University of Twente},
             %addressline={},
             %city={},
             %postcode={},
             %state={},
             {The Netherlands}}
\address[label2]{organization={King Abdullah University of Science and Technology (KAUST)},
             %addressline={},
             %city={},
             %postcode={},
             %state={},
             {Kingdom of Saudi Arabia}}
\address[label3]{organization={University of Pavia},
             %addressline={},
             %city={},
             %postcode={},
             %tate={},
             {Italy}}
\address[label4]{organization={H.S. College, Munger University},
             %addressline={},
             %city={},
             %postcode={},
             %tate={},
             {India}}

\begin{abstract}
In this article, we propose a data-driven reduced basis (RB) method for the approximation of parametric eigenvalue problems. The method is based on the offline and online paradigms. In the offline stage, we generate snapshots and construct the basis of the reduced space, using a POD approach. Gaussian process regressions (GPR) are used for approximating the eigenvalues and projection coefficients of the eigenvectors in the reduced space. All the GPR corresponding to the eigenvalues and projection coefficients are trained in the offline stage, using the data generated in the offline stage. The output corresponding to new parameters can be obtained in the online stage using the trained GPR. The proposed algorithm is used to solve affine and non-affine parameter-dependent eigenvalue problems. The numerical results demonstrate the robustness of the proposed non-intrusive method.
\end{abstract}

\begin{keyword}
Reduced basis method\sep Gaussian process regression\sep Eigenvalue problem\sep Proper orthogonal decomposition\sep non-intrusive method
\end{keyword}

\end{frontmatter}

%\maketitle
\section{Introduction}
Eigenvalue problems occur naturally in many areas of science and engineering. Just to name a few applications, let us mention fluid dynamics~\cite{BoffiCodina16,Akhan18}, structural mechanics~\cite{Bertrand19,Zheng16}, electromagnetism~\cite{Boffi19}, nuclear reactor physics \cite{Buchanetal13,reactorref1,GermanRagusa19}, photonics~\cite{Photonic1,photonic2,Gilbert19}, optics~\cite{optics01,optics04} for optimal control and optimal design,  population dynamics~\cite{Pdynamics}, in acoustic problem ~\cite{acoustic13} and solid-state physics~\cite{Pau07a} for calculating band structure. In nuclear physics, the eigenvalue problem represents the energetic diffusion of the neutron, where the first eigenvalue represents the criticality of the reactor and the corresponding eigenvalue represents neutron flux. For photonic crystal fibers, one has to solve two eigenvalue problems in order to calculate the photonic band gap~\cite{Gilbert19,photonic2}.

There is an active area of research in solving eigenvalue problems using the finite element method (FEM) and adaptive FEM. This includes dealing with eigenvalue problems in mixed form~\cite{Boffi10, Boffietal13, Boffietal97, Boffietal00}.
But there are very few works related to parametric eigenvalue problems~\cite{Buchanetal13, Fumagallietal16}. The parameters in the eigenvalue problems represent various quantities like material properties, physical domain, and initial/boundary conditions. In the many-query context, we need to solve the eigenvalue problem for different parameters and it becomes much more expensive. In the past few decades, reduced order modeling (ROM) has been used successfully to reduce the dimensionality of the problem~\cite{Peterson89,ItoRavindran98,KunischVolkwein01,KunischVolkwein02,Prudetal02a,Prudetal02b,WillcoxPeraire02,Veroyetal03,Benneretal03,Nguyenetal05,Huynhetal07,QuarteroniRozza07,Hesthavenetal16}. The idea is to reduce the high-fidelity (e.g. finite element) problem into a problem with much smaller dimension by projection or any other means. Reduced basis methods are efficient for solving parametric PDEs. The basis of the reduced space can be constructed in two ways, namely proper orthogonal decomposition (POD) and greedy approach. In the POD approach, the high-fidelity (HF) problem is solved for a certain number of parameters, and putting them into a matrix column-wise we form a snapshot matrix. The first few left-singular values of the snapshot matrix are chosen as the basis of reduced space. In the greedy approach first, the parameters are chosen based on some optimal criterion, and the corresponding solutions of the high-fidelity problem are used as a basis after orthogonalization. In general, the reduced solution is a Galerkin projection onto the reduced space that approximates the high-fidelity solution of the given PDE.

Parametric eigenvalue problems have been considered in this context, in particular, in~\cite{Machielsetal00,GermanRagusa19,Buchanetal13,ourrom22,Fumagallietal16,Horgeretal17,Pau07a,Pau07b}. All of these works solve the affine parameter dependent problem using projection based Galerkin method. The POD-based approach is used in~\cite{Machielsetal00,GermanRagusa19,Buchanetal13}.
In~\cite{Machielsetal00} a reduced basis method for calculating the first eigenvalues was proposed and an a posteriori estimate for the first eigenvalues was discussed.  In \cite{Buchanetal13} a fictitious time parameter was introduced and a POD approach was used for the first eigenpair of the parameter-dependent eigenvalue problem. The approach of introducing a fictitious time parameter was extended for finding later eigenpairs in~\cite{ourrom22} and the behavior of the approximation of the eigenvalues on different meshes was investigated. In~\cite{GermanRagusa19}, a system of parametric eigenvalue problems has been considered.

An a posteriori error estimate for the first eigenvalue and eigenvector was derived in~\cite{Fumagallietal16} and used for the selection of the parameters in the greedy approach. An a posteriori estimate for the repeated eigenvalues was established in~\cite{Horgeretal17} and the estimate was used for the selection of the parameters of the greedy algorithm and numerical results for an elasticity eigenvalue problem are presented. Finding multiple eigenvalues and eigenvectors simultaneously is the focus of~\cite{Pau07a,Pau07b} where a projection based Galerkin method is used and the greedy approach is adopted for the construction of the RB-basis. 

The affine dependence of the problem is the main key behind the success of the RB method, as it helps to decouple the problem into the online-offline stage. For non-affine and non-linear problems it is well known that this decoupling is not possible or quite difficult depending on the problem. The empirical interpolation method (EIM) and its discretized version have been developed to tackle these issues, but these methods are problem dependent and difficult to generalize. So, in this paper, we propose a data-driven model for solving the eigenvalue problems by adapting data-driven techniques used in~\cite{mguo18,mguo19} for parametric PDEs that will be applicable for both affine and non-affine parameter dependent problem.

The main tool will be the Gaussian process regression (GPR). In the offline stage, we solve the high-fidelity (HF) eigenvalue problem for some selected parameters. Then using a standard POD technique we construct the basis of the reduced space. Using the change-of-basis matrix, we get the reduced solutions at the sample parameters. A GP-based regression function is constructed using the selected parameters as input and the eigenvalues as output. For the approximation of the eigenfunctions, we construct one GP-based regression for each coefficient of the reduced eigenvector. Our proposed algorithm is capable to calculate multiple eigenvectors simultaneously as well.
In the online stage, for a new parameter, we get the approximation of the eigenvalue of interest, by evaluating the GPR corresponding to the eigenvalues; the projection coefficients are obtained by evaluating other GPRs. Multiplying the reduced vector obtained by the evaluation of the GPRs corresponding to the reduced coefficients, with the change-of-basis matrix, we get an approximation of the eigenvector of the HF problem.

The paper is organized as follows: in Section~\ref{sec:GPR} we discuss the Gaussian process regression and its training. In Section~\ref{sec:EVP_seeting} we discuss the POD approach, the generation of training data, and the proposed data-driven reduced order model. Numerical results for different eigenvalue problems are discussed in Section~\ref{sec:result}.
%%%%%%%%%%%%%%%%%%%%%%%%%%%%%%%%%%%%%%%%%%%%%%%%%%%%%
%%%%%%.  SECTION
%%%%%%%%%%%%%%%%%%%%%%%%%%%%%%%%%%%%%%%%%%%%%%%%%%%%
\section{Gaussian process regression}\label{sec:GPR}
A Gaussian process (GP) is a collection of random variables $x$ with the property that any finite number of them follows the joint Gaussian distribution. The Gaussian process is specified by its mean function $m(x)$ and covariance function $\kappa(x,x')$.
We write: $$f \sim GP(m,\kappa)$$ and say that the function $f$ is distributed as a GP with the mean function $m$ and covariance function $\kappa$. %For example, 
%let $f \sim GP(m,\kappa)$ with $m(x)=\frac{1}{4}x^2$ and $k(x,x')=exp(-\frac{1}{2}(x-x')^2)$.

Regression provides the approximation of an unknown output function from a set of input-output data. In Gaussian process regression (GPR), we approximate the unknown function by the mean of a GP. In GPR we assume a priori that the input-output pairs are generated from a GP, involving some parameters in the mean and variance function. Then using the training data, we find the parameters of the mean and the variance of the GP and get the posterior GP. The posterior GP is used to make predictions for new inputs. One advantage of using the GPR is that we are able to quantify the approximate function in terms of the variance, thus obtaining automatic error bounds. In our case, the approximate function will be a function from the parameter space $\mathcal{P}$ to the set of real numbers. The output of the functions will be either the eigenvalues or the projection coefficients of the eigenvectors onto the RB space. First, we will discuss the GPR of single valued parameter and assume that there are no errors in the output data. 
%The GPR for vector valued parameters are straight forward and can be found in ~\cite{mguo19}.

Let the output data be $y_n=U({\mu_n})$ for $n=1,2,\dots, n_s$. We want to use these data to approximate the function $U: \mathcal{P} \to \mathbb{R}$ by a GPR.
Let $X=[U({\mu_1}),\dots,U({\mu_{n_s}})]$ and $\pmb{y}=\{y_1,y_2,\dots,y_{n_s}\}$. We want to predict the value $y_{\ast}$ of $U$ at the point ${\mu}_{\ast}$.
Let $U$ be a Gaussian random field with mean $m:\mathscr{P} \to \mathbb{R} $ and covariance function $C: \mathscr{P}\times\mathscr{P} \to \mathbb{R}$, the random vector $(\pmb{y},y_{\ast})$ follows a Gaussian distribution with mean $(m(X),m({\mu}_{\ast}))$ and covariance matrix
\begin{center}
$\begin{bmatrix}
C(X,X)& C(X,{\mu}_{\ast})^\top\\
C(X,{\mu_*})& C({\mu}_{\ast},{\mu}_{\ast})
\end{bmatrix}.$
\end{center}
Then (see~\cite{Sullivan15}, for instance) the conditional distribution of $U({\mu}_{\ast})$ with given observation $U(X)=\pmb{y}$ follows a Gaussian distribution with mean 
$$m^{\star}({\mu}_{\ast})=m({\mu}_{\ast})+C({\mu}_{\ast},X)C(X,X)^{-1}(\pmb{y}-m(X))$$
and variance
$$K^{\star}({\mu}_{\ast})=C({\mu}_{\ast},{\mu}_{\ast})-C({\mu}_{\ast},X)C(X,X)^{-1}C(X,{\mu}_{\ast}). $$
For example, let $\mathcal{P}$ be the unit interval $[0,1]$ and $U$ follow a Gaussian process with mean $m(x)=ax+b$ and variance
$$C(x,x'):=\exp\left(-\frac{(x-x')^2}{2 l^2}\right),$$
where $l>0$ is the correlation length. Then our goal is to find a suitable hyperparameter $\theta=(a,b,l)$ by using the training data and by maximizing the log-likelihood function.
\subsection{Gaussian process regression with vector valued input and noisy output}
%%%%%% 
Let us now consider the case of a vector valued input parameter and noisy output data. Let $\mathcal{D}=\{ (\pmb{x}_i,y_i): i=1,2,\dots,n_s\}$ be the training set where $\pmb{x}_i\in \mathscr{P}\subset \mathbb{R}^d$ and $y_i$ are the corresponding outputs. Let the input-output pairs follow some regression function $U: \mathscr{P} \to \mathbb{R}$, whose prior is defined as a GP. 
Let the regression function be a GP corrupted by independent Gaussian noise~\cite{mguo19}, that is 
\[
y=U(\pmb{x})+\epsilon, \quad \epsilon \sim \mathcal{N}(0,\sigma^2), \quad U(\pmb{x})\sim GP(m(\pmb{x}),\kappa(\pmb{x},\pmb{x}')),
\]
where $m(\pmb{x})=\theta_1 H_1(\pmb{x})+\dots+\theta_k H_k(\pmb{x})$ for some basis $\{H_1(\pmb{x}),\dots,H_k(\pmb{x}) \}$ defined on $\mathscr{P}$. Suppose that for our numerical result we use a linear basis, that is $m(\pmb{x})=\theta_1+\theta_2 \pmb{x}$. One possible choice for the kernels is automatic relevance determination squared exponential
$$\kappa(\pmb{x},\pmb{x}')=\sigma_1^2 exp\Big( -\frac{1}{2} \sum\limits_{i=1}^d \frac{(x_i-x_i')^2}{l_i^2}\Big),$$
where $l_i > 0$ are the correlation length in the coordinate directions.
Given a finite number of points in the input domain, the outputs follow a joint Gaussian/Normal distribution, that is
\begin{equation}\label{ynormal}
   \pmb{y}|\pmb{X} \sim \mathcal{N} (m(\pmb{X}),\pmb{K}_y),\quad\pmb{K}_y=cov(\pmb{y}|\pmb{X})=\kappa(\pmb{X},\pmb{X})+\sigma ^2\mathbb{I}_{n_s},
\end{equation}
where $\pmb{y}=\{y_1,y_2,\dots,y_{n_s} \}^\top$, $\pmb{X}=[\pmb{x}_1|\pmb{x}_1|\dots|\pmb{x}_{N_s}]$, and $\mathbb{I}_{n_s}$ is the unit matrix of dimension $n_s$.
Let $\pmb{y}^{\ast}$ be a set of functional values corresponding to the test set inputs $\pmb{X}^{\ast}$, then 
the goal is to predict the output $\pmb{y}^{\ast}(\pmb{X}^{\ast})$ for a new test input set $\pmb{X}^{\ast}$. From the properties of Gaussian distribution the random vector $(\pmb{y}, \pmb{y}^{\ast})^\top$ follow the following joint normal distribution
\[
\begin{bmatrix}
\pmb{y}\\ \pmb{y}^{\ast}
\end{bmatrix}
\sim \mathcal{N}\Bigg(\begin{bmatrix}
m(\pmb{X})\\ m(\pmb{X}^\ast)
\end{bmatrix},\begin{bmatrix}
\pmb{K}_y & \kappa(\pmb{X},\pmb{X}^\ast)\\ \kappa(\pmb{X}^\ast,\pmb{X})  &\kappa(\pmb{X}^\ast,\pmb{X}^\ast)\\
\end{bmatrix}\Bigg).
\]
%%%%%%%%%%%%%%%%%%%%%%%%
Since we know the values of $\pmb{y}$ at the training set, we are interested in the conditional distribution of $\pmb{y}^{\ast}$ given $\pmb{y}$ which is expressed as:
$$\pmb{y}^{\ast}|\pmb{y}\sim \mathcal{N}\bigg(m(\pmb{X}^\ast)+\kappa(\pmb{X}^\ast, \pmb{X}) \pmb{K}_y^{-1} (\pmb{y}-m(\pmb{X})),\kappa(\pmb{X}^\ast,\pmb{X}^\ast)-\kappa(\pmb{X}^\ast,\pmb{X}) K_y^{-1} \kappa(\pmb{X},\pmb{X}^\ast)\bigg)$$
This is the posterior distribution for a specific set of test cases.
Thus the posterior function follow a new GP~\cite{Gprpaper}: 
$$U^{\ast}(\pmb{\mu})\sim GP(m^{\ast}(\pmb{\mu})), \kappa^{\ast}(\pmb{\mu},\pmb{\mu}')),$$
where
\begin{gather*}
m^{\ast}(\pmb{\mu})= m(\pmb{\mu})+\kappa(\pmb{\mu}, \pmb{X}) \pmb{K}_y^{-1} (\pmb{y}-m(\pmb{X}))\\
\kappa^{\ast}(\pmb{\mu},\pmb{\mu}')=\kappa(\pmb{\mu},\pmb{\mu}')-\kappa(\pmb{\mu},\pmb{X}) K_y^{-1} \kappa(\pmb{X},\pmb{\mu}').
\end{gather*}
Note that the posterior variance is smaller than the prior variance which is achieved using the data.

The hyperparameter $\pmb{\theta}=(\theta_1,\dots,\theta_k, l_1,\dots,l_d,\sigma, \sigma_1)$ plays an important role for the prediction. The random vector $\pmb{y}|\pmb{X}$ follows the mutivariate Gaussian distribution, so from~\eqref{ynormal} the likelihood function is given by:
\[
p( \pmb{y}|\pmb{X} ,\pmb{\theta} )=\frac{1}{\sqrt{(2\pi)^{N_s}|\pmb{K}_y(\pmb{\theta})|}}\exp\Big(-\frac{1}{2} (\pmb{y}-m(\pmb{X}))^\top\pmb{K}_y^{-1}(\pmb{\theta})(\pmb{y}-m(\pmb{X}))\Big),
\]
where $|\pmb{K}_y(\pmb{\theta})|$ denotes the determinant of the matrix $\pmb{K}_y$.
The optimal hyperparameter is obtained my maximizing the log likelihood function, that is
\begin{align*}
\pmb{\theta}_{opt}&=\arg\max\limits_{\pmb{\theta}}\log p(\pmb{y}|\pmb{X},\pmb{\theta} ) \\
 &=\arg\max\limits_{\pmb{\theta}} \Big\{-\frac{1}{2} (\pmb{y}-m(\pmb{X}))^\top\pmb{K}_y^{-1} (\pmb{y}-m(\pmb{X}))-\frac{1}{2}\log|\pmb{K}_y| -\frac{N_s}{2}\log(2\pi)\Big\}.
\end{align*}
%[see. \cite{Bishop07}]
For more details on GPR we refer the readers to~\cite{Gprbk,Gprpaper,Sullivan15,mguo19}.
%%%%%%%%%%%%%%%%%%%%%%%%%%%%%%%%%%%%%%%%%%%%%%%%%%%%%%%%%%%%%%%%%%%%
\section{Problem settings and data-driven reduced order model}\label{sec:EVP_seeting}
Let us consider a parametric PDE eigenvalue problem in variational form: given $\mu\in\mathcal{P}$ find $\lambda(\mu)\in \RE$ and $u(\mu)\in V$ with $u(\mu)\ne0$ such that
\[
a(u(\mu),v; \mu)=\lambda(\mu) b(u(\mu),v;\mu) \quad \forall v \in V,
\]
where $a,b: V \times V \to \RE$ are symmetric bilinear forms and $V$ is a suitable Hilbert space.
Let $V_h \subset V$ be a finite element space so that $\dim(V_h)=N_h$. The discrete problem reads: for any given $\mu \in \mathcal{P}$ find $\lambda_h(\mu)\in\RE$ and $u_h(\mu)\in V_h$ with $u_h(\mu)\ne0$ such that
\begin{equation}\label{Vardis}
a(u_h(\mu),v_h; \mu)=\lambda_h(\mu) b(u_h(\mu),v_h; \mu) \quad \forall v_h \in V_h.
\end{equation}
Let $\{\phi^j\}_{j=1}^{N_h}$ denote the basis of $V_h$ and $u_h(\mu)=\sum\limits_{j=1}^{N_h} u^j_h(\mu) \phi^j$, then the problem is equivalent to the solution of the generalized eigenvalue problem
\begin{equation}\label{hfls}
\mathbb{A}_h(\mu)\pmb{u}_h(\mu) =\lambda_h \mathbb{B}_h(\mu)\pmb{u}_h(\mu),
\end{equation}
where $\pmb{u}_h(\mu)=(u^1_h(\mu),\dots,u^{N_h}_h(\mu))^\top\in\mathbb{R}^{N_h}$ is the vector containing the nodal values and
\[(\mathbb{A}_h(\mu))_{i,j}=a( \phi^j,\phi^i;\mu),\quad(\mathbb{B}_h(\mu))_{i,j}=b( \phi^j,\phi^i;\mu),\quad 1\leq i,j\leq N_h.
\]
\subsection{Projection based reduced order modeling}
Let $V_N\subset V_h$ be the reduced space of dimension $N$ and $V_N=span\{\zeta_1, \dots, \zeta_N \} $. There are two popular methods to obtain a reduced basis, namely POD and greedy approach. In this paper, we use the POD approach which is discussed in details later in Subsection~\ref{sec:POD}. The reduced order problem reads~\cite{Quarteronietal16}: for any given $\mu \in \mathcal{P}$ find $\lambda_N(\mu)\in\RE$ and $u_N(\mu)\in V_N$ with $u_N(\mu)\ne0$ such that
\[
a(u_{N}(\mu),v_{N}; \mu)=\lambda_N(\mu) b(u_{N}(\mu),v_{N}; \mu) \quad \forall v_{N} \in V_N.
\]
Let $u_N(\mu)=\sum\limits_{j=1}^{N} u^j_N(\mu) \zeta_j$. We are then led to the following generalized eigenvalue problem for the reduced model
\begin{equation}\label{rbls}
\mathbb{A}_N({\mu})\pmb{u}_N({\mu}) =\lambda_N(\mu)\mathbb{B}_N({\mu})\pmb{u}_N({\mu})
\end{equation}
where $\pmb{u}_N(\mu)=(u^1_N(\mu),\dots,u^{N}_N(\mu))^\top\in\mathbb{R}^{N}$ is the vector containing projection coefficients of the eigenvector $u_h(\mu)$ and
\[
(\mathbb{A}_N({\mu}))_{n,m}=a(\zeta_m,\zeta_n; {\mu}),\quad(\mathbb{B}_h(\mu))_{m,n}=b(\zeta_m,\zeta_n;\mu)\quad 1\leq m,n\leq N.
\]
If we assume the affine parametric dependence of the bilinear forms then it is standard to write $\mathbb{A}_N({\mu})$ in terms of $\mathbb{A}_h({\mu})$ and $\mathbb{B}_N({\mu})$ in terms of $\mathbb{B}_h({\mu})$. Let consider the case when the bilinear forms are affinely dependent, that is they can be written as follows:
\[
\aligned
a(w,v;{\mu})&=\sum\limits_{q=1}^{Q_a}\theta_a^q({\mu})a_q(w,v)\quad \forall v,w\in V,\ {\mu}\in\mathcal{P}\\
b(w,v;{\mu})&=\sum\limits_{q=1}^{Q_b}\theta_b^q({\mu})b_q(w,v)\quad \forall v,w\in V,\ {\mu}\in\mathcal{P},
\endaligned
\]
where the bilinear forms $a_q(\cdot,\cdot)$ ($q=1,2,\dots,Q_a$) and $b_{q'}(\cdot)$ ($q'=1,2,\dots,Q_b$) are parameter independent. That is, we are able to write the bilinear forms as a sum of products of few parameter dependent functions $\theta_q$ and parameter independent bilinear forms $a_q$.
Then the matrix $\mathbb{A}_h$ and the vector $\mathbb{B}_h$ of the HF problem can be written as
$$\mathbb{A}_h({\mu})= \sum\limits_{q=1}^{Q_a} \theta_a^q({\mu})\mathbb{A}_h^q, \quad 
\mathbb{B}_h({\mu})= \sum\limits_{q=1}^{Q_b} \theta_b^q({\mu})\mathbb{B}_h^q,$$
and the matrix $\mathbb{A}_N$ and $\mathbb{B}_N$ of the RB problem can be written as
$$\mathbb{A}_N({\mu})= \sum\limits_{q=1}^{Q_a} \theta_a^q({\mu})\mathbb{A}_N^q, \quad
\mathbb{B}_N({\mu})= \sum\limits_{q=1}^{Q_b} \theta_b^q({\mu})\mathbb{B}_N^q.$$

Let $\mathbb{V}\in \mathbb{R}^{N_h\times N}$ be the transformation matrix of the basis of $V_h$ and $V_N$. More precisely, 
let us write the basis of $V_N$ in terms of the basis of $V_h$ as
\[
\zeta_m=\sum\limits_{i=1}^{N_h} \zeta_m^{(i)}\phi^i, \quad 1\leq m\leq N.
\]
 Define the transformation matrix $\mathbb{V}\in \mathbb{R}^{N_h\times N}$ as $\mathbb{V}=[\pmb{\zeta}_1| \dots| \pmb{\zeta}_N]$ where $\pmb{\zeta}_m=(\zeta_m^{(1)},\zeta_m^{(2)},\dots, \zeta_m^{(N_h)})^\top$.
Thus the matrix $\mathbb{A}_N$ and $\mathbb{B}_N$ of the RB problem can be written as~\cite{Quarteronietal16}
\begin{equation}\label{affromse}
\mathbb{A}_N({\mu})=\sum\limits_{q=1}^{Q_a} \theta_a^q({\mu}) \mathbb{V}^T \mathbb{A}_h^q \mathbb{V}, \quad \mathbb{B}_N({\mu})=\sum\limits_{q=1}^{Q_b} \theta_b^q({\mu}) \mathbb{V}^T \mathbb{B}_h^q \mathbb{V}.
\end{equation}

In the offline stage we calculate and store the matrices $\mathbb{V}^\top \mathbb{A}_h^q \mathbb{V}$ and $\mathbb{V}^\top \mathbb{B}_h^q \mathbb{V}$. In the online stage, for a new parameter ${\mu}$ we only evaluate the functions $\theta_a^q$ and $\theta_b^q$ and then with the help of the stored matrices we will get the linear system~\eqref{affromse} for the ROM.

It turns out that if the bilinear forms are not affine dependent, that is if we cannot separate the parameter and space dependent functions, then we cannot use the matrices of the HF problem the to get the system of RB method, as they will be parameter dependent. We propose a data-driven model which is discussed in Section~\ref{sec:DD} and which can be a viable solution to face with non affine problems.
\subsection{Proper orthogonal decomposition (POD)}\label{sec:POD}
Let $\mathbb{S}$ be the matrix of the $n_s$ snapshots, containing the eigenvectors at the sample points $\mu_1,\dots,\mu_{n_s}$. The sample points can be selected using various strategies. For instance we can use uniform tensorial sampling, Latin hypercube sampling, or Smolyak sparse sampling. The SVD gives
$$\mathbb{S}=\mathbb{U}\Sigma \mathbb{Z}^\top$$
with $$\mathbb{U}=[\pmb{\zeta}_1,\dots,\pmb{\zeta}_{N_h}] \in \mathbb{R}^{N_h\times N_h} \quad \text{and} \quad \mathbb{Z}=[\pmb{\psi}_1,\cdots,\pmb{\psi}_{n_s}] \in \mathbb{R}^{n_s\times n_s}$$
orthogonal matrices and $\Sigma=\mathrm{diag}(\sigma_1,\dots,\sigma_r,0,\dots,0)\in \mathbb{R}^{N_h \times n_s}$ with $\sigma_1\geq \sigma_2\geq \cdots \geq \sigma_r>0$, where $r$ is the rank of the matrix $\mathbb{S}$. The first $N$ left eigenvectors will be the basis, i.e. $\{\pmb{\zeta}_1,\dots,\pmb{\zeta}_{N}\}$ will be the basis of the low-dimensional subspace $V_N$. We will not use the SVD directly to find the left singular vectors unless $n_s \geq N_h$. Next, we describe the process of calculating the POD basis.
We can write
\begin{align}
\mathbb{S} \pmb{\psi}_i=\sigma_i\pmb{\zeta}_i \quad \text{and}\quad \mathbb{S}^\top \pmb{\zeta}_i =\sigma_i\pmb{\psi}_i, \quad i=1,\dots,r
\end{align}
or, equivalently,
\begin{align}
\mathbb{S}^\top\mathbb{S} \pmb{\psi}_i=\sigma^2_i\pmb{\psi}_i \quad \text{and}\quad \mathbb{S}\mathbb{S}^\top \pmb{\zeta}_i =\sigma^2_i\pmb{\zeta}_i, \quad i=1,\dots,r.
\end{align}
The matrix $\mathbb{C}=\mathbb{S}^T\mathbb{S}$ is called correlation matrix. The POD basis $\mathbb{V}$ of dimension $N\leq n_s$ is defined as the set of the first $N$ left singular vectors $\{\pmb{\zeta}_1,\dots,\pmb{\zeta}_N\}$ of $\mathbb{S}$ or, equivalently, the set of vectors
\begin{align}
\pmb{\zeta}_j=\frac{1}{\sigma_j} \mathbb{S} \pmb{\psi}_j\qquad 1\leq j\leq N
\end{align}
obtained from the first $N$ eigenvectors $\{\pmb{\psi}_1,\dots,\pmb{\psi}_N\}$ of the correlation matrix $\mathbb{C}$. So, instead of using the SVD of the matrix $\mathbb{S}$ we can solve the eigenvalue problem $\mathbb{C} \pmb{\psi}=\sigma^2 \pmb{\psi}$ which is of size $n_s$; notice that $n_s$ is usually much smaller than $N_h$.

Following a common practice, we choose the POD dimension $N$ as the smallest positive integer such that
\begin{equation}\label{criterion}
I(N)=\frac{\sum\limits_{i=1}^N\sigma_i^2}{\sum\limits_{i=1}^r\sigma_i^2} \geq 1-\varepsilon,
\end{equation}
where $\varepsilon$ is a prescribed tolerance. For the numerical computation of this paper we have used $\varepsilon=10^{-8}$.
%%%%%%%%%%%%%%%%%%%%%%%%%%%%%%%
\subsection{Data-driven model}\label{sec:DD}
In this subsection, we describe our data-driven model for calculating the first eigenpair. The approach can be extended to other eigenpairs also, by considering the eigenvector of interest in the snapshot matrix. In the next subsection we will show how to deal with the problem of finding multiple eigenvectors simultaneously.

First, select the sample parameters $\mu_1, \mu_2, \dots, \mu_{n_s}$ using sampling techniques like tensorial sampling, Monte Carlo sampling, Latin hypercube sampling, Smolyak sparse sampling with Clenshaw-Curits points from the parameter space. We calculate the first eigenpairs $\big(\pmb{u}_h(\mu_1),\lambda_1(\mu_1)\big)$, $\big(\pmb{u}_h(\mu_2),\lambda_1(\mu_2)\big)$,\dots, $\big(\pmb{u}_h(\mu_{n_s}),\lambda_1(\mu_{n_s})\big)$ of the high fidelity eigenvalue problem~\eqref{Vardis} in the finite dimensional space $V_h$  of dimension $N_h$. These eigenvectors are the snapshots corresponding to the selected parameters $\mu_1$, $\mu_2$,\dots, $\mu_{n_s}$. Consider the snapshot matrix
\[
\mathbb{S}=[\pmb{u}_h(\mu_1)|\cdots |\pmb{u}_h(\mu_{n_s})]
\]
containing the HF solutions as columns. Applying the SVD or using POD approach to the snapshot matrix $\mathbb{S}$, we construct the transformation matrix $\mathbb{V}$ by taking first few dominating left singular vectors of $\mathbb{S}$. The columns of $\mathbb{V}$ are the POD basis. If $\pmb{u}_h$ is the vector containing the nodal values of the solution $u_h\in V_h$ and $\pmb{u}_N$ is the vector containing the projection coefficients of the corresponding reduced solution $u_N \in V_N$, then the HF solutions and the RB solutions are related as:
\[
\textbf{u}_N=\mathbb{V}^\top \textbf{u}_h,\quad\pmb{u}_h =\mathbb{V}\pmb{u}_N.
\]
Thus the columns of
\[
\mathbb{S}_N:=\mathbb{V}^\top \mathbb{S}=[\mathbb{V}^\top\pmb{u}_h(\mu_1)| \cdots |\mathbb{V}^\top\pmb{u}_h(\mu_{n_s}) ]\approx[\pmb{u}_N(\mu_1)| \cdots |\pmb{u}_N(\mu_{n_s}) ] \in \mathbb{R}^{N\times n_s}
\]
are solutions of the reduced problem at $\mu_1$, $\mu_2$,\dots, $\mu_{n_s}$. Hence the $i^{th}$ row of the matrix $\mathbb{S}_N$ contains the $i^{th}$ coefficients of the reduced solutions at $\mu_1$, $\mu_2$,\dots, $\mu_{n_s}$. So, $i^{th}$ row of $\mathbb{S}_N$ serves as the outputs of the $i^{th}$ coefficient at the input points $\mu_1$, $\mu_2$,\dots, $\mu_{n_s}$.
 
Let $P_{\lambda}: \mathcal{P} \to \mathbb{R}$ be the mapping defined by $\mu \mapsto \lambda(\mu)$. We want to approximate the map $P_{\lambda}$ by a nonlinear GP-based regression $\hat{P}_{\lambda}$ using the training data  $\{\big(\mu_i, \lambda_1(\mu_i)\big): i=1,2,\dots, n_s \}$.% For a given new parameter $\mu_{new}$ the approximate eigenvalue will be $\hat{P}_{\lambda}(\mu_{new})$.

Similarly, let $P_{N}: \mathcal{P} \to \mathbb{R}^N$ be the mapping defined by $\mu \mapsto \pmb{u}_N(\mu)$. We want to approximate the map $P_{N}$ by a nonlinear GP-based regression $\hat{P}_{N}$ using the training data  $\{\big(\mu_i, \mathbb{V}^\top\pmb{u}_h(\mu_i)\big): i=1,2,\dots, n_s \}$. The GP-based regression models are constructed during the offline stage. In the online stage the eigenvalue at a new parameter $\mu_{\ast}$ is obtained by evaluating the regression $\hat{P}_{\lambda}$ at that point and the corresponding eigenvector is obtained as
$u_N(\mu_{\ast})= \sum\limits_{j=1}^{N}\hat{P}^j_N(\mu_{\ast}) \zeta_j$. Using the notation of \eqref{rbls}, we assume that $\hat{P}_{N}=(\hat{P}^1_{N},\dots,\hat{P}^N_{N})$, where $\hat{P}^j_{N}$ is the regression corresponding to the mapping $P^j_{N}: \mathcal{P} \to \mathbb{R}$ defined by $\mu \mapsto u^j_N(\mu)$ and the corresponding training data:
\[
\{\big(\mu_i, {u}_N^j(\mu_i)\big): i=1,2,\dots, n_s,  j=1,2,\dots,N.
\]
  
% Using these $N+1$ set of input-output data we construct $N+1$ Gaussian process regression model, one for the eigenvalue and $N$ regression model for the coefficients of the reduced eigenvector. In the online stage when a new parameter $\mu$ is given, we get all the coefficients of the corresponding RB solution and the eigenvalue using the above $N+1$ regression models.
Let $\pmb{u}_N(\mu_{\ast})=(\hat{P}^1_N(\mu_{\ast}),\dots,\hat{P}^N_N(\mu_{\ast}))^T$ be the vector containing the projection coefficients obtained using the GPR, then pre-multiplying it by $\mathbb{V}$, we get the corresponding approximate eigenvector $\pmb{u}_h(\mu_{\ast})$ of the eigenvalue problem~\eqref{Vardis}. In a similar way, considering other eigenpairs, we can construct the corresponding GPR in the offline phase and, evaluating those GPR in online phase, we can get the other eigenpairs for a new parameter. For the numerical implementation of the GPR we have used the MATLAB command \texttt{fitrgp}.

\subsection{Data-driven model for finding multiple eigenvectors simultaneously}\label{sec:simul}
We now consider the case when we are interested in finding $n_e$ eigenvalues and eigenvectors for a given eigenvalue problem. Let $\mu_1,\dots,\mu_{n_s}$ be sample points, then we need to find all the $n_e$ eigenpairs at these sample points. For each eigenvalue we need to train one GPR using the training data $\{\big(\mu_i, \lambda_j(\mu_i)\big): i=1,2,\dots, n_s \}$, for $j=1,2,\dots,n_e$. We can find the eigenvectors simultaneously by considering the following snapshot matrix
$$\mathbb{S}_s=[\pmb{s}_1| \dots|\pmb{s}_{n_s}]\in \mathbb{R}^{n_eN_h\times n_s},$$
where 
\begin{align}
    \pmb{s}_j &= \begin{bmatrix}
           \pmb{u}_{1,h}(\mu_j) \\
           \pmb{u}_{2,h}(\mu_j) \\
           \vdots \\
            \pmb{u}_{n_e,h}(\mu_j)
         \end{bmatrix}
\end{align}
Then we proceed as earlier with the modified snapshot matrix; we will get a transformation matrix $\mathbb{V}_s\in \mathbb{R}^{n_eN_h\times N}$, where $N$ is  the number of POD basis corresponding to the new snapshot matrix, obtained using the criterion expressed in~\eqref{criterion}. Then we get the matrix $\mathbb{S}_{s,N}=\mathbb{V}_s^\top\mathbb{S}_s\in\mathbb{R}^{N\times n_s}$, each row will serve as the training data. Let $P^s_{N}: \mathcal{P} \to \mathbb{R}^N$ be the mapping from $\mu$ to the reduced coefficients. We want to approximate the map $P^s_{N}$ by a nonlinear GP-based regression $\hat{P}^s_{N}$ using the training data $\{\big(\mu_i, \mathbb{V}_s^\top\pmb{s}(\mu_i)\big): i=1,2,\dots, n_s \}$. The GP-based regression models are constructed during the offline stage. In the online stage all the GPR $\hat{P}^s_{N}$ are evaluated at the parameter $\mu_{\ast}$, and the collection of all the coefficients forms a vector ${\pmb{\tilde{u}}_N}=(\hat{P}^s_{1}(\mu_{\ast}),\dots,\hat{P}^s_{N}(\mu_{\ast}))^\top.$ Then the vector $\mathbb{V}_s {\pmb{\tilde{u}}_N}$ contains the first $n_e$ eigenvectors at $\mu_{\ast}$. We can extract all the eigenvectors from this vector.

The advantage of this method is that we need to train only a limited number $N$ of GPRs. Hence, in order to find $n_e$ eigenpairs, we need to train $n_e+N$ GPRs, while if we do it separately we would need $n_e(N+1)$ GPRs.

%%%%%%%%%%%%%%%%%%%%%%%. SECTION-2%%%%%%%%%% STARTS %% HERE
\section{Numerical result}\label{sec:result}
First, we test our data-driven model on some affine parametric dependent eigenvalue problems.
After we have checked the good performances of our method, we move to eigenvalue problem where the dependence on the parameters is nonlinear and nonaffine.
Finally, we test our model on the the case when the problem presents crossing eigenvalues. This is a very challenging situation, as observed in~\cite{eccomas,moataz} and the references therein.

In all this section we compare the computations of our DD reduced order method with the reference high fidelity FEM solution computed with standard P1 elements.

\subsection{Results for affine parameter dependent eigenvalue problems}
In this subsection, we considered a simple harmonic oscillator eigenvalue problem in one and two dimensions. The analytical solutions of these problems are known. Note that the bilinear forms for these problems can be written affinely and we could use RB methods. Here we will use oue new approach in order to verify that the data-driven model is working well for this case.
\subsubsection{Simple harmonic oscillator in 1D}
Let us consider the eigenvalue problem defined on $\Omega=(-L,L) \subset \mathbb{R}$
\begin{align}\label{model}
-\frac{1}{2} \frac{d^2 u}{dx^2}+\frac{1}{2}\mu^2 x^2 u=\lambda u
\end{align}
with boundary conditions $u(-L)=u(L)=0$.

The eigensolutions of this problem are
\begin{align*}
\lambda_n(\mu)&=\left(n+\frac{1}{2}\right)\mu\quad n=0,1,2,3,\dots\\
u_n(x;\mu)&=\frac{\pi^{-\frac{1}{4}}}{\sqrt{2^n n!}} e^{-\mu x^2/2} H_n(\sqrt{\mu}x),
\end{align*}
where $H_n$ is the Hermite polynomial of order $n$. First few Hermite polynomials are $H_0(x)=1$, $H_1(x)=2x$, $H_2(x)=2x^2-2$, $H_3(x)=8x^3-12x$, $x_4(x)=16x^4-48x^2+12$, $H_5(x)=32x^5-160x^3+120x$.

The variational formulation of Problem~\eqref{model} is
\[
a(u(\mu),v;\mu)=\lambda(\mu) b(u(\mu),v)\quad \forall v \in H_0^1(\Omega)
\]
with
$$a(u(\mu),v;\mu)=a_1(u,v)+\mu^2 a_2(u,v),$$
where $$a_1(u,v)=\frac{1}{2} \int_{\Omega}\frac{du}{dx} \frac{dv}{dx}\,dx, \quad a_2(u,v)=\frac{1}{2} \int_{\Omega}x^2 u v\,dx, \quad b(u,v)=\int_{\Omega} uv\,dx. $$

For the numerical test, we choose $L=10$ and the parameter space $\mathcal{P}= [1,9]$. We have calculated the results of full order model (FEM) with mesh size $h=0.05$, $0.02$, and $0.01$, respectively. We selected $41$ uniformly distributed sample parameters from the parameter domain $\mathcal{P}$. The first eigenvectors are calculated at the selected sample points and put into a matrix columnwise to form the snapshot matrix.  We choose $\mu=2.5$, $3.5$,\dots, $8.5$ as test points and calculated the first eigenvalues and eigenvectors at these test points using data-driven (DD)  reduced order model. The first eigenvalues of FEM and the DD model at the test points are reported in Table~\ref{table1}. Note that the sample points used for snapshot matrix and the test points used for prediction are different. From the table one can see that the results are matching with the FEM results.
\begin{table}
\footnotesize
 	 	\centering
 	\begin{tabular}{|c|c|c|c|c|c|c|c|c|c|c|c|} 
\hline
$h$& $\mu \to$&2.5&	3.5	&4.5&	5.5&	6.5&	7.5&	8.5\\	

\hline
0.05 &DD &1.2508	&1.7528	&2.2548&	2.7568&	3.2588&	3.7608&	4.2628	\\
&FEM&1.2511&	1.7522&	2.2537&	2.7555&	3.2577&	3.7603&	4.2632\\	
\hline
0.02&DD&1.2501& 1.7504& 2.2508 &2.7511 &3.2514& 3.7517& 4.2521\\
&FEM&1.2502&	1.7504&	2.2506&	2.7509&	3.2512&	3.7516&	4.2521\\
\hline
0.01&DD&1.2500&	1.7501&	2.2502&	2.7503&	3.2504&	3.7504&	4.2505\\	
& FEM&1.2500&	1.7501&	2.2501&	2.7502&	3.2503&	3.7504&	4.2505\\	
\hline
 	\end{tabular}
	\caption{First eigenvalues of~\eqref{model} corresponding to $\mu=2.5$, $3.5$, $4.5$, $5.5$, $6.5$, $7.5$, $8.5$, respectively, using Data-Driven (DD) and full order model (FEM) with different mesh size of the domain when we consider $41$ uniformly distributed sample points.}
 	\label{table1}
\end{table}
 
Next, we choose $21$ uniformly distributed points from the sample space $\mathcal{P}$ and calculated the first eigenpairs at these sample points. Using the eigenvalues at these points we trained a GPR model. Using the GPR model, we predict the eigenvalues at the test points $\mu=2.5$, $3.5$, $4.5$, $5.5$, $6.5$, $7.5$, $8.5$. The predicted eigenvalues are reported in Table~\ref{table2} and compared with the first eigenvalues obtained by FEM. Also in this case the eigenvalues are matching pretty well with the FEM eigenvalues. If we compare these results with the results of Table~\ref{table1} then one can see that for coarse mesh the DD results in Table~\ref{table1} are better than the results in Table~\ref{table2}. Thus if we take more points in the parameter domain then the results of the DD model is better even if the mesh is not fine enough.
%For the numerical test we took $(-10,10)$ and the parameter $\mu\in [1,10].$
\begin{table}
\footnotesize
 	 	\centering
 	\begin{tabular}{|c|c|c|c|c|c|c|c|c|c|c|c|} 
	\hline
$h$& $\mu \to$&2.5&	3.5	&4.5&	5.5&	6.5&	7.5&	8.5\\	
 \hline
0.05 &DD& 1.2511&	1.7529&	2.2547&	2.7565&	3.2584&	3.7602&	4.2620	\\
&FEM&1.2511&	1.7522&	2.2537&	2.7555&	3.2577&	3.7603&	4.2632\\	
\hline
0.02&DD&1.2502&	1.7505&	2.2508&	2.7510&	3.2513&	3.7516&	4.2519\\	
&FEM&1.2502&	1.7504&	2.2506&	2.7509&	3.2512&	3.7516&	4.2521\\
\hline
0.01&DD&1.2500&	1.7501&	2.2502&	2.7503&	3.2503&	3.7504&	4.2505	\\
& FEM&1.2500&	1.7501&	2.2501&	2.7502&	3.2503&	3.7504&	4.2505\\	
\hline
 	\end{tabular}
	\caption{ First eigenvalues of~\eqref{model} corresponding to $\mu=2.5$, $3.5$, $4.5$, $5.5$, $6.5$, $7.5$, $8.5$, respectively, using DD and FEM with different mesh size of the domain when we consider $21$ uniformly distributed sample points.}
 	\label{table2}
\end{table}

\begin{figure}
\centering
\includegraphics[width=6.5cm]{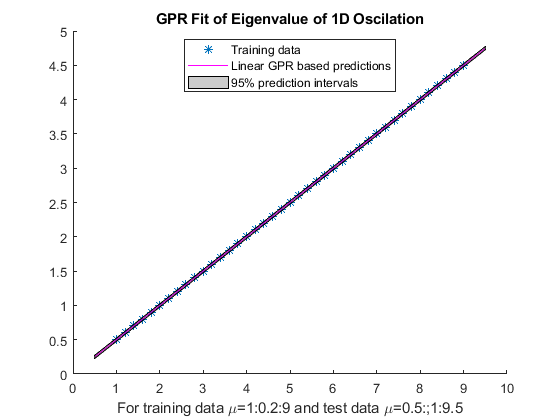}
\includegraphics[width=6.5cm]{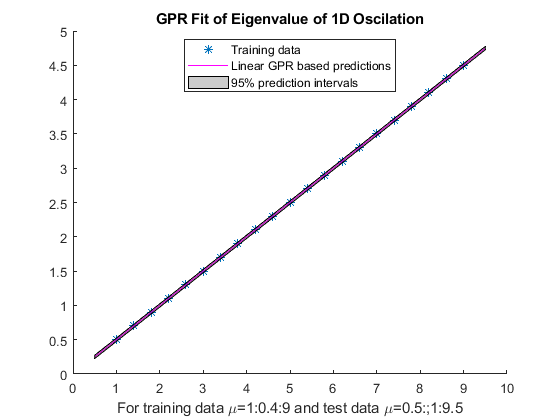}
\caption{First eigenvalue plot with $95\%$ confidence interval with mesh size $h=0.05$ for the problem \eqref{model} corresponding to 41 and 21 sample points respectively.}
\label{fig0}       % Give a unique label
\end{figure}

In Figure~\ref{fig0}, we plot the first eigenvalue for the DD model with $95$ percentage confidence interval for both the cases with $41$ and $21$ points, respectively. We can see that the confidence interval is very narrow that means that the predictions by the DD model is very accurate.

\begin{figure}
\centering
\includegraphics[width=4.5cm]{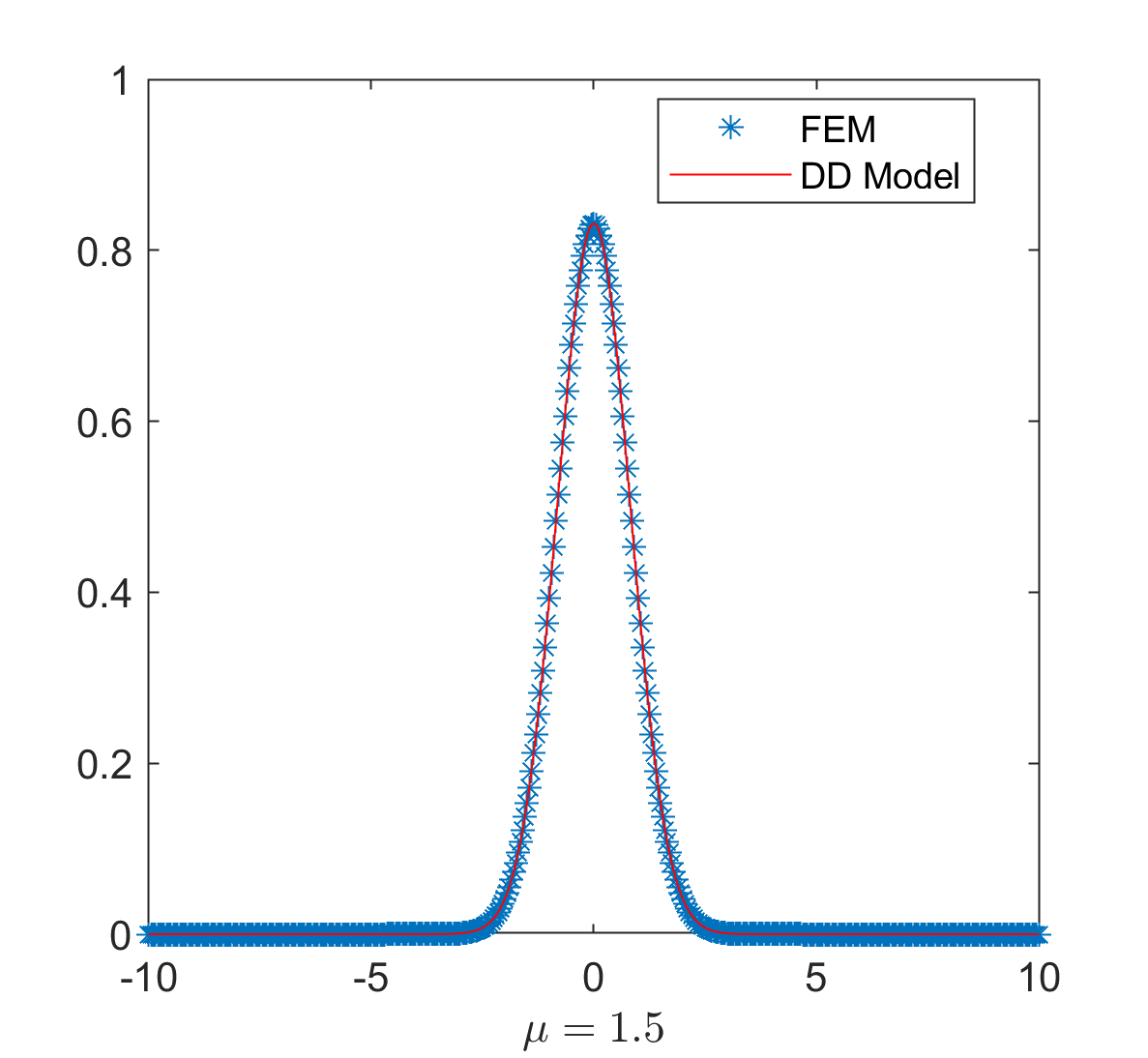}
\includegraphics[width=4.5cm]{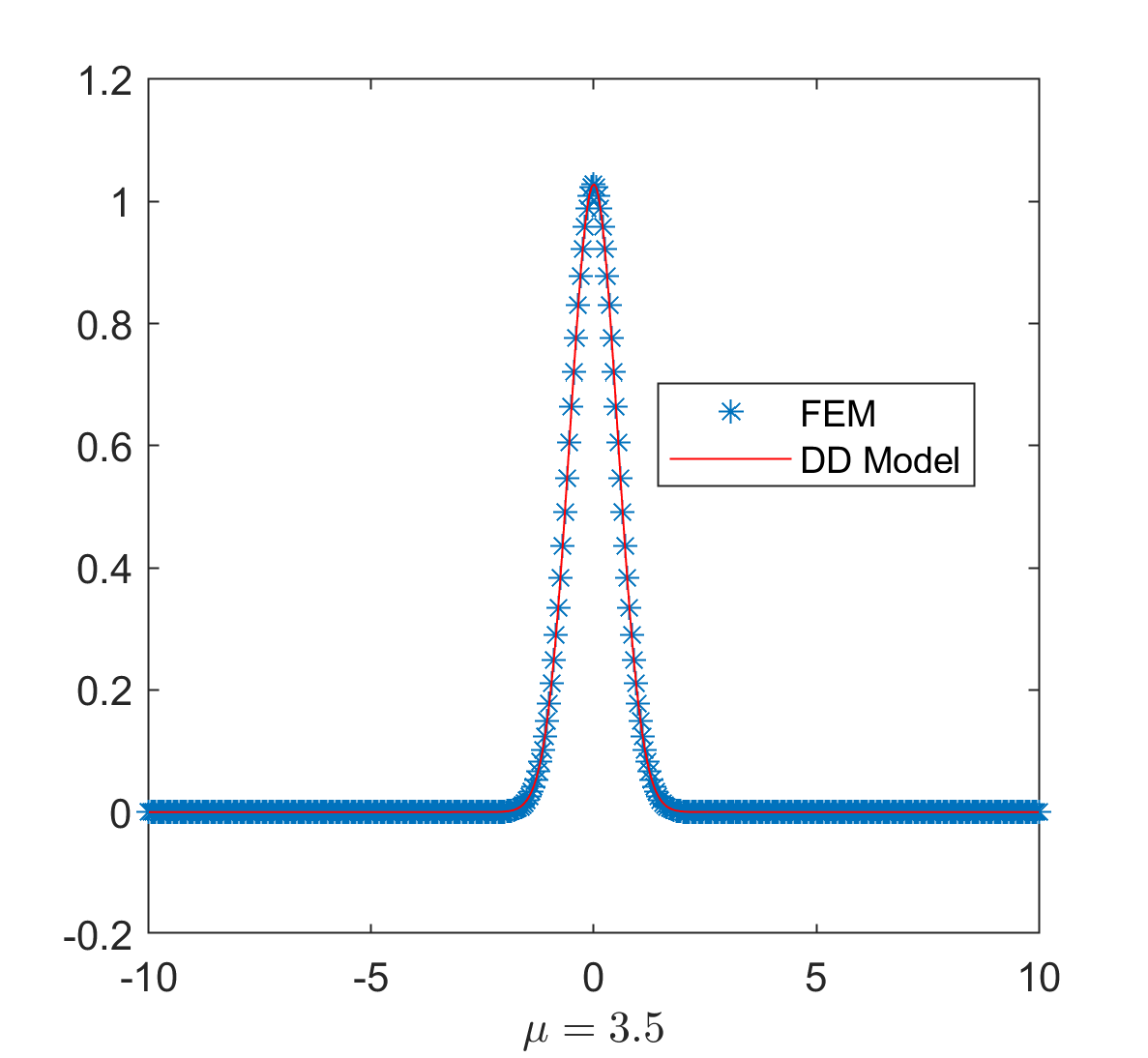}
\includegraphics[width=4.5cm]{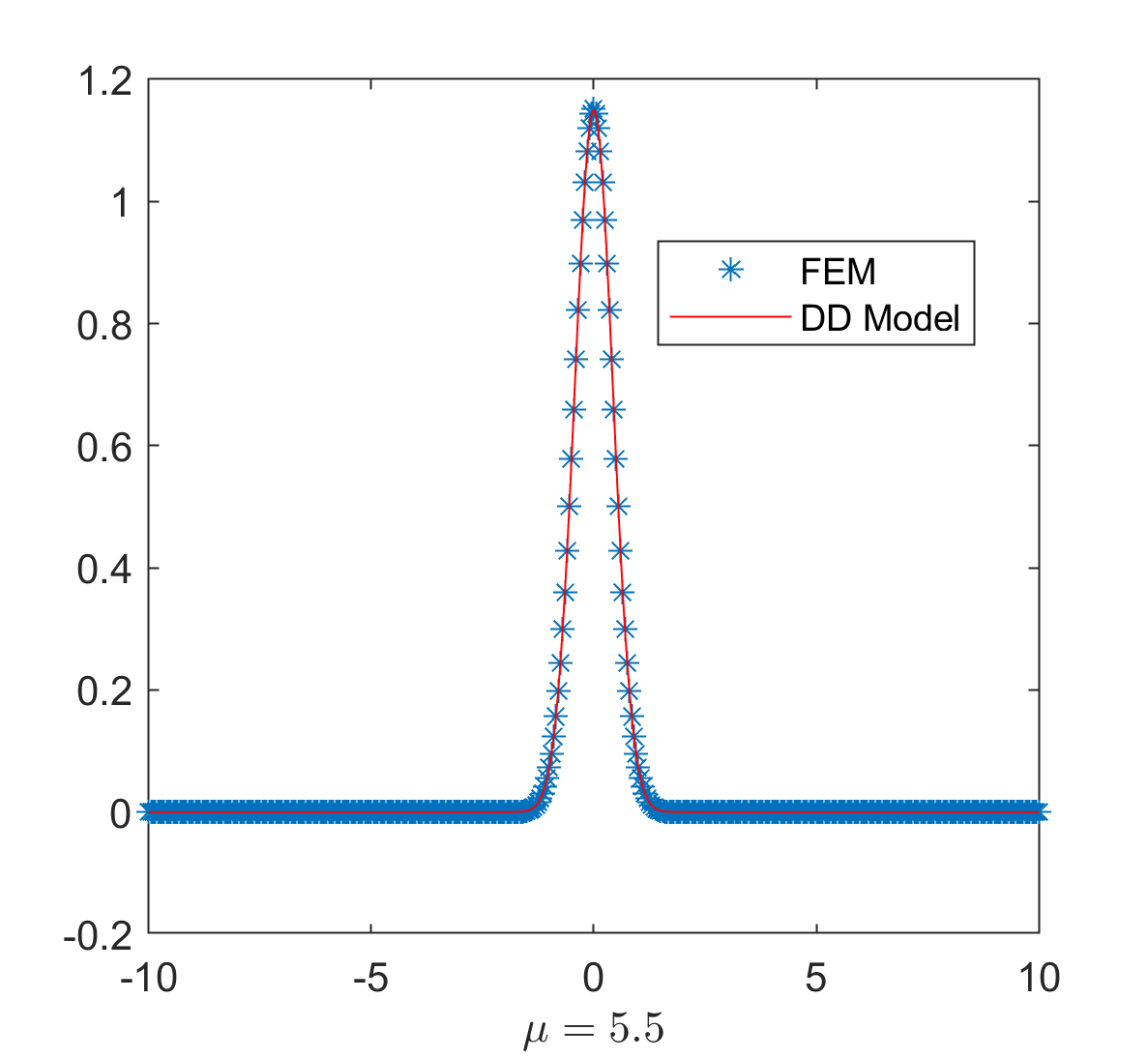}
\caption{First eigenvectors of the Problem~\eqref{model} at $\mu=1.5$, $3.5$, $5.5$ respectively for the case of $21$ sample points and mesh size $h=0.05$.}
     \label{fig:evct_evp1d}       % Give a unique label
\end{figure}

\begin{figure}
\centering
\includegraphics[width=4.5cm]{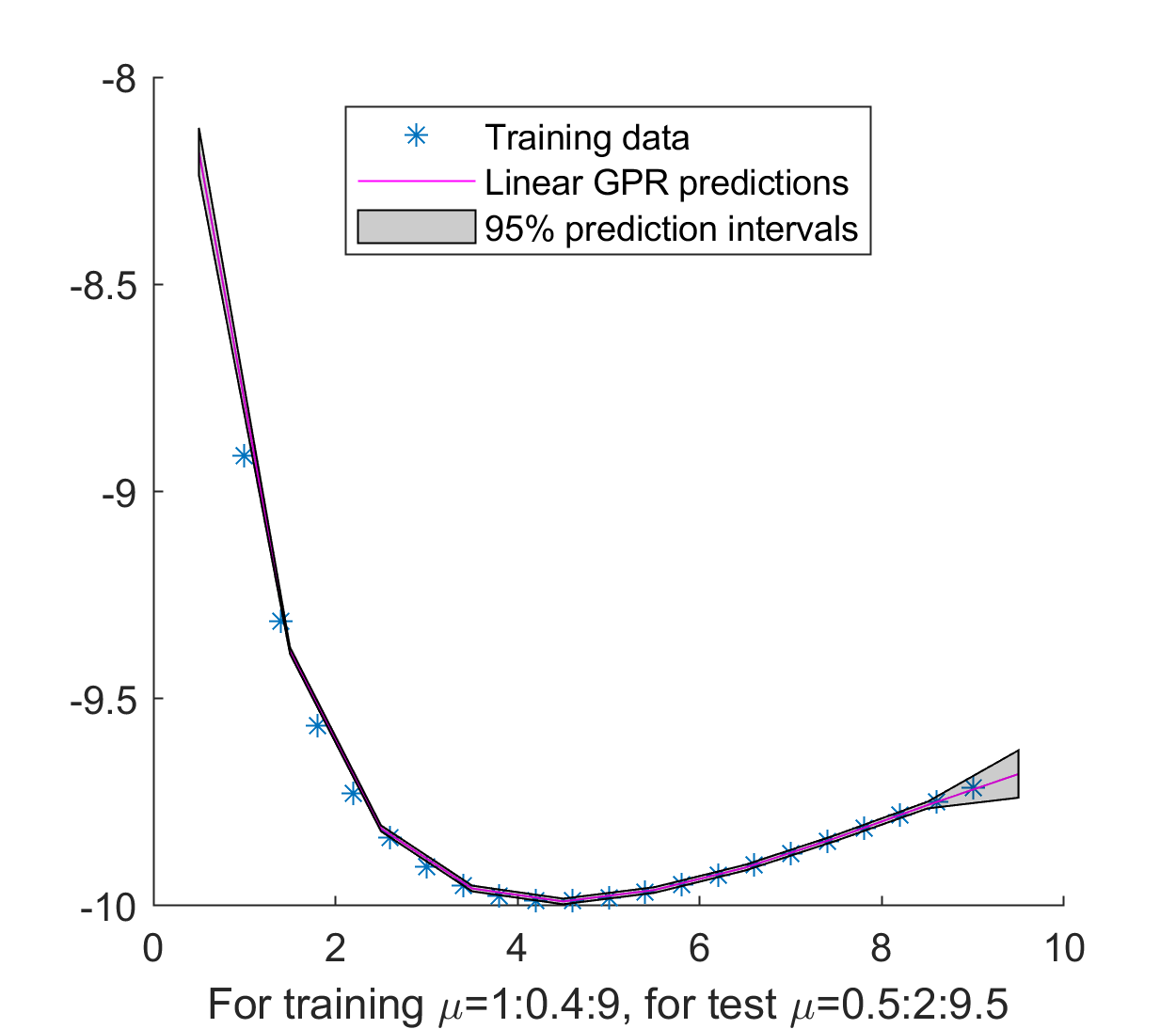}
\includegraphics[width=4.5cm]{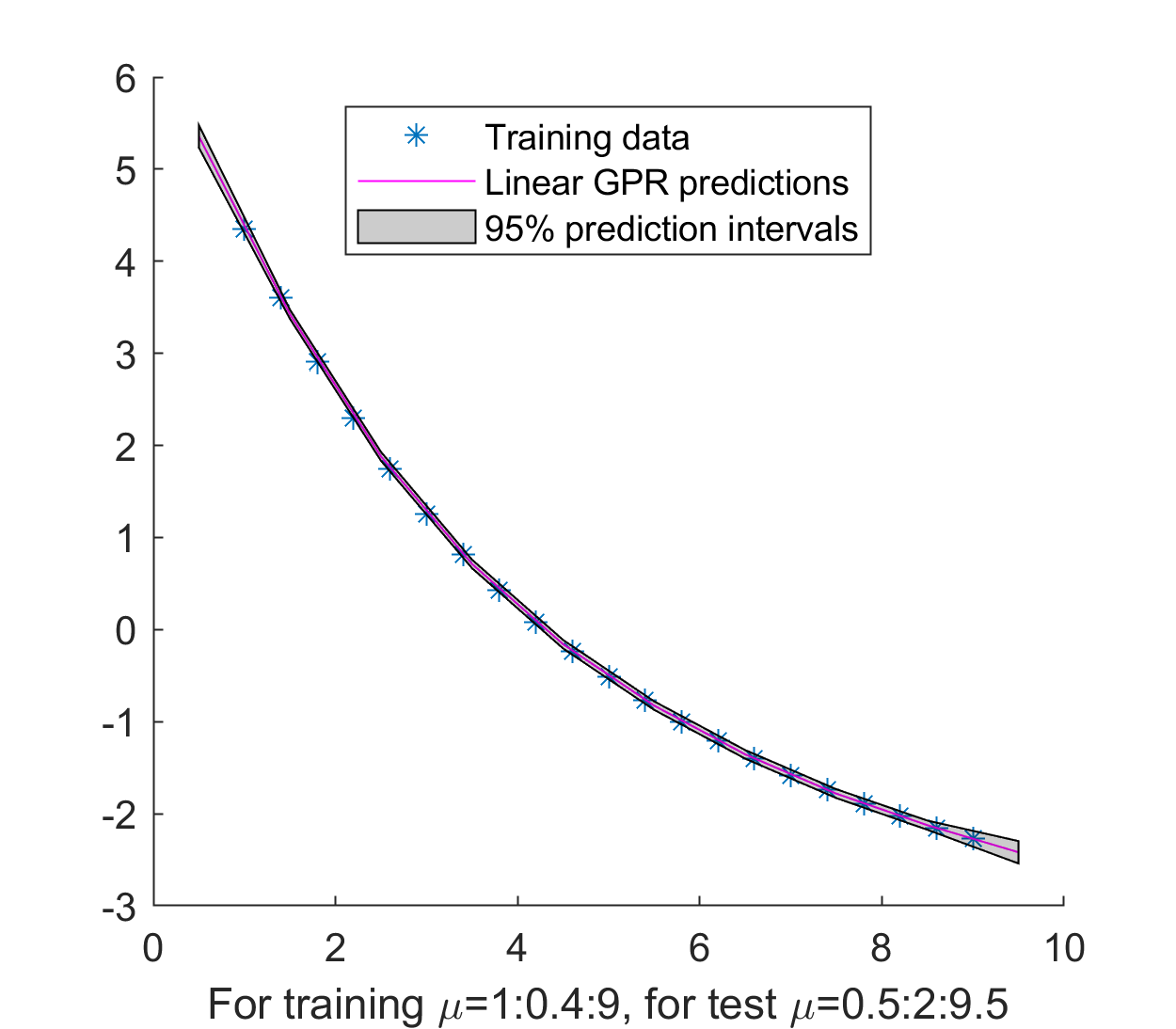}
\includegraphics[width=4.5cm]{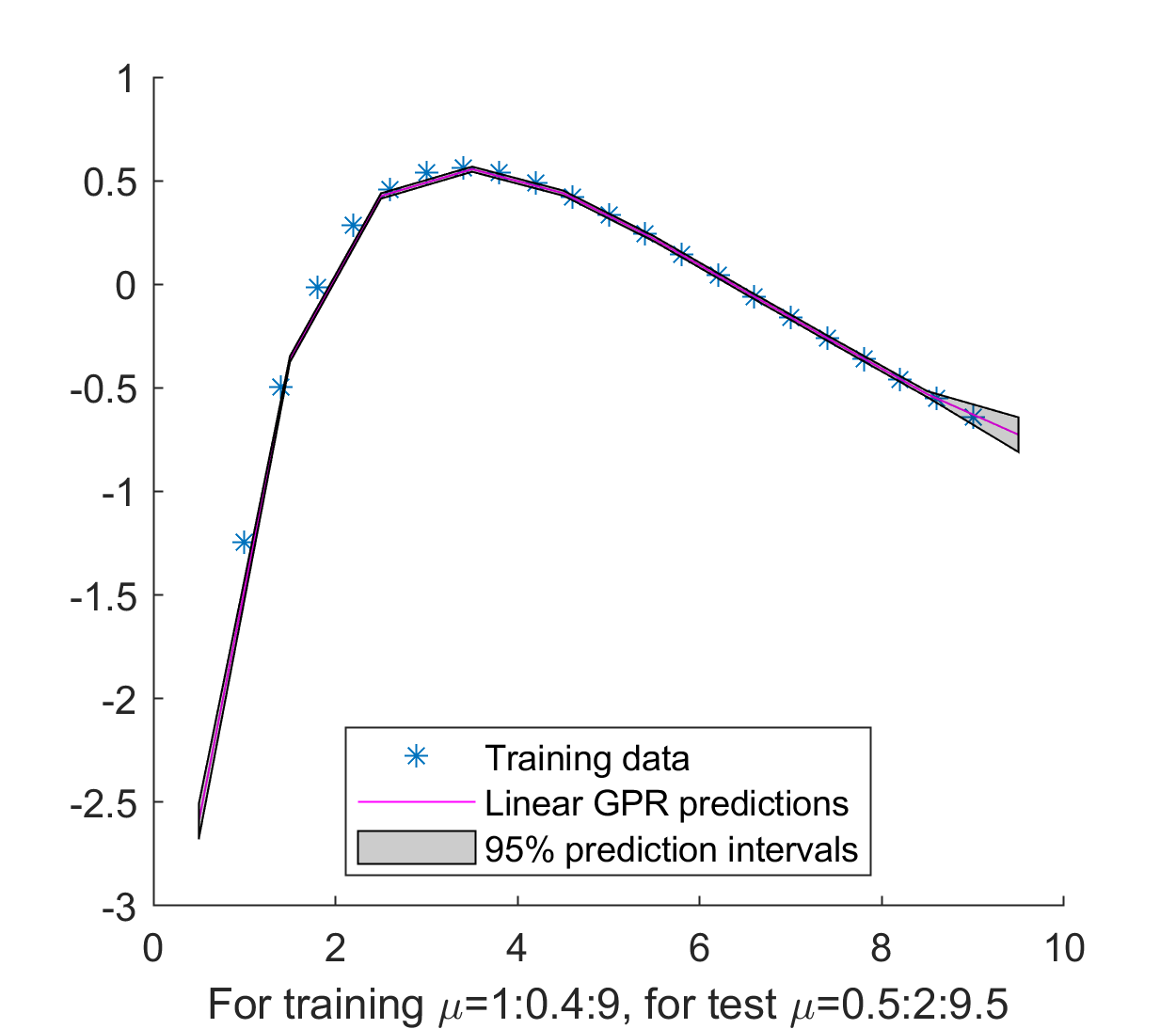}

\includegraphics[width=4.5cm]{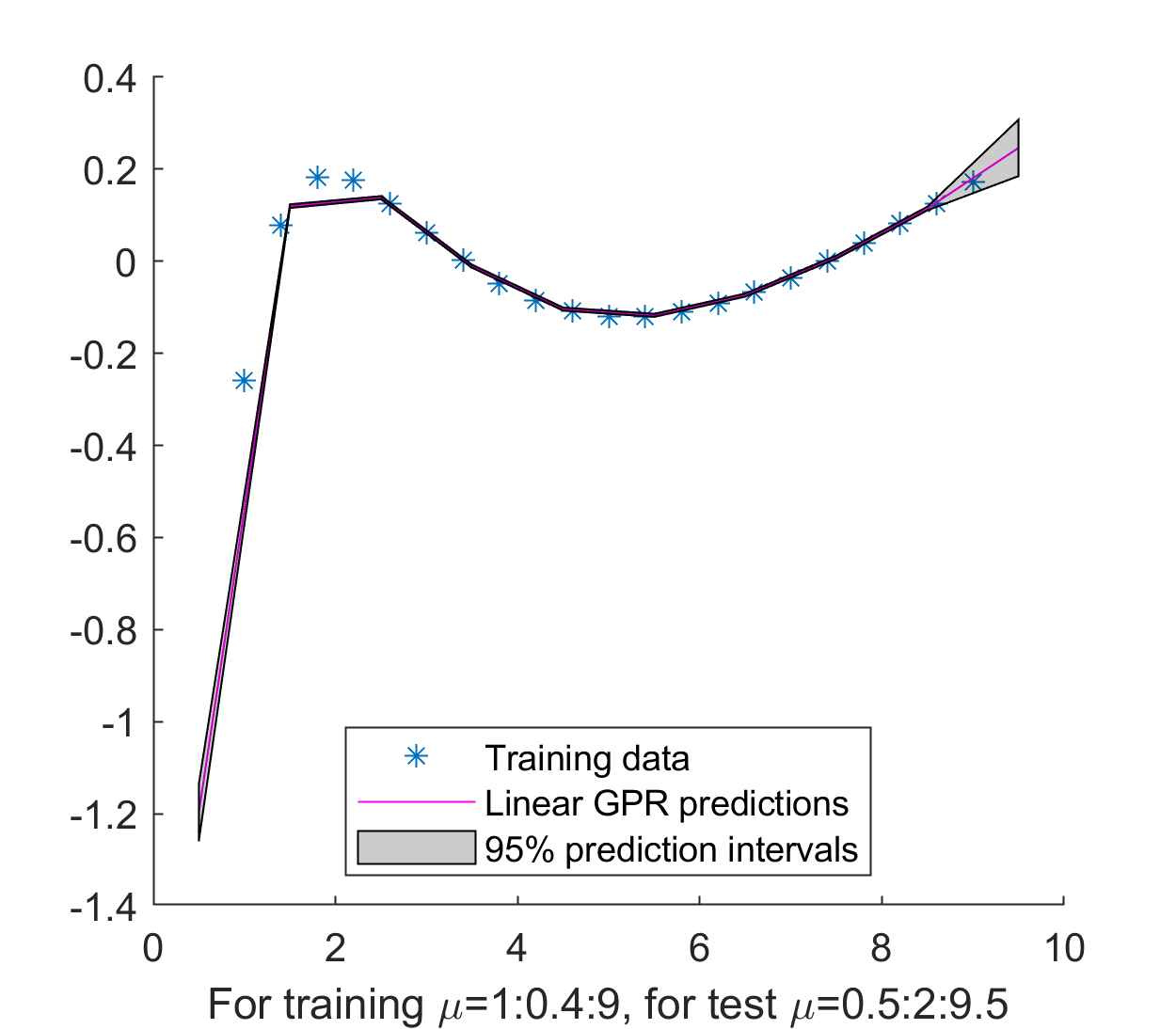}
\includegraphics[width=4.5cm]{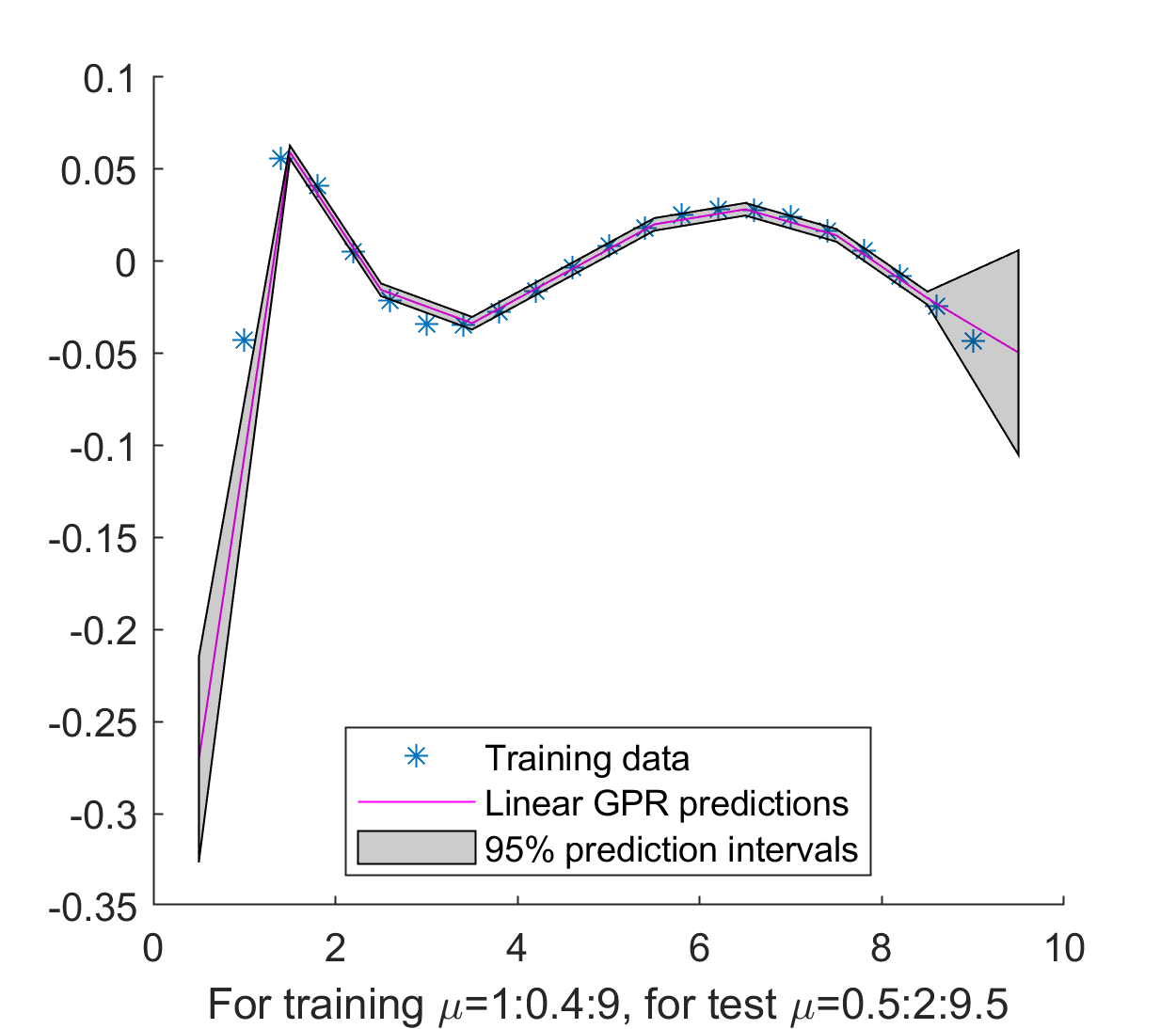}
\includegraphics[width=4.5cm]{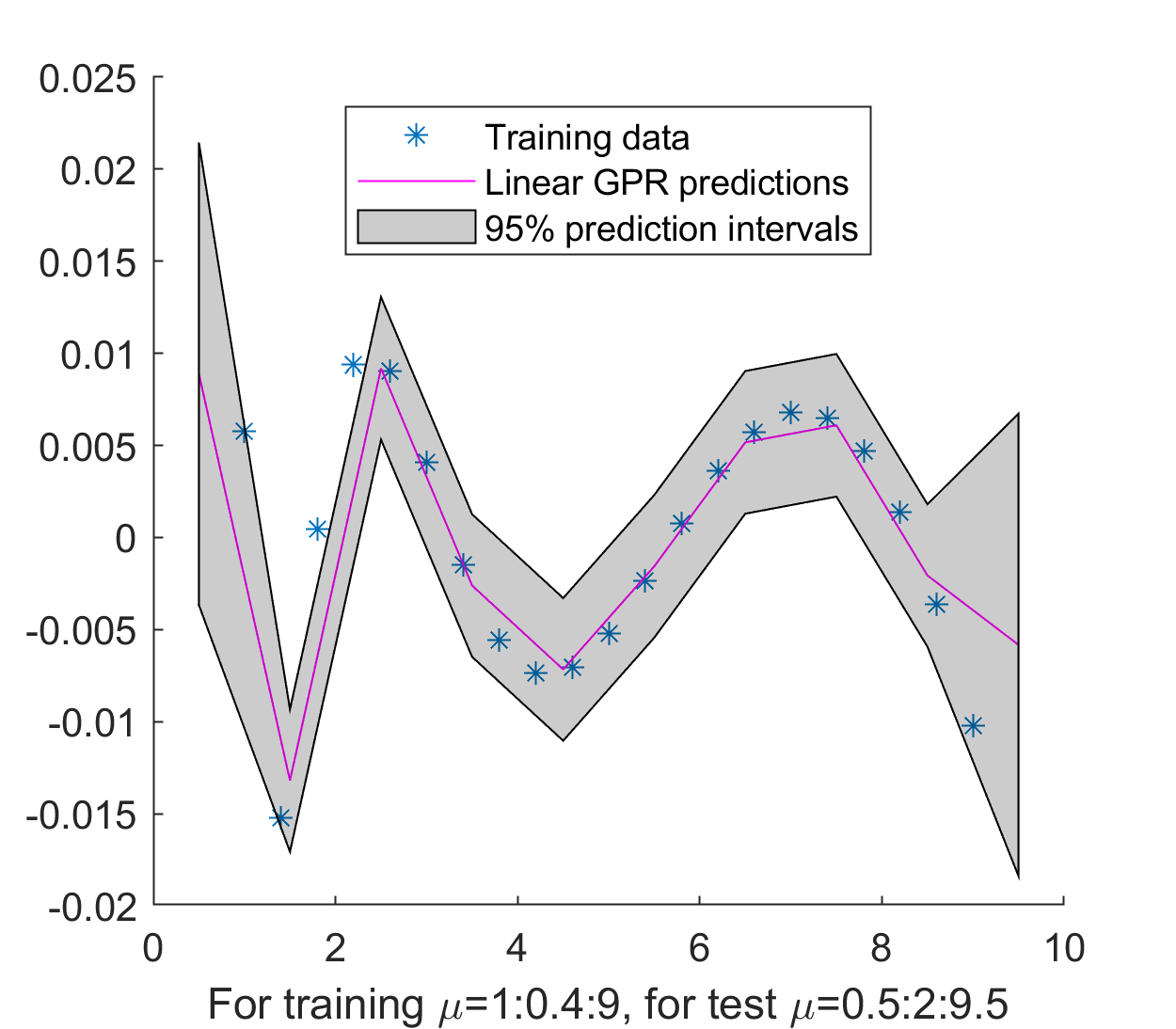}
\caption{First six RB coefficients with $95\%$ confidence interval for the case of 21 sample points and mesh size $h=0.05$.}
     \label{fig:coeff_evp1d}       % Give a unique label
\end{figure}

In Figure~\ref{fig:evct_evp1d}, we show the first eigenvector at test points $\mu=1.5$, $3.5$, $5.5$, respectively, obtained by our DD model with $21$ sample points and compare it with the corresponding FEM eigenvectors. It is evident that the DD eigenvectors match well with the FEM eigenvectors.
 
In Figure~\ref{fig:coeff_evp1d}, we show the GPR plot of the first six projection coefficients of the RB solutions. We show the test data (blue star) and we also plot the GPR mean function in magenta color. For each GPR, we plot the $95$ percentage confidence interval. Note that in the test set we take two points, namely $\mu=0.5$, $9.5$, which lie outside of the parameter interval so at the end points the confidence interval is large. Our DD model will work best for all parameter within the parameter space.
%%%%%%%%%%%%%%%%%%%%%%%%%%%%%%%%%%%%%%%%%%%%%%%%%%%%%%%%%%%%%%%%%%%%%%%%%%%%%%%%%%%%%%%%%%%%%%%
\subsubsection{Simple harmonic oscillator in 2D}
Let us consider the eigenvalue problem defined on $\Omega=\left(-\frac{\pi}{2},\frac{\pi}{2}\right)^2\subset\mathbb{R}^2$
\begin{equation}\label{model1}
    \left\{
\aligned
&-\frac{1}{2}\Delta u(x,y)+\frac{1}{2}\mu^2 (x^2+y^2) u(x,y)=\lambda u(x,y)&&\text{in }\Omega\\
&u=0&&\text{on }\partial\Omega
\endaligned
\right.
\end{equation}
The analytic eigenvalues and eigenvectors of this problem are given by
\begin{gather*}
\lambda_{m,n}(\mu)=(m+n+1)\mu\\
u_{m,n}(x,y;\mu)=\frac{1}{\sqrt{2^{m+n} m! n! \pi}} e^{-\mu (x^2+y^2)/2} H_m(\sqrt{\mu}x)H_n(\sqrt{\mu}y)
\end{gather*}
for $ m,n=0,1,2,3,\dots$, where $H_n$ is the Hermite polynomial of order $n$.

The variational form of the Problem~\eqref{model1} is 
\[
a(u(\mu),v;\mu)=\lambda(\mu) b(u(\mu),v) \quad \forall v \in H_0^1(\Omega)
\]
with
$$a(u,v;\mu)=a_1(u,v)+\mu^2 a_2(u,v)$$
where
\[
\aligned
&a_1(u,v)=\frac{1}{2} \int_{\Omega}\nabla u\cdot \nabla v\,dx\\
&a_2(u,v)=\frac{1}{2} \int_{\Omega}(x^2+y^2) uv\,dx\\
&b(u,v)=\int_{\Omega} uv\,dx.
\endaligned
\]

As in the one dimensional case we choose the parameter space as $\mathcal{P}=[1,9]$. We select two sets of uniformly distributed parameters from the parameter space with parameter spacing $\delta\mu=0.2$, $0.4$, respectively. For $\delta\mu=0.2$ we have $41$ uniformly distributed points and we get $21$ uniformly distributed points when $\delta\mu=0.4$.

\begin{table}
\footnotesize
 	 	\centering
 	\begin{tabular}{|c|c|c|c|c|c|c|c|c|c|c|c|c|} 
 		\hline
	h& {\begin{tabular}[c]{@{}c@{}} $\delta \mu$ \end{tabular}} &
	{\begin{tabular}[c]{@{}c@{}} Model \end{tabular}} &
 		  {\begin{tabular}[c]{@{}c@{}} $\mu=1.5$ \end{tabular}} &
		 {\begin{tabular}[c]{@{}c@{}} $2.5$ \end{tabular}} &
 		{\begin{tabular}[c]{@{}c@{}}  $3.5$ \end{tabular}}& 
 		{\begin{tabular}[c]{@{}c@{}} $4.5$ \end{tabular}} & 
 		{\begin{tabular}[c]{@{}c@{}} $5.5$ \end{tabular}}&
		{\begin{tabular}[c]{@{}c@{}} $6.5$ \end{tabular}} & 
		{\begin{tabular}[c]{@{}c@{}} $7.5$ \end{tabular}} \\
		%{\begin{tabular}[c]{@{}c@{}} $\mu=2.5$ \end{tabular}}  \\	
	\hline
0.1& 0.4&DD&1.6683& 2.5374& 3.5040& 4.5205& 5.5275& 6.5297& 7.5457\\	
		& & FEM&1.6490& 2.5317 &3.5136& 4.5167& 5.5241& 6.5334& 7.5441 \\
			\hline
0.1& 0.2&DD&1.6633& 2.5318& 3.5075& 4.5214& 5.5248& 6.5302& 7.5461 \\
		& & FEM&1.6490& 2.5317 &3.5136& 4.5167& 5.5241& 6.5334& 7.5441 \\
	\hline
0.05& 0.4&DD& 1.6668& 2.5334& 3.4970& 4.5084 &5.5088& 6.5044 &7.5125\\
& & FEM&1.6475& 2.5279& 3.5062& 4.5043& 5.5058& 6.5081& 7.5107\\	
		\hline
		0.05& 0.2& DD&1.6618& 2.5279& 3.5004& 4.5090& 5.5062& 6.5052& 7.5128 \\
		& & FEM&1.6475& 2.5279& 3.5062& 4.5043& 5.5058& 6.5081& 7.5107\\
\hline
 	\end{tabular}
	\caption{First eigenvalues corresponding to different values of $\mu$ for the 2D oscillation problem. Note that the exact first eigenvalue is equal to $\mu$.}
	\label{table:2dosc}
\end{table}

\begin{figure}
\centering
%\begin{subfigure}{0.45\textwidth}
%\includegraphics[height=7cm,width=8cm]{evp_2d_sho/1stev_evp2d_h1_dmu4.png}
% \caption{$h=0.1,\delta \mu=0.4$}
%\end{subfigure}
%\begin{subfigure}{0.45\textwidth}
%\includegraphics[height=7cm,width=8cm]{evp_2d_sho/1stev_evp2d_h1_dmu2.png}
%  \caption{$h=0.1,\delta \mu=0.2$}
%\end{subfigure}\\
%\begin{subfigure}{0.45\textwidth}
%\includegraphics[height=7cm,width=8cm]{evp_2d_sho/1stev_evp2d_h05_dmu4.png}
%  \caption{$h=0.05,\delta \mu=0.4$}
%\end{subfigure}
%\begin{subfigure}{0.45\textwidth}
%\includegraphics[height=7cm,width=8cm]{evp_2d_sho/1stev_evp2d_h05_dmu2.png}
%  \caption{$h=0.05,\delta \mu=0.2$}
%  \end{subfigure}
\subcaptionbox{$h=0.1$, $\delta \mu=0.4$}{
\includegraphics[width=6cm]{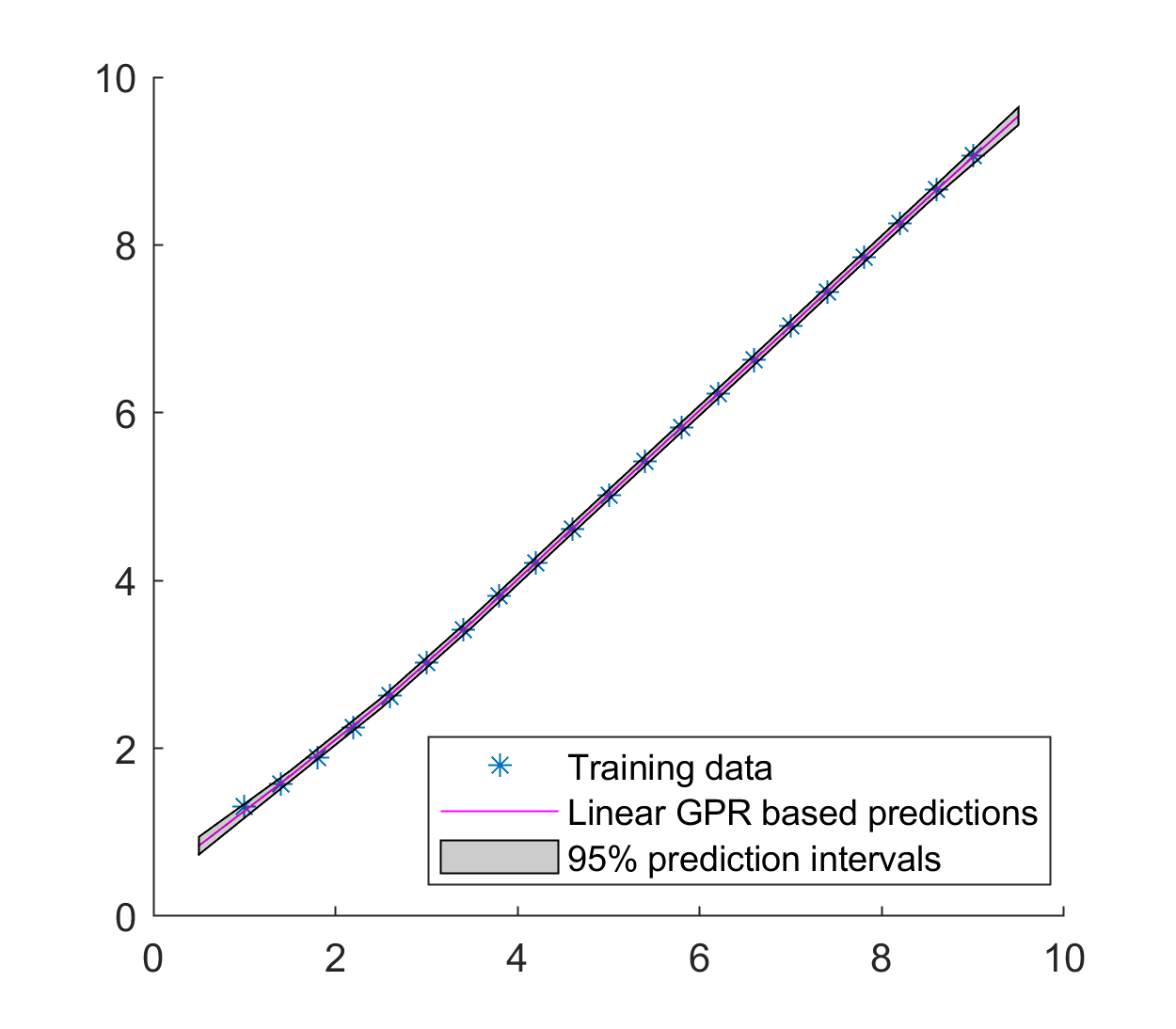}}
\subcaptionbox{$h=0.1$, $\delta \mu=0.2$}{
\includegraphics[width=6cm]{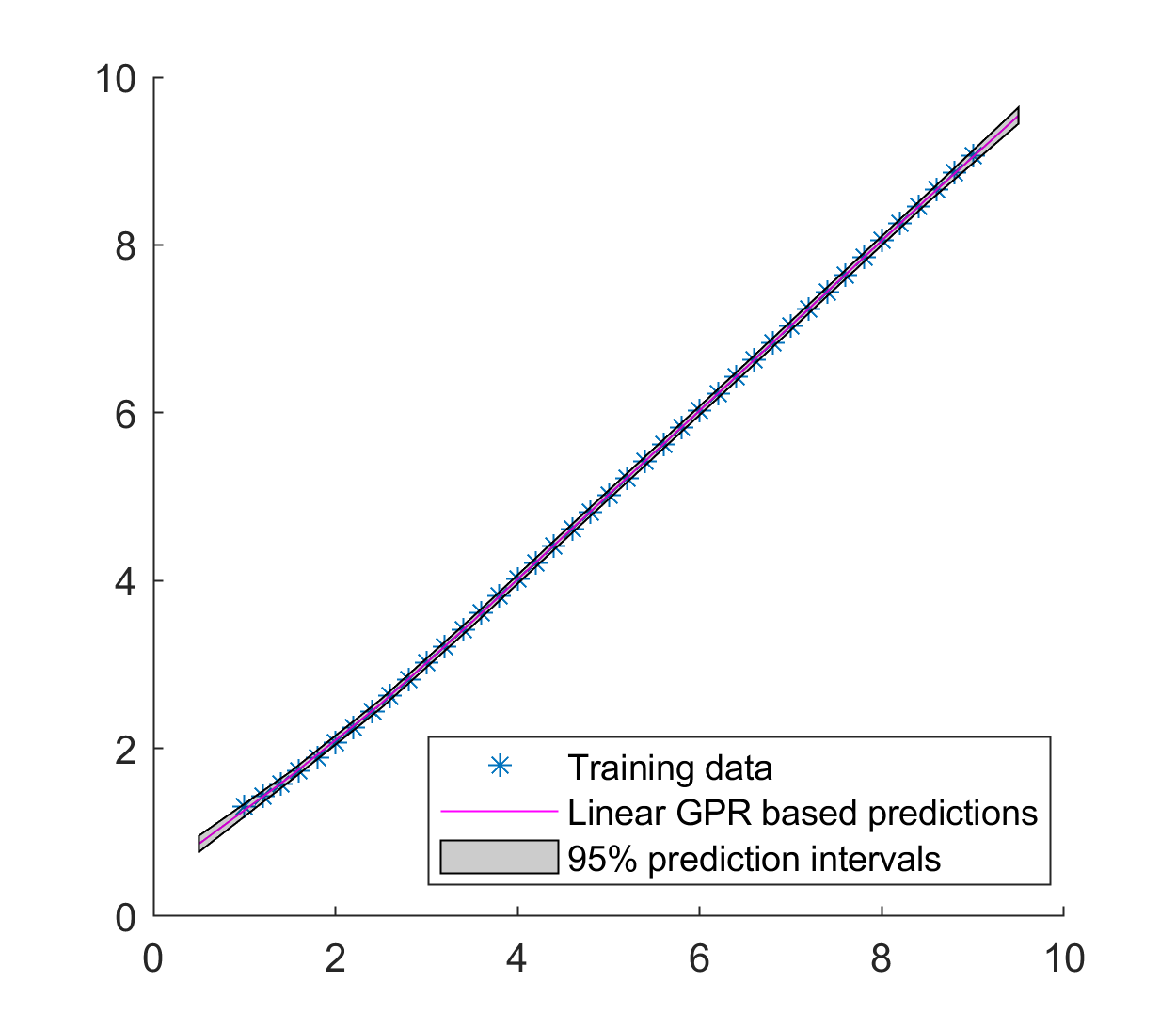}}

\subcaptionbox{$h=0.05$, $\delta \mu=0.4$}{
\includegraphics[width=6cm]{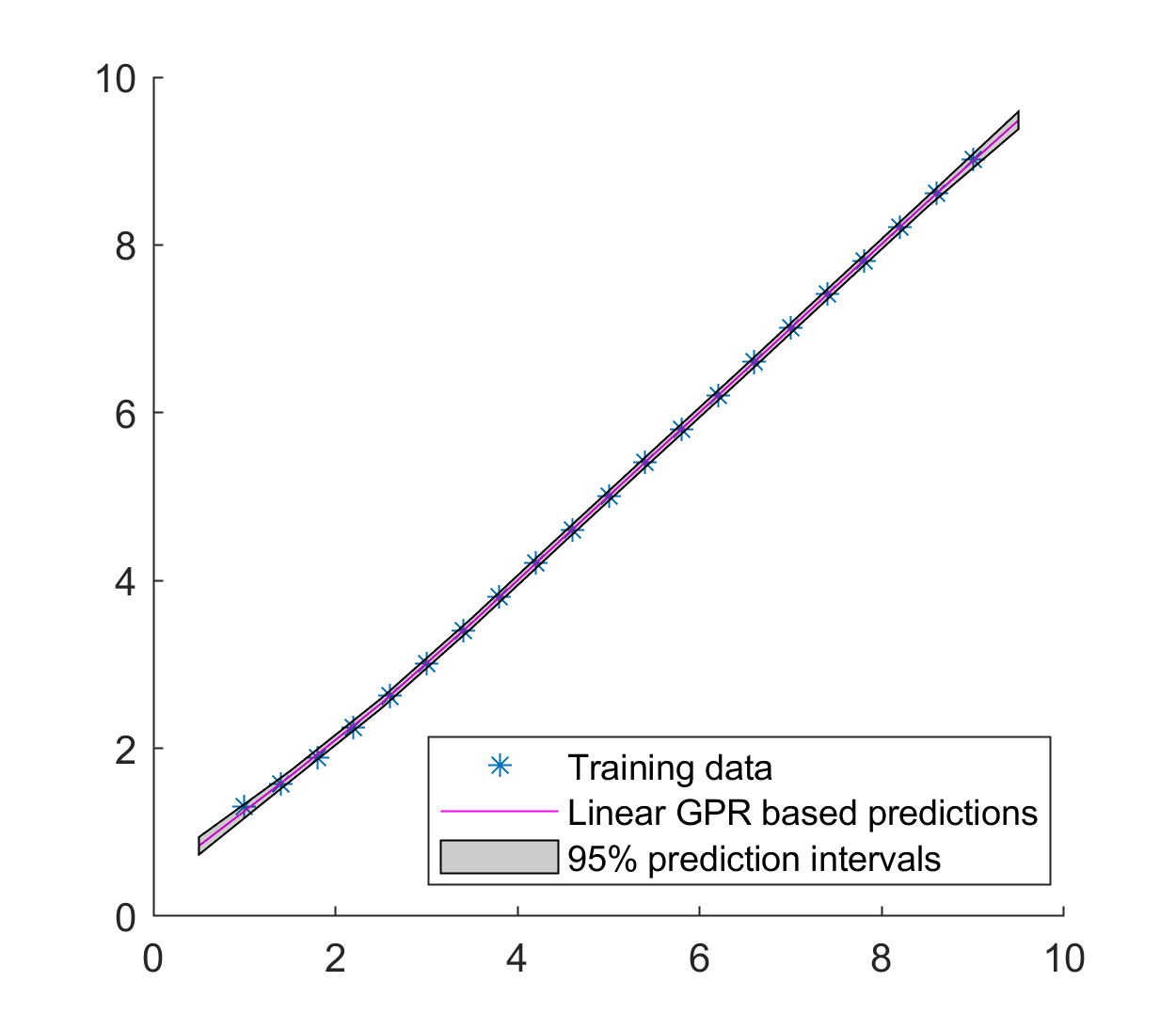}}
\subcaptionbox{$h=0.05$, $\delta \mu=0.2$}{
\includegraphics[width=6cm]{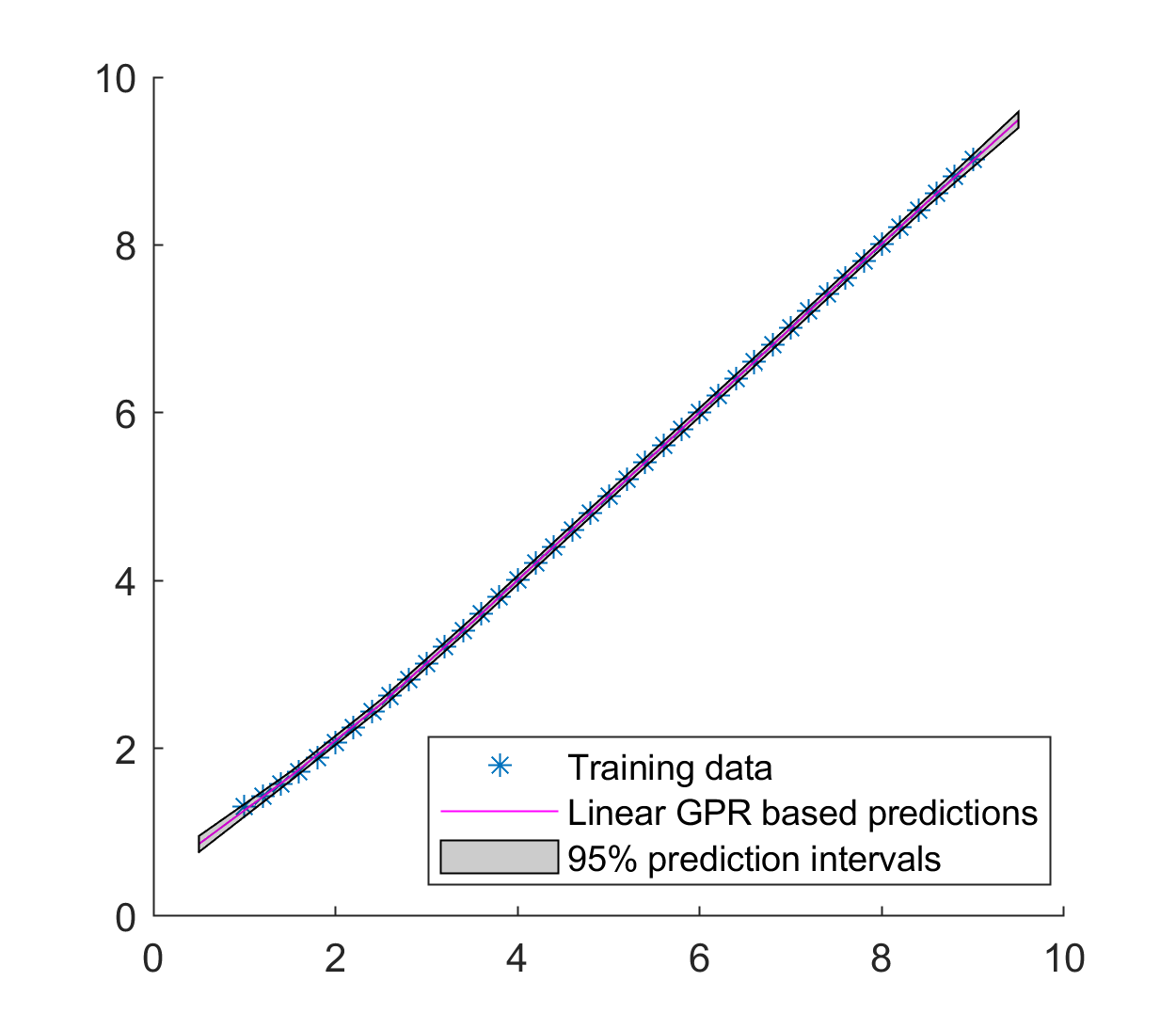}}
   \caption{The GPR plot for the first eigenvalues with $95\%$ confidence interval for Problem~\eqref{model1}.}
     \label{fig:ev_2dosc}       % Give a unique label
\end{figure}

In Table \ref{table:2dosc}, we have reported the first eigenvalue using the data-driven model corresponding to the mesh size $h=0.1$ and $0.05$ at test points $\mu=1.5$, $2.5$, $3.5$, $4.5$, $5.5$, $6.5$, $7.5$ and compared them with the results obtained by FEM. Note that the analytic first  eigenvalue is equal to $\mu$. From Table~\ref{table:2dosc} one can see that the results match pretty well with the result of FEM. In Figure~\ref{fig:ev_2dosc}, we have shown the eigenvalues used for the training and the predicted eigenvalues for all the four cases corresponding to $h=0.1$ and $0.05$ and $\delta\mu=0.2$ and $0.4$. We have also indicated the confidence interval for the predicted eigenvalues. The confidence interval is very narrow.

\begin{figure}
\centering
\includegraphics[width=4.5cm]{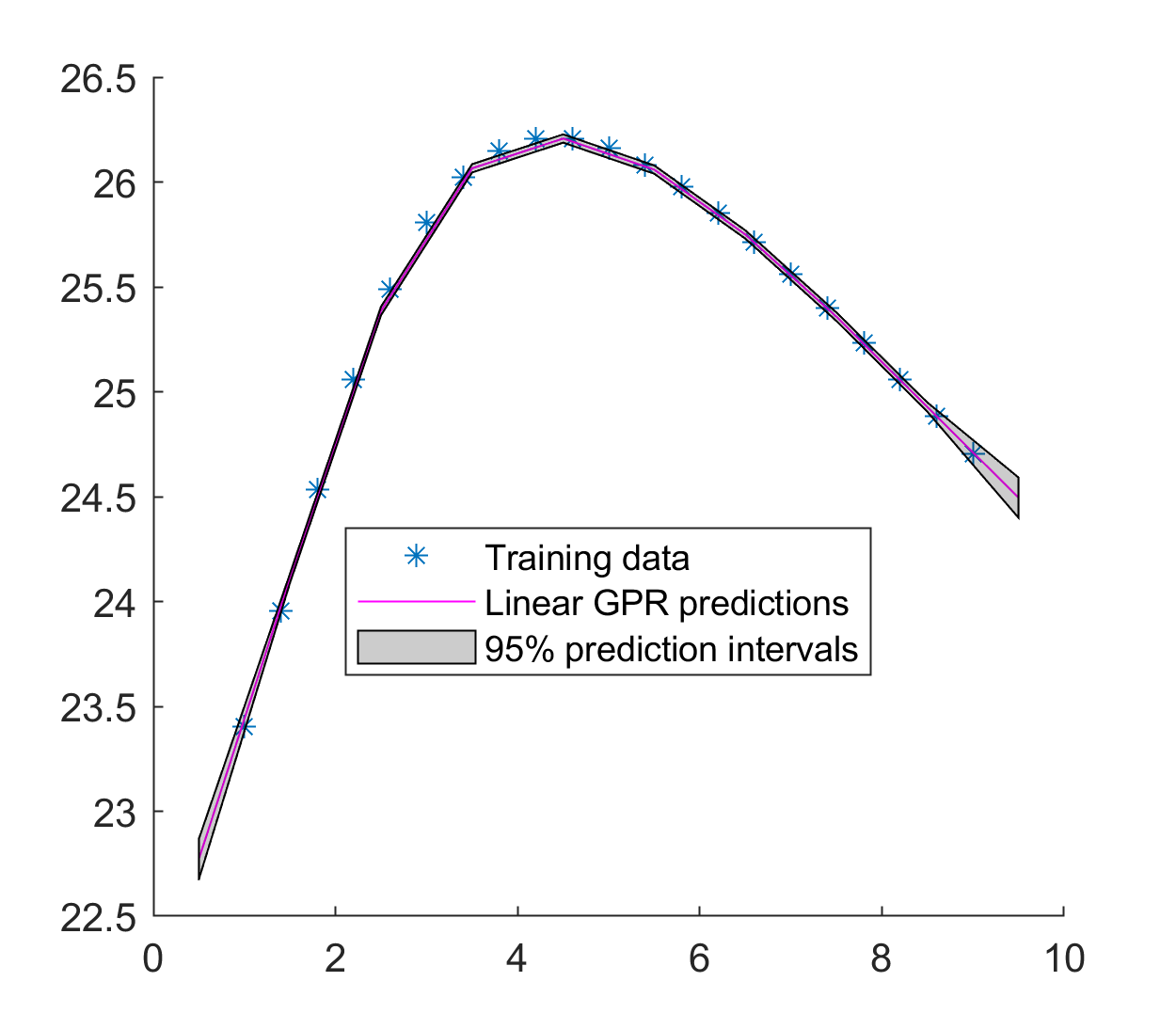}
\includegraphics[width=4.5cm]{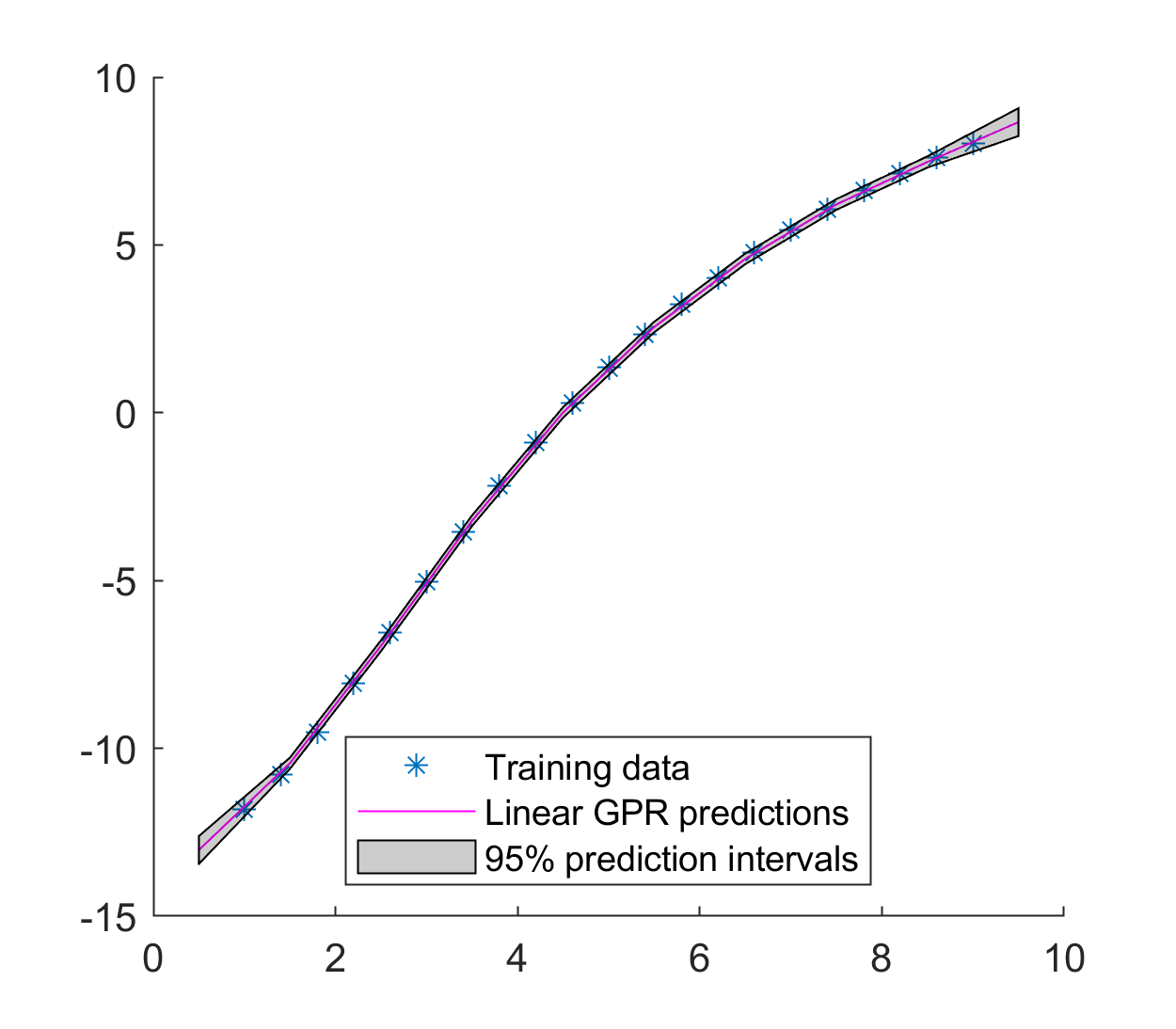}
\includegraphics[width=4.5cm]{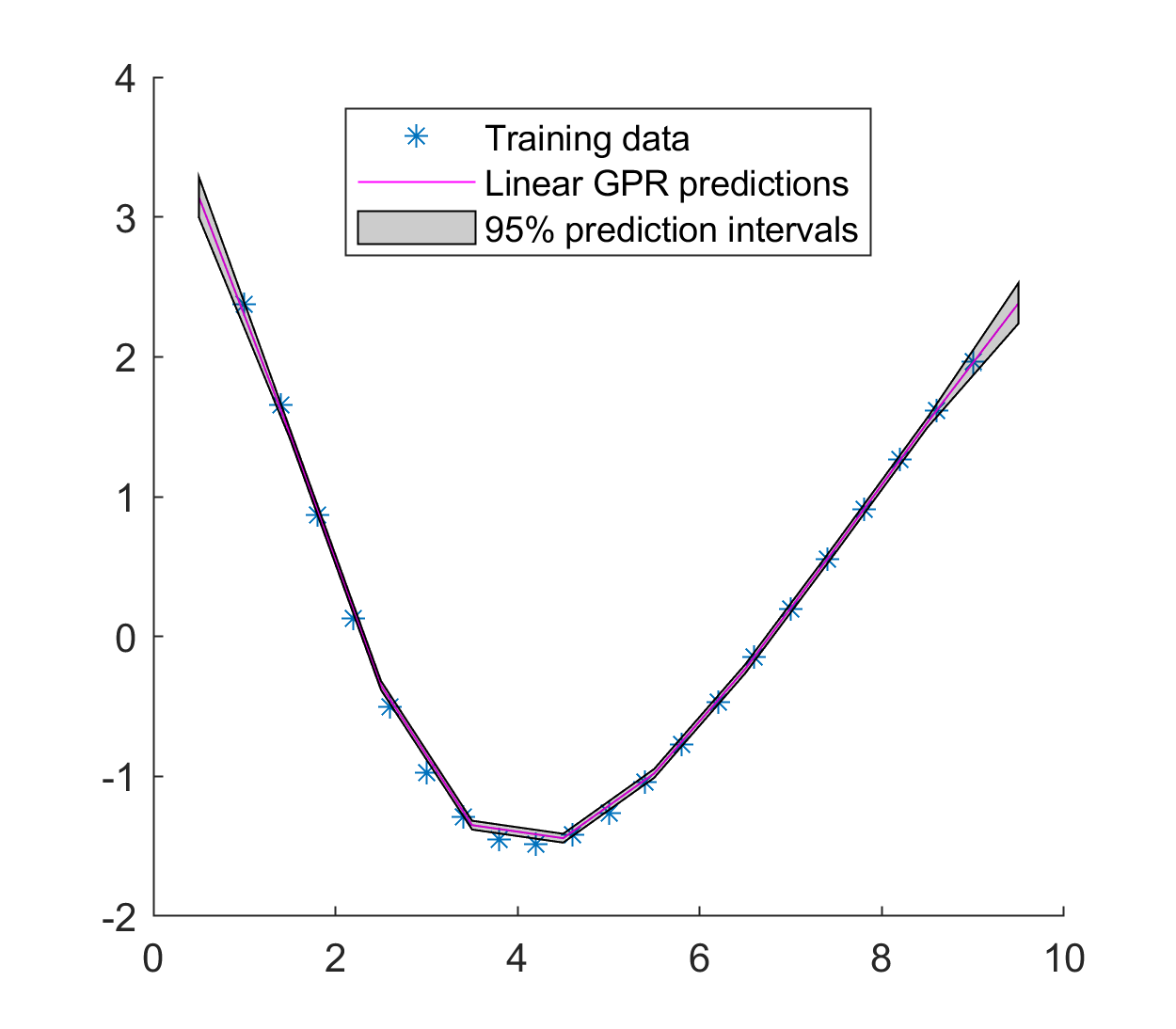}

\includegraphics[width=4.5cm]{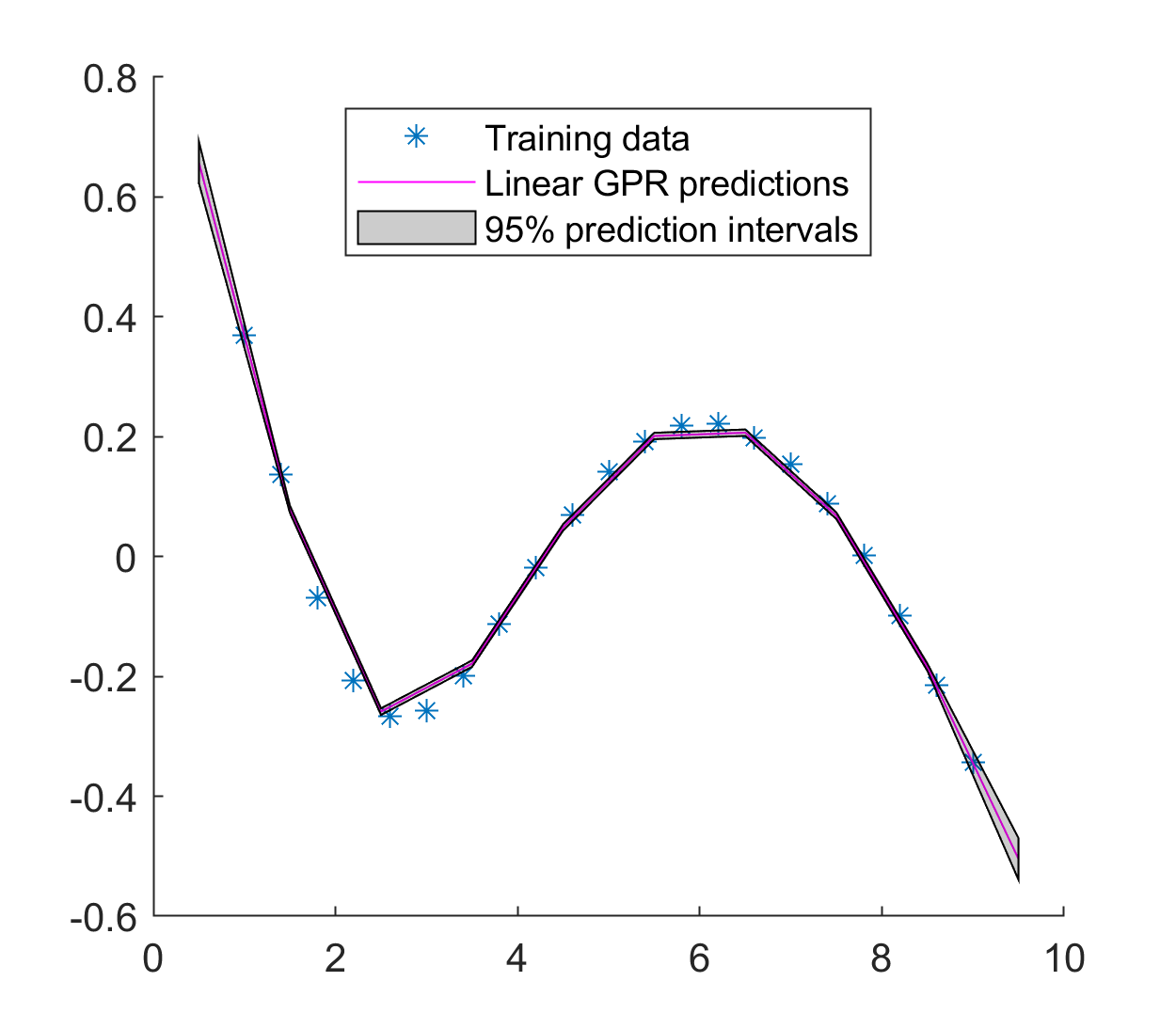}
\includegraphics[width=4.5cm]{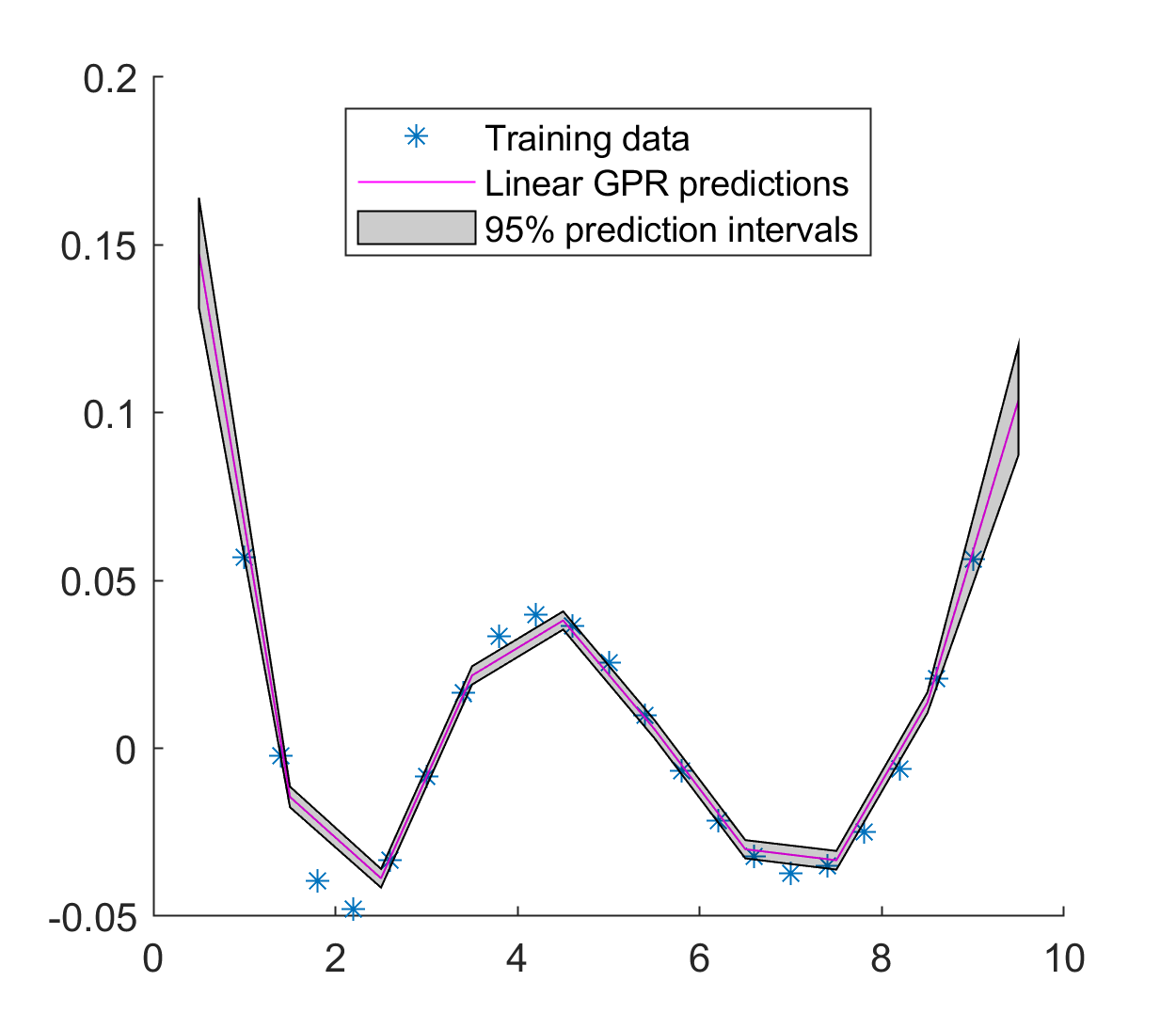}
\includegraphics[width=4.5cm]{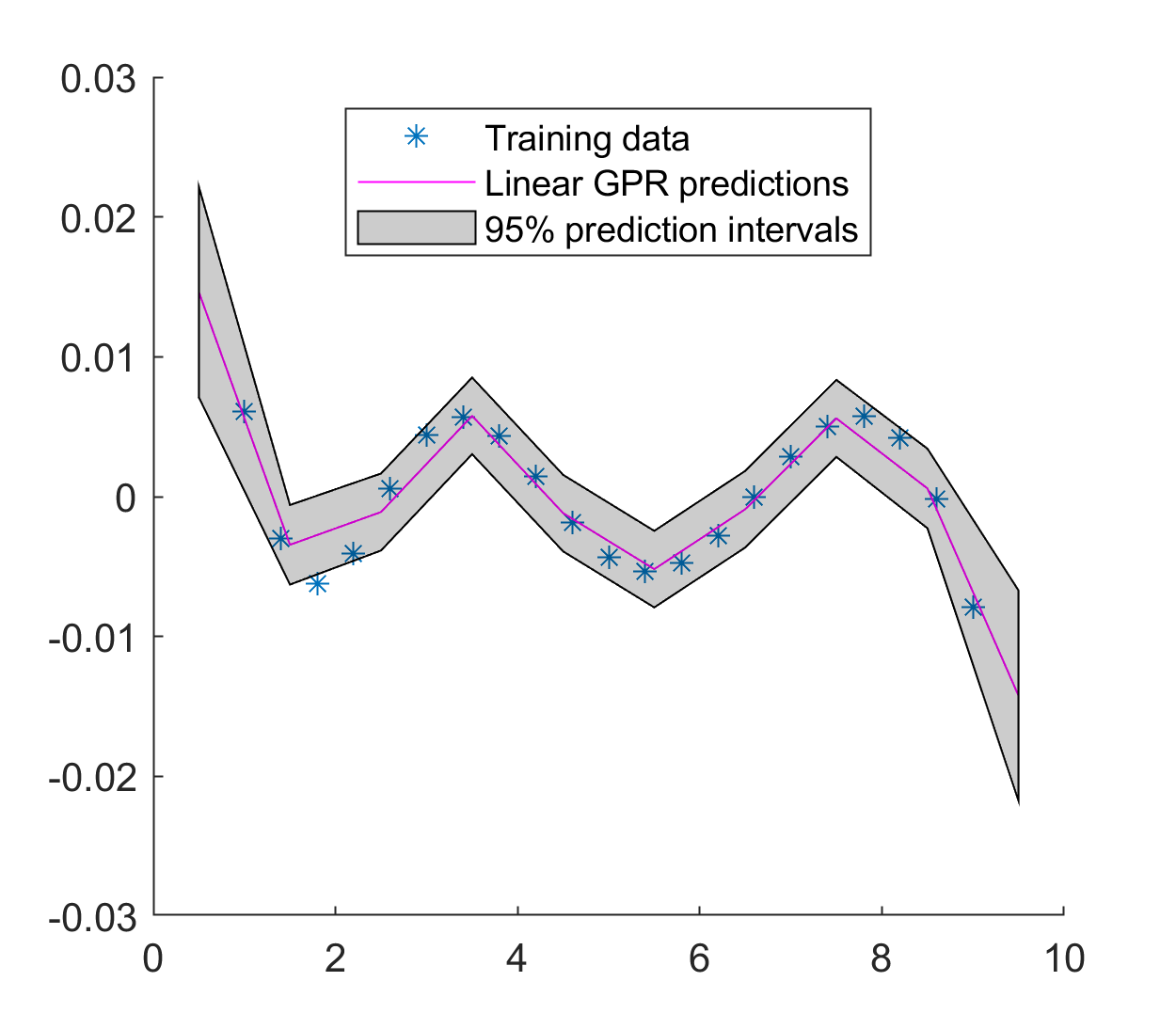}
   \caption{The GPR corresponding to the six RB coefficients for Problem~\eqref{model1} with $21$ sample points and $h=0.05$.}
     \label{fig:2dosc}       % Give a unique label
\end{figure}

In Figure~\ref{fig:2dosc}, we have presented the training set and the predicted results for the GPR corresponding to the six projected coefficients corresponding to first eigenvectors for the 2D oscillator problem for the case $h=0.05$ and $\delta\mu =0.4$. The training sets are plotted as star and the predicted values in magenta line and the 95 percentage of confidence interval is shown as shaded region. The test points are a uniform subdivision of $(0.5,9.5)$ with stepsize $0.1$. The endpoints of the testset lie outside the parameter interval, so the confidence interval is little wider at the endpoints.

\begin{figure}
\centering
\includegraphics[width=4.5cm]{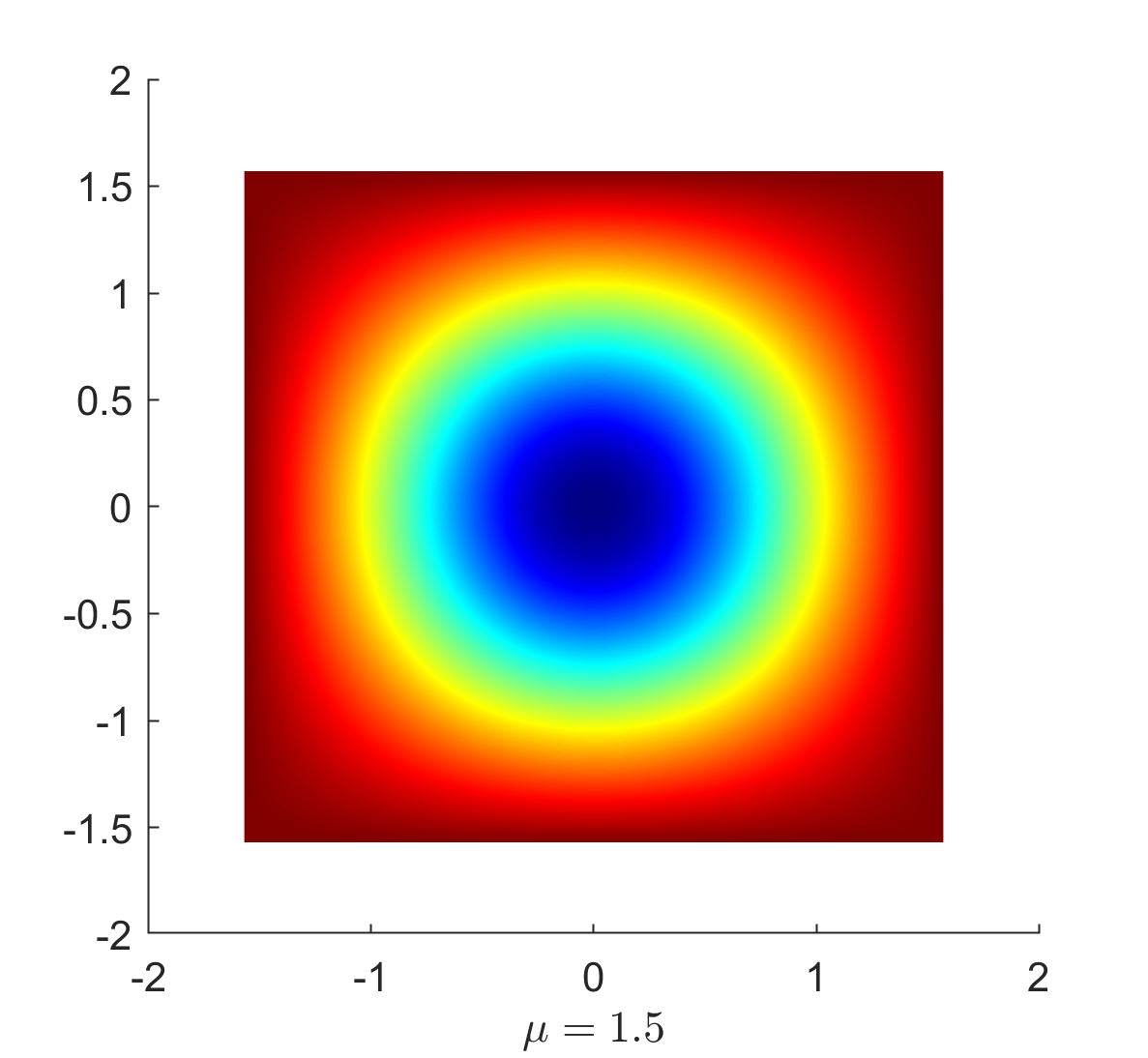}
\includegraphics[width=4.5cm]{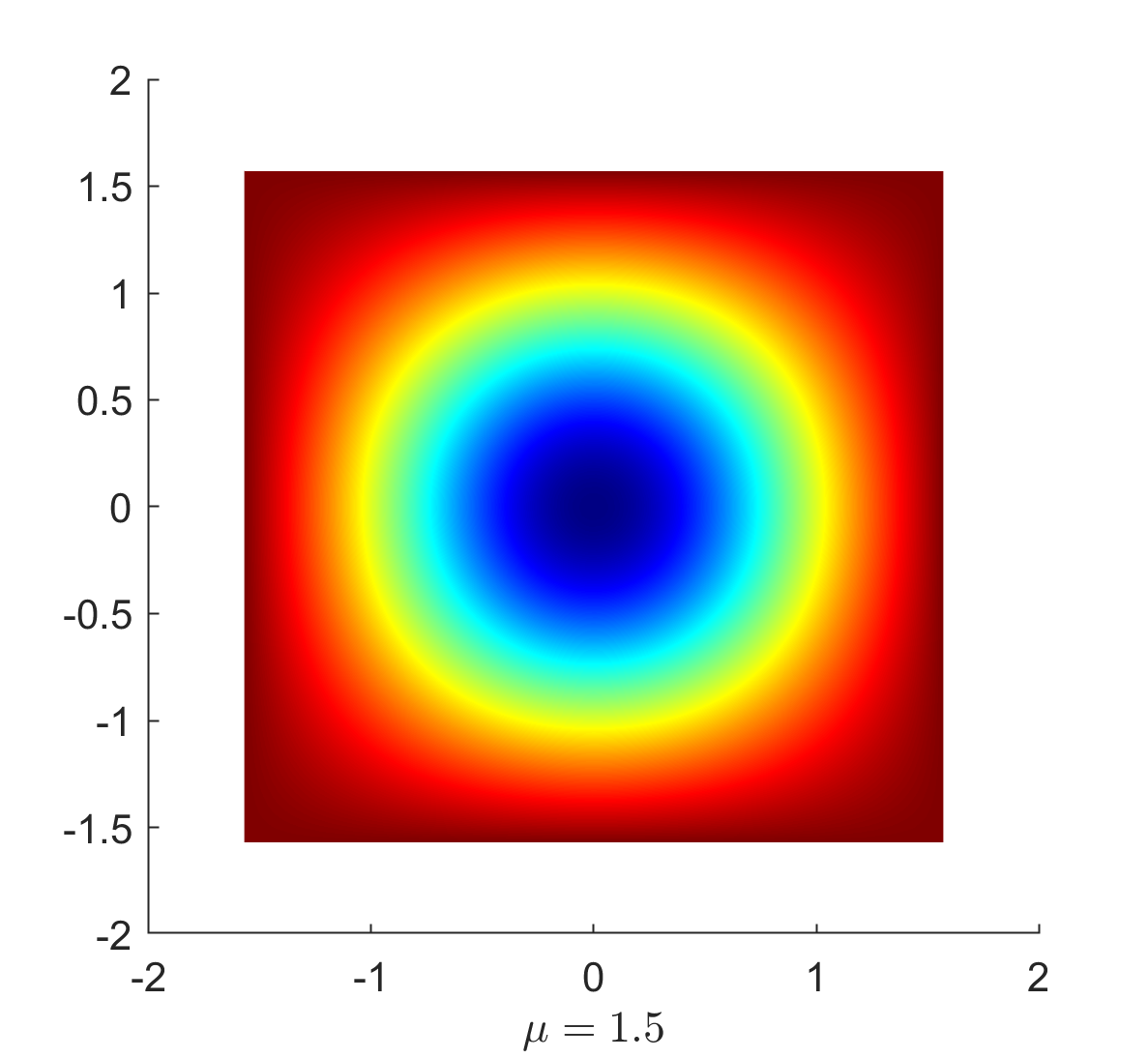}
\includegraphics[width=4.5cm]{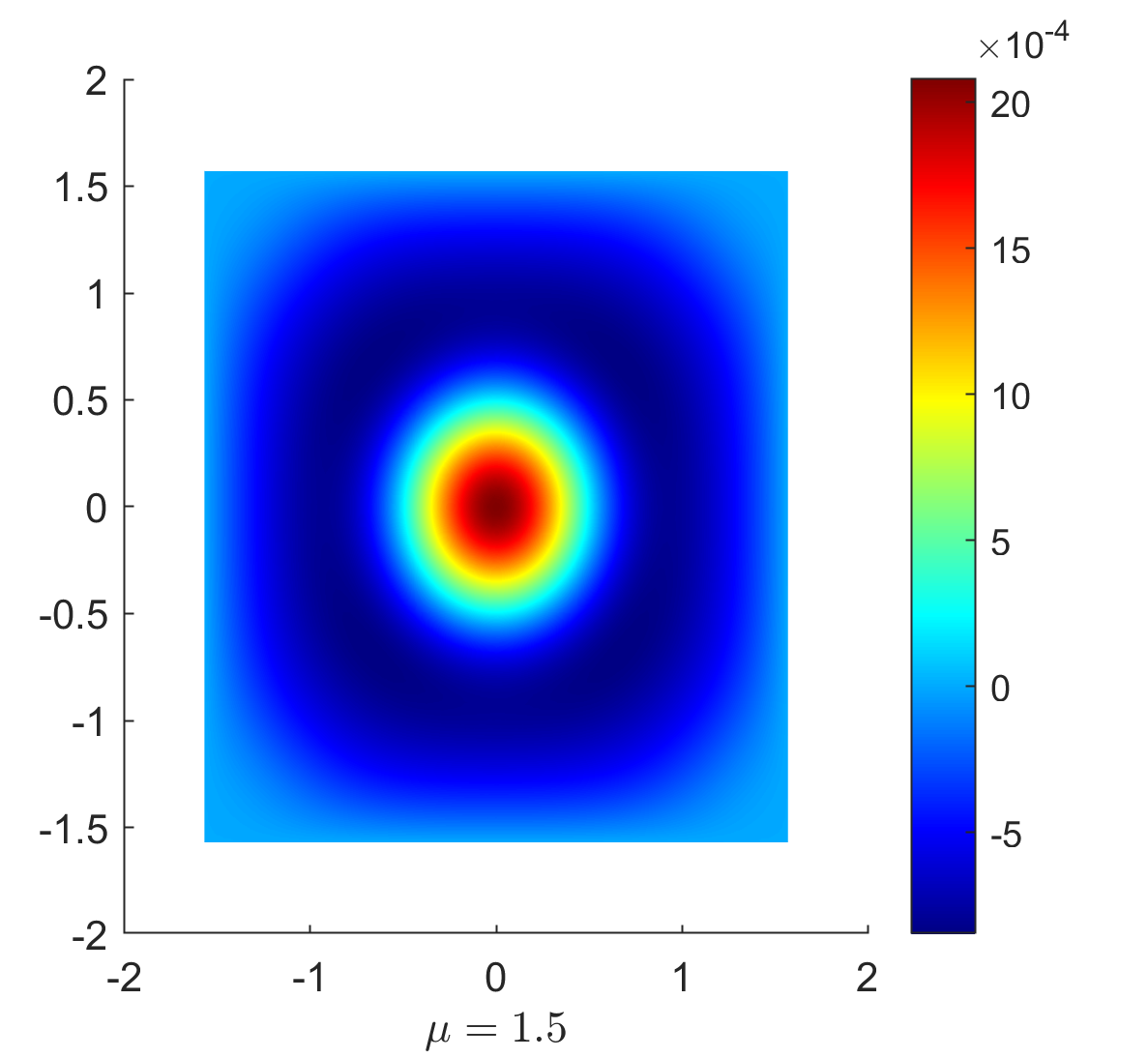}

\includegraphics[width=4.5cm]{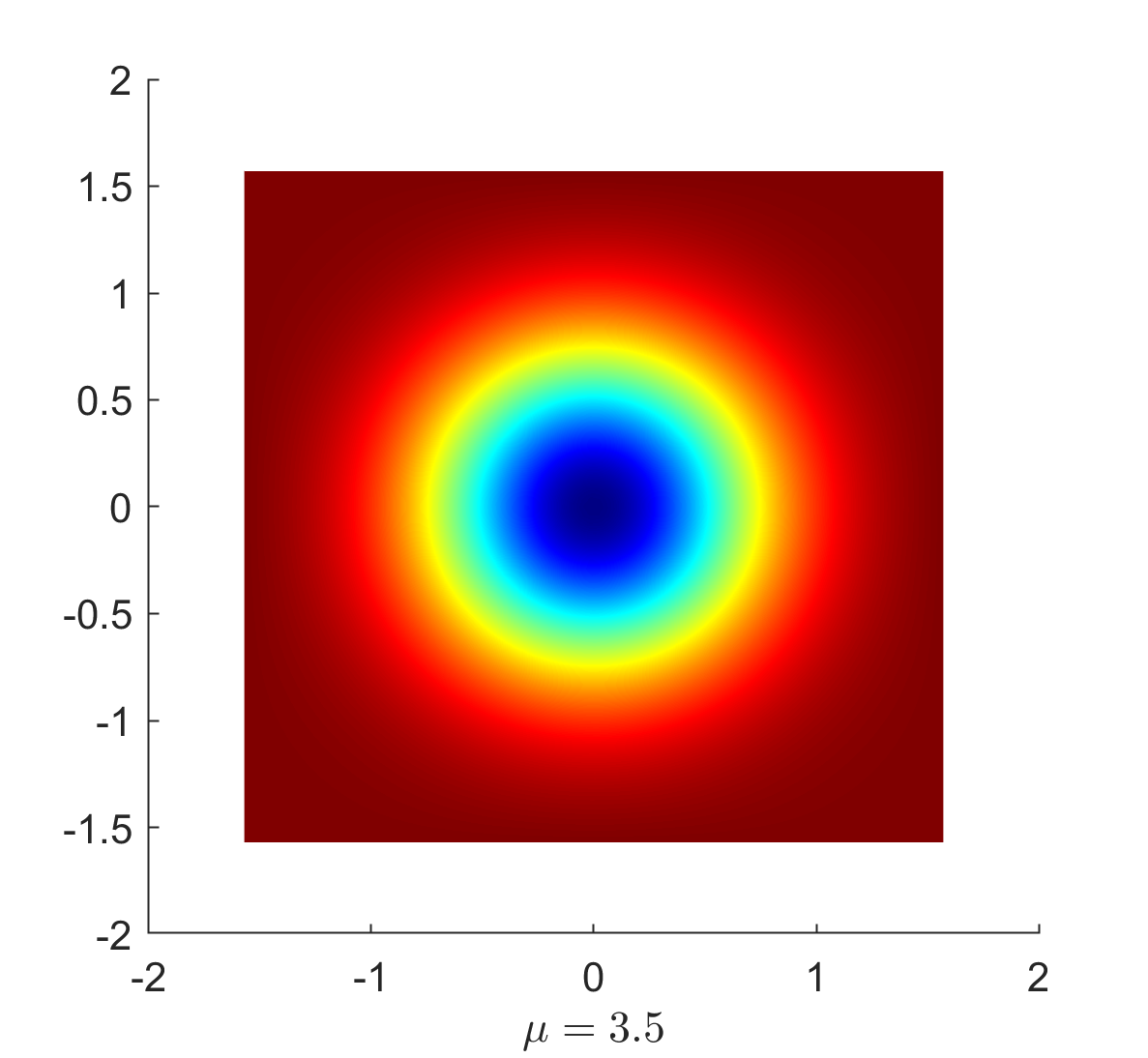}
\includegraphics[width=4.5cm]{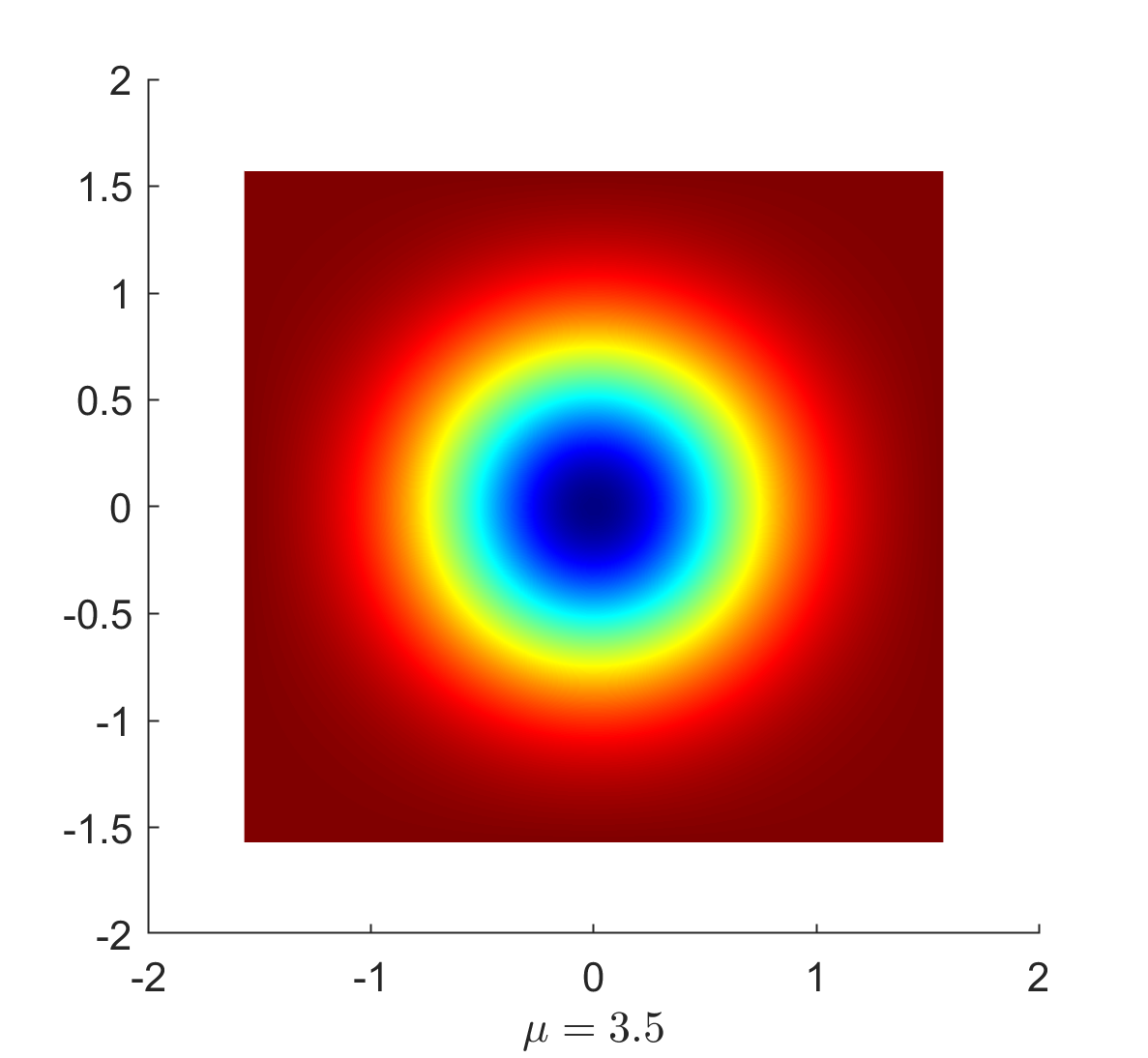}
\includegraphics[width=4.5cm]{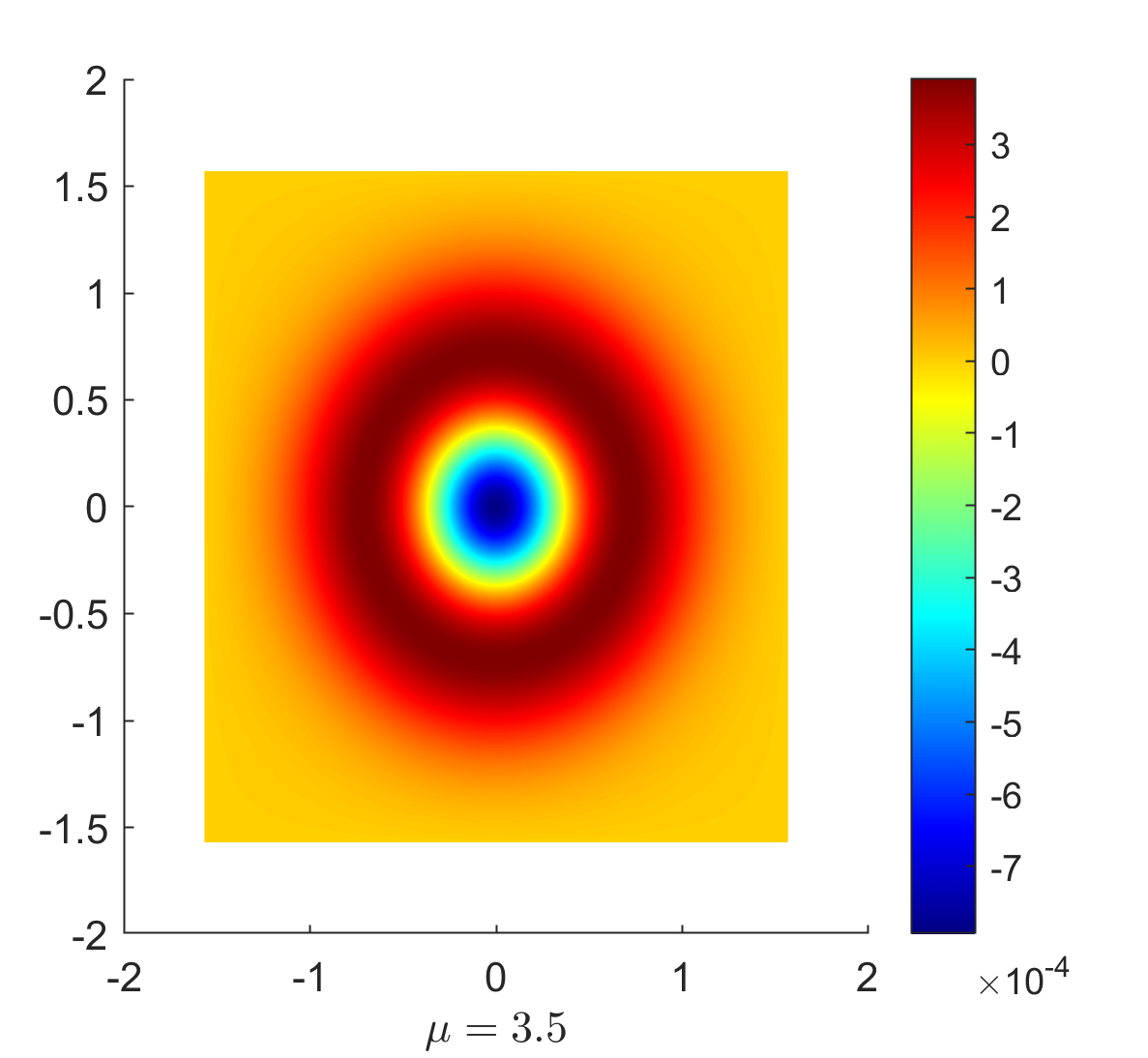}

\includegraphics[width=4.5cm]{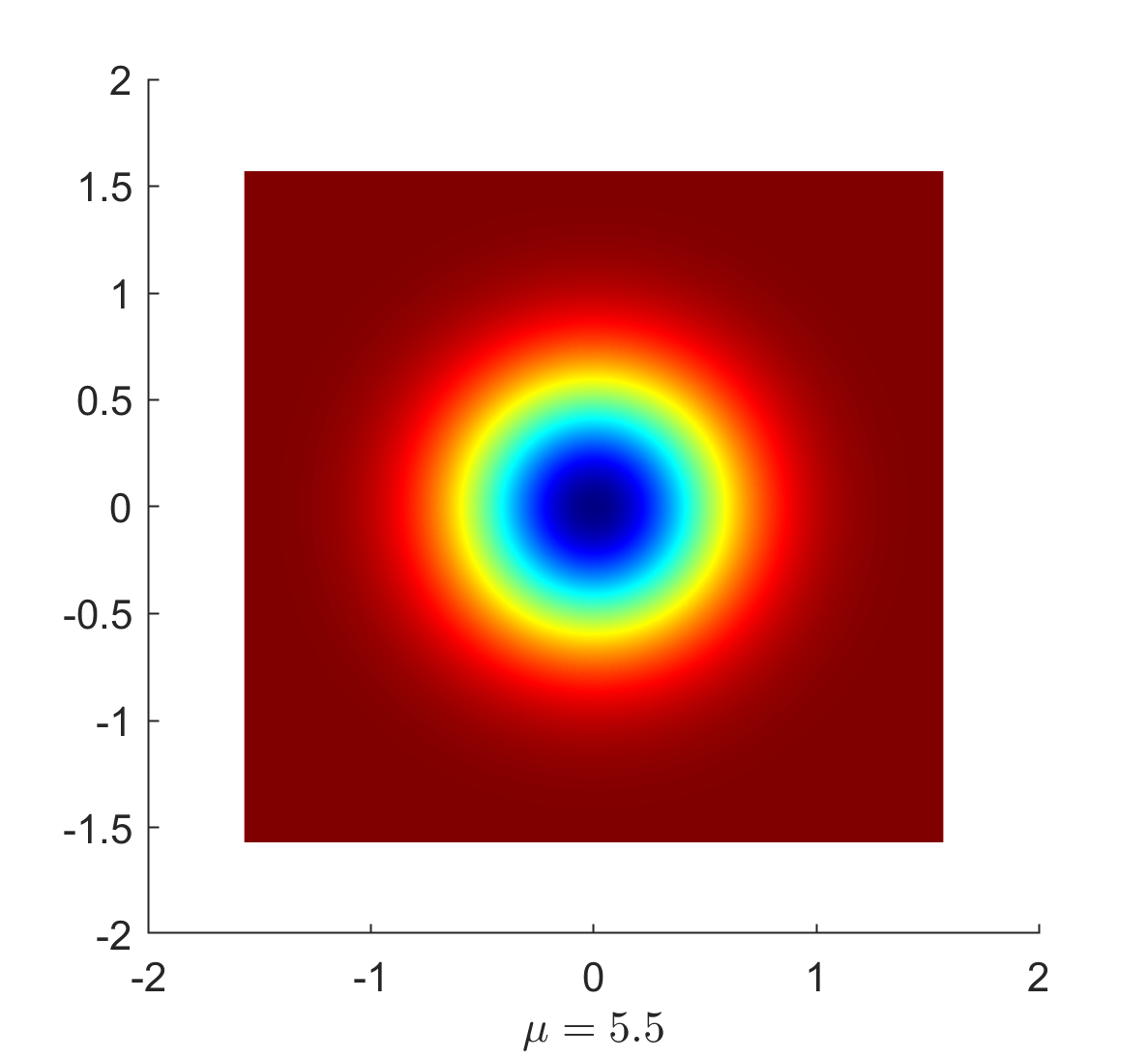}
\includegraphics[width=4.5cm]{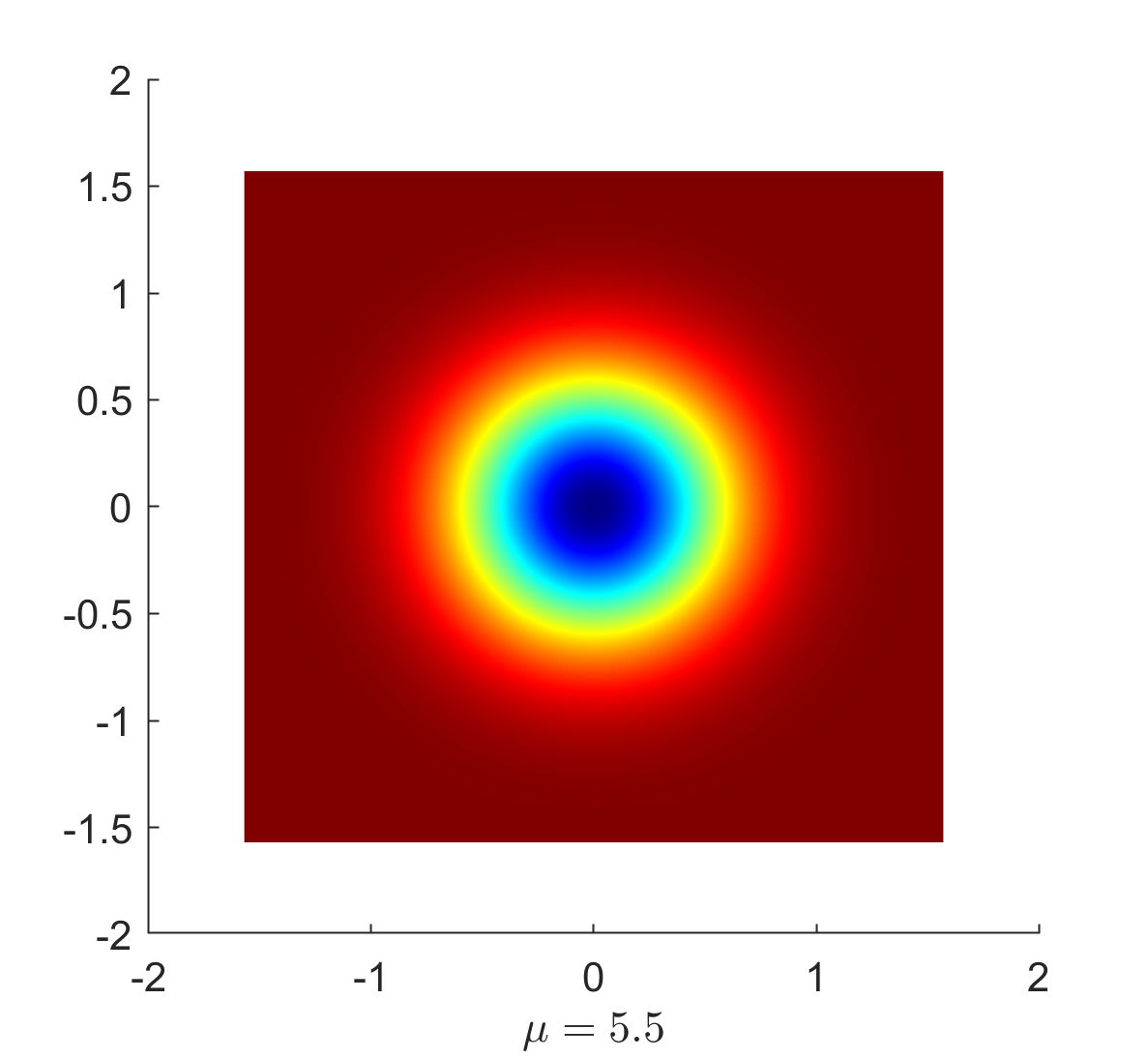}
\includegraphics[width=4.5cm]{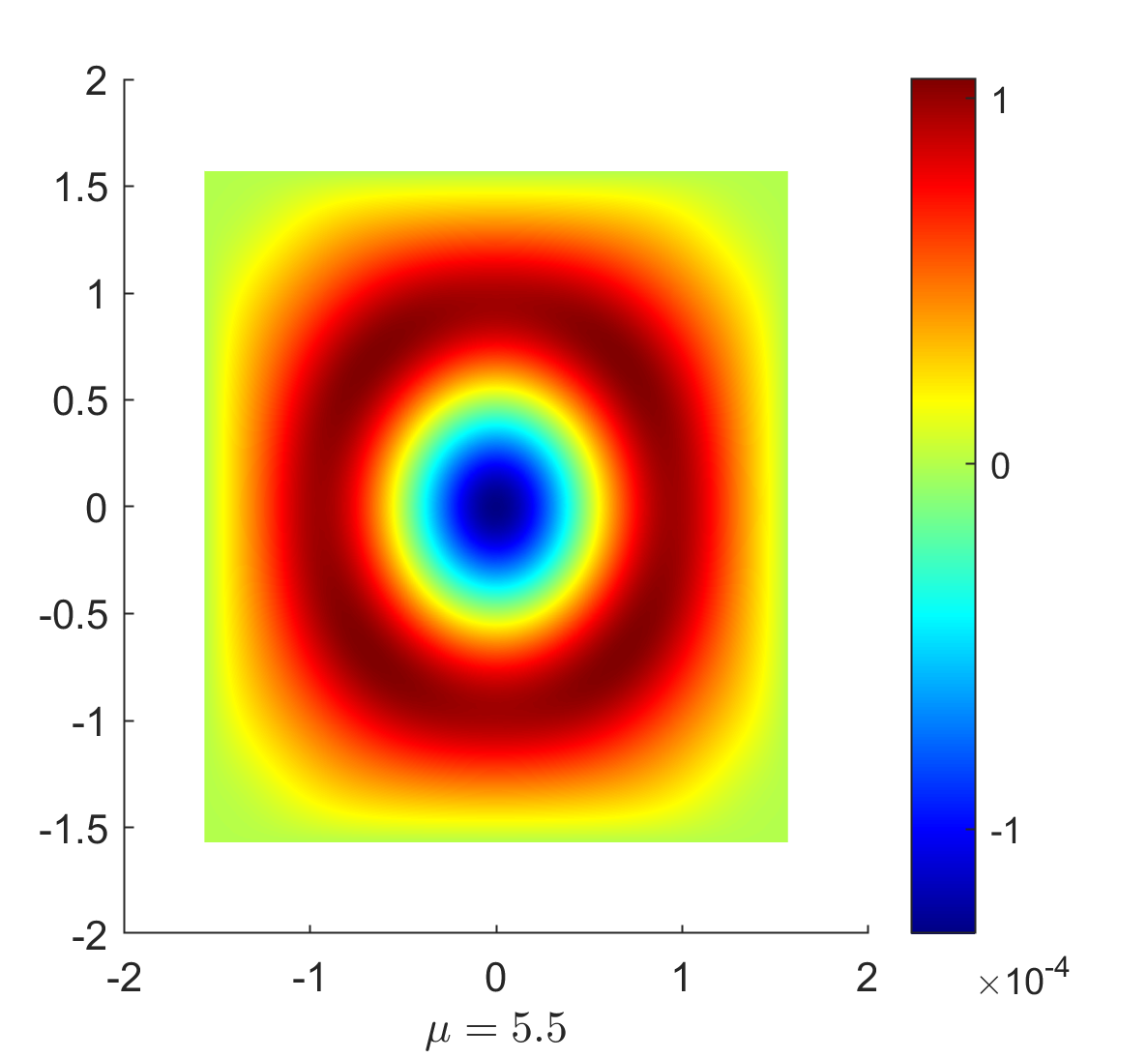}

\includegraphics[width=4.5cm]{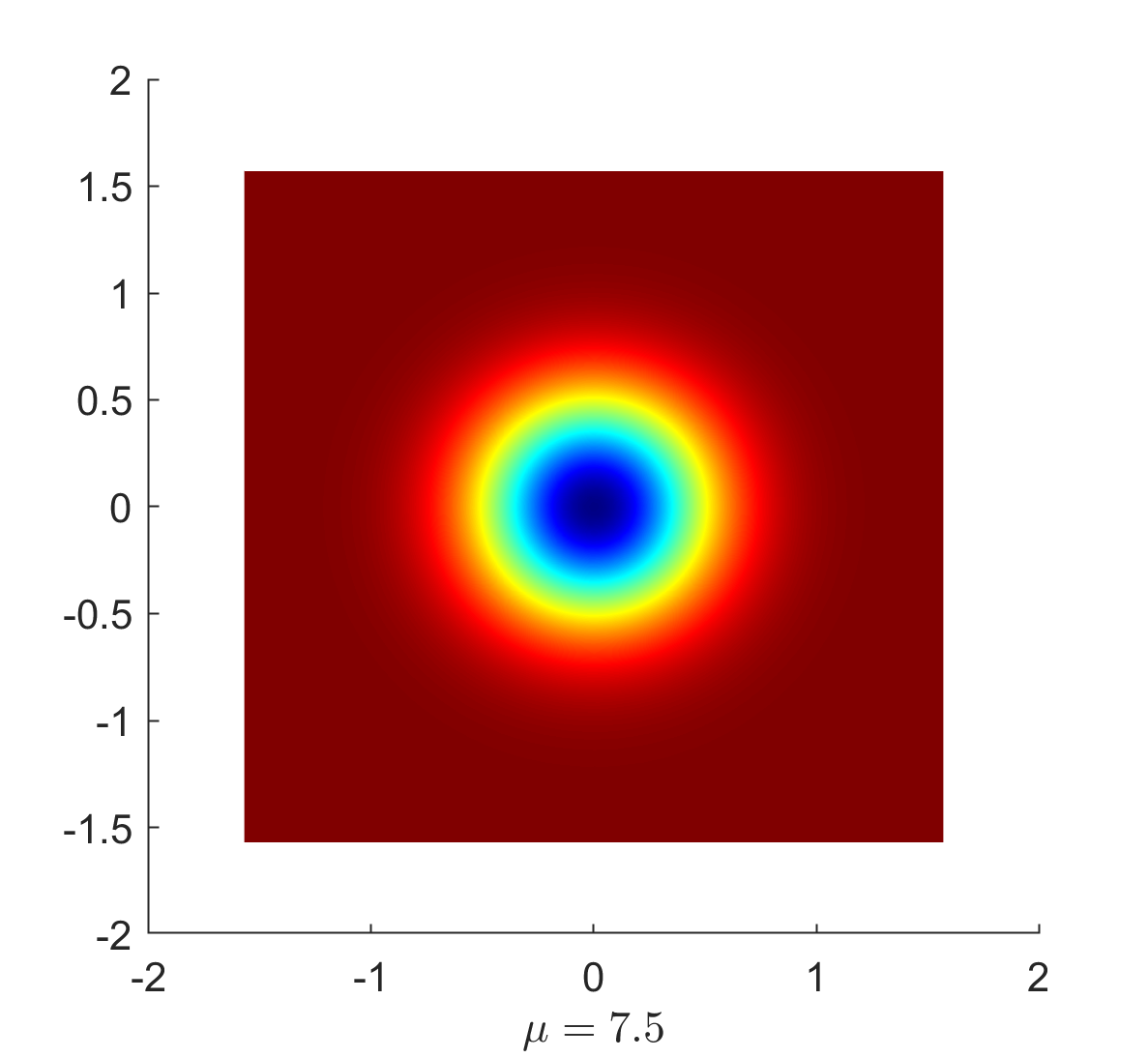}
\includegraphics[width=4.5cm]{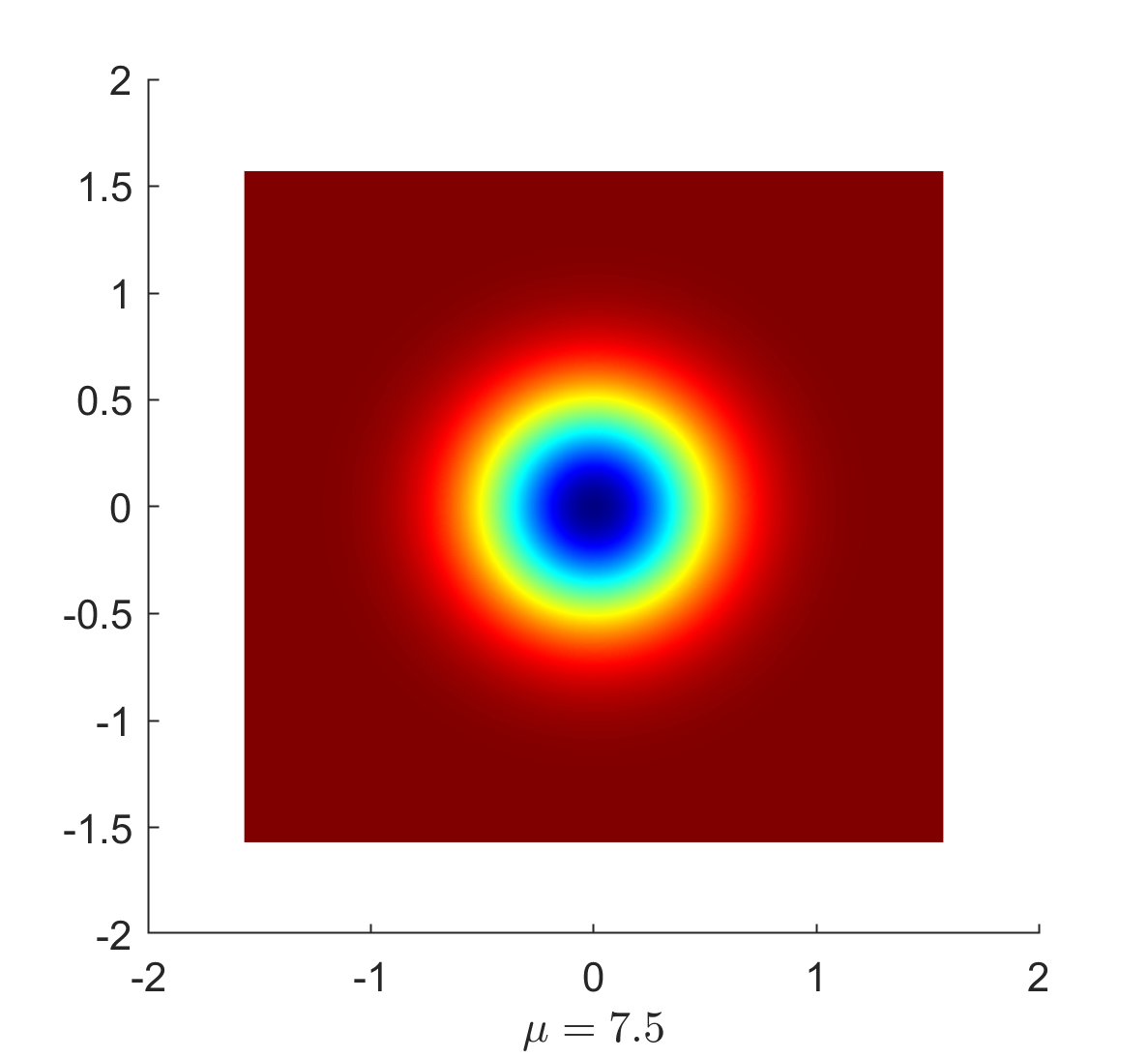}
\includegraphics[width=4.5cm]{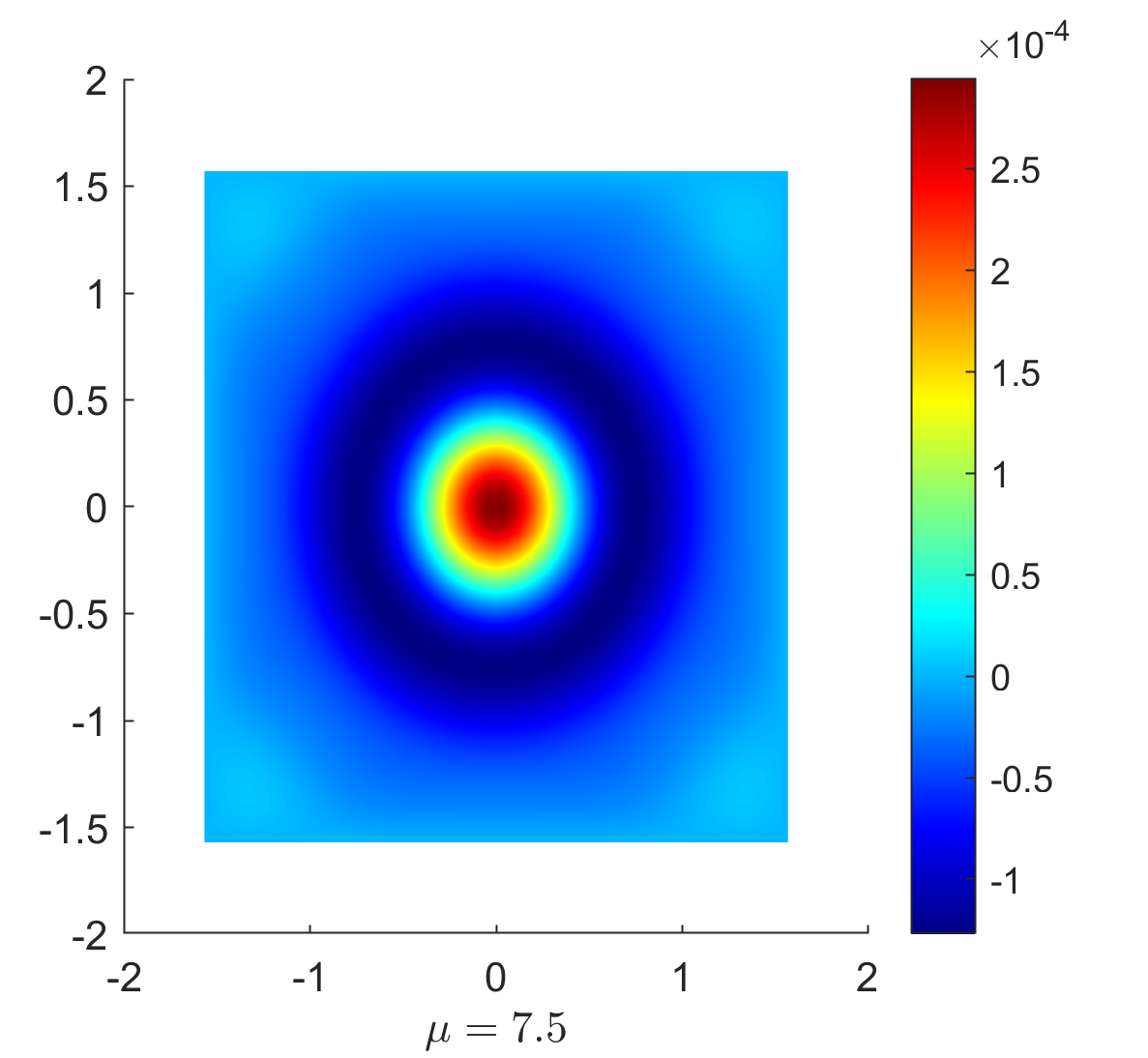}
\caption{Results for Problem~\eqref{model1} at the test points $\mu=1.5$, $3.5$, $5.5$, $7.5$, respectively, with $21$ sample points and $h=0.05$. First eigenfunction $u_{1,h}$ of FEM (left), first eigenfunction $\tilde{u}_{1,h}$ of DD (middle), error $(u_{1,h}-\tilde{u}_{1,h}$ (right)}
      \label{fig:evp2d}       % Give a unique label
 \end{figure}

The first eigenvectors obtained by our GPR-based DD model for the case  $h=0.05$ and $\delta \mu=0.4$, at the test points $\mu=1.5$, $3.5$, $5.5$, $7.5$ are presented in Figure~\ref{fig:evp2d} and compared with the eigenvectors obtained using FEM. Also, we have presented the error between the FEM eigenvectors and the DD-based eigenvectors. We can see that the error is of order $10^{-4}$.

% 	 \begin{table}[t]
%  	 	\centering
%  	\begin{tabular}{|c|c|c|c|c|c|c|c|c|c|c|c|c|} 
%  		\hline
% 		n& {\begin{tabular}[c]{@{}c@{}}Bound \end{tabular}} &
%  		  {\begin{tabular}[c]{@{}c@{}} $\mu=2.5$ \end{tabular}} &
% 		 {\begin{tabular}[c]{@{}c@{}} $3.5$ \end{tabular}} &
%  		{\begin{tabular}[c]{@{}c@{}}  $4.5$ \end{tabular}}& 
%  		{\begin{tabular}[c]{@{}c@{}} $5.5$ \end{tabular}} & 
%  		{\begin{tabular}[c]{@{}c@{}} $6.5$ \end{tabular}}&
% 		{\begin{tabular}[c]{@{}c@{}} $7.5$ \end{tabular}} & 
% 		{\begin{tabular}[c]{@{}c@{}} $8.5$ \end{tabular}} \\
% 		%{\begin{tabular}[c]{@{}c@{}} $\mu=2.5$ \end{tabular}}  \\	
% 	\hline
% 	32&Lower&2.4837&	3.4633&	4.4838&	5.4953&	6.5088&	7.5356&	8.5581	\\	
% &Upper&2.5901&	3.5686&	4.5890&	5.6004&	6.6141&	7.6420&	8.6669		\\	
% 	\hline
% 64&Lower&2.4765&	3.4504&	4.4608&	5.4601&	6.4611&	7.4710&	8.4746	\\	
% &Upper&2.5819&	3.5548&	4.5651&	5.5643&	6.5656&	7.5764&	8.5823	\\	
% \hline
% 128&Lower& 2.4747&	3.4471&	4.4550&	5.4513&	6.4492&	7.4548&	8.4537\\	
% &Upper&2.5799&	3.5514&	4.5591&	5.5553&	6.5534&	7.5600&	8.5612\\
% \hline
%  	\end{tabular}
% 	\caption{Upper and lower bound of the Eigenvalues obtained using DD Model.}
%  	\label{table:2dbound}
%  \end{table}
%%%%%%%%%%%%%%%%%%%%%%%%%%%%%%%%%%%%%%%%%%%%%%%%%%%%%%%%%%
\subsection{Eigenvalue problems with nonlinear parametric dependence}
%%%%%%%%%%%%%%%%%%%%%%%%%%%%%%%%%%%%%%%%%%%%%%%%%%%%%%%%%%%%%%%
Consider the nonlinear eigenvalue problem of the following form: given any $\mu \in \mathcal{P}$ find $\lambda(\mu)\in \mathbb{R}$ and $u(\mu)\in V$ with $u(\mu)\ne0$ satisfying
\begin{equation}\label{meqn}
\left\{
\aligned
&a(u(\mu),v)+\mu^2\int_{\Omega}g(u(\mu))v=\lambda(\mu) m(u(\mu),v), \quad \forall v\in V\\
& (u(\mu),u(\mu))=1,
\endaligned
\right.
\end{equation}
where $g$ is a nonlinear function.
In order to solve the problem we apply Newton's method as follows
\[
u^{k+1}=u^k+\delta u, \quad \lambda^{k+1}=\lambda^k+\delta \lambda
\]
so that~\eqref{meqn} reduces to
\[
\left\{
\aligned
&a(u^k+\delta u,v)+\mu^2\int_{\Omega}g(u^k+\delta u)v=(\lambda^k+\delta \lambda) (u^k+\delta u,v), \quad \forall v\in V\\
& (u^k+\delta u,u^k+\delta u)=1.
\endaligned
\right.
\]
Using Taylor's expansion of the function $g$ and neglecting the higher order terms we get
\[
\left\{
\aligned
&a(u^k,v)+a(\delta u,v)+\mu^2\int_{\Omega}g(u^k)v+\mu^2\int_{\Omega}\delta u g'(u^k)v\\
&\hspace{5cm}=\lambda^k (u^k,v)+\delta \lambda(u^k,v)+\lambda^k(\delta u,v)&&\forall v\in V\\
&(u^k+,u^k)+2(\delta u,u^k)=1.
\endaligned
\right.
\]
This leads to
\[
\left\{
\aligned
&a(\delta u,v)+\mu^2\int_{\Omega}\delta u g'(u^k)v-\lambda^k(\delta u,v)-\delta \lambda(u^k,v)\\
&\hspace{4.5cm}=\lambda^k (u^k,v)-a(u^k,v)-\mu^2\int_{\Omega}g(u^k)v&&\forall v\in V\\
&2(\delta u,u^k)=1-(u^k,u^k).
\endaligned
\right.
\]
Let $V_h \subset V$ be the finite dimensional space. Then the Galerkin FE method corresponding to the nonlinear problem reads: given $(u_h^k,\lambda^k)$ find $(u_h^{k+1},\lambda^{k+1})=(u_h^k+\delta u_h,\lambda^k+\delta \lambda_h)$ such that
\[
\left\{
\aligned
&a(\delta u_h,v_h)+\mu^2\int_{\Omega}\delta u_h g'(u_h^k)v_h-\lambda^k(\delta u_h,v_h)-\delta \lambda_h (u_h^k,v_h)=\\
&\hspace{4cm}\lambda^k (u_h^k,v_h)-a(u_h^k,v_h)-\mu^2\int_{\Omega}g(u_h^k)v_h&&\forall v_h\in V_h\\
&2(\delta u_h,u_h^k)=1-(u_h^k,u_h^k),
\endaligned
\right.
\]
Let $\{ \varphi_j: j=1,\dots, N_h\}$ be the basis of $V_h$. Now expressing $u_h^k=\sum\limits_{j=1}^{N_h}u_{j,h}^k \varphi_j$ and $\delta u_h=\sum\limits_{j=1}^{N_h}d_j \varphi_j$
 we get the matrix form of the above discretized equation as
\begin{align*}
\begin{bmatrix}
A & -\pmb{c}\\
2 \pmb{c}^T & 0\\
\end{bmatrix}
\begin{bmatrix}
\pmb{d} \\ \delta \lambda_h
\end{bmatrix} = \begin{bmatrix}
\pmb{b} \\ 1-s
\end{bmatrix},
\end{align*}
where we have $A=A_1+\mu^2 M_{g'} -\lambda^k M$ with
\[
A_1(i,j)=a(\varphi_j,\varphi_i),\quad M(i,j)=(\varphi_j,\varphi_i), \quad M_{g'}(i,j)=\int_{\Omega}g'(u^k) \varphi_j\varphi_i,
\]
and $\pmb{b}, \pmb{c}$ are the vectors satisfying $\pmb{c}=M\pmb{u}^k_h(\mu)$, $\pmb{b}=\lambda^k M \pmb{u}_h^k-A_1\pmb{u}_h^k-\mu^2 \pmb{g}$ with $\pmb{g}=(g_1,\dots,g_{N_h})^\top$ and $g_j=\int_{\Omega}g(u^k_h)\varphi_j$ and the scalar $s=\pmb{u}_h^k(\mu)^\top M\pmb{u}_h^k(\mu)$ with $\pmb{u}_h^k=(u_{1,h}^k,\dots,u_{N_h,h}^k)^\top$.

Solving the above equation we will get $\delta u_h$ and $\delta \lambda $. Thus the updated eigenvalue is $\lambda^{k+1}=\lambda^k+\delta \lambda$ and updated eigenvectors $u^{k+1}(\mu)=u^k(\mu)+\delta u_h(\mu)$. For nonlinear eigenvalue problems we only focus on the first eigenvalue and first eigenvector.

For our numerical tests we consider $\Omega=(0,1)\subset \mathbb{R}$ and the parameter space is $\mathcal{P}=[1,9]$. For given $\mu \in \mathcal{P}$, we solve Problem~\eqref{meqn} for the following particular choice of data:
\[
a(w,v)=\int_{\Omega} \frac{dw}{dx}\frac{dv}{dx}, \quad m(w,v)=\int_{\Omega} wvdx, \quad g(w)=|w|^{7/3}w.
\]
We select two sets of sample parameters in order to calculate the snapshots. First we take $41$ uniformly distributed points and then $21$ uniform points in the parameter space $\mathcal{P}$, that is the parameter stepsize will be $\delta\mu=0.2$ and $0.4$, respectively. First eigenvalues at test parameters uniformly distributed in $\mu=(1.5,6.5)$ with stepsize $1$, obtained by our GPR-based DD model, are reported in Table~\ref{table5} and compared with the first eigenvalues obtained using FEM. We have reported the results with mesh size $h=0.05$ and $h=0.01$. The eigenvalues obtained by FEM and the DD model are very close to each other. Note that because of the nonlinearity, the FEM methods is iterative in nature and it takes few iterations to get the solution to converge but in the data-driven model we are getting directly the solution by evaluating the regression functions.

\begin{table}
\footnotesize
 	 	\centering
 	\begin{tabular}{|c|c|c|c|c|c|c|c|c|} 
 		\hline
 		h & {\begin{tabular}[c]{@{}c@{}} $\delta \mu$ \end{tabular}} &
		  {\begin{tabular}[c]{@{}c@{}} Method \end{tabular}} &
 		 {\begin{tabular}[c]{@{}c@{}} 1.5 \end{tabular}} & 
 		 {\begin{tabular}[c]{@{}c@{}} 2.5 \end{tabular}} & 
 		 {\begin{tabular}[c]{@{}c@{}} 3.5 \end{tabular}} & 
 		 {\begin{tabular}[c]{@{}c@{}} 4.5 \end{tabular}} & 
 		 {\begin{tabular}[c]{@{}c@{}} 5.5 \end{tabular}} & 
 		{\begin{tabular}[c]{@{}c@{}} 6.5 \end{tabular}}\\
 \hline
   0.05&0.4&FEM&13.4620&19.5299&28.1425&39.0260&52.0268&67.0749\\
   & & DD&13.4926& 19.6184& 28.1278& 38.9666 &51.9936& 67.0948\\
   \hline
   0.05&0.2&FEM&13.4620&19.5299&28.1425&39.0260&52.0268&67.0749\\
   & & DD&13.4807& 19.5773& 28.1241& 38.9963& 52.0197& 67.0857\\
   \hline
     0.01&0.4&FEM&13.4553&19.5367&28.1596&39.0516&52.0624&67.1239\\
      & & DD&13.4863& 19.6256& 28.1447& 38.9918 &52.0289& 67.1440\\ 
  \hline
  0.01&0.2&FEM&13.4553&19.5367&28.1596&39.0516&52.0624&67.1239\\
      & & DD&13.4743& 19.5842 &28.1409& 39.0218& 52.0552 &67.1348\\ 
  \hline
 	\end{tabular}
	\caption{1st eigenvalues of the nonlinear problem~\eqref{meqn} obtained using DD model with mesh size $h=0.05,0.01$.}
	 	\label{table5}
\end{table}

\begin{figure}
\centering
\subcaptionbox{$\mu_{tr}=1:0.2:9$}{
\includegraphics[width=6.5cm]{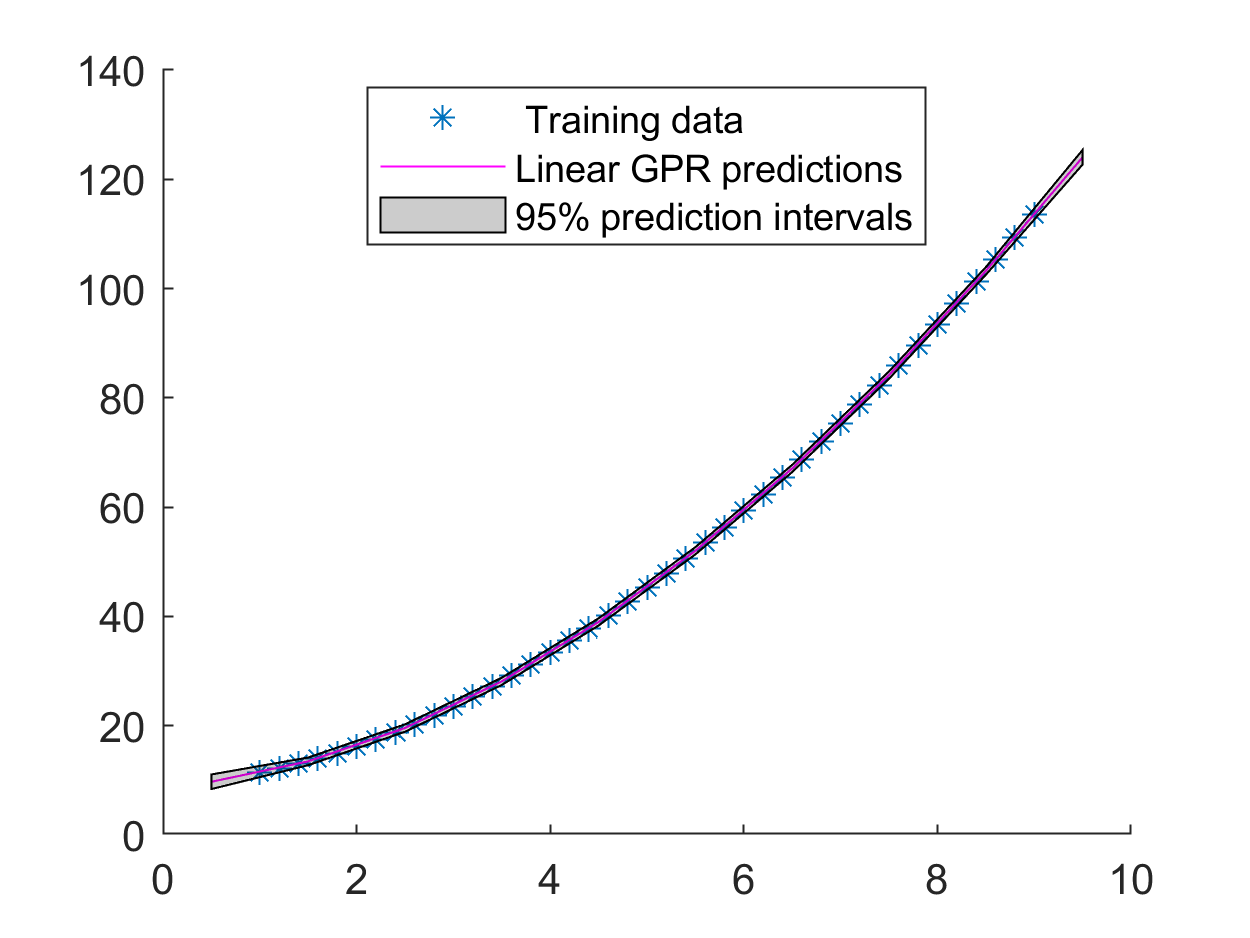}}
\subcaptionbox{$\mu_{tr}=1:0.4:9$}{
\includegraphics[width=6.5cm]{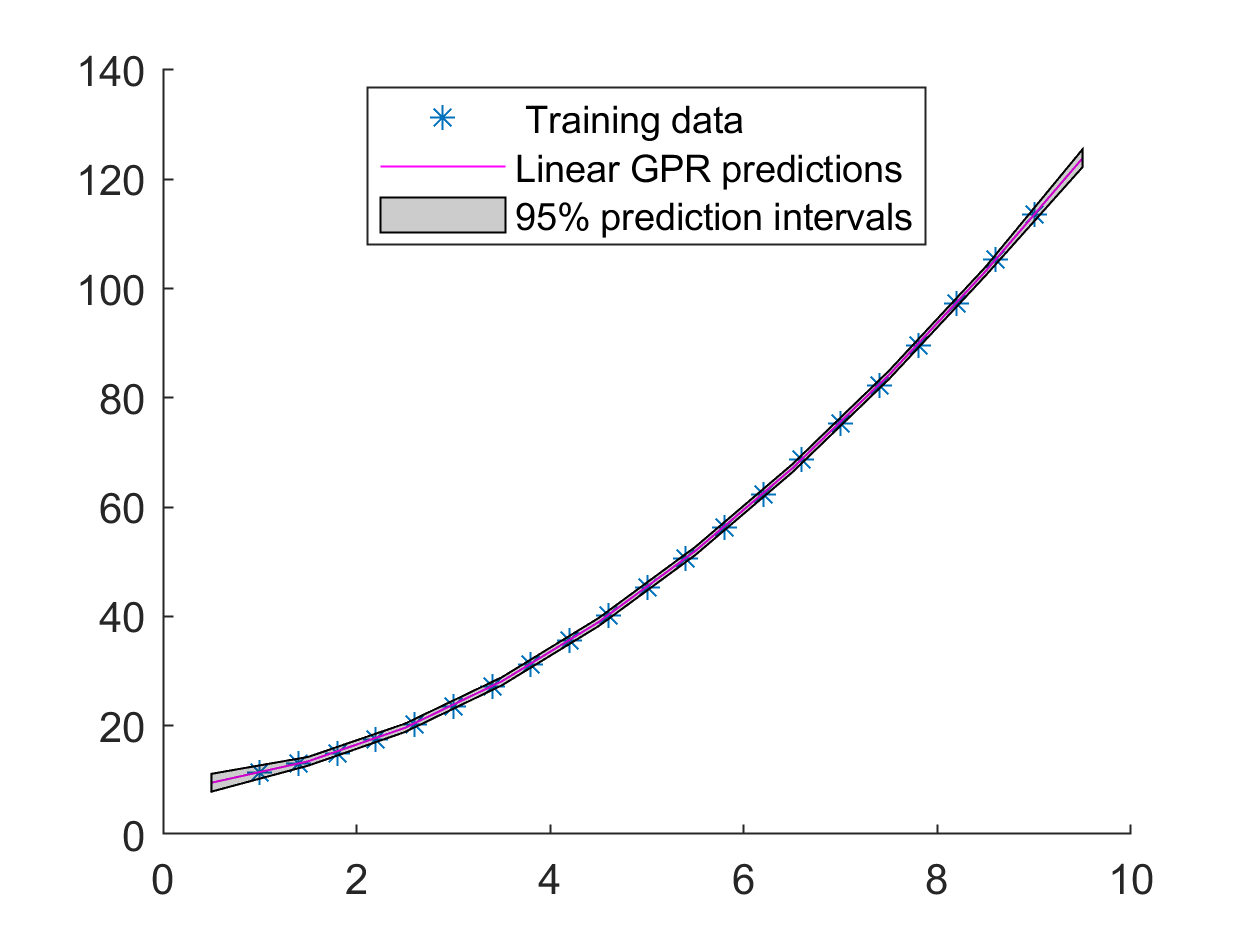}}
\caption{The GPR for the first eigenvalues of the nonlinear problem~\eqref{meqn} with different number of snapshots with $h=0.05$.}
\label{fig1_nonlin}       % Give a unique label
\end{figure}

\begin{figure}
\centering
\includegraphics[width=6.5cm]{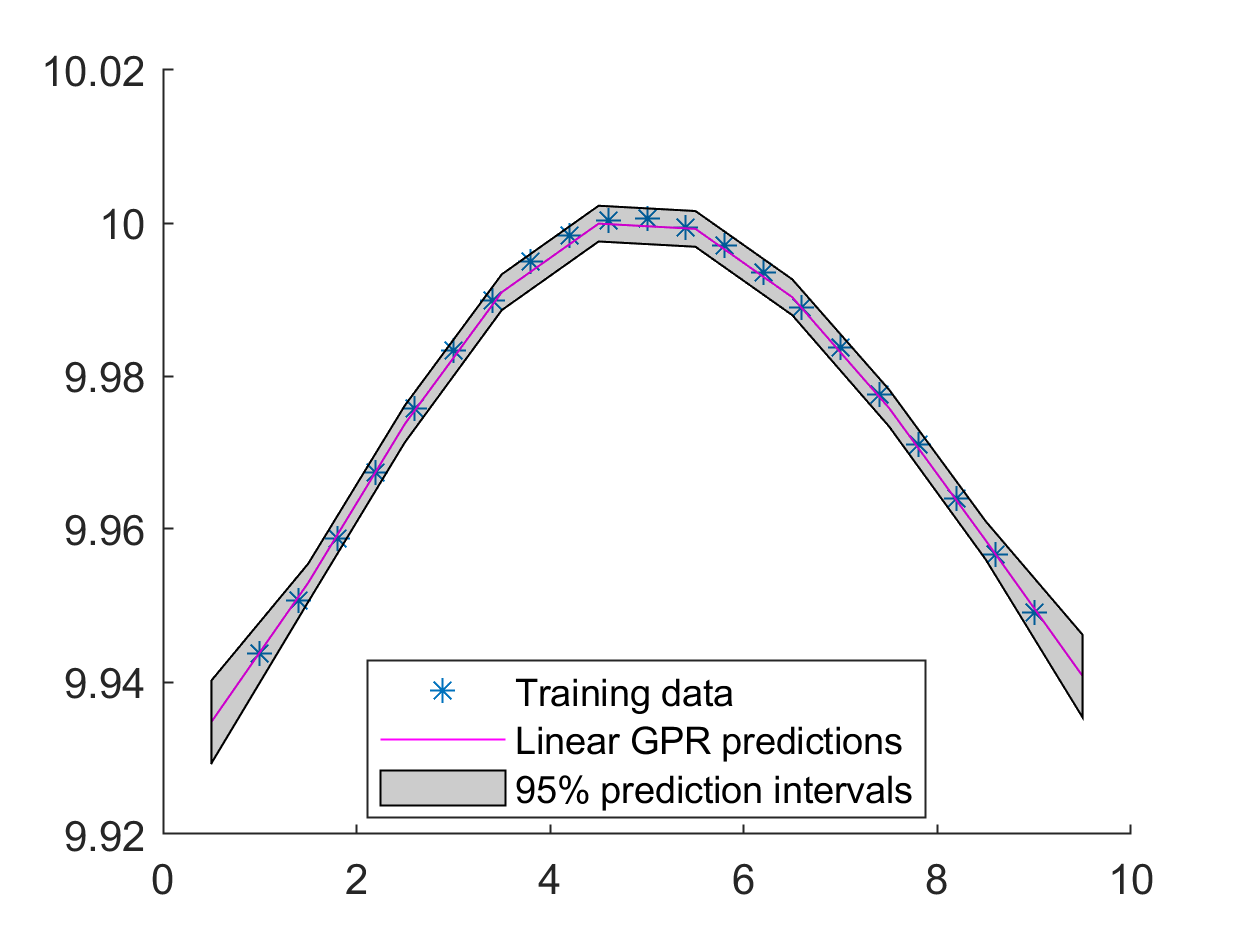}
\includegraphics[width=6.5cm]{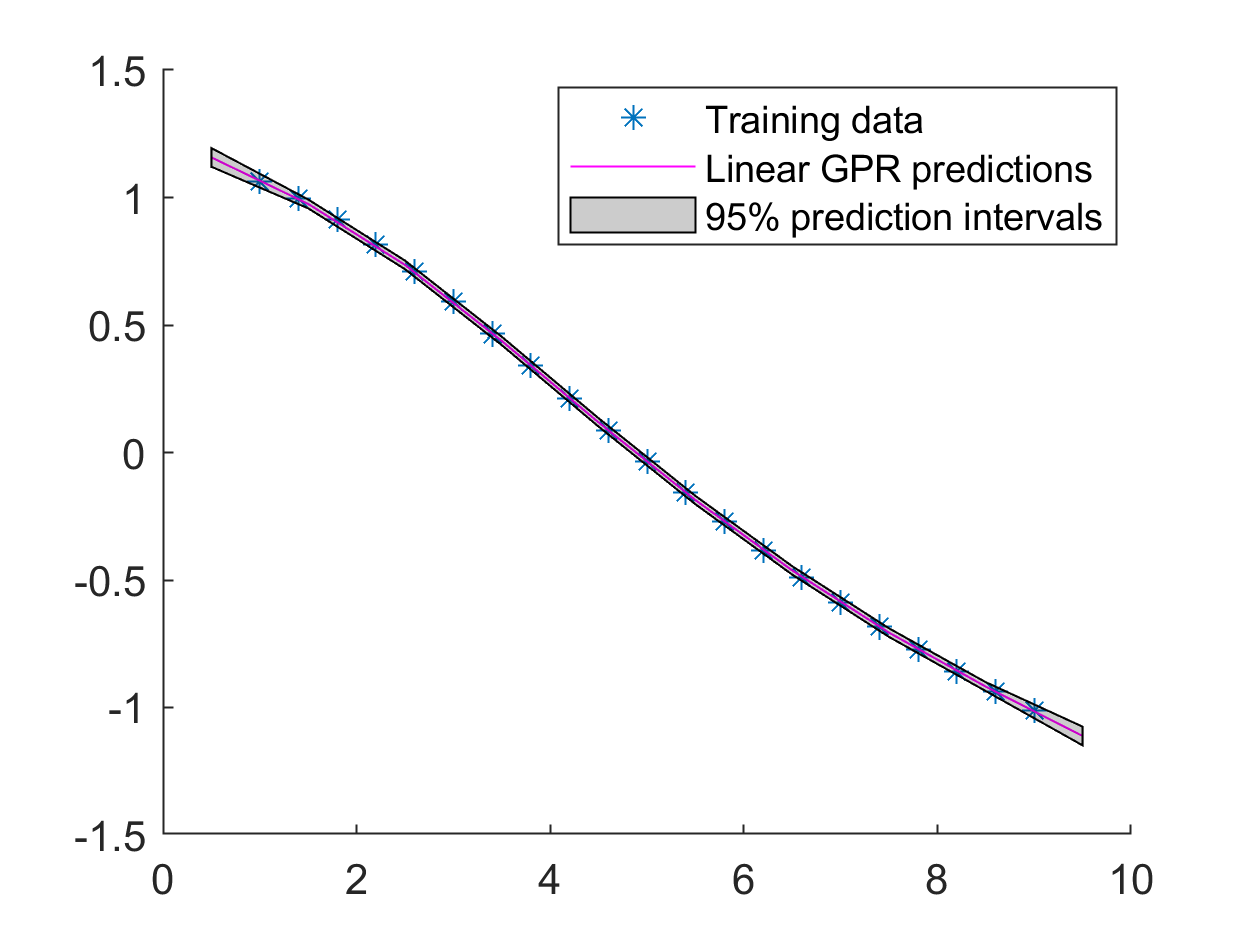}

\includegraphics[width=6.5cm]{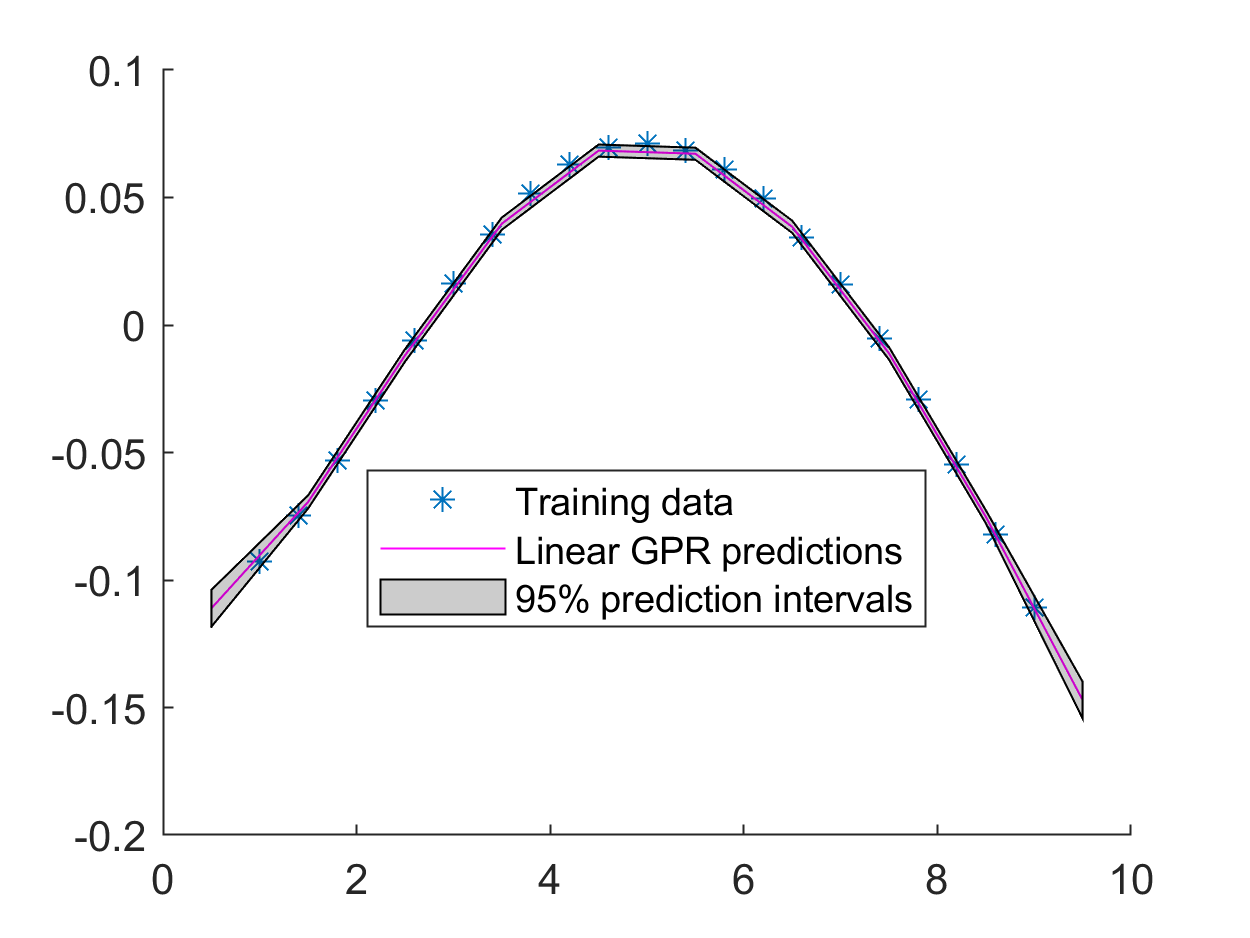}
\includegraphics[width=6.5cm]{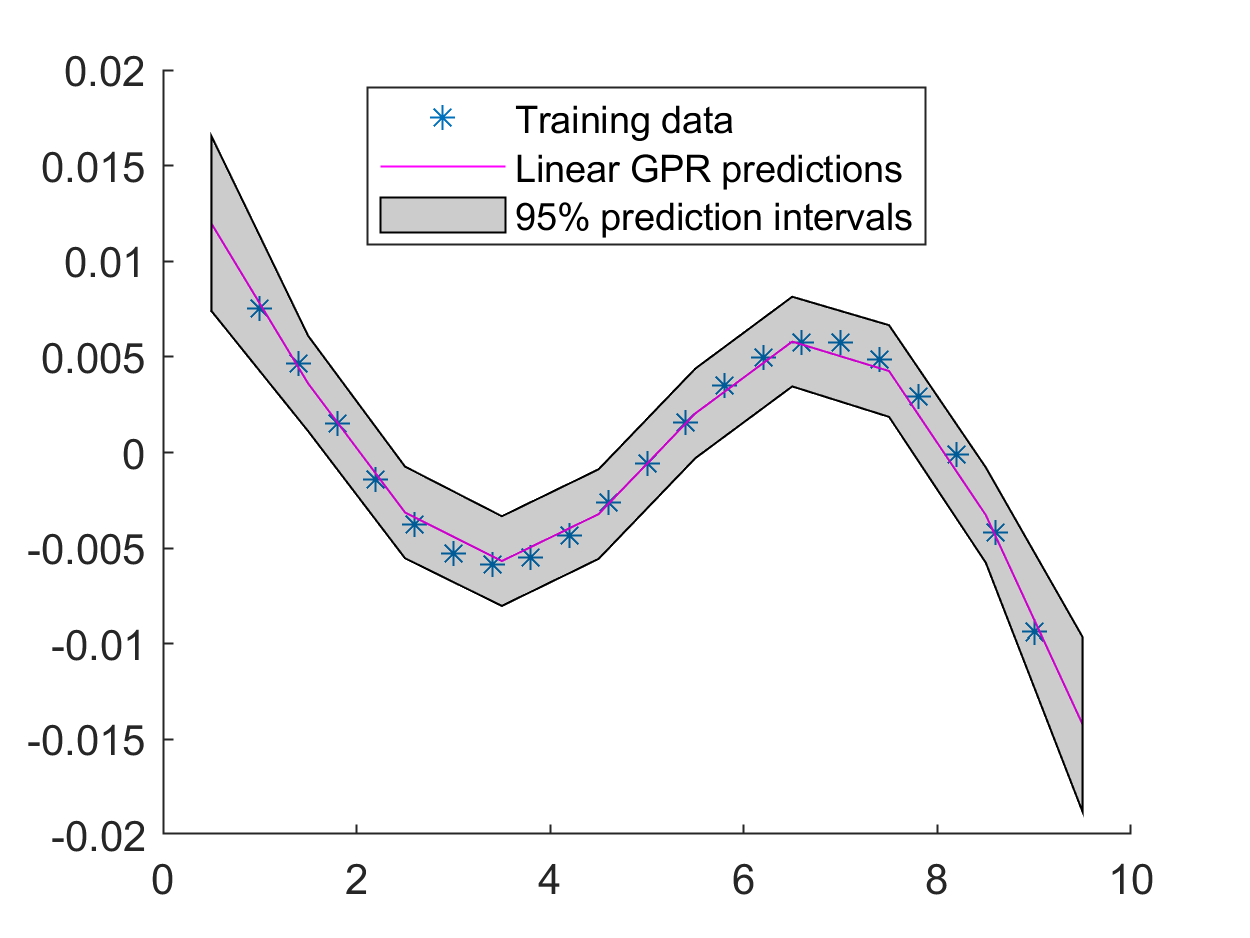}
\caption{The four GPR corresponding to the four coefficents of the first reduced eigenvector and their confidence interval of Problem~\eqref{meqn}, when the snapshots are calculated at $\mu=1:0.4:9$ and mesh size $h=0.05$.}
\label{fig2_nonlin}       % Give a unique label
\end{figure} 

In Figure~\ref{fig1_nonlin}, we show the first eigenvalue obtained using FEM and DD model when the mesh parameter is $h=0.05$ and $\delta\mu=0.2$ and $0.4$, respectively. Also, we indicate the confidence interval for the GPR corresponding to the first eigenvalue. The confidence interval is very narrow, that is the eigenvalues obtained by FEM and DD model match very well. That is the eigenvalues we get using the GPR-based DD model are very accurate.
Also, in Figure~\ref{fig2_nonlin}, we plot the mean function of the GPR corresponding to the four projected coefficients of the RB solutions with $95$ percent confidence of interval for the choice $h=0.05$; the sample points are uniformly distributed in the interval $[1,9]$ with stepsize $0.4$, that is the snapshot matrix contains all the first eigenvectors at these sample points.

\begin{figure}
\centering
\subcaptionbox{$\mu=1.5$}{
\includegraphics[width=4.4cm]{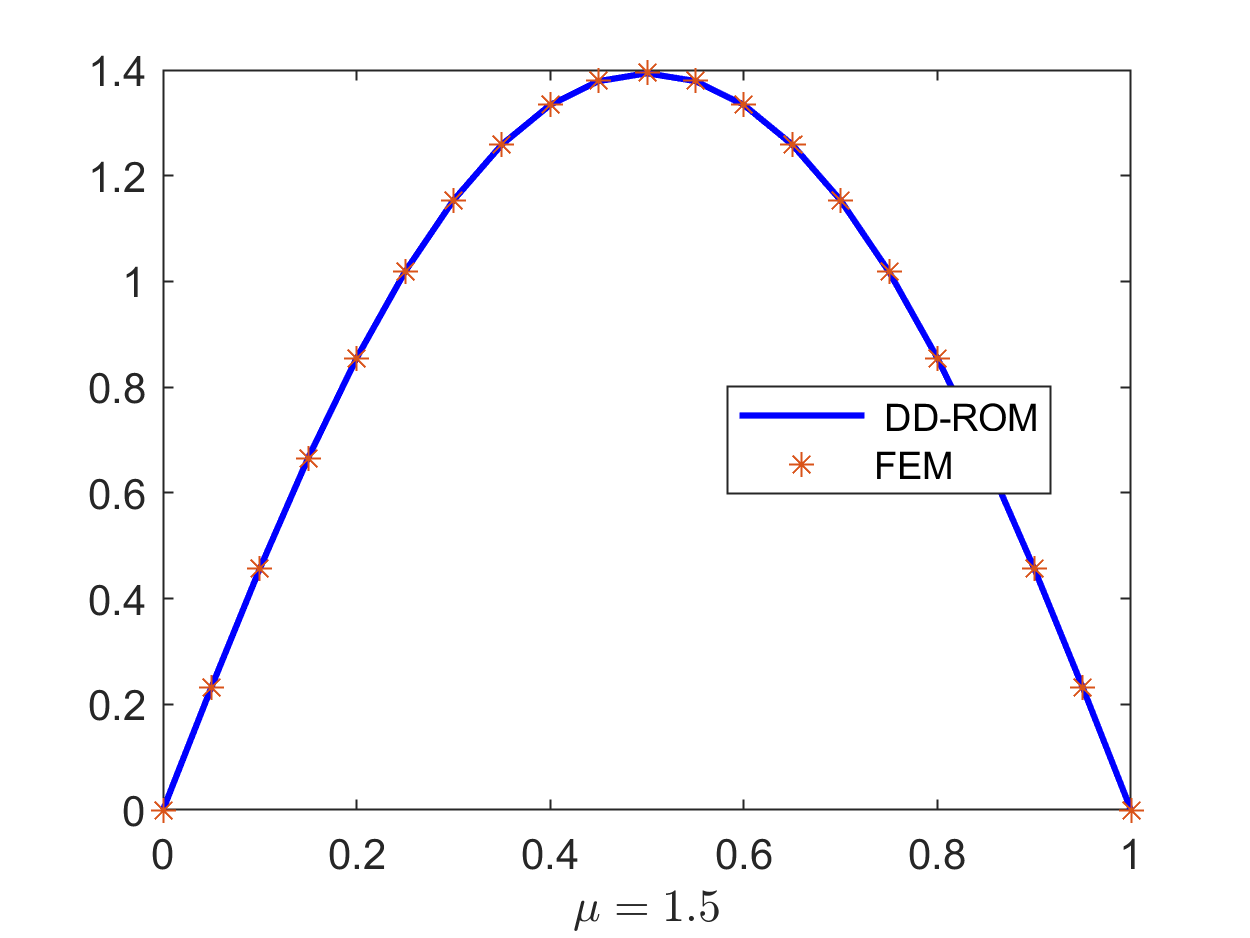}}
\subcaptionbox{$\mu=2.5$}{
\includegraphics[width=4.4cm]{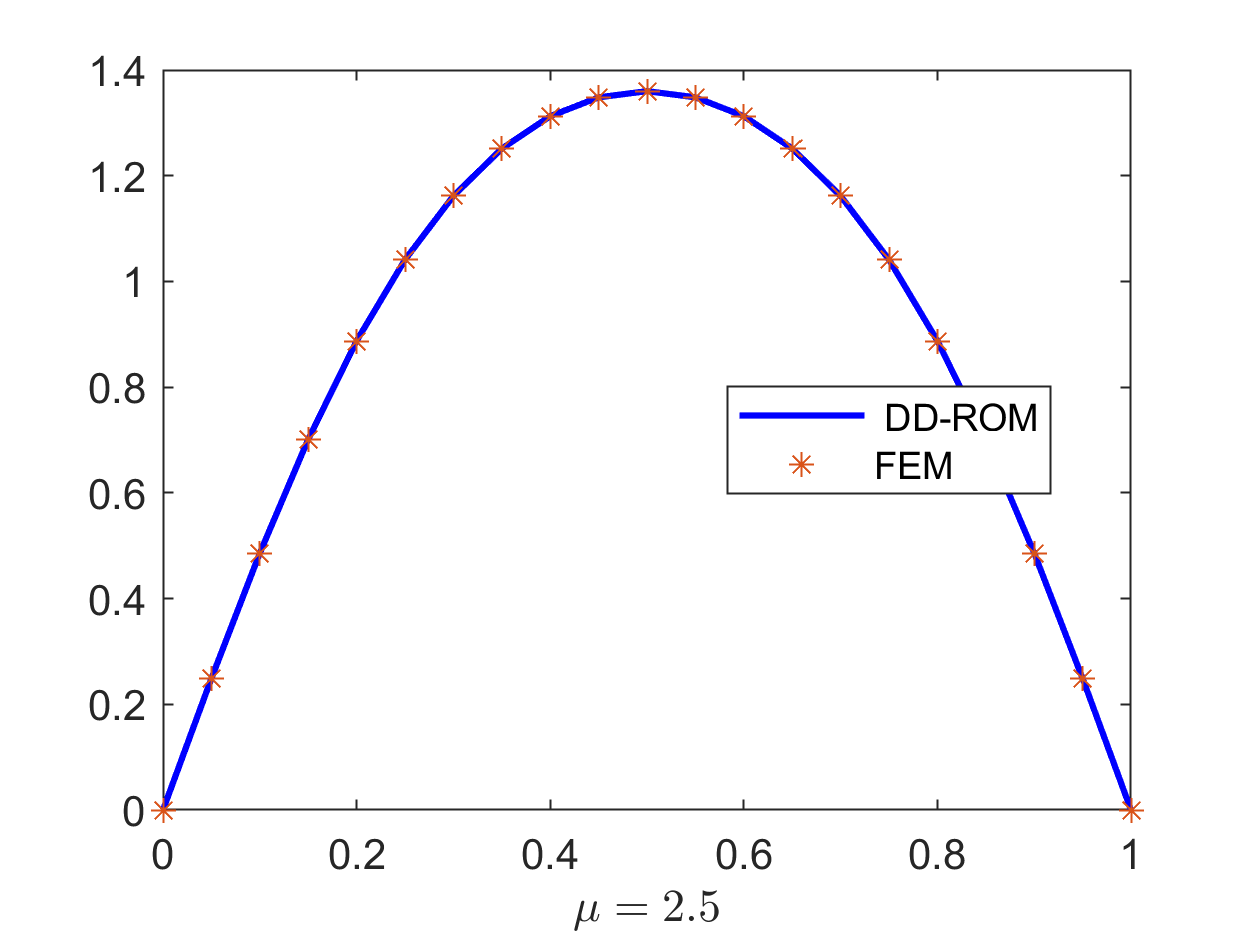}}
\subcaptionbox{$\mu=3.5$}{
\includegraphics[width=4.4cm]{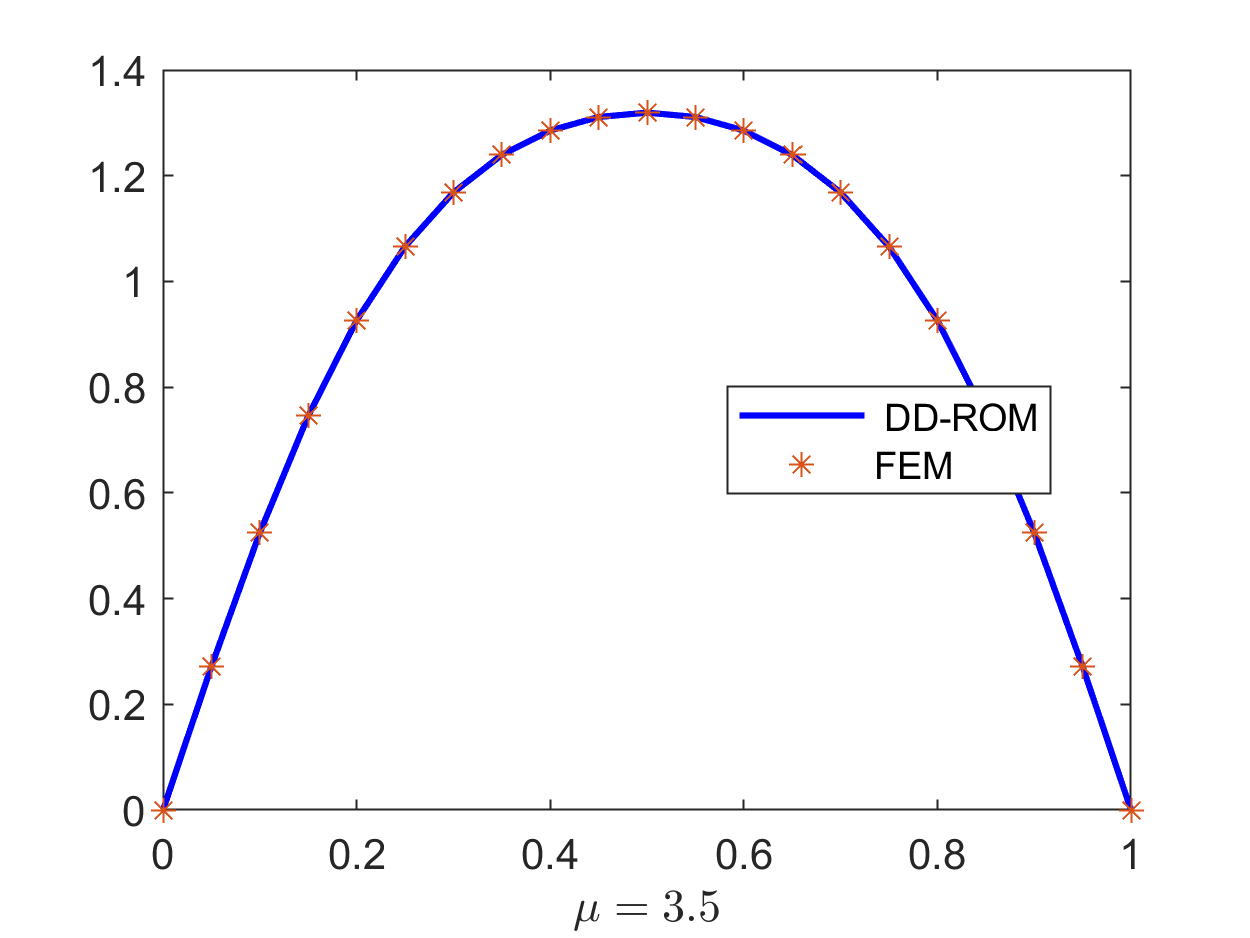}}

\subcaptionbox{$\mu=4.5$}{
\includegraphics[width=4.4cm]{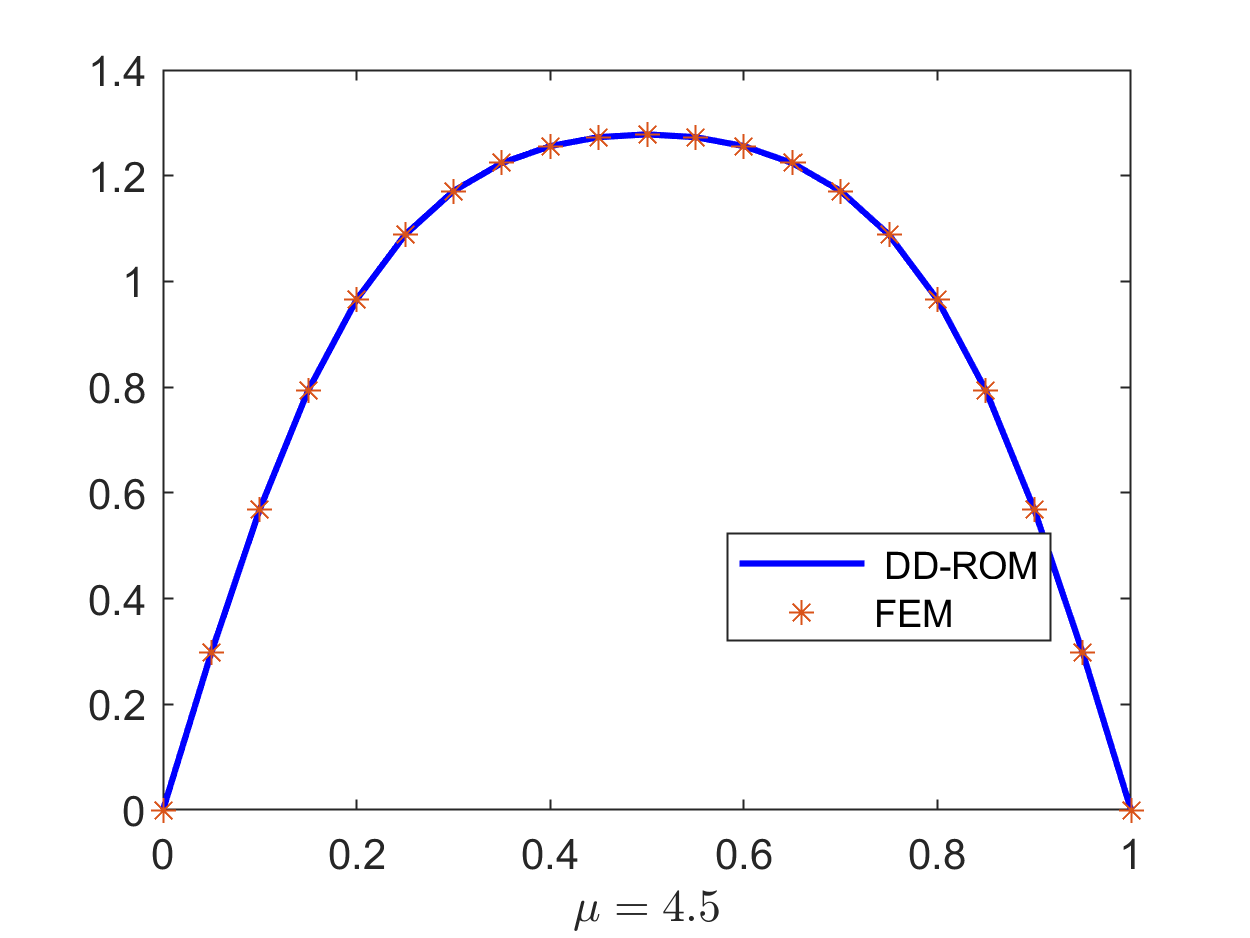}}
\subcaptionbox{$\mu=5.5$}{
\includegraphics[width=4.4cm]{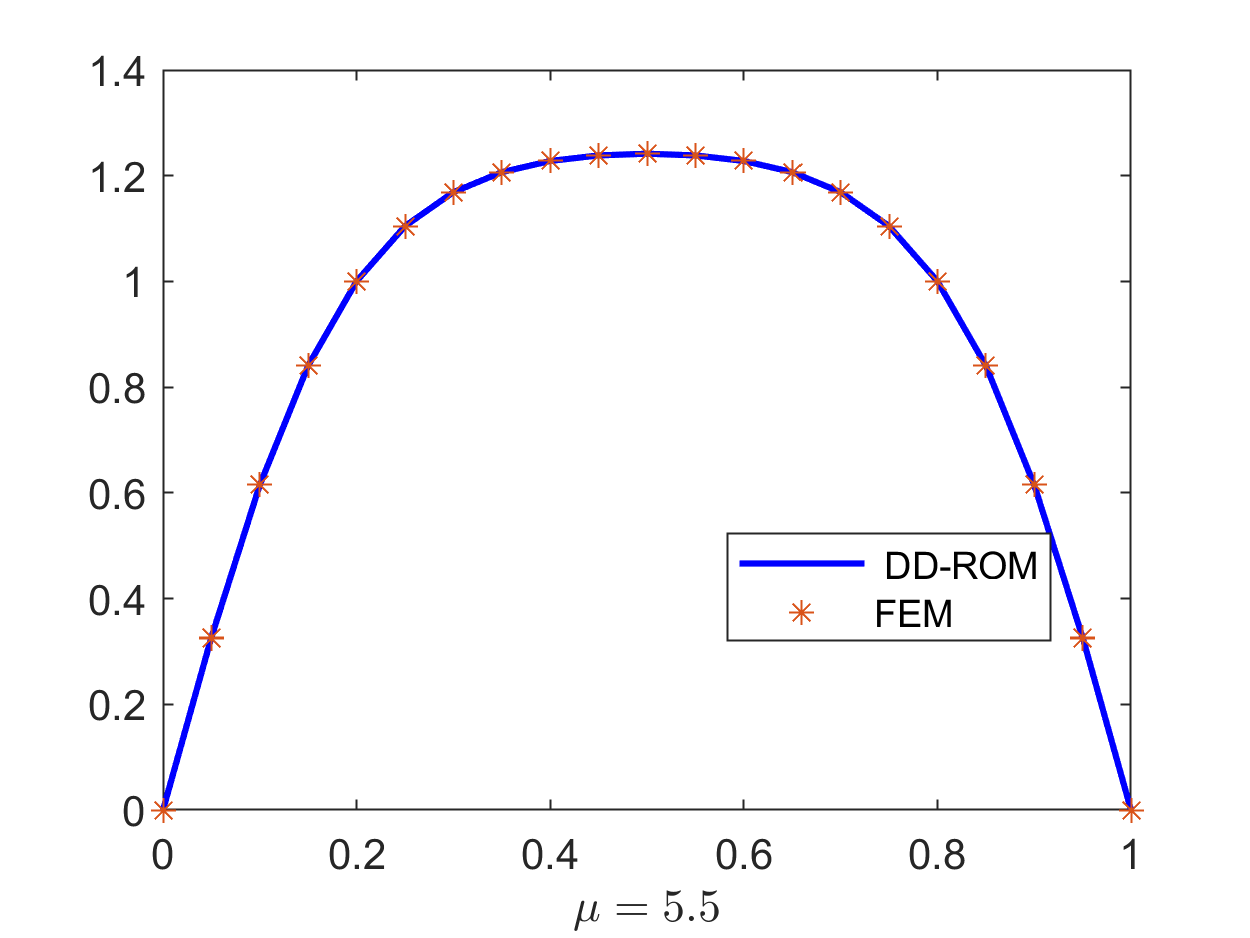}}
\subcaptionbox{$\mu=6.5$}{
\includegraphics[width=4.4cm]{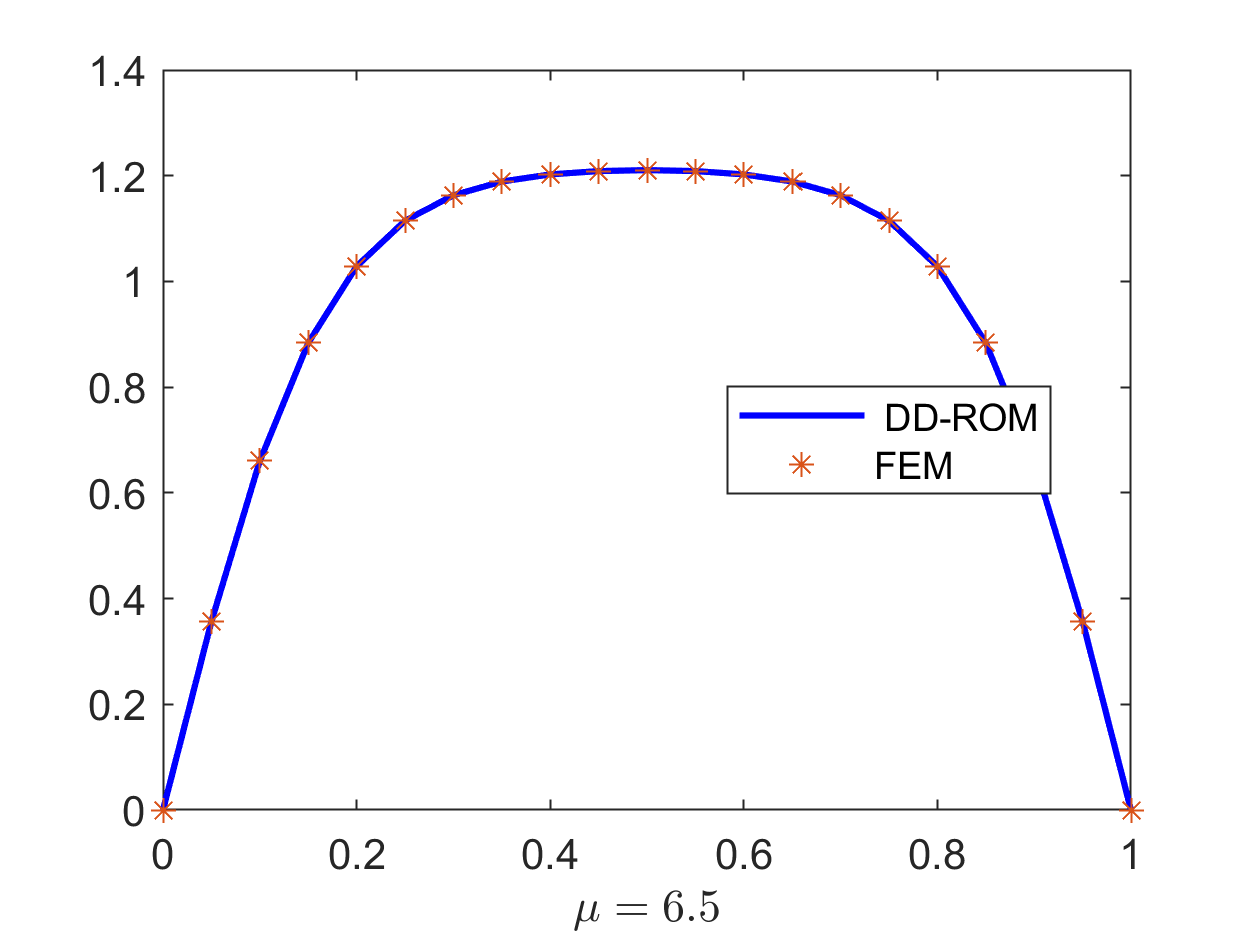}}
\caption{First eigenvectors of the nonlinear problem~\eqref{model} using DD model and FEM with $h=0.05$ and $\delta\mu=0.4$.}
\label{fig:nonlin1d}
\end{figure}
In Figure~\ref{fig:nonlin1d}, we show the first eigenvector obtained by our DD model at the test points in the interval $[1.5,6.5]$ with stepsize $1$; and also plot the corresponding eigenvector using FEM for the case $h=0.05$ and $\delta\mu=0.4$. From the Figure~\ref{fig:nonlin1d}, one can see that the eigenvectors using FEM and DD model are coinciding to each other at all the test points.

%%%%%%%%%%%%%%%%%%%%%%%%%%%%%%%%%%%%%%%%%%%%%%%%%%%%%%%%%%%%%%%%%%%%%%
\subsection{Non-affine parameter dependent eigenvalue problem}  
In this section, we consider eigenvalue problems that are not affine parameter dependent. For these problem, we cannot apply the standard projection based reduced method. On the other hand, since our model is based on the data, we can easily get the reduced solution using the trained GPR.
\subsubsection{Non-affine eigenvalue problem with single parameter}
Let us consider the following eigenvalue problem
\begin{equation}
\label{1par_nonaffn}
    \left\{
\aligned
&-\Delta u(\mu)=\lambda(\mu)e^{-\mu(x^2+y^2)}u(\mu)&&\text{in }\Omega=(0,1)^2\\
&u(\mu)=0&&\text{on }\partial\Omega,
\endaligned
\right.
\end{equation}
where the parameter $\mu$ belongs to the set $\mathcal{P}=[1,8]$. Note that the right hand side of the equation is non-affine. For the training our GPR model, we have calculated the eigensolutions of the EVP at sample points  with stepsize $0.4$ and $0.2$, respectively. We used as mesh size $h=0.05$ and $0.01$. Using the technique described in Section~\ref{sec:DD}, we construct one GPR for each projection coefficients of the RB solution and one GPR for the eigenvalue.
We report the eigenvalues obtained using our model at the test points $\mu=1.5$, $3.5$, $5.5$, $7.5$ in Table~\ref{table:nonafnpar1} for all the four cases corresponding to $h=0.05$, $0.01$ and $\delta\mu=0.4$, and $0.2$. We compare the eigenvalues with the FEM based eigenvalues. Note that the test points are different from the sample points. Also in this case the DD eigenvalues match pretty well with the FEM ones.

\begin{table}
\footnotesize
 	 	\centering
 	\begin{tabular}{|c|c|c|c|c|c|c|c|} 
 		\hline
 	h&  {\begin{tabular}[c]{@{}c@{}} $\delta \mu$ \end{tabular}} &
		  {\begin{tabular}[c]{@{}c@{}}Method \end{tabular}} &
 		 {\begin{tabular}[c]{@{}c@{}} $\mu=1.5$ \end{tabular}} & 
 		{\begin{tabular}[c]{@{}c@{}} $\mu=3.5$   \end{tabular}}&
 		{\begin{tabular}[c]{@{}c@{}}$\mu=5.5$   \end{tabular}}&
 		{\begin{tabular}[c]{@{}c@{}} $\mu=7.5$   \end{tabular}}\\
  \hline
0.05 &0.4 & FEM&39.73473805	&72.62248873&108.48892187&145.54749238\\
     & & DD&39.81667386	&	72.57625001&	108.50118905&	145.49823357\\	
\hline
0.05 &0.2 & FEM&39.73473805	&72.62248873&108.48892187&145.54749238\\
     & & DD&39.79566310	&72.57602363&	108.50603170&	145.54219385\\	
\hline
0.01 &0.4 & FEM&39.62169633&72.32111463&107.85749755&	144.42856653\\
     & & DD&39.70506978	&	72.27383282	&107.87538316&	144.37709615\\	
\hline
0.01 &0.2 & FEM&39.62169633&72.32111463&107.85749755&	144.42856653\\
     & & DD&39.68283223&72.27478421	&107.87527246&	144.42228776\\
\hline
 	\end{tabular}
	\caption{Comparison of 1st eigenvalues of the EVP \eqref{1par_nonaffn} using FEM and DD Model.}
	 	\label{table:nonafnpar1}
\end{table}

\begin{figure}
\centering
\subcaptionbox{$\mu_{tr}=1:0.4:8$}{
\includegraphics[width=6.5cm]{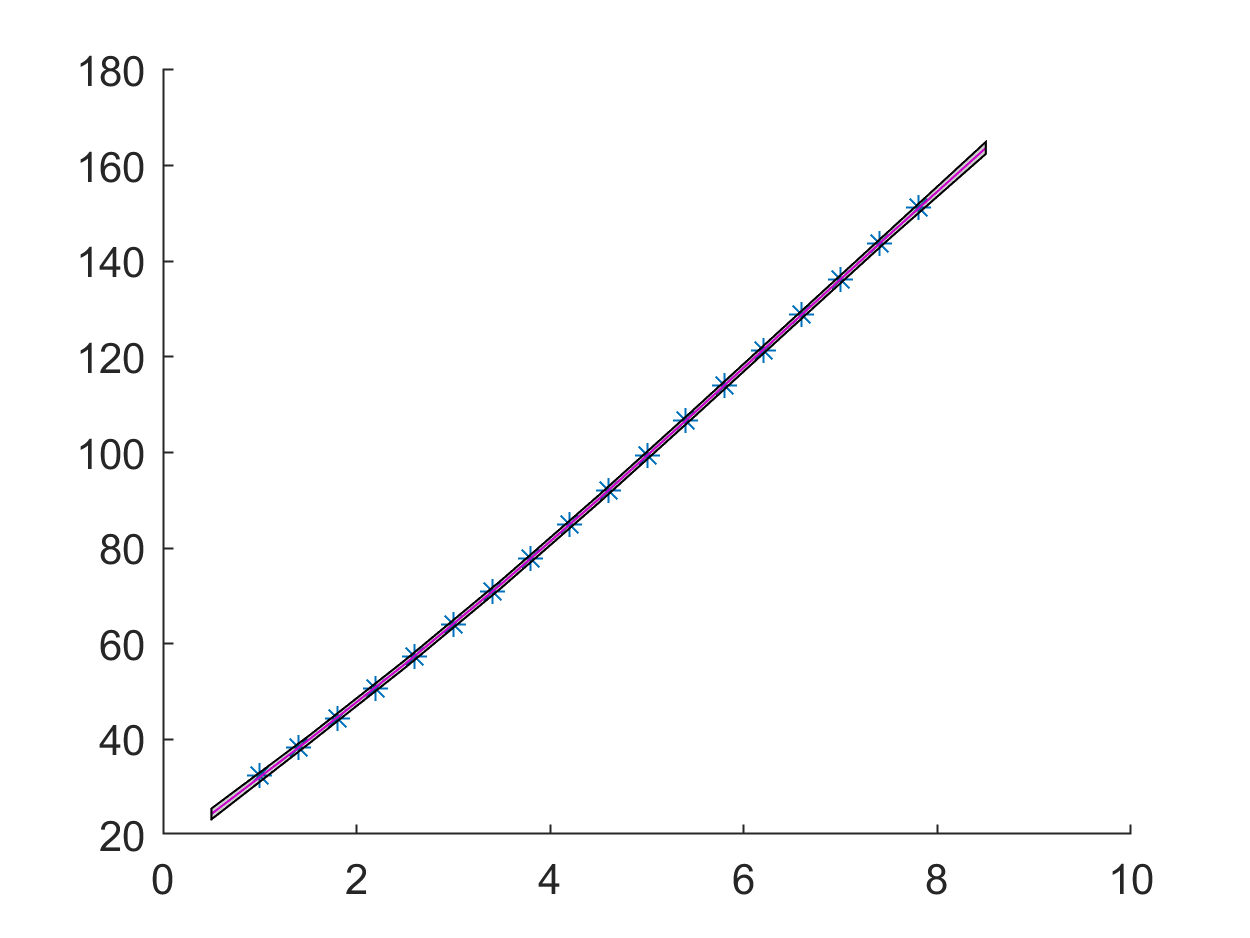}}
\subcaptionbox{$\mu_{tr}=1:0.2:8$}{
\includegraphics[width=6.5cm]{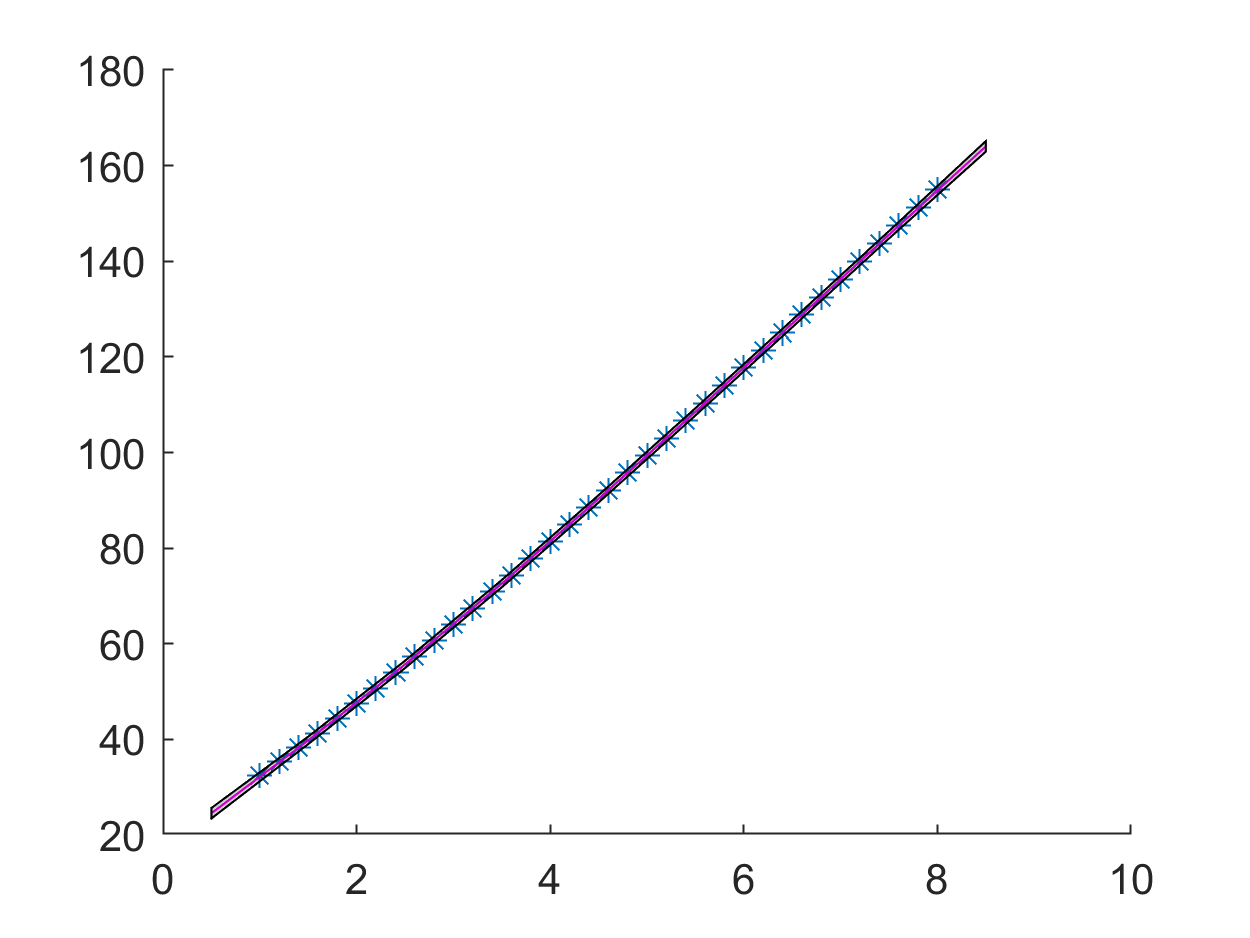}}
\caption{The GPR for the first eigenvalues of the nonlinear problem \eqref{1par_nonaffn} with different number of snapshots with $h=0.05$.}
\label{fig1_nonafnpar1}       % Give a unique label
\end{figure}

In Figure~\ref{fig1_nonafnpar1} we show the first eigenvalue using FEM and DD model with mesh parameter $h=0.05$ and $\delta\mu=0.2$, $0.4$, respectively. We show the confidence interval for the GPR of the eigenvalues. The confidence interval is very narrow, that is the eigenvalues obtained by FEM and data-driven model match pretty well. That is the eigenvalues we get using the GPR-based DD model is very accurate.

\begin{figure}
\centering
\includegraphics[width=4.5cm]{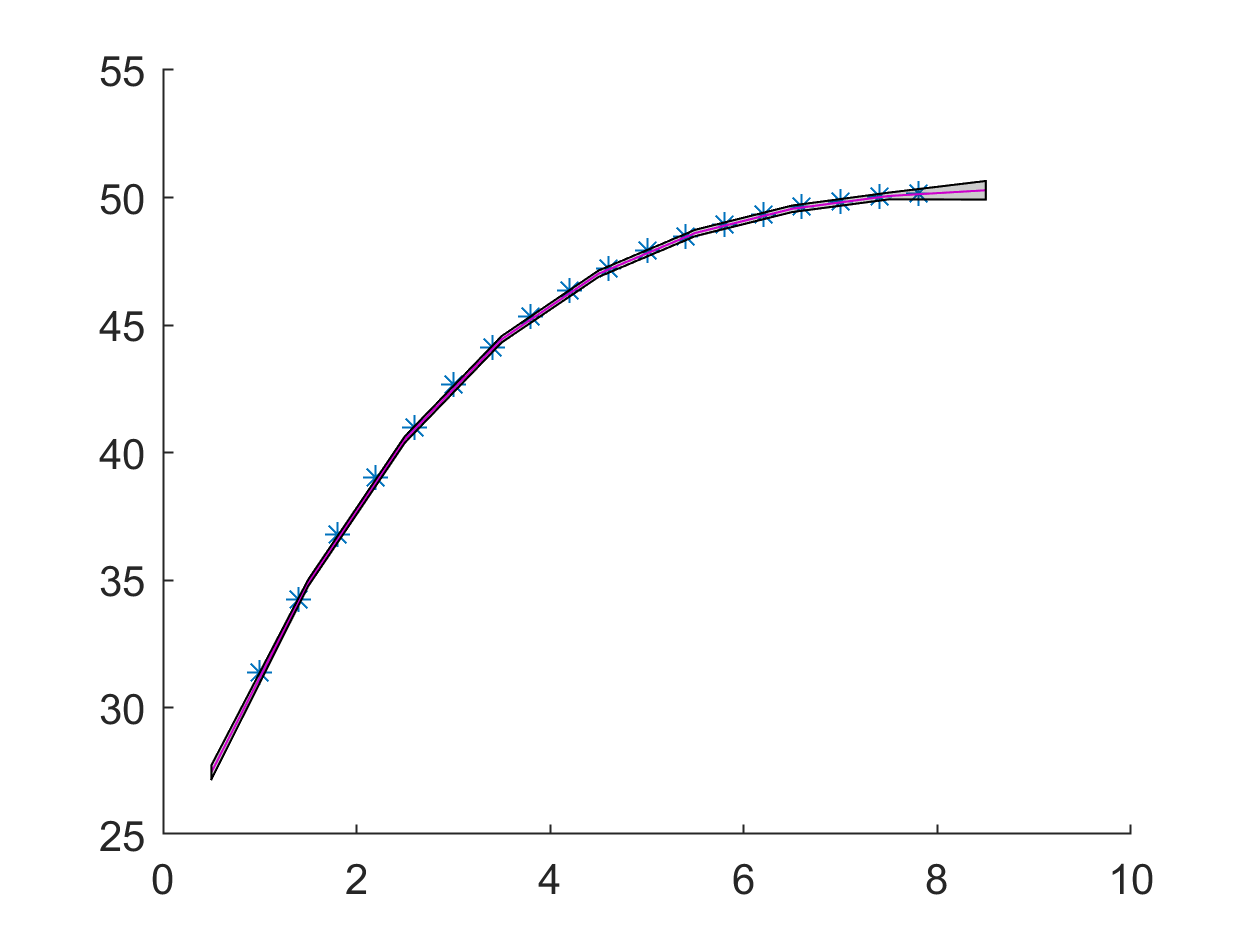}
\includegraphics[width=4.5cm]{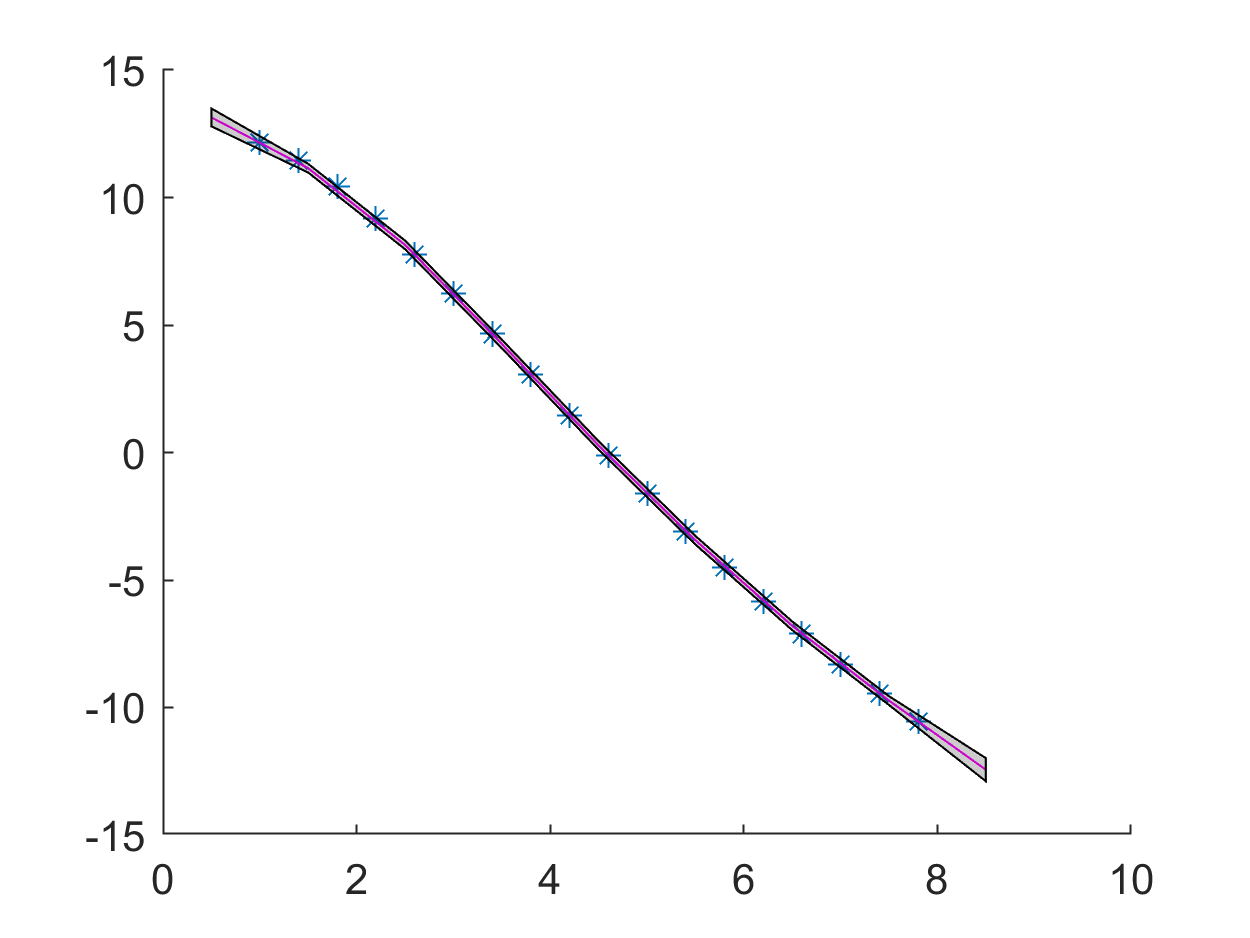}
\includegraphics[width=4.5cm]{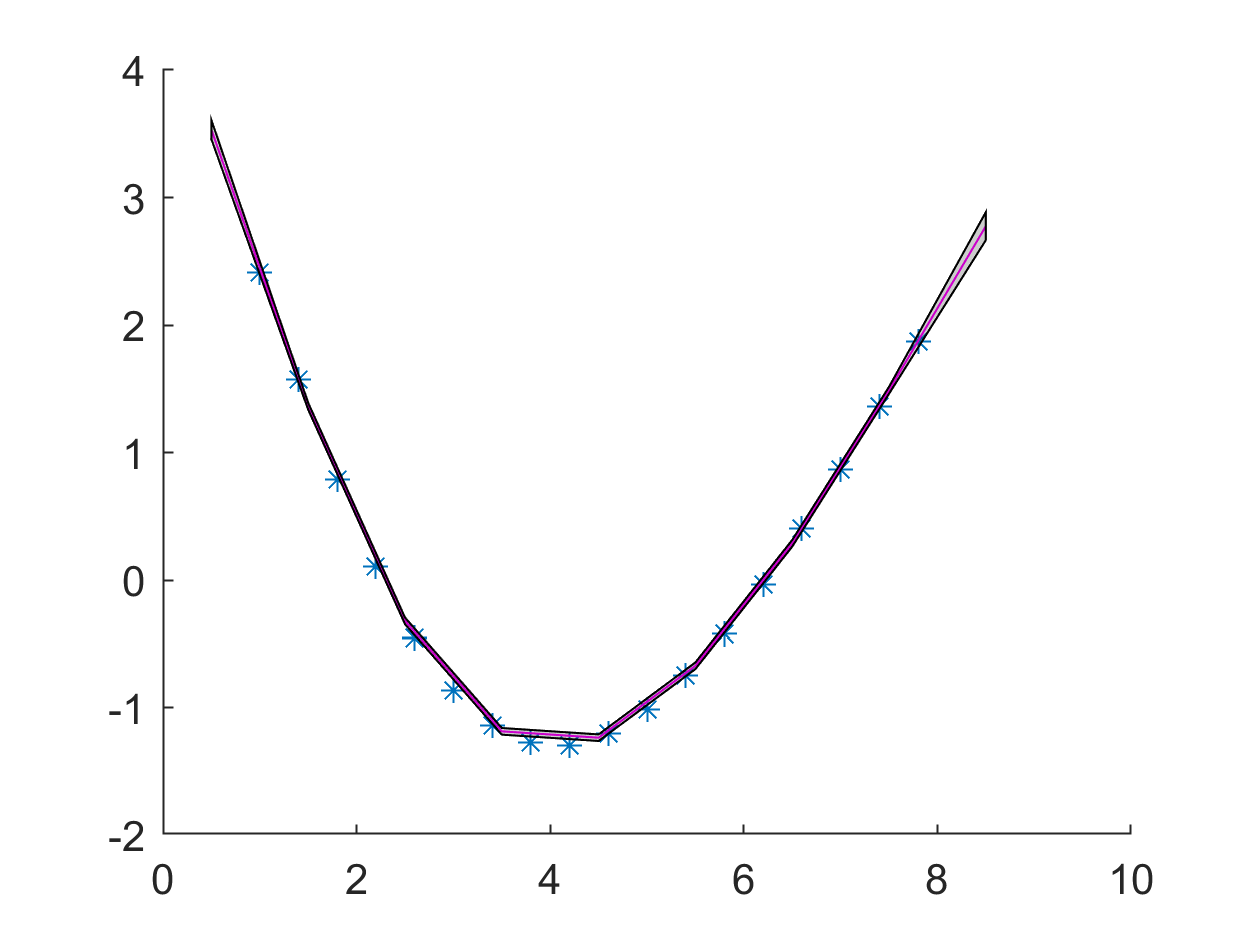}

\includegraphics[width=4.5cm]{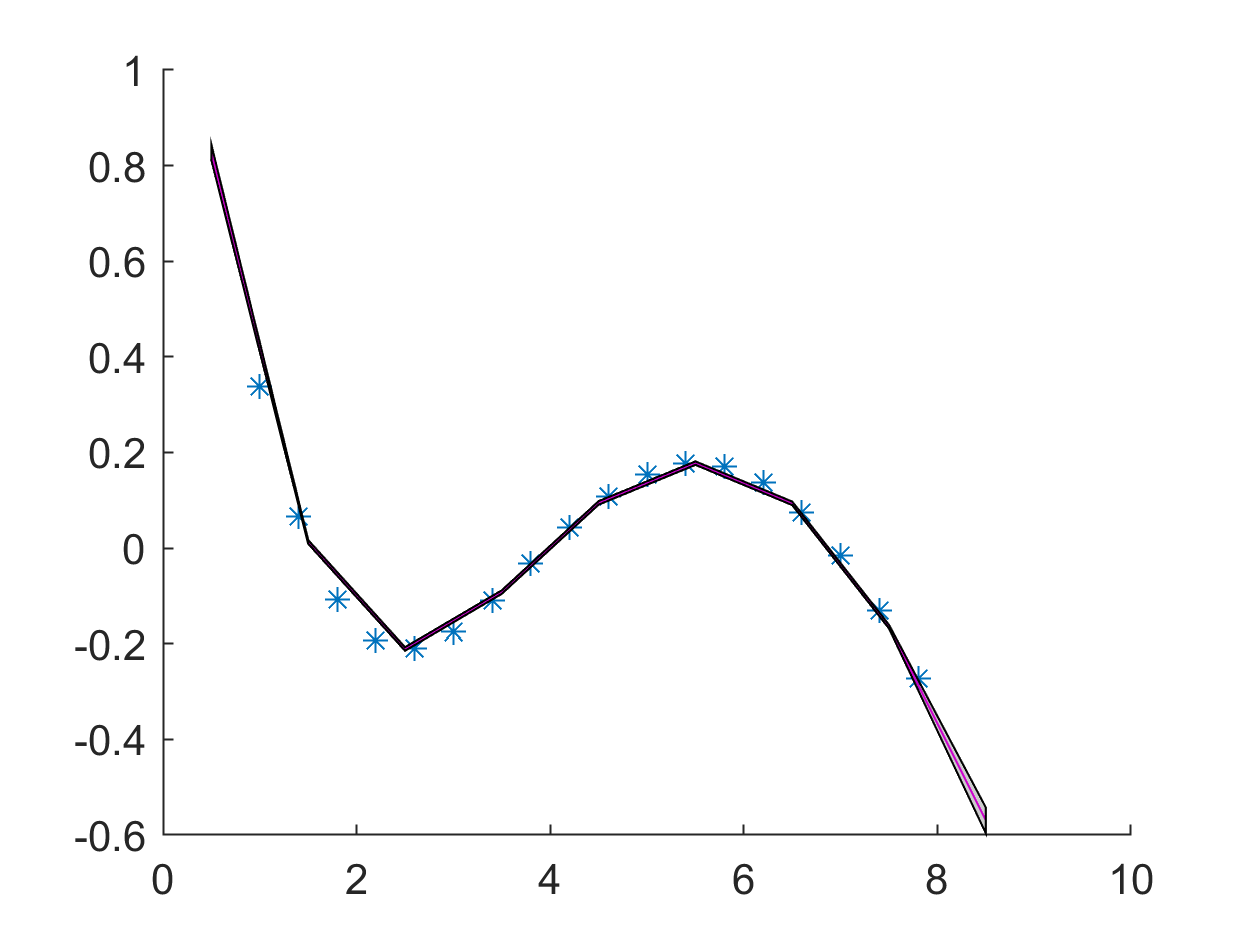}
\includegraphics[width=4.5cm]{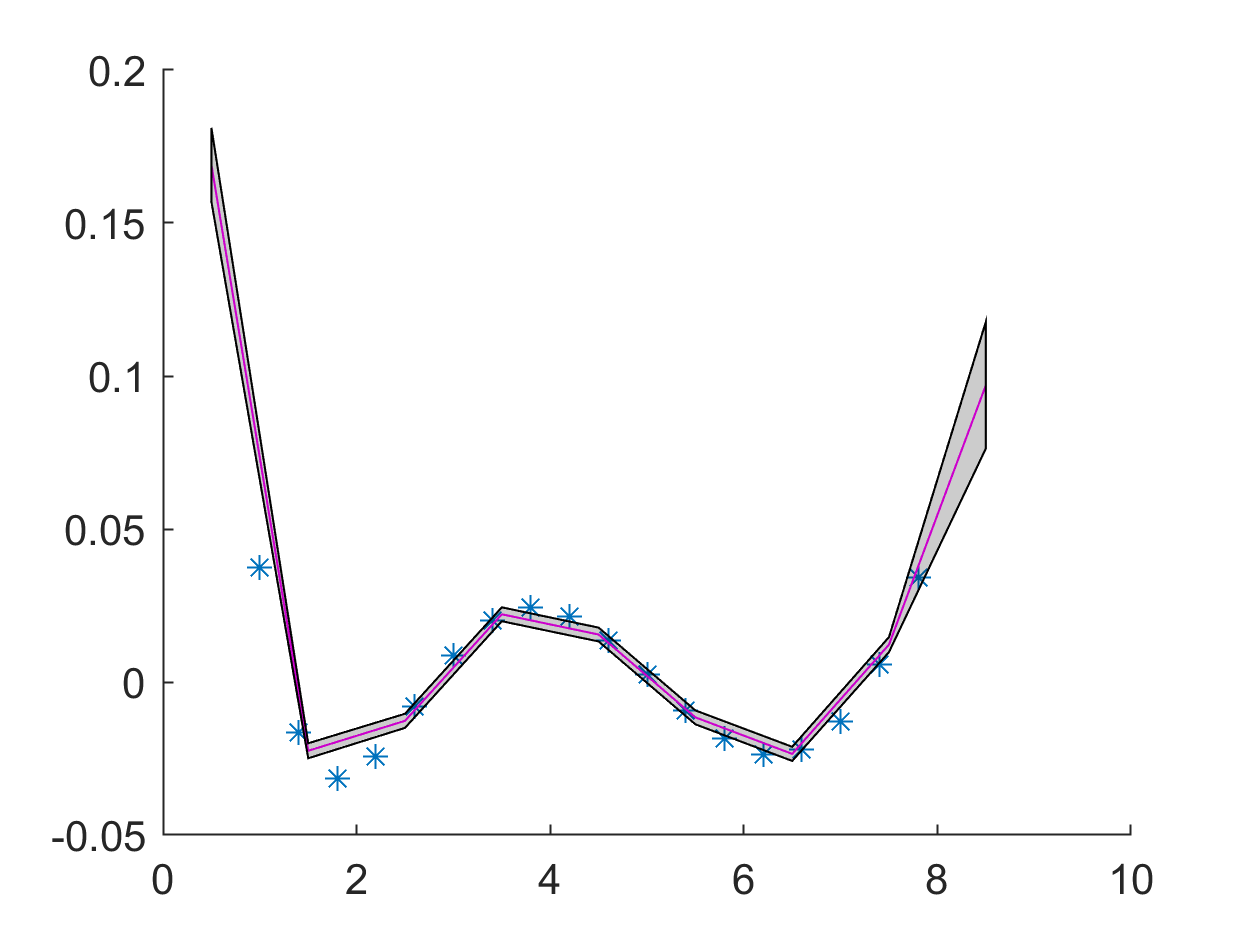}
\caption{The GPR plot for the first five coefficients of the eigenfunction for Problem~\eqref{1par_nonaffn} with $h=0.05$, $\delta \mu=0.4$, and test points $\mu=0.5:1:8.5$.}
\label{fig2:nonafnpar1}
\end{figure}

In Figure~\ref{fig2:nonafnpar1} we plot all the five coefficients of the RB solutions with $95$ percent confidence of interval for the choice $h=0.05$ when the sample points are taken in the interval $[1,8]$ with stepsize $0.4$, that is, the snapshot matrix contains all the first eigenvectors at these sample points. The test points are taken in the interval $[0.5,8.5]$ with stepsize $1$, that is, the end points are outside of the parameter space. we notice that the confidence interval at the end points is a little larger.

\begin{figure}
\centering
\includegraphics[width=4.5cm]{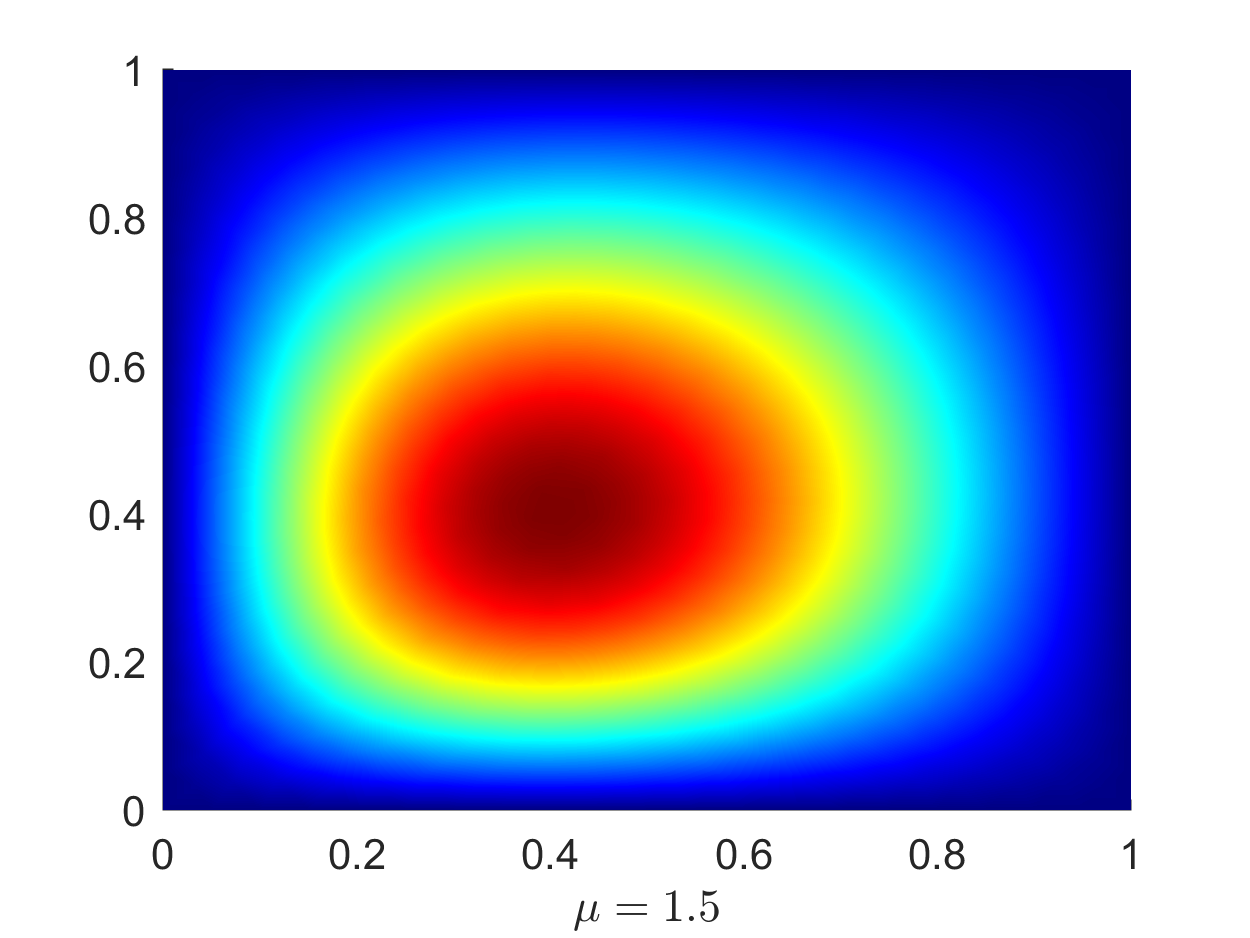}
\includegraphics[width=4.5cm]{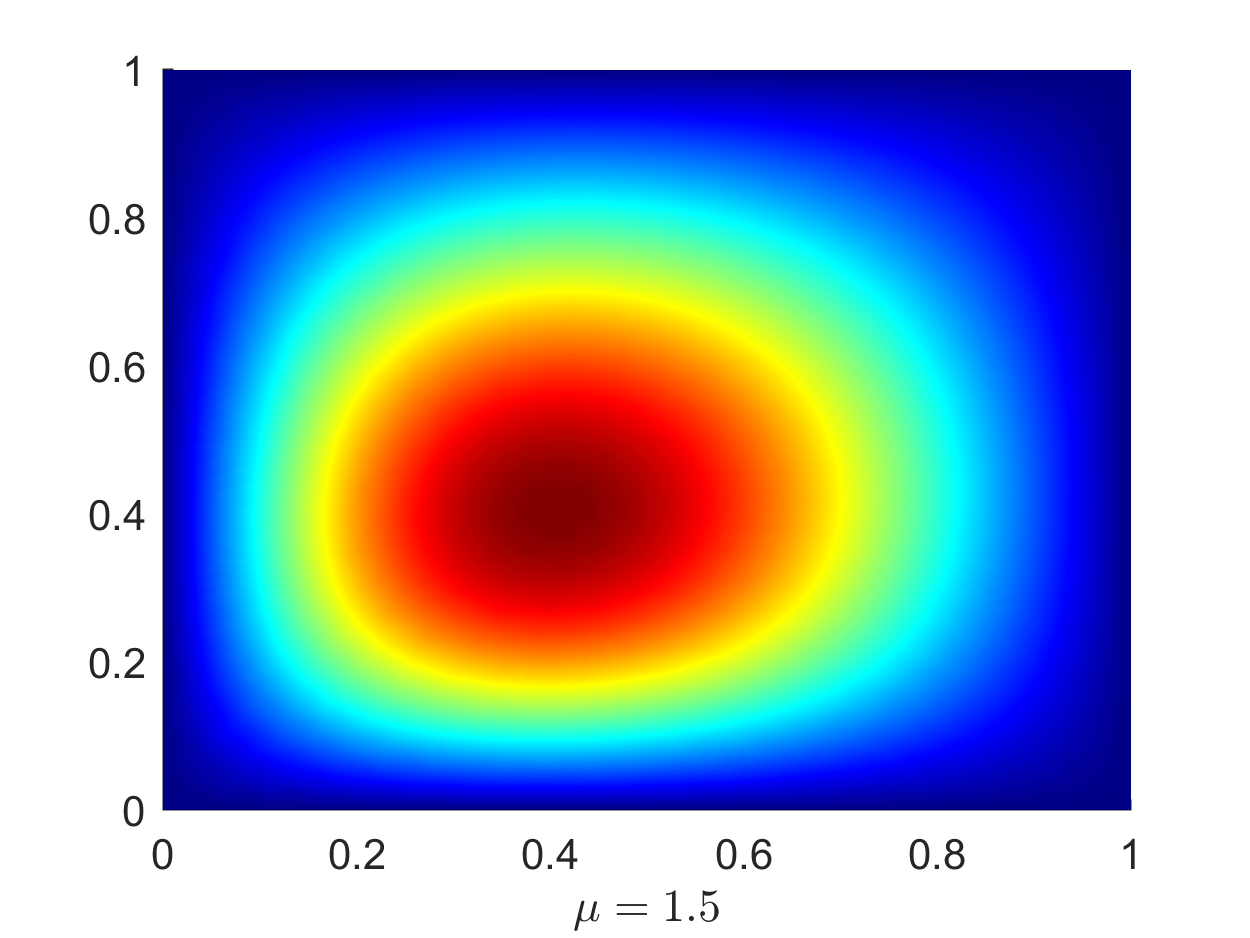}
\includegraphics[width=4.5cm]{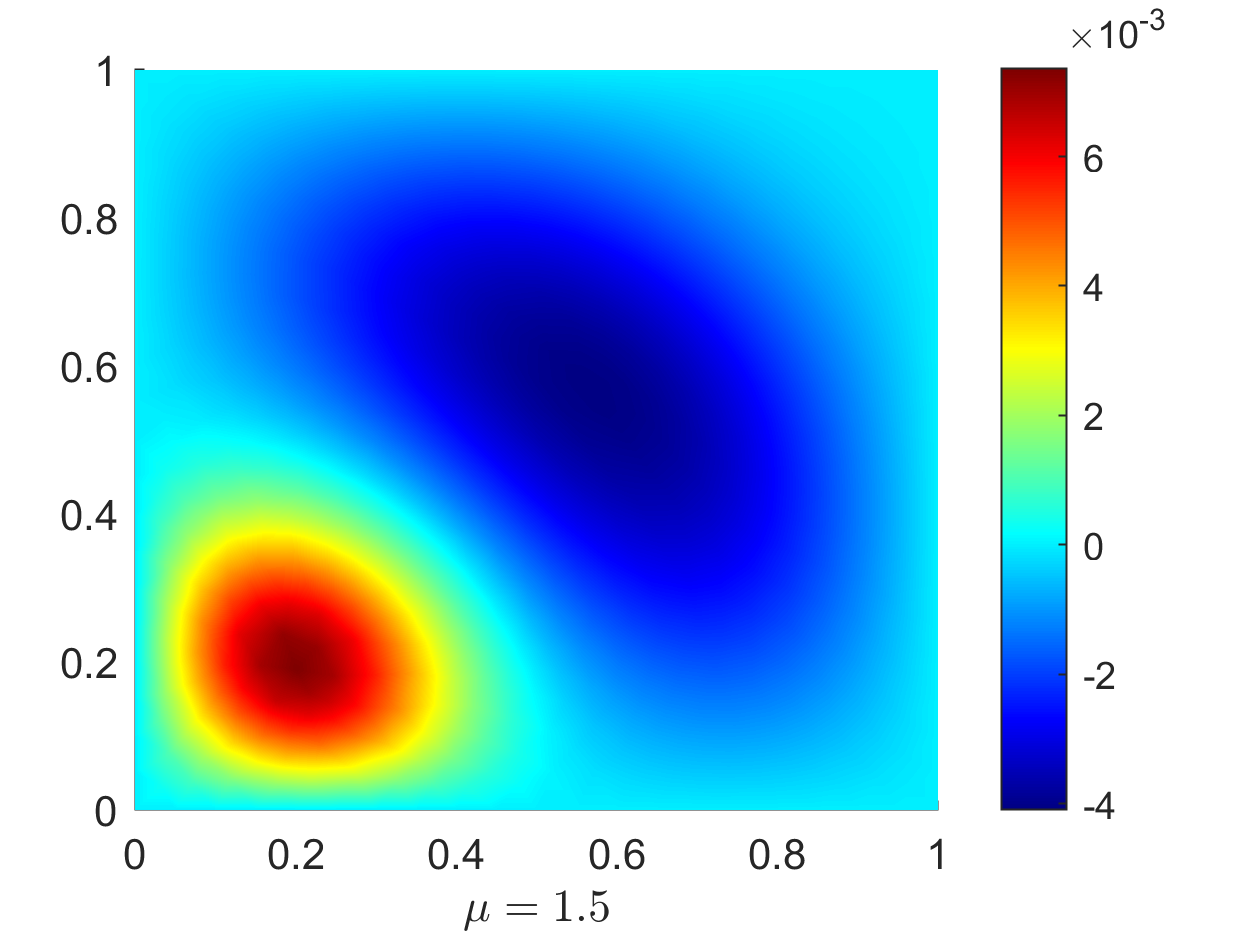}

\includegraphics[width=4.5cm]{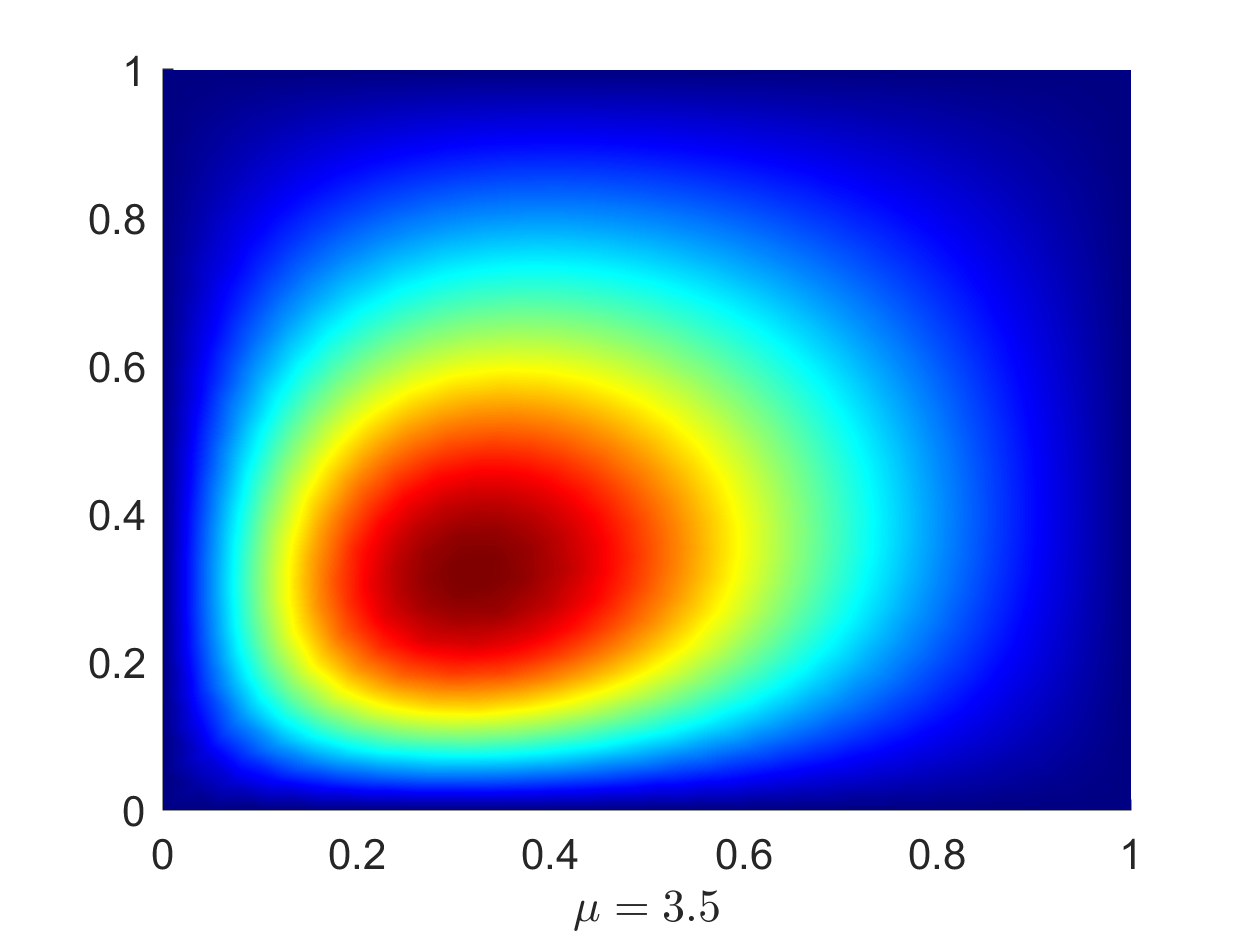}
\includegraphics[width=4.5cm]{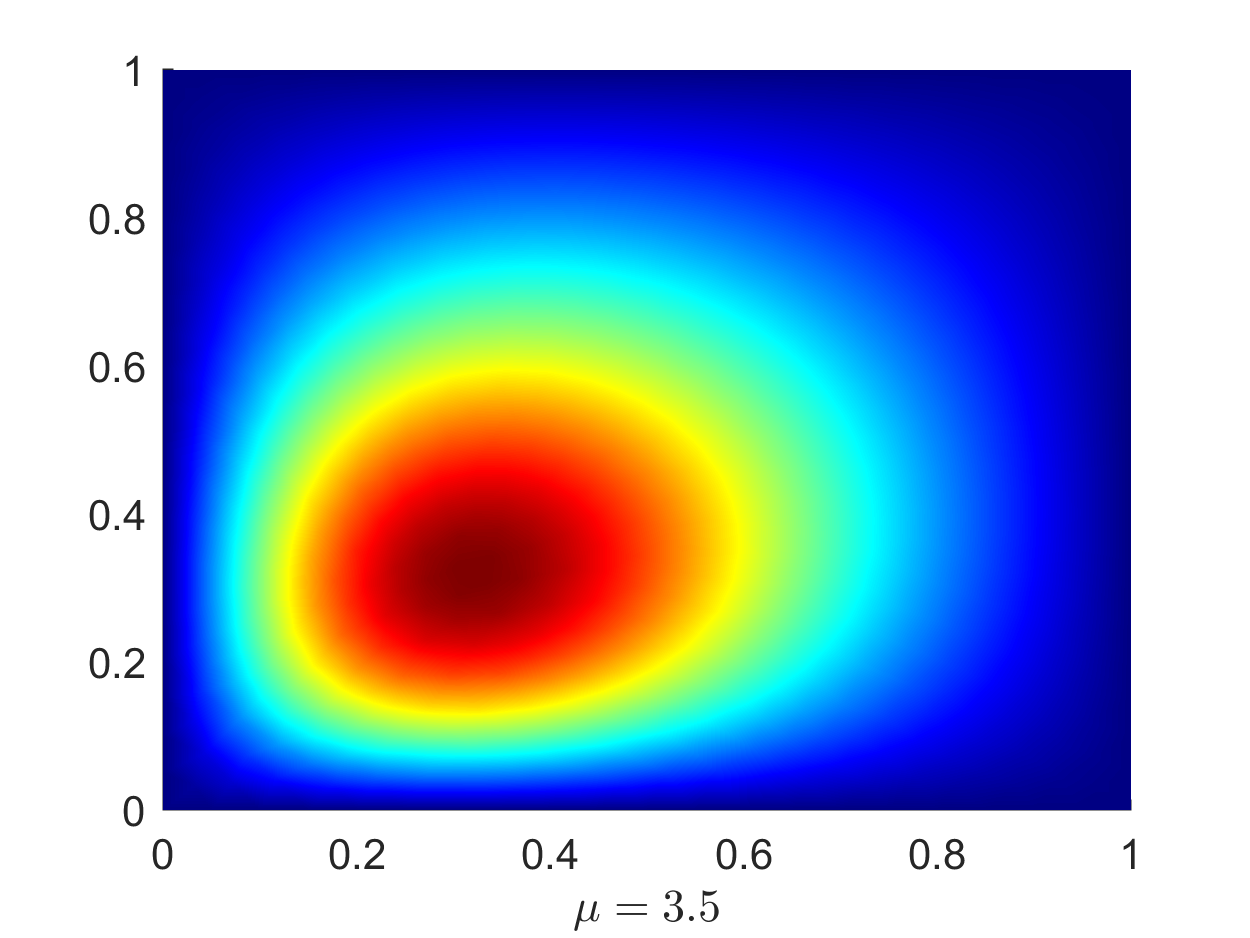}
\includegraphics[width=4.5cm]{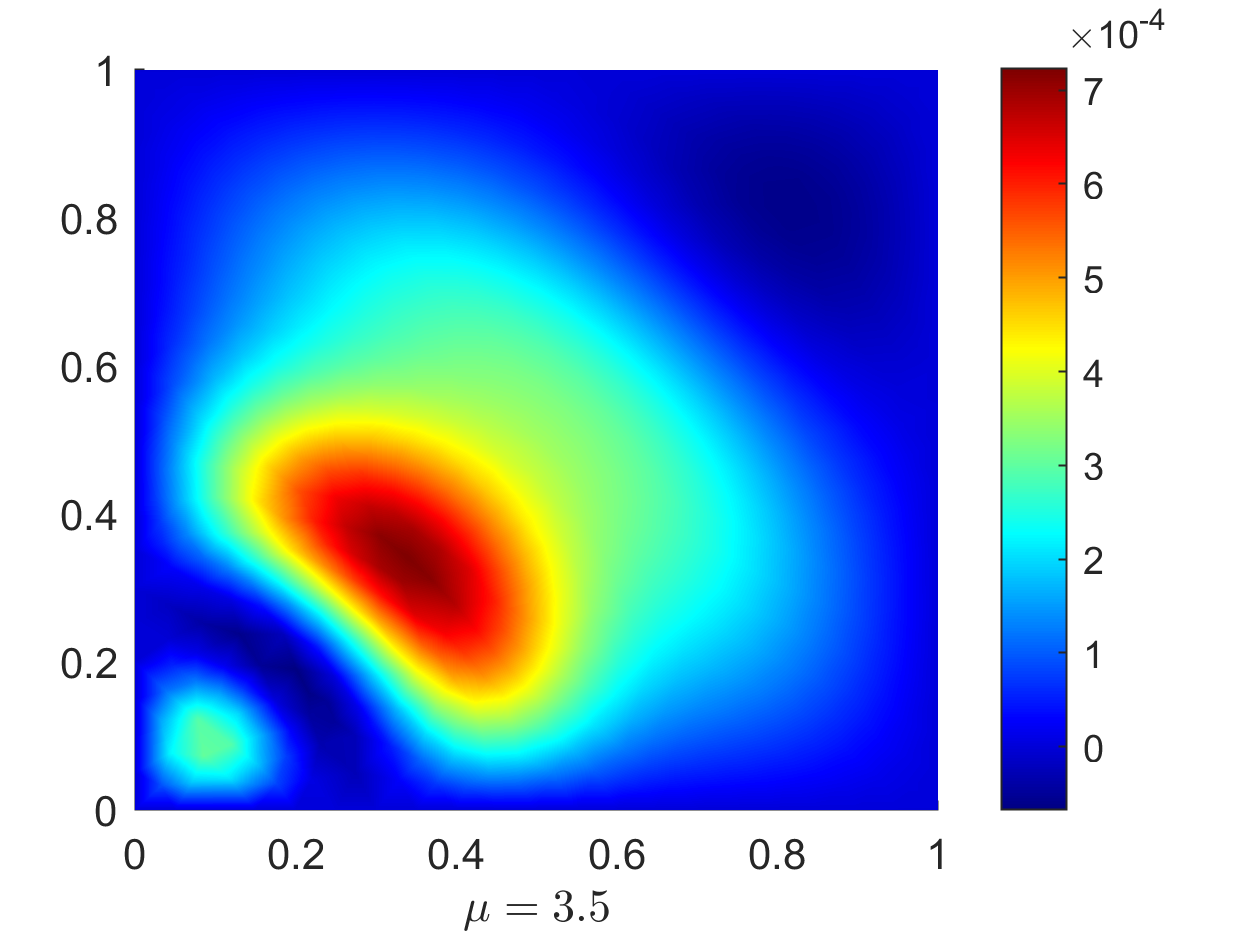}

\includegraphics[width=4.5cm]{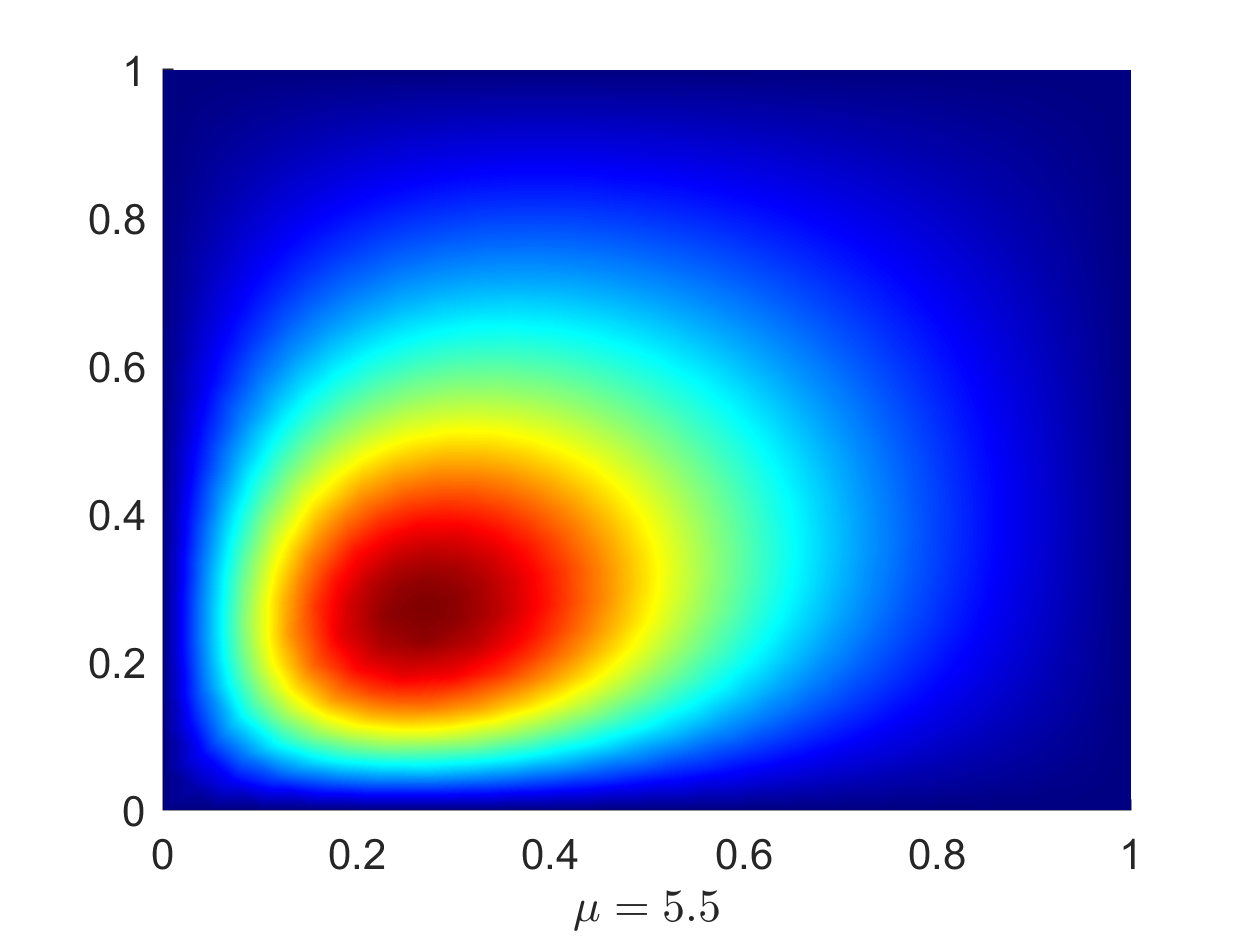}
\includegraphics[width=4.5cm]{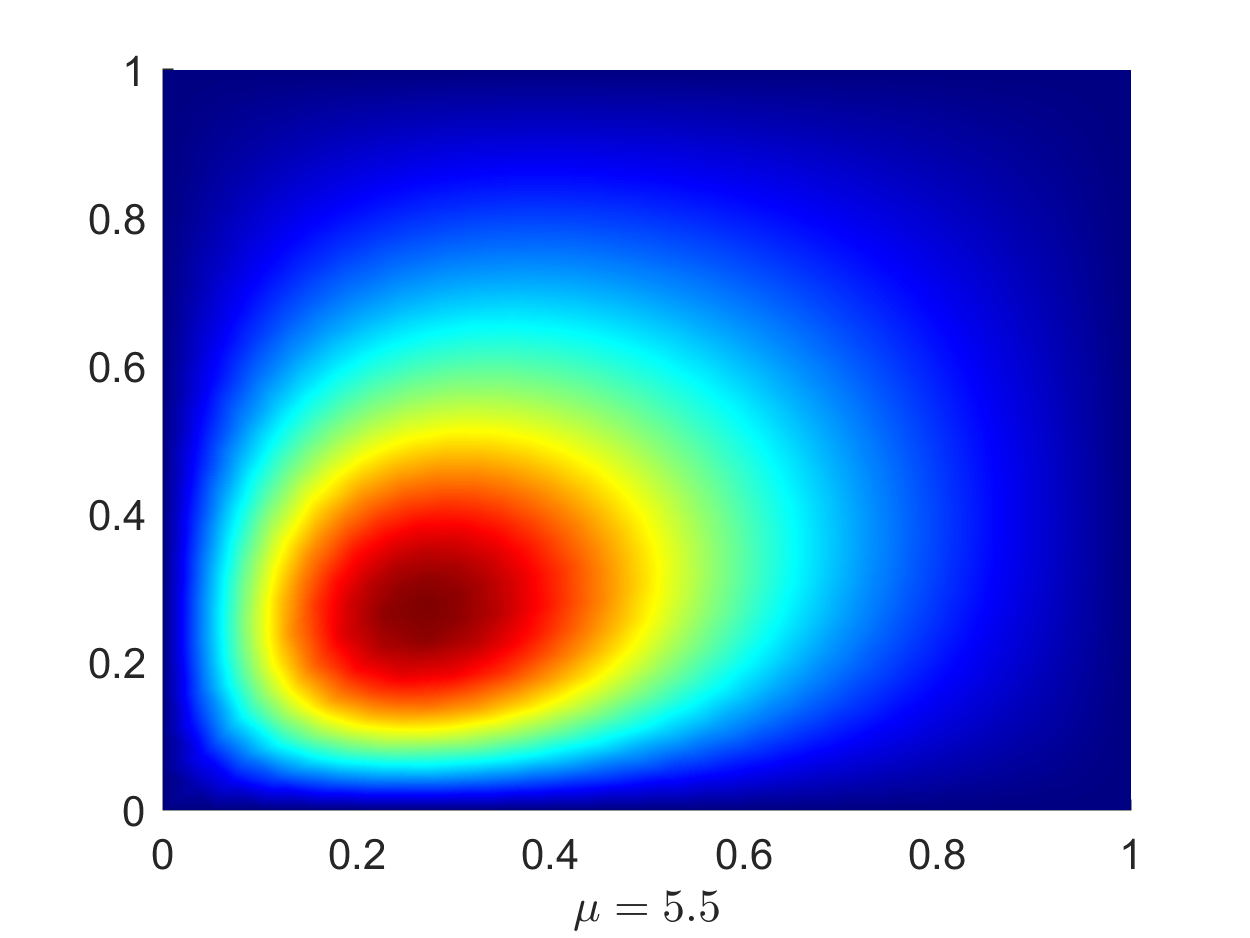}
\includegraphics[width=4.5cm]{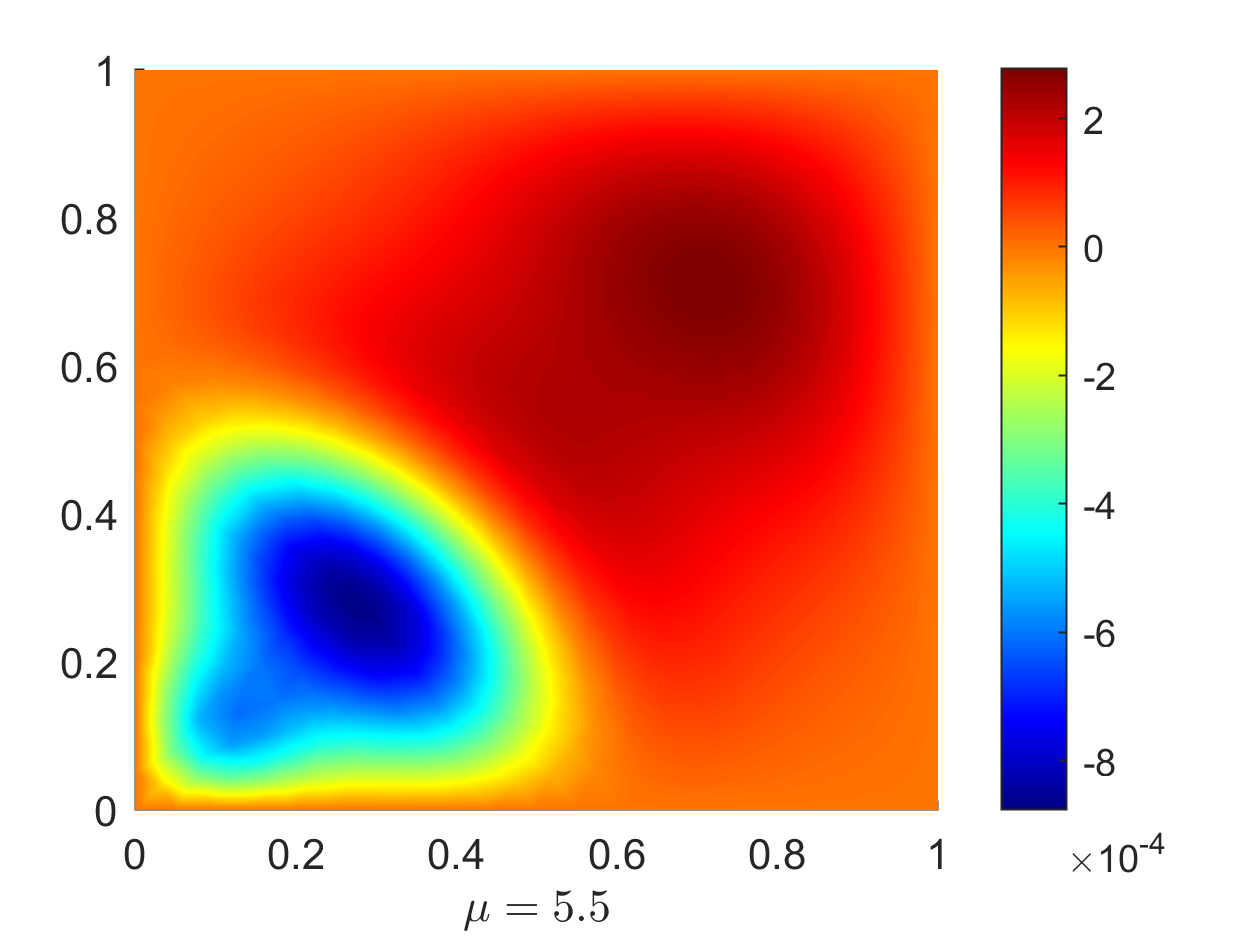}

\includegraphics[width=4.5cm]{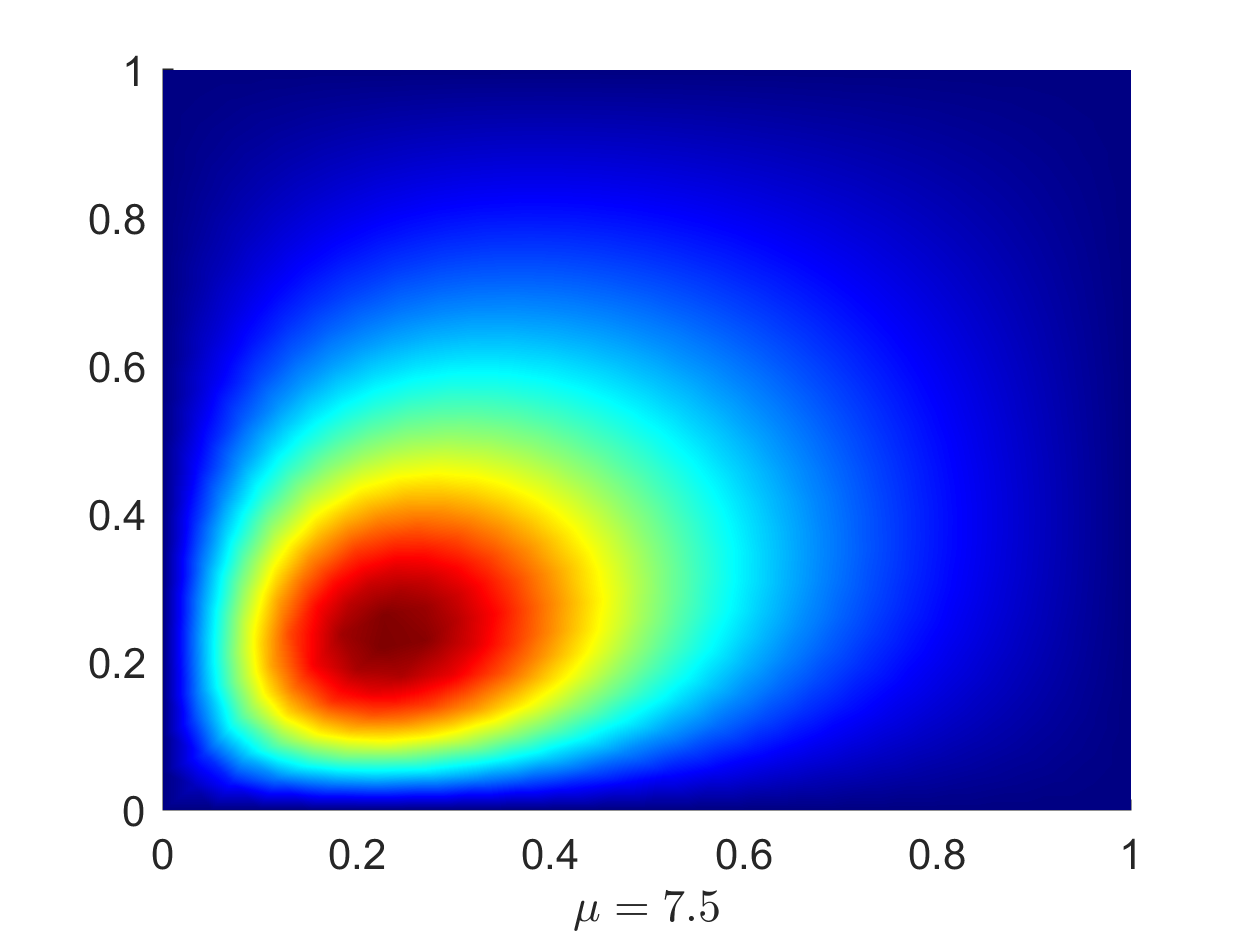}
\includegraphics[width=4.5cm]{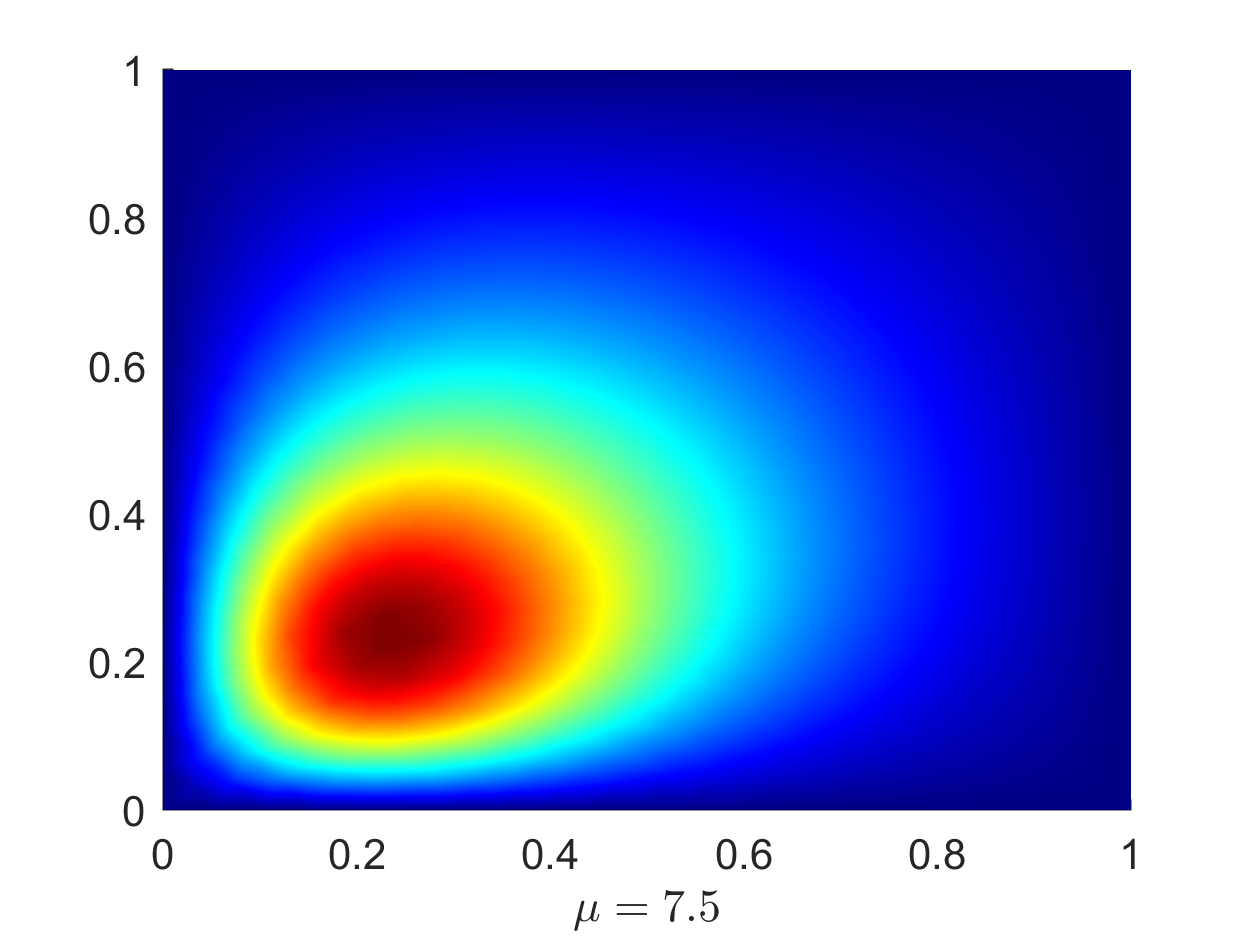}
\includegraphics[width=4.5cm]{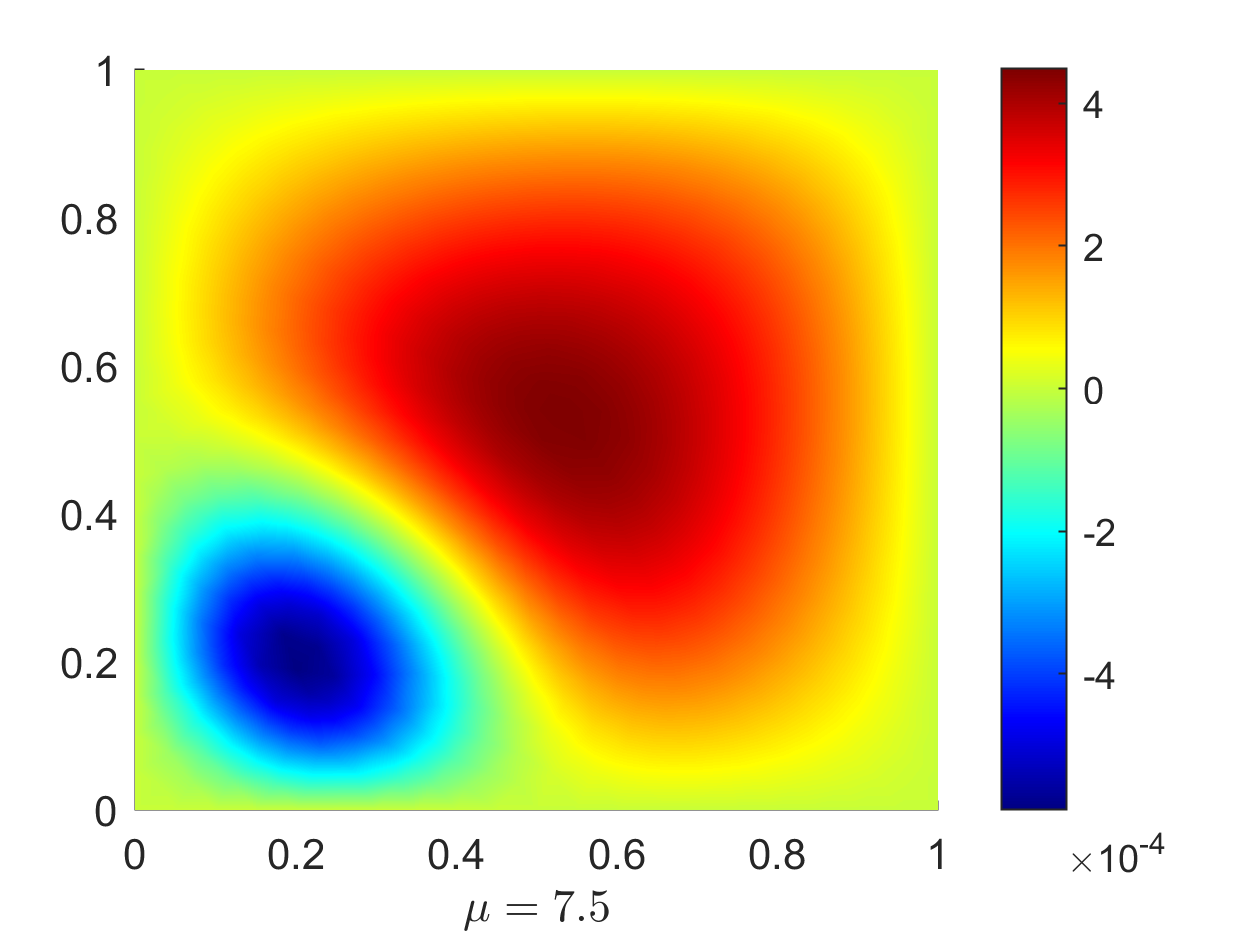}
\caption{1st eigenvectors of Problem~\eqref{1par_nonaffn} using FEM and DD model with $h=0.05$ and $\delta \mu=0.4$. We compare the FEM solution (left) with the DD solution (middle) and report on their difference (right)}
\label{fig:nonafnpar1}
\end{figure}

In Figure~\ref{fig:nonafnpar1}, we display the 1st eigenvectors obtained by using FEM and our DD model at the test points. The eigenvectors look very similar to each other. We show the difference between them and we can see that the error is of order $10^{-4}$.

%%%%%%%%%%%%%%%%%%%%%%%%%%%%%%%%%%%%%%%%%%
\subsubsection{Non-affine eigenvalue problem with multiple parameters}   
\begin{figure}
\centering
\subcaptionbox{$\mu_1=0.2,\mu_2=1$}{
\includegraphics[width=6.5cm]{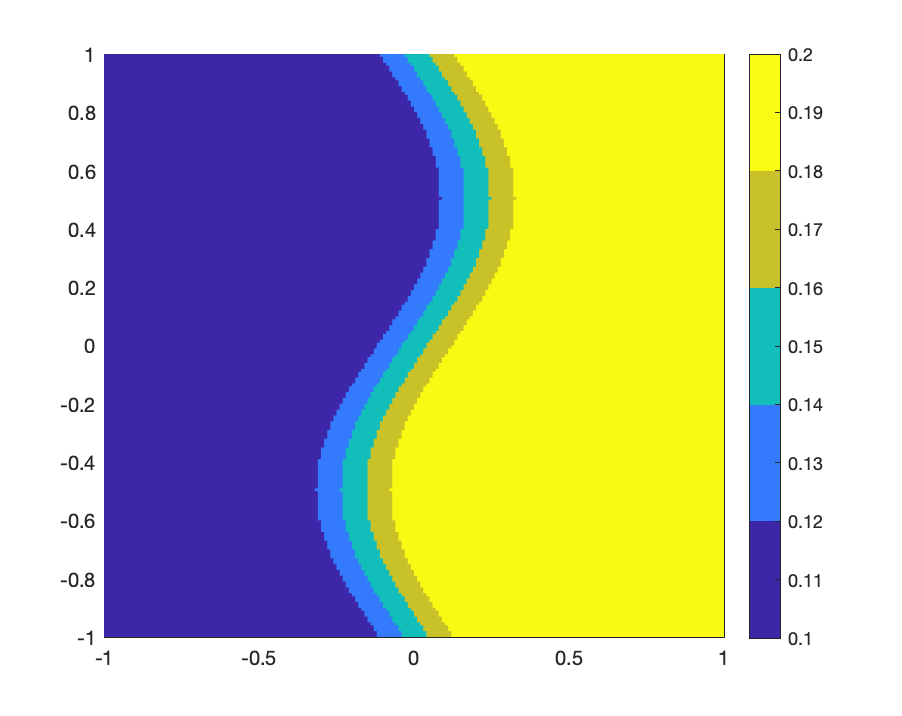}}
\subcaptionbox{$\mu_1=0.2,\mu_2=8$}{
\includegraphics[width=6.5cm]{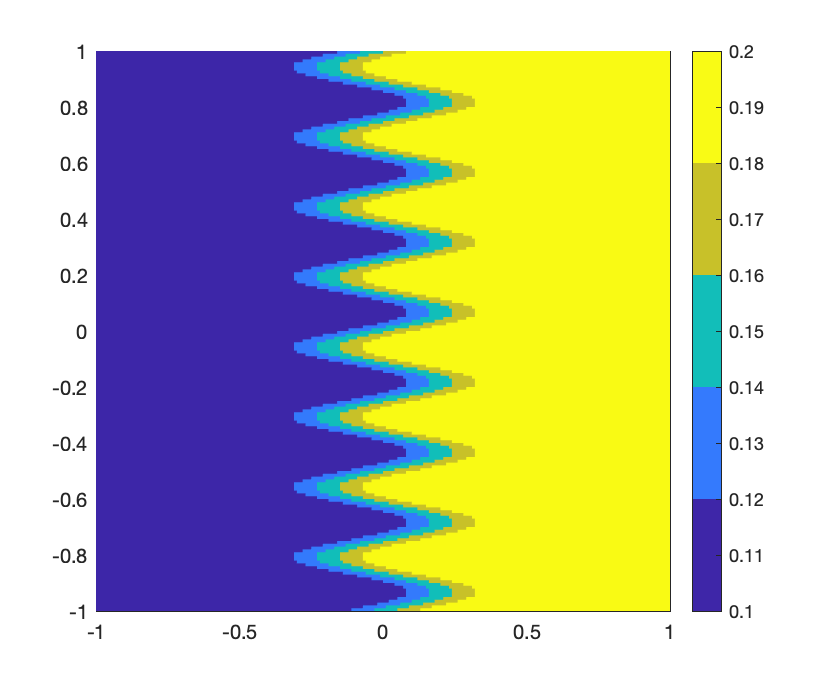}}
\caption{The weight function $\varepsilon(\pmb{\mu})$ in~\eqref{2par_nonaffn} for two different values of $\pmb{\mu}$.}   
\label{fig1_nonaff}  
\end{figure}

We now consider the eigenvalue problem: given $\pmb{\mu}=(\mu_1,\mu_2)\in\mathcal{P}$, find $\lambda(\pmb{\mu})\in\RE$ and a nonvanishing $u(\pmb{\mu})$ such that
\begin{equation}
\label{2par_nonaffn}
\left\{
\aligned
&-\Delta u=\lambda(\pmb{\mu})\varepsilon(x;\pmb{\mu})u&&\text{in }\Omega\\
&u(\mu)=0&&\text{on }\partial\Omega,
\endaligned
\right.
\end{equation}
where
\[
\varepsilon(x;\pmb{\mu})=\epsilon_1 \chi_{\Omega_1(\pmb{\mu})}(x)+\epsilon_2 \chi_{\Omega_2(\pmb{\mu})}(x)
\]
with $\Omega_1(\pmb{\mu})\cap\Omega_2(\pmb{\mu})=\emptyset$ and $\Omega_1(\pmb{\mu})\cup\Omega_2(\pmb{\mu})=\Omega=(-1,1)^2$. In our numerical result we take $\epsilon_1=0.1$, $\epsilon_2=0.2$ and the interface separating $\Omega_1(\pmb{\mu})$ and $\Omega_2(\pmb{\mu})$ is
$$ x_1=\mu_1\sin(\mu_2\pi x_2).$$
The parameter $\pmb{\mu}=(\mu_1,\mu_2)$ varies in $\mathcal{P}=[0.1,0.2]\times [1,8]$. We consider a regularized problem, where we introduce a linear transition between $\epsilon_1$ and $\epsilon_2$ near the interface, like it is done in~\cite{Fumagallietal16}. In particular, the function $\epsilon$ is presented in Figure~\ref{fig1_nonaff} for the parameters choice $\mu_1=0.2$, and $\mu_2=1$ and $8$, respectively.

In this example we can not write the right hand side, and hence the corresponding bilinear form, as an affinely parametric dependent form because of the nature of the function $\varepsilon(\pmb{\mu})$. We solve this eigenvalue problem using our GPR based DD model.

For this example we have taken two uniformly distributed in the parameter space $\mu_1$ and $36$ uniformly distributed points in the parameter space of $\mu_2$, that is in the interval $[1,8]$ with stepsize $0.2$. We take their tensor product and get $72$ sample points in the parametric space $\mathcal{P}$. The snapshot matrix is constructed by taking the first eigenvector at those sample points as columns. As a test set we choose $\mu_1=0.15$ and  $\mu_2=1$, $2.5$, $4$, $4.5$, $6.5$, $8$, respectively. These are different values from those taken to generate snapshots. 

\begin{table}
\footnotesize
 	 	\centering
 	\begin{tabular}{|c|c|c|c|c|c|c|c|c|c|c|c|} 
 		\hline
		$(\mu_1,\mu_2)$&  {\begin{tabular}[c]{@{}c@{}} Method \end{tabular}} &
 		  {\begin{tabular}[c]{@{}c@{}} Eigenvalue-I \end{tabular}} &
		 {\begin{tabular}[c]{@{}c@{}} Eigenvalue-II \end{tabular}} &
 		{\begin{tabular}[c]{@{}c@{}}   Eigenvalue-III \end{tabular}}& 
 		{\begin{tabular}[c]{@{}c@{}}  Eigenvalue-IV \end{tabular}} \\
		%{\begin{tabular}[c]{@{}c@{}} $\mu=2.5$ \end{tabular}}  \\	
		\hline
(0.15,1.5)
& FEM&31.48699165&75.22816969 &89.16775007 &143.71355581 	\\
&DD  &31.49811030 &75.28078739	&89.10207780 &142.38589826 	\\
\hline
(0.15,3.5)&
FEM &31.54972648& 75.32810691& 89.03858404& 144.08480023	\\
&DD&31.56475170 & 75.35994556&88.96400162&143.35827613 	\\
\hline
(0.15,5.5)&
FEM&31.56473098  &  75.38499108&  89.05984549 &  144.34148451 	\\
&DD&31.57095675& 75.40394116 & 89.00349397 &143.69597428\\	
\hline
(0.15,7.5)&
FEM& 31.58159550&  75.42566021& 89.01447975&144.45513191	\\
&DD&  31.58010472&75.42981927&89.00210188& 143.80005269	\\
\hline
\end{tabular}
\caption{ First four eigenvalues using FEM and DD model at the test points of Problem~\eqref{2par_nonaffn} with mesh size $h=0.1$.}
\label{table6}
\end{table}

In Table~\ref{table6} we report the first four eigenvalues obtained by using FEM and our DD model for the test parameter mentioned in the last paragraph. The predicted eigenvalues by the DD model are matching with the eigenvalues obtained by using FEM. Note that for the construction of GPR we do not use any information of the snapshot matrix, so we can predict more than one eigenvalues using the data. 

\begin{figure}
\centering
\includegraphics[width=4.5cm]{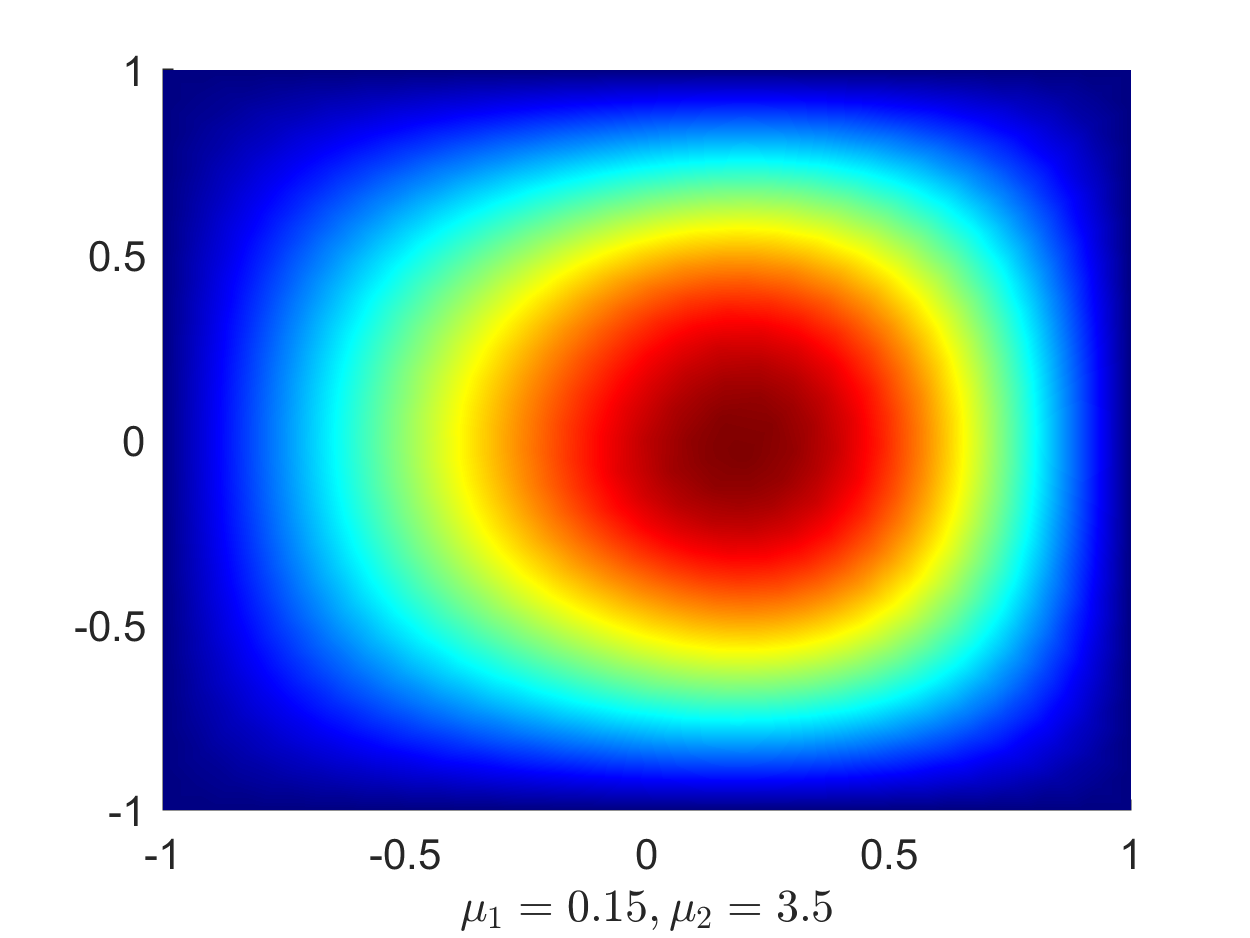}
\includegraphics[width=4.5cm]{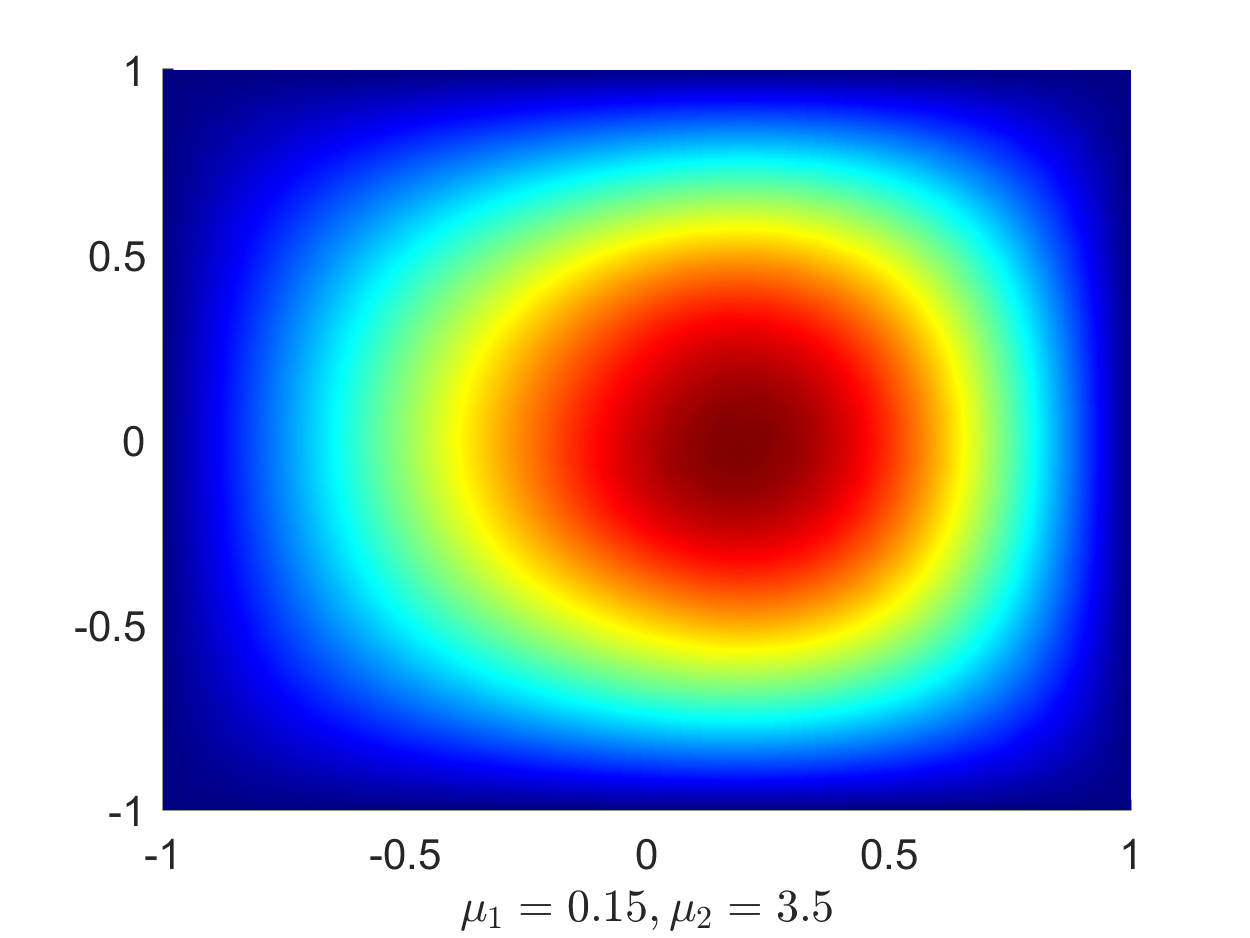}
\includegraphics[width=4.5cm]{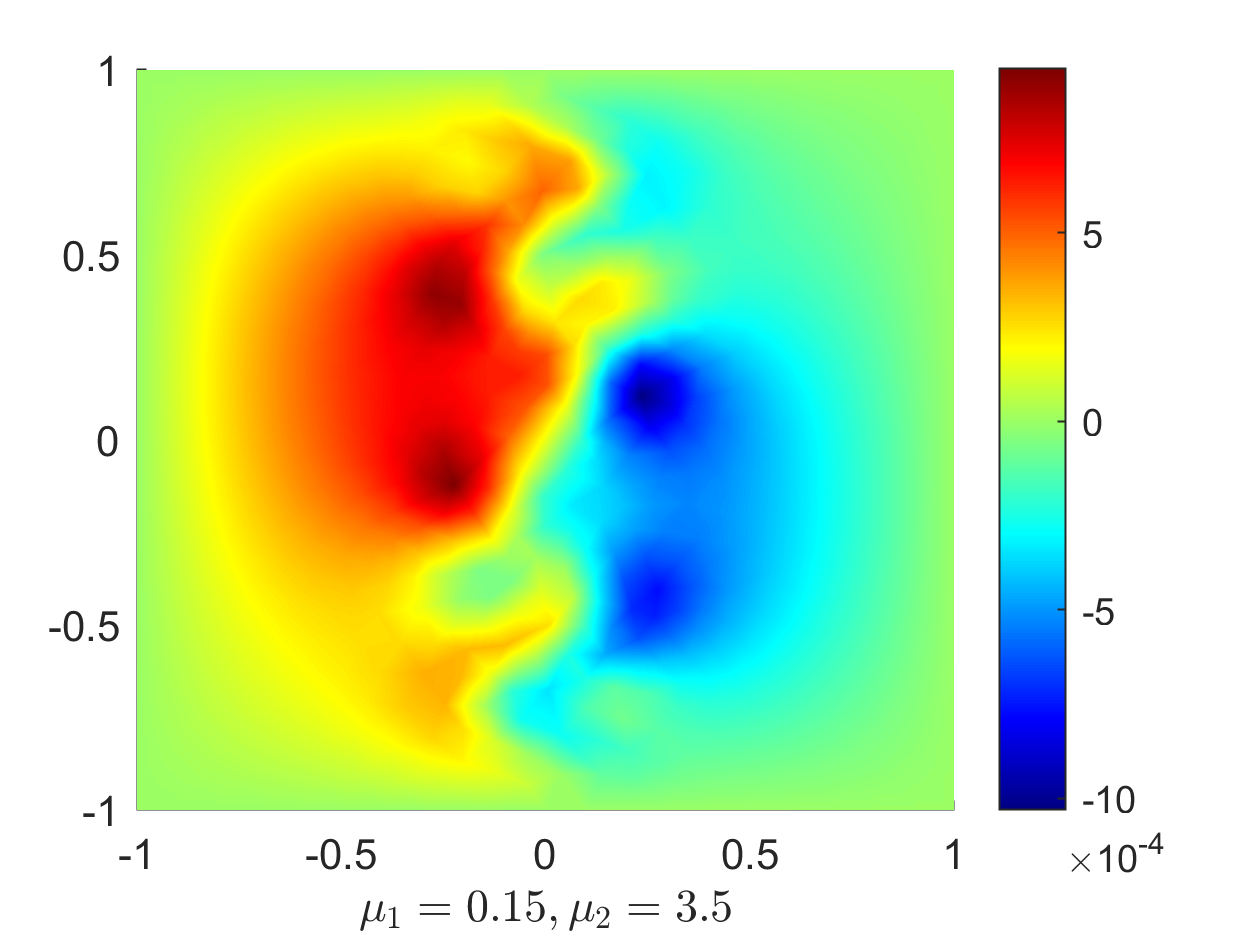}

\includegraphics[width=4.5cm]{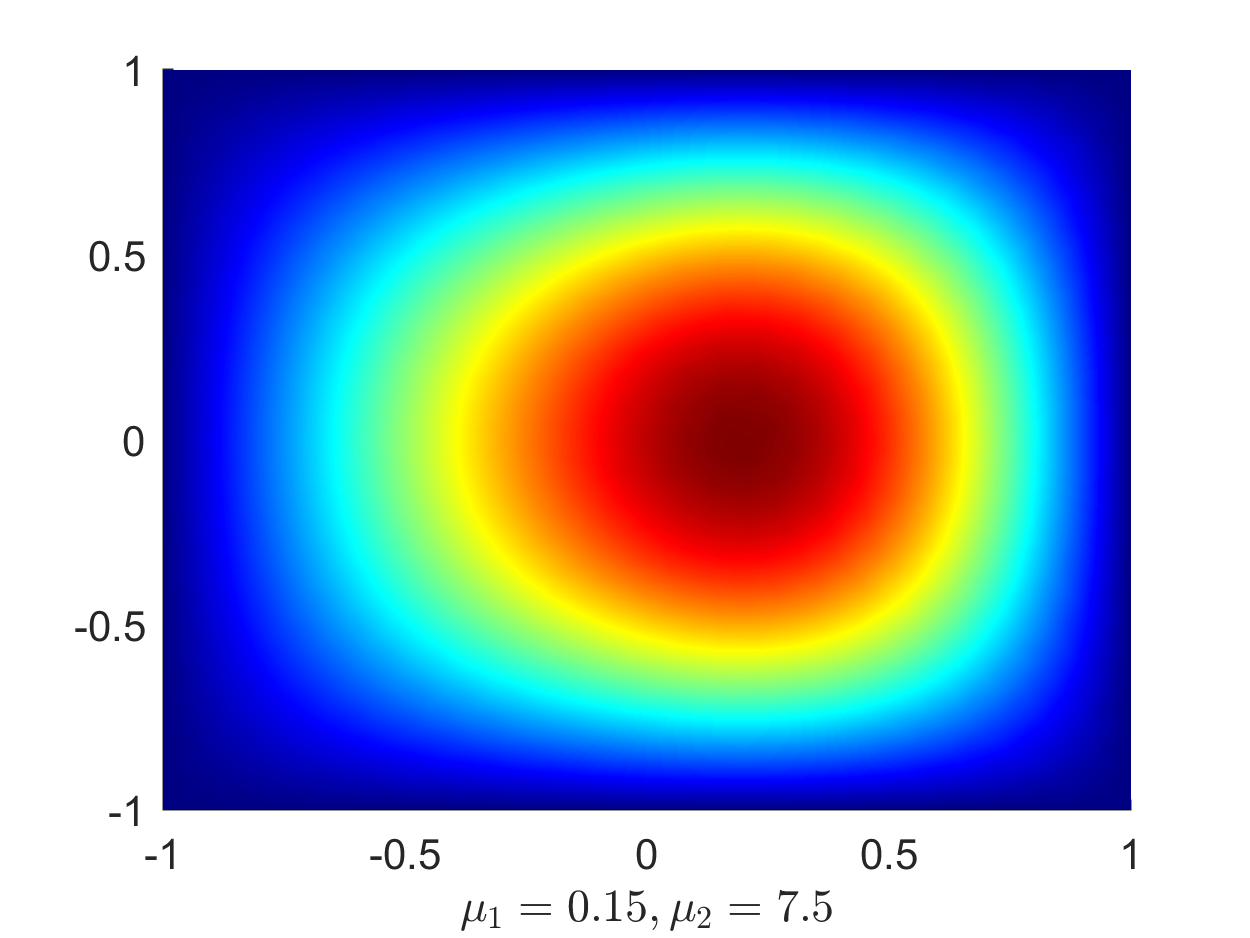}
%         \caption{FEM}
\includegraphics[width=4.5cm]{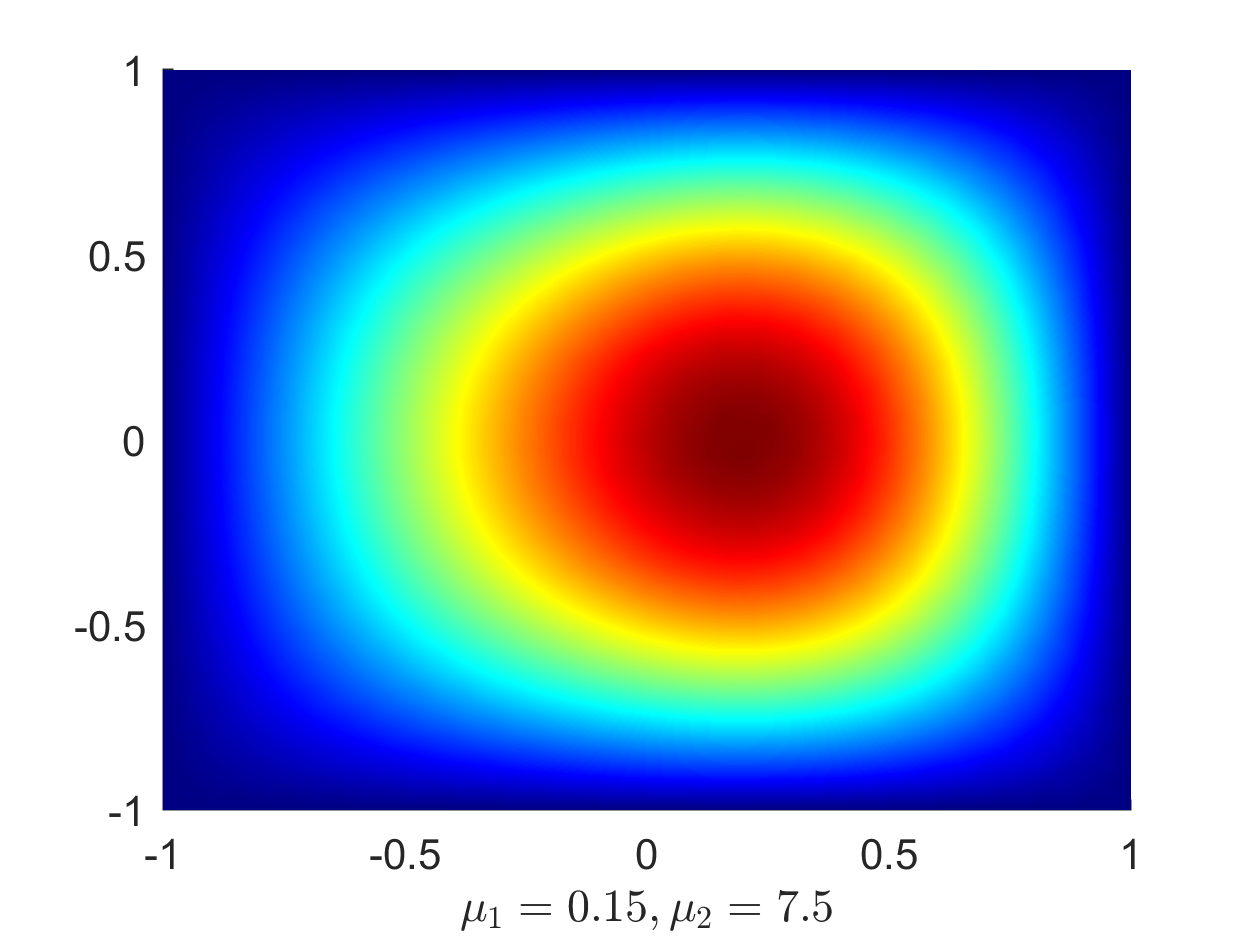}
%          \caption{DD ROM}
\includegraphics[width=4.5cm]{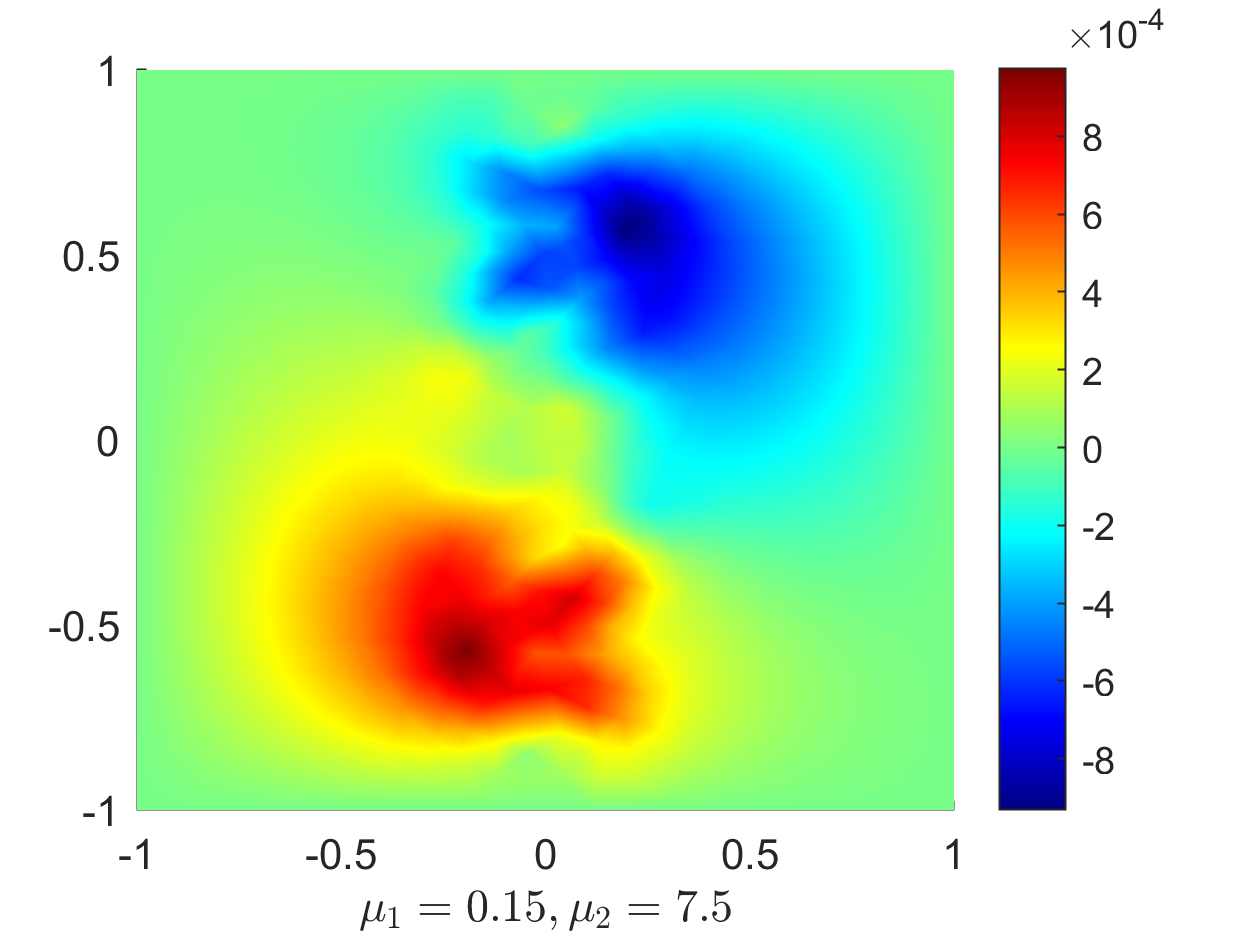}
%          \caption{Error}
\caption{First eigenvectors of Problem~\eqref{2par_nonaffn} at parameters $\pmb{\mu}=(0.15,3.5)$ and $\pmb{\mu}=(0.15,7.5)$ respectively, using FEM (left) and DD model (middle). The difference is reported on the right.}
\label{fig2_nonaff}       % Give a unique label
\end{figure} 

In Figure~\ref{fig2_nonaff} we show the first eigenvectors obtained by FEM and DD for the parameters $\pmb{\mu}=(0.15,3.5)$ and $\pmb{\mu}=(0.15,7.5)$. We also show the relative error between them which is of order $10^{-4}$.
%%%%%%%%%%%%%%%%%%%%%%%%%%%%%%%%%%%%%%%%%%%%%%%%%%%%%%%%%%%%%%%%%%%%%%%%
\subsection{Some results on crossing eigenvalues}
\begin{figure}
\centering
\includegraphics[width=12cm]{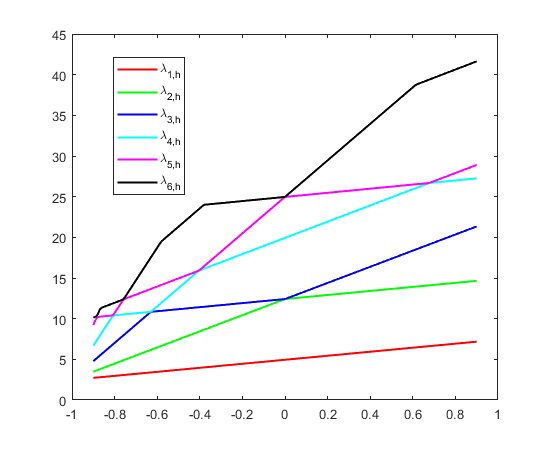}
\caption{First six sorted eigenvalues of Problem~\eqref{mdl:crossing}.}
\label{fig:evplot}       % Give a unique label
\end{figure}

Let us consider the following eigenvalue problem 
\begin{equation}
\label{mdl:crossing}
    \left\{
\aligned
&-\nabla \cdot(A\nabla u)=\lambda u&&\text{in }\Omega=[-1,1]^2\\
&u(\mu)=0&&\text{on }\partial\Omega,
\endaligned
\right.
\end{equation}
where the diffusion matrix is chosen as
$$A=\begin{bmatrix}
1 &0 \\
0 & 1+\mu
\end{bmatrix}$$
with the parameter $\mu$ varying in the space $\mathcal{P}=[-0.9,0.9]$.
The analytic eigenvalues are
\begin{align*}
\lambda_{m,n}(\mu)=\frac{\pi^2}{4}(m^2+(1+\mu)n^2), \quad m,n=1,2,\dots
\end{align*}
and the eigenvectors are
\begin{align*}
u_{m,n}(\mu)=\cos\left(\frac{m\pi}{2}x\right)\cos\left(\frac{n\pi}{2}y\right).
\end{align*}
This problem, already considered in~\cite{eccomas,moataz}, is particularly simple because the eigenfunctions do not depend on $\mu$, but are already complicated enough because of the eigenvalues crossings.
In Figure~\ref{fig:evplot} we show the first six eigenvalues obtained with FEM. We can see that the eigenvalues other than the first eigenvalues are intersecting. We will investigate how our model works for the approximation of the intersecting eigenvalues and eigenvectors. We present the numerical results for the first and the third eigenvalues and eigenvectors for this problem computed with the DD model. In order to calculate the third eigenvector, we consider the snapshot matrix consisting of the third eigenvectors at the sample parameters $\mu_1$, $\mu_2$,\dots, $\mu_{n_s}$ as the columns.

In our computations we took $19$ uniformly distributed sample points on the parameter space $[-0.9,0.9]$, that is the parameter stepsize is $0.1$. Then, using the eigenvalues and eigenvectors at these points, we construct the GPR following the idea described in Section~\ref{sec:DD}. We present the eigenvalues and eigenvectors at four test points $\mu=-0.75$, $-0.25$, $0.25$, $0.75$.

\begin{table}
\footnotesize
 	 	\centering
 	\begin{tabular}{|c|c|c|c|c|c|c|c|} 
 		\hline
 		h &  {\begin{tabular}[c]{@{}c@{}} Method \end{tabular}} &
		  {\begin{tabular}[c]{@{}c@{}} $\mu=-0.75$ \end{tabular}} &
 		 {\begin{tabular}[c]{@{}c@{}} $\mu=-0.25$ \end{tabular}} & 
 		 {\begin{tabular}[c]{@{}c@{}} $\mu=0.25$ \end{tabular}} &
 		{\begin{tabular}[c]{@{}c@{}}   $\mu=0.75$  \end{tabular}}\\
 \hline
 0.1&FEM&3.09172930&4.32853369&5.56526834 &6.80197424\\
    & DD&3.09175362&4.32850239&5.56525117&6.80199994\\
  \hline
  0.05& FEM&3.08606437&4.32052203 &5.55496589 &6.78940395\\
  &DD&3.08606913&4.32051579 &5.55496244 &6.78940910\\
   \hline
    0.01&FEM&3.08432204 &4.31805168 &5.55178071&6.78550948\\
    &DD&3.08432225&4.31805140&5.55178056&6.78550971\\         
  \hline
 	\end{tabular}
\caption{First eigenvalues of \eqref{mdl:crossing} for different mesh using DD model with sample points $\mu_{tr}=-0.9:0.1:0.9$ and FEM model.}
\label{crossing:1stev}
\end{table}

\begin{figure}
\includegraphics[width=4.5cm]{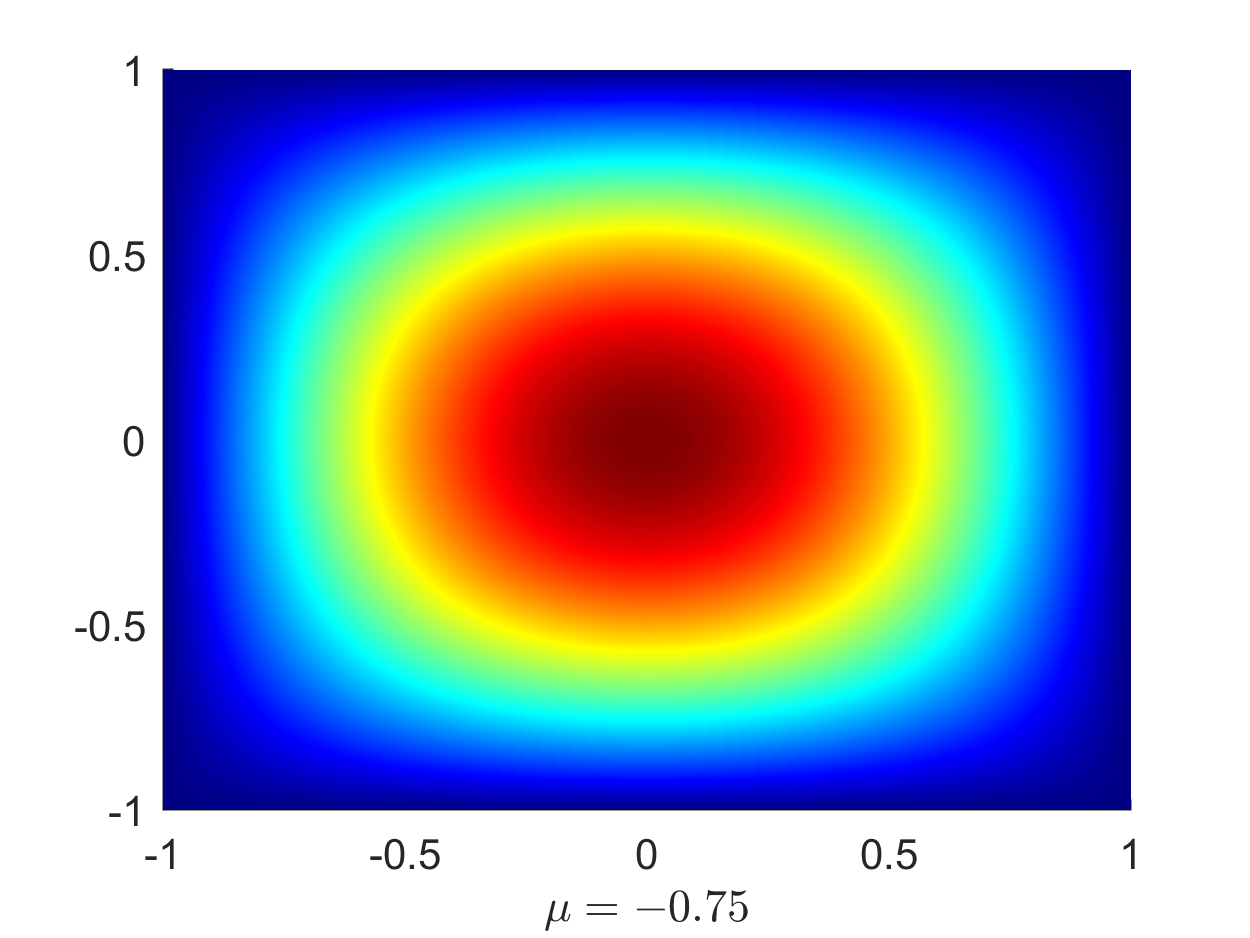}
        % \caption{$\mu=-0.75$,FEM}
\includegraphics[width=4.5cm]{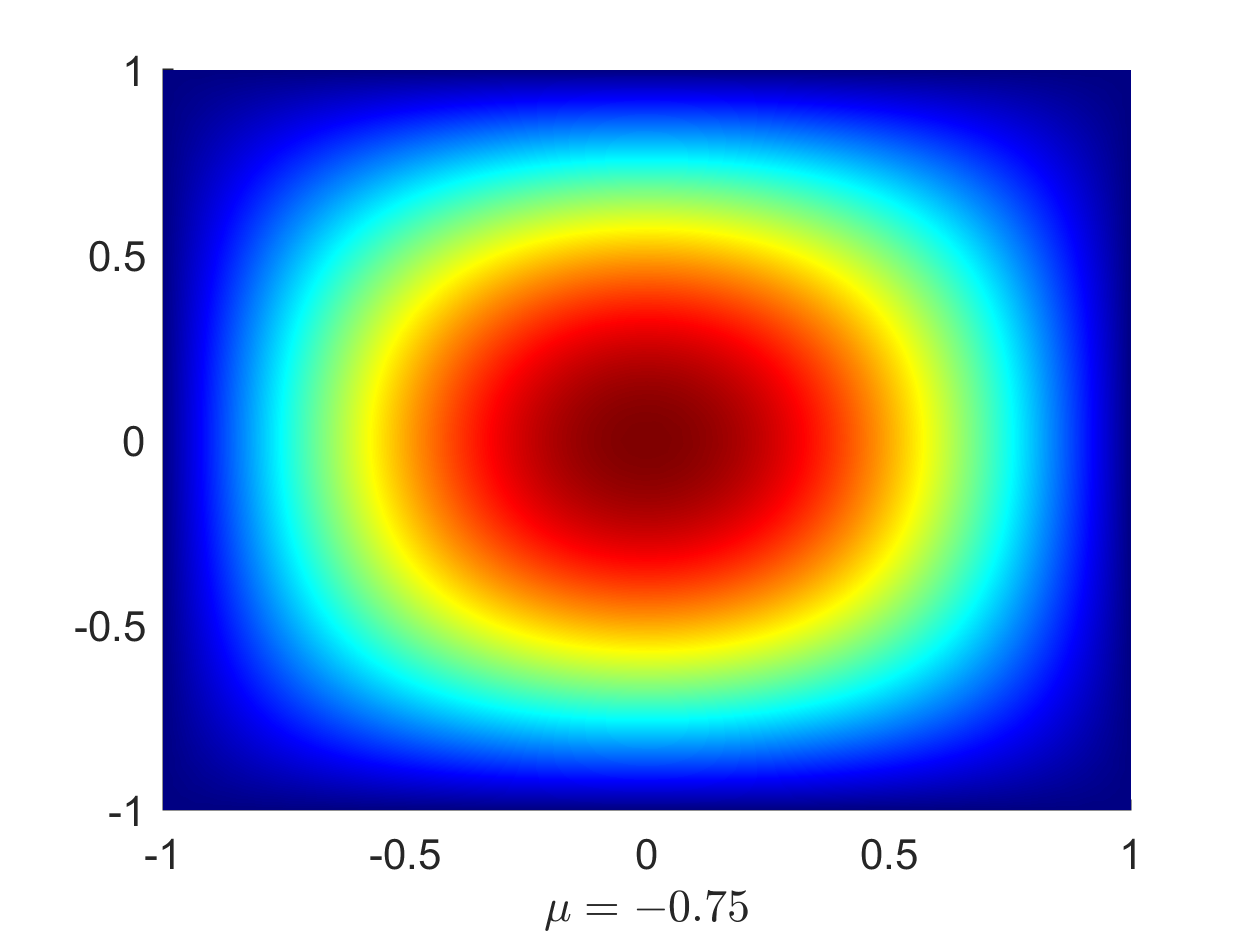}
       %  \caption{$\mu=-0.75$,DD}
\includegraphics[width=4.5cm]{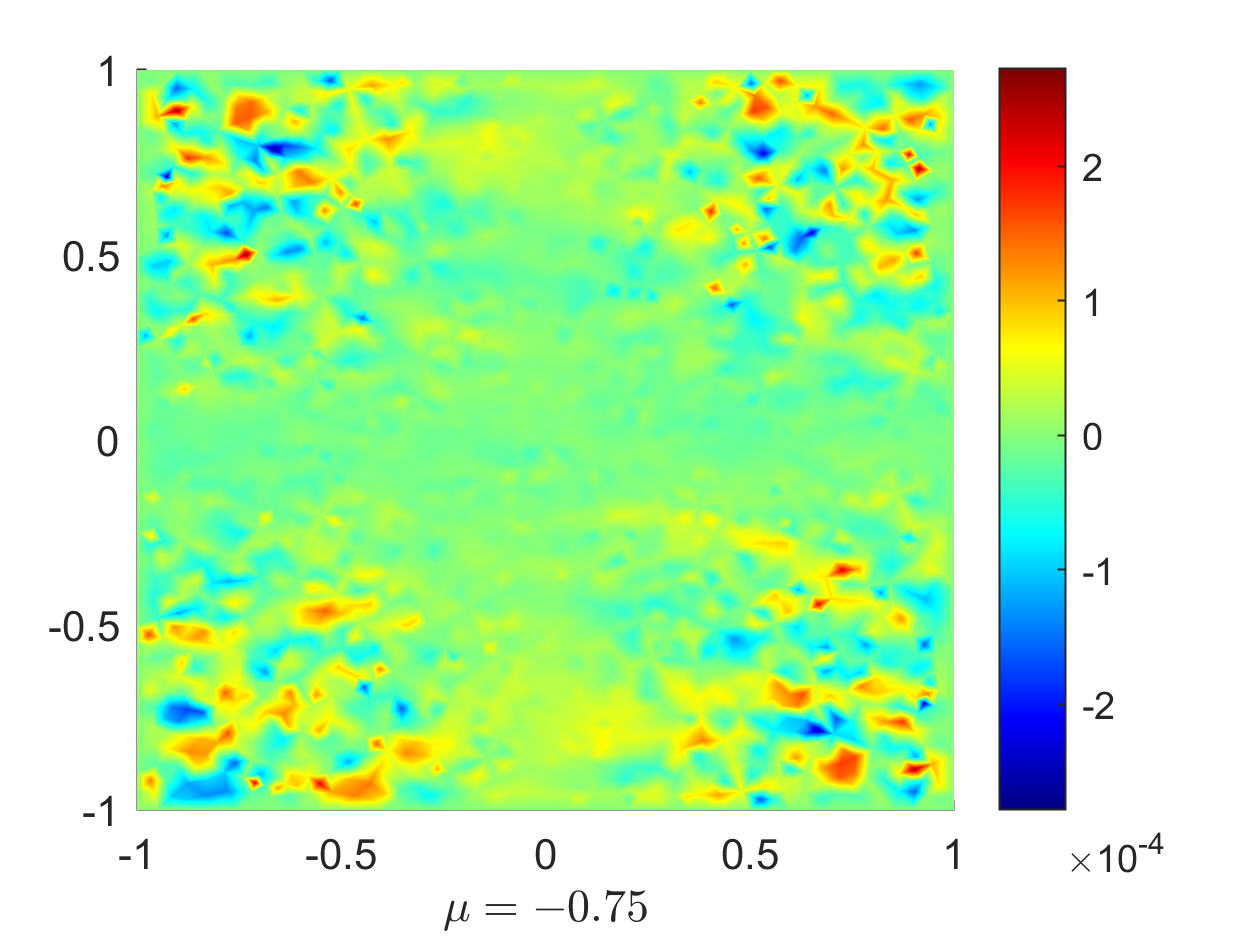}
       %  \caption{$\mu=-0.75$,DD}
       
\includegraphics[width=4.5cm]{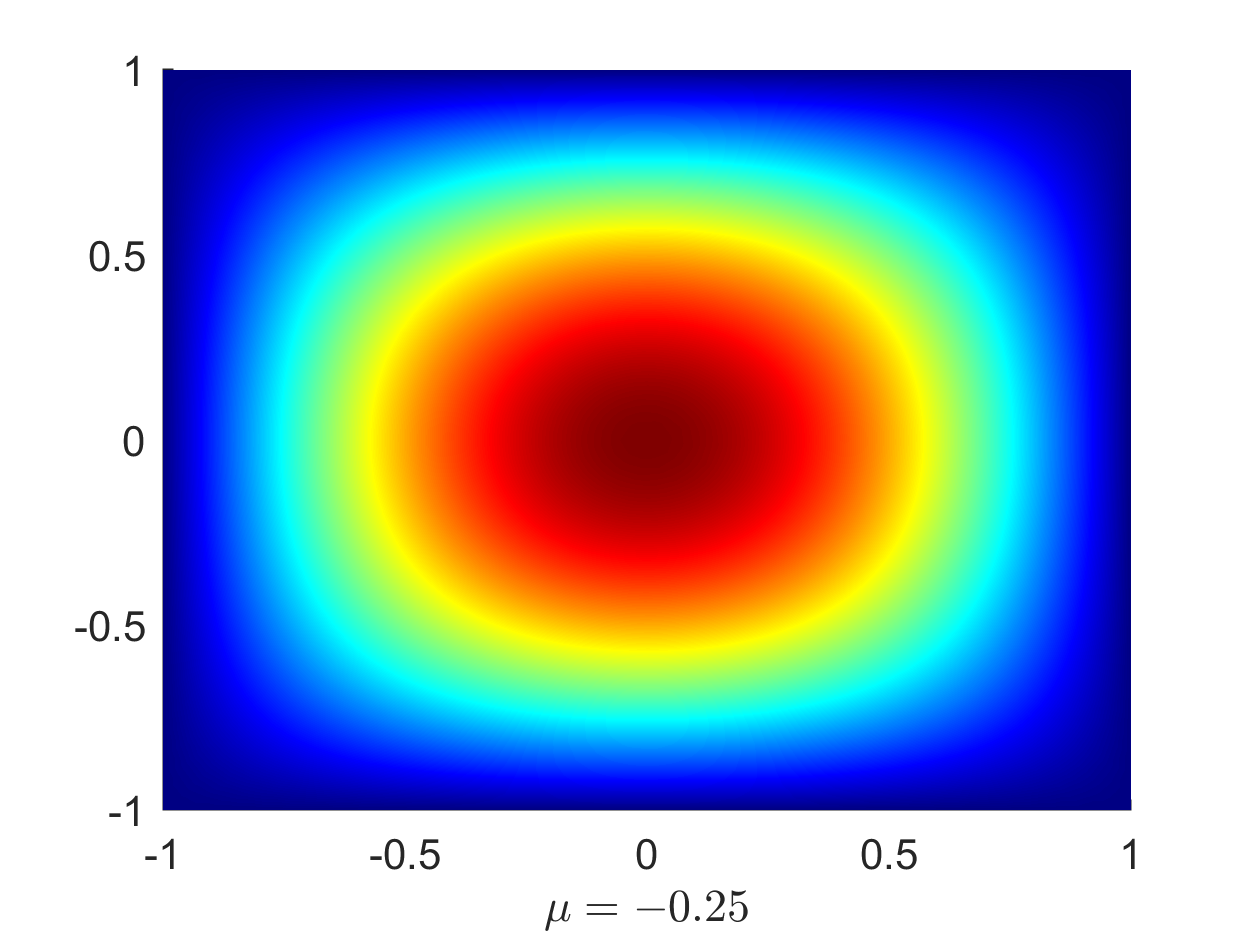}
       %  \caption{$\mu=-0.25$,FEM}
\includegraphics[width=4.5cm]{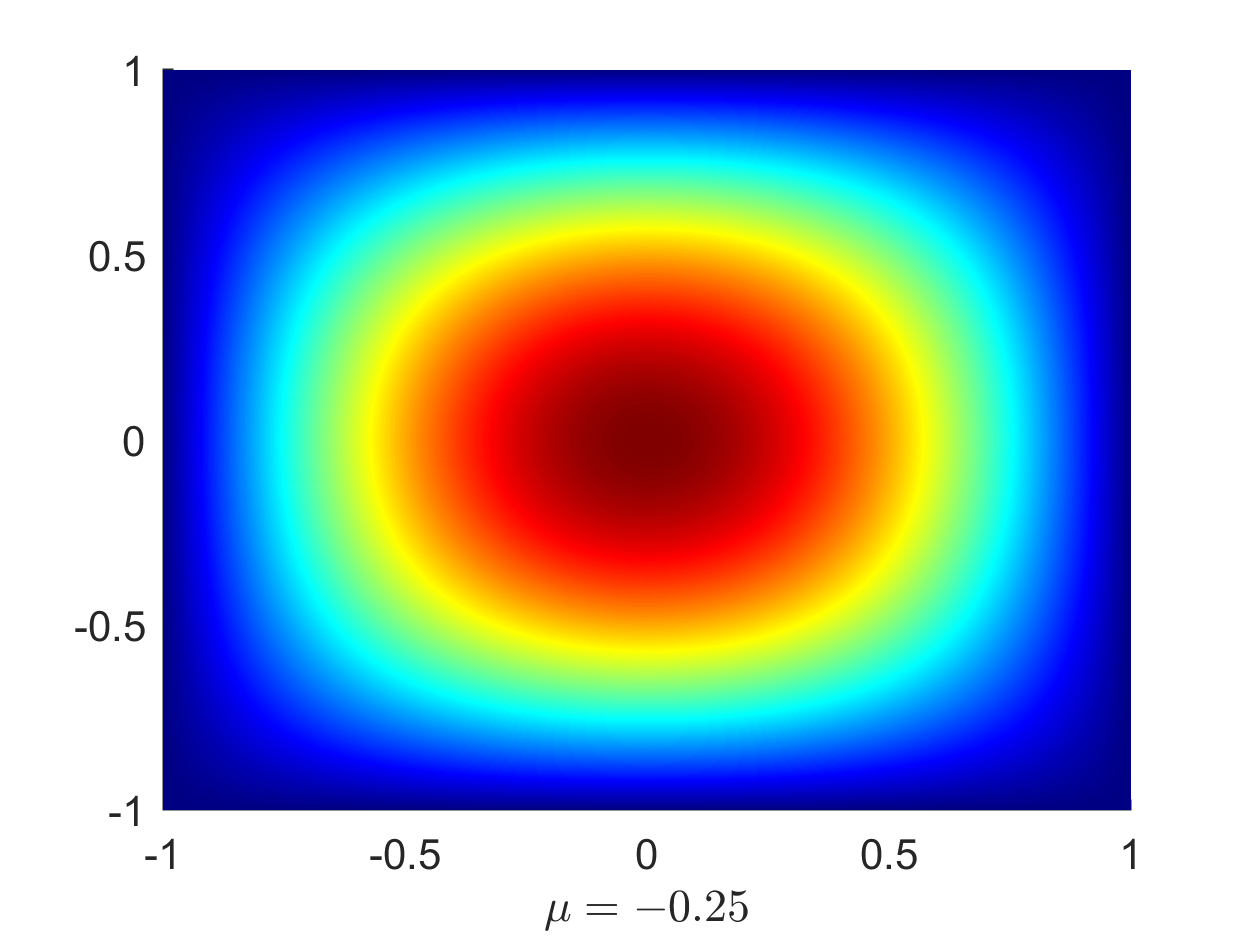}
      %   \caption{$\mu=-0.25$,DD}
\includegraphics[width=4.5cm]{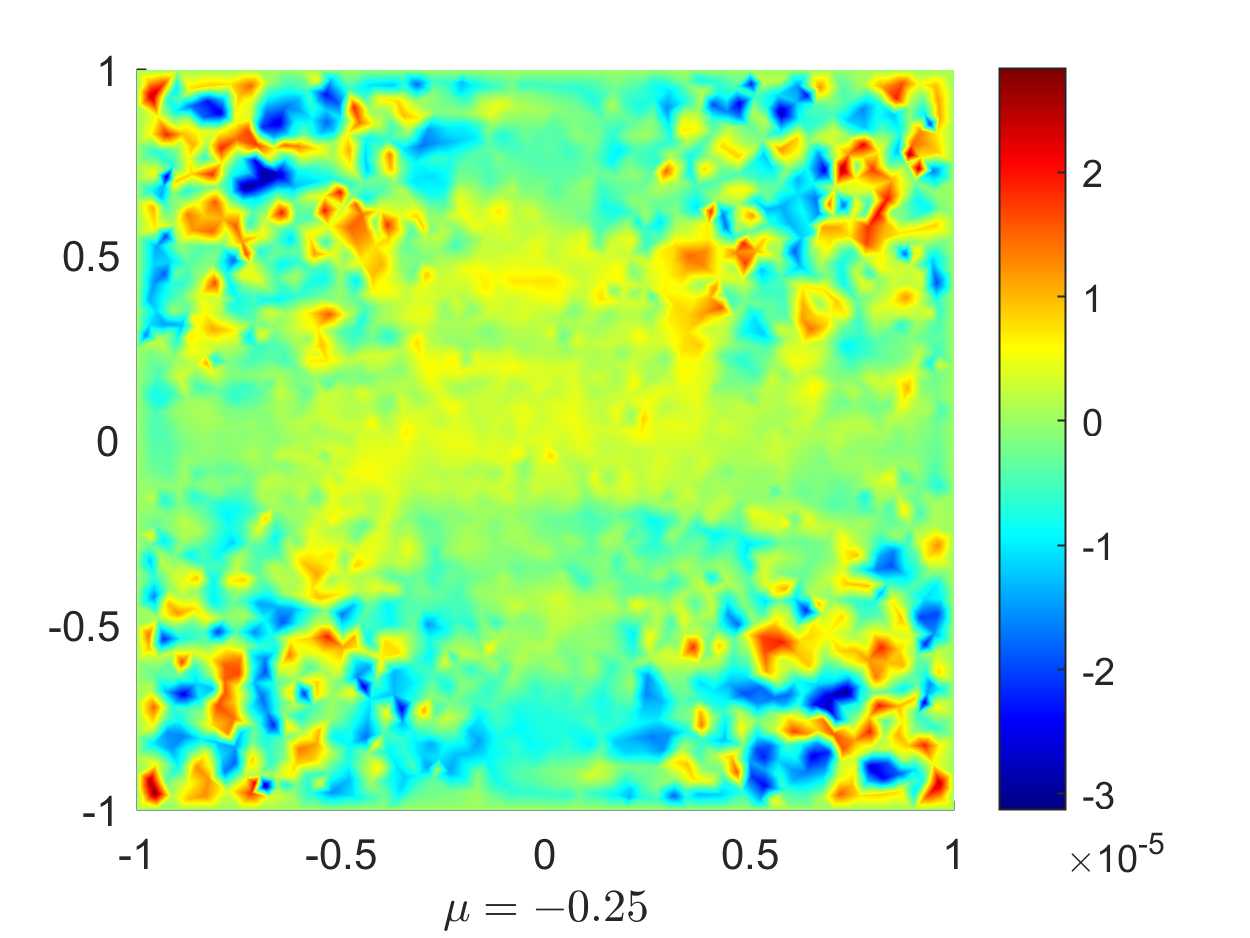}
       %  \caption{$\mu=-0.75$,DD}

\includegraphics[width=4.5cm]{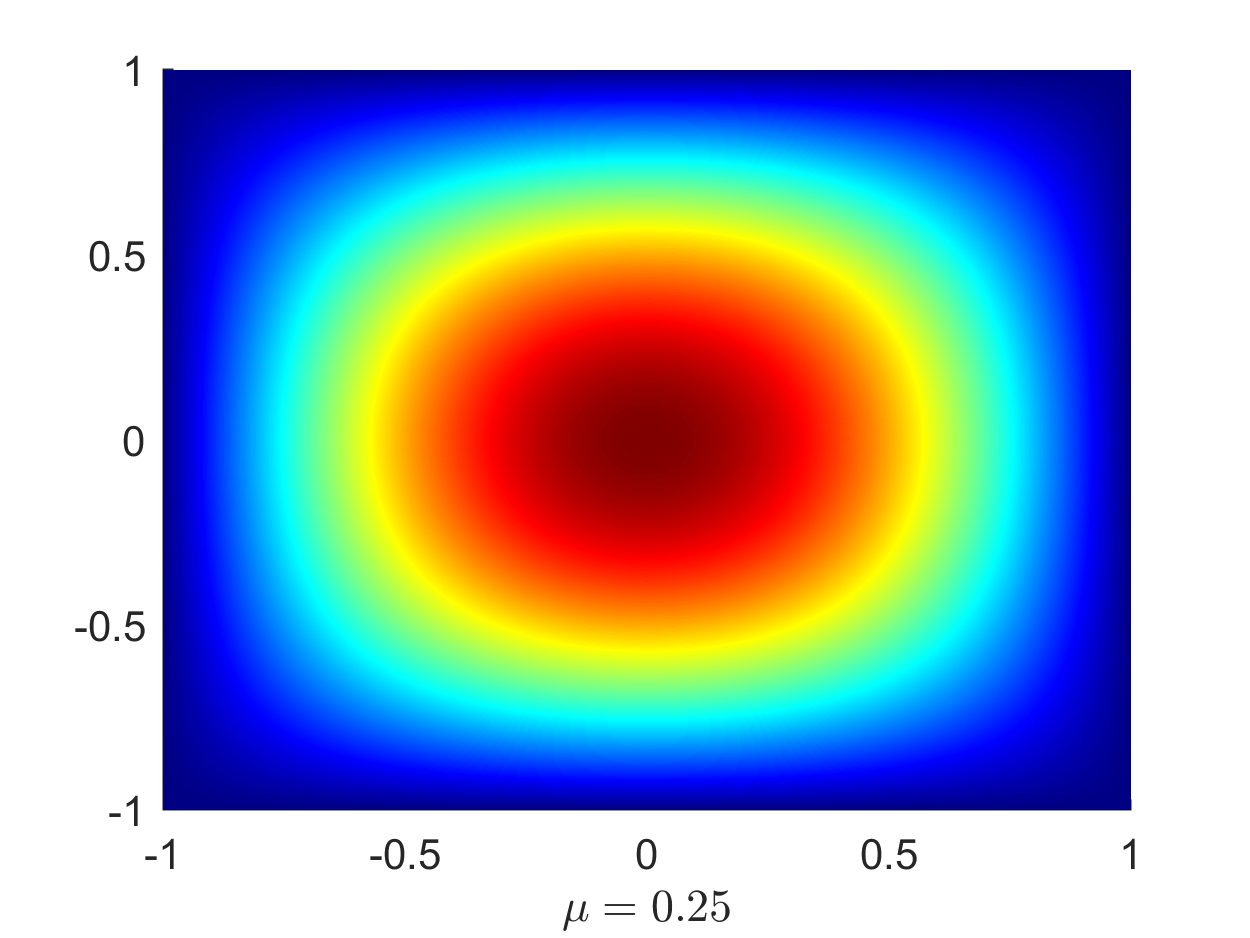}
      %   \caption{$\mu=0.25$,FEM}
\includegraphics[width=4.5cm]{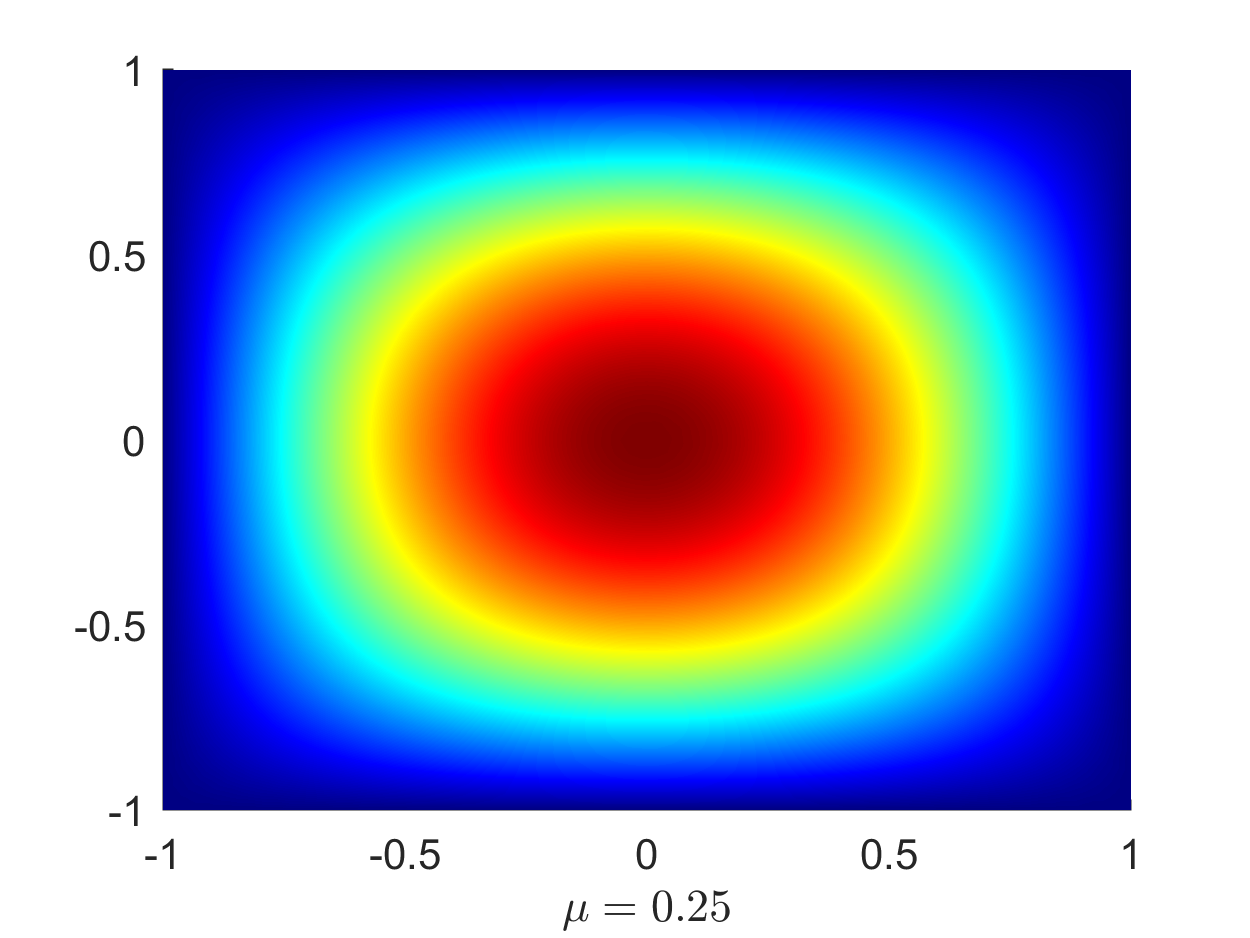}
       %  \caption{$\mu=0.25$,DD}
\includegraphics[width=4.5cm]{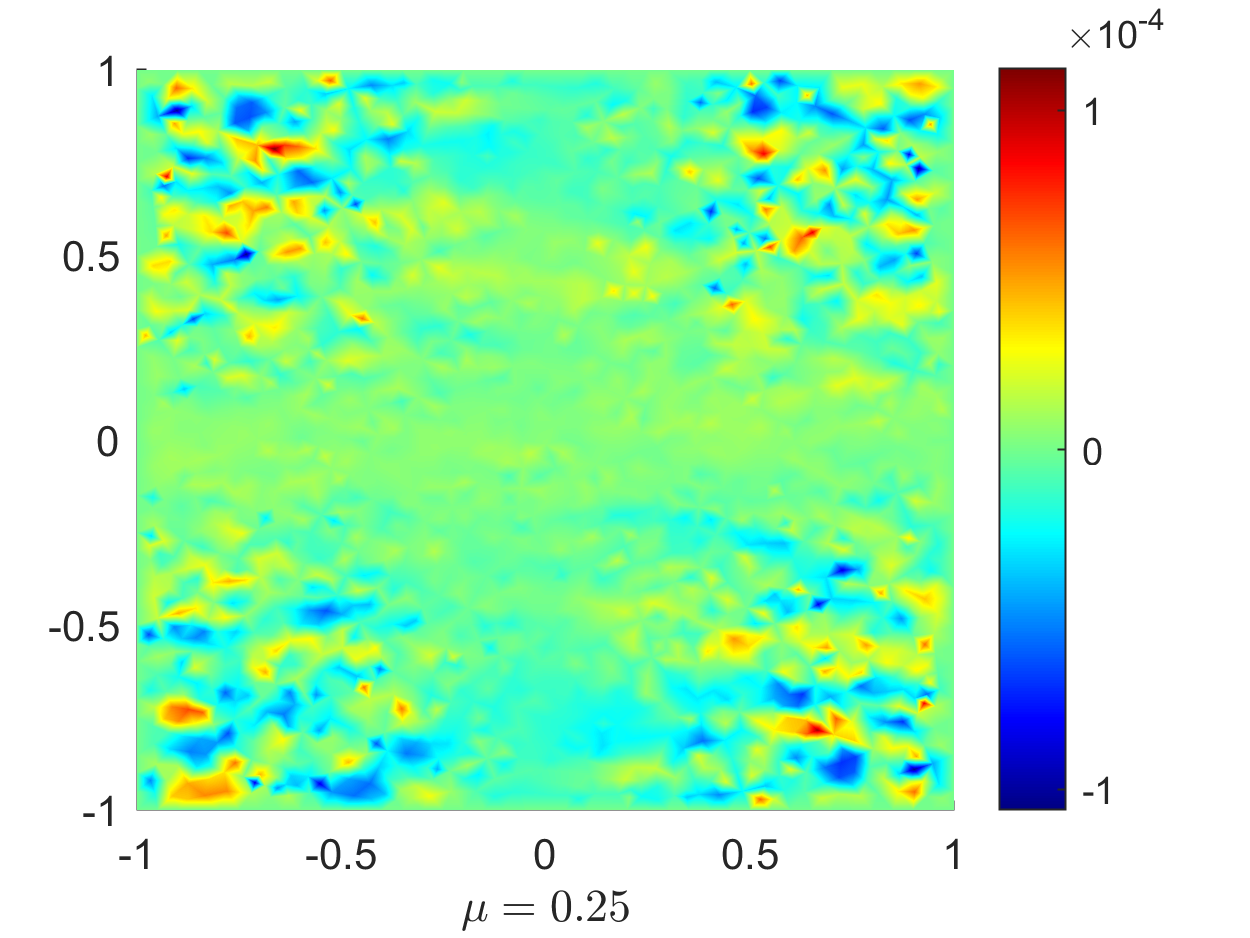}
       %  \caption{$\mu=-0.75$,DD}
\includegraphics[width=4.5cm]{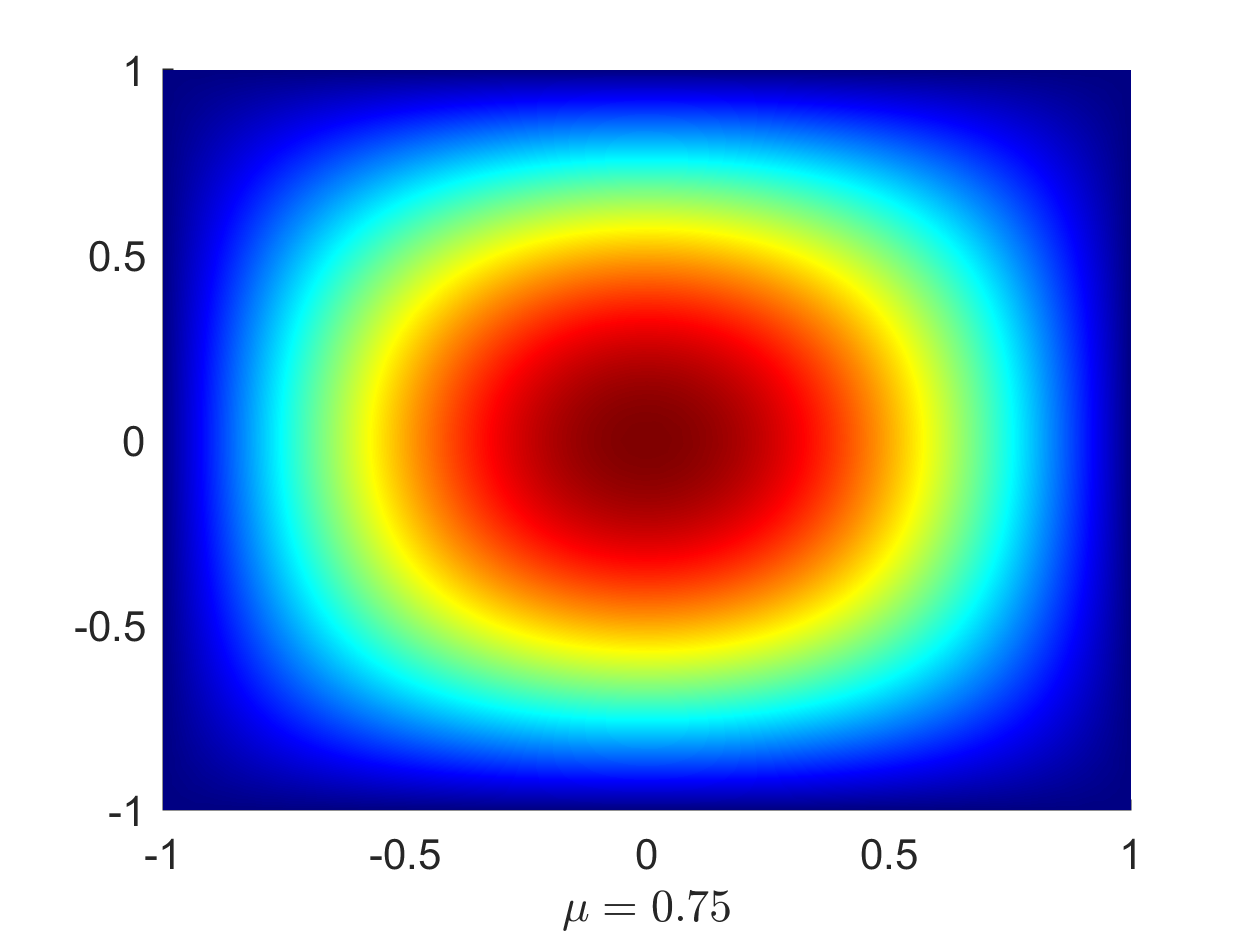}
%         \caption{FEM}
\includegraphics[width=4.5cm]{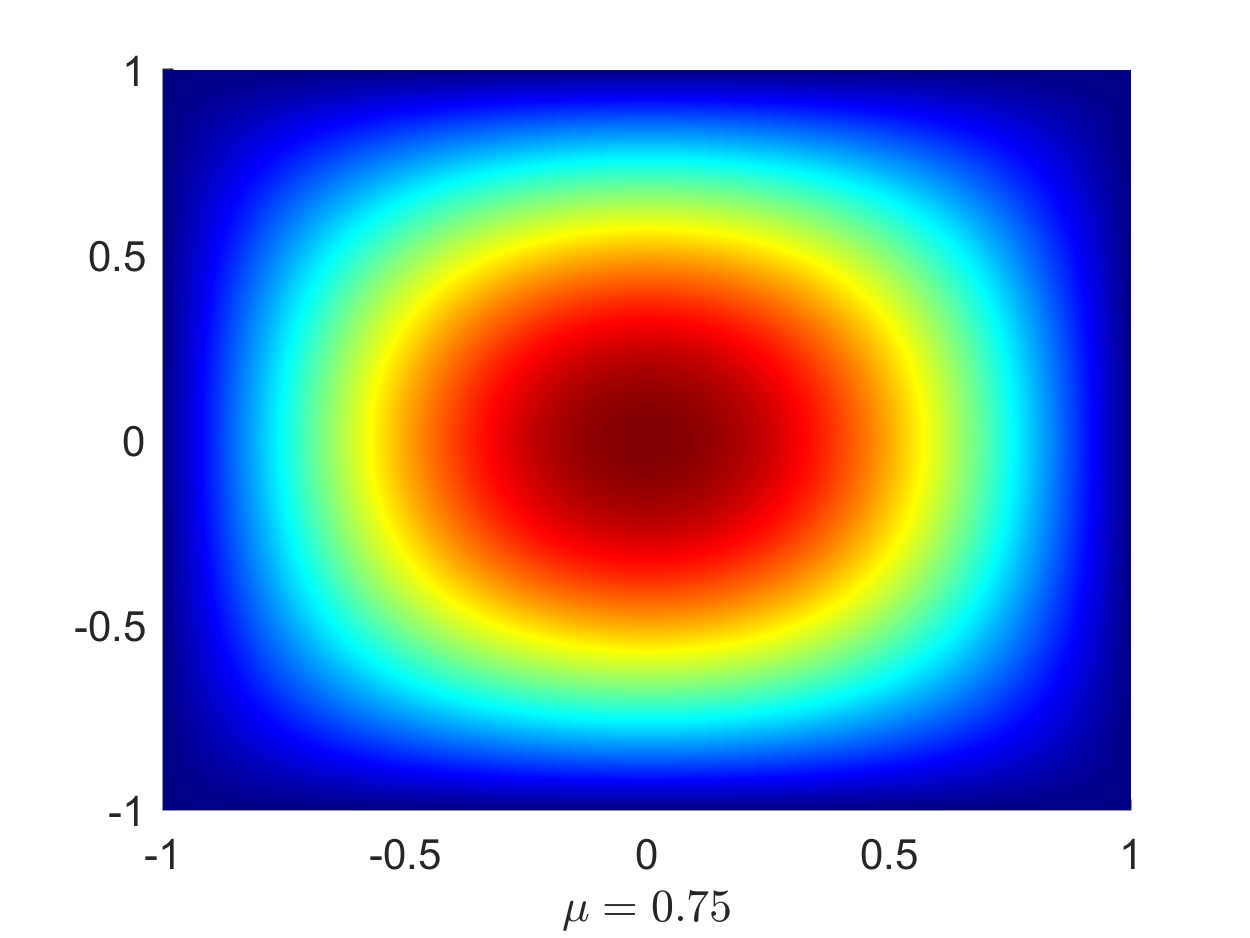}
%         \caption{DD}
\includegraphics[width=4.5cm]{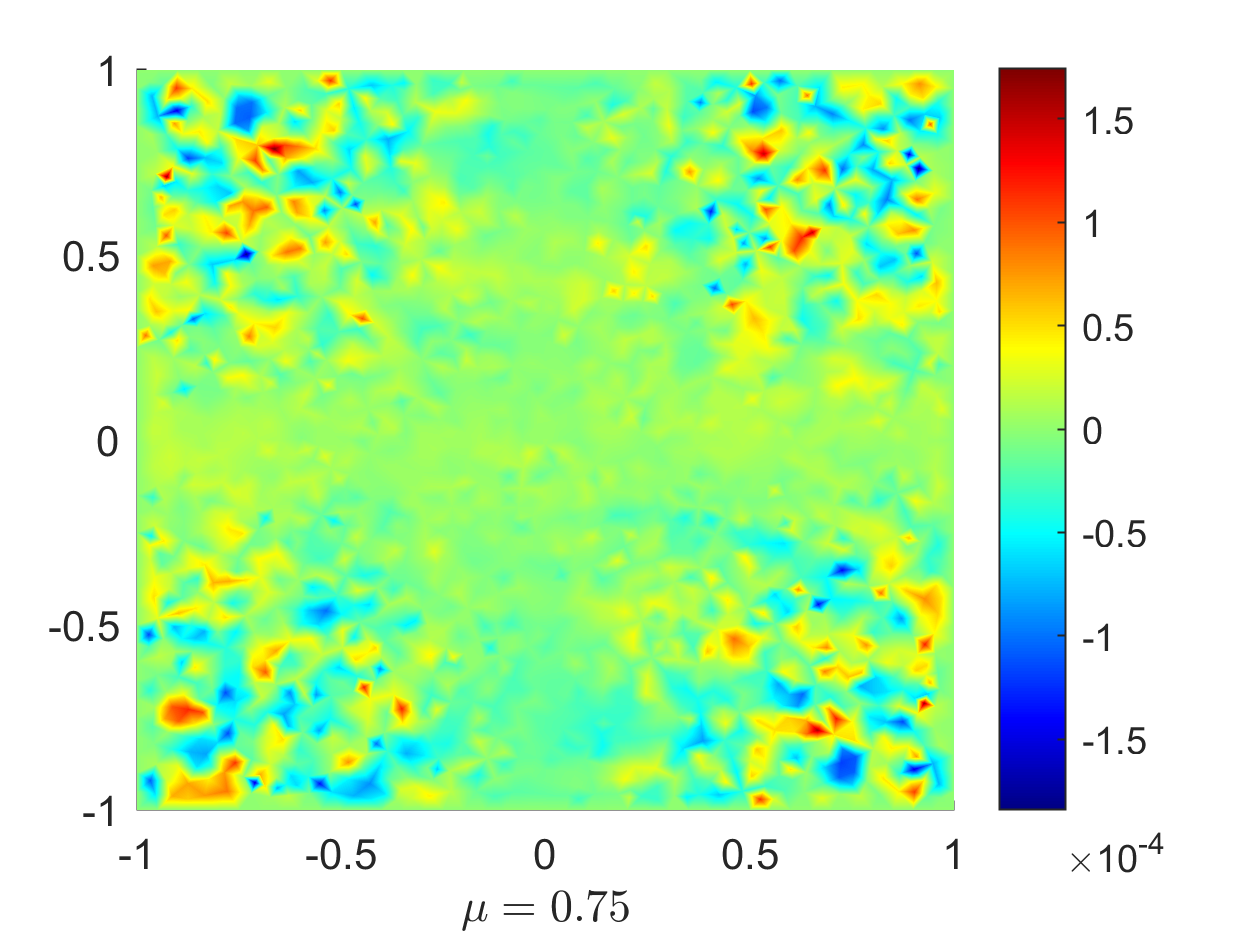}
%         \caption{Error}
\caption{1st eigenvectors corresponding to mesh size $h=0.05$ using FEM (left) and  DD model with $1$ POD basis (middle) when the training set is $\mu_{tr}=-0.9:0.1:0.9$. The difference is reported on the right column.}
\label{fig:crossing1stev}
\end{figure}

In Table~\ref{crossing:1stev}, we report the first eigenvalues obtained by our DD model at these four points for the mesh size $h=0.1$, $0.05$, and $0.01$, respectively. For comparison we report the eigenvalues using FEM. We can see that the eigenvalues of the DD prediction are matching with the FEM ones at least up to four decimal places. In Figure~\ref{fig:crossing1stev} we have display the first eigenvector at those four parameters obtained by FEM and DD. We also display the difference between the twp. The error is of order $10^{-4}$. Note that the eigenvector is independent of the parameter $\mu$; for this reason the criterion expressed in~\eqref{criterion} selects only one POD basis for this example.

\begin{figure}
\centering
\includegraphics[width=4.5cm]{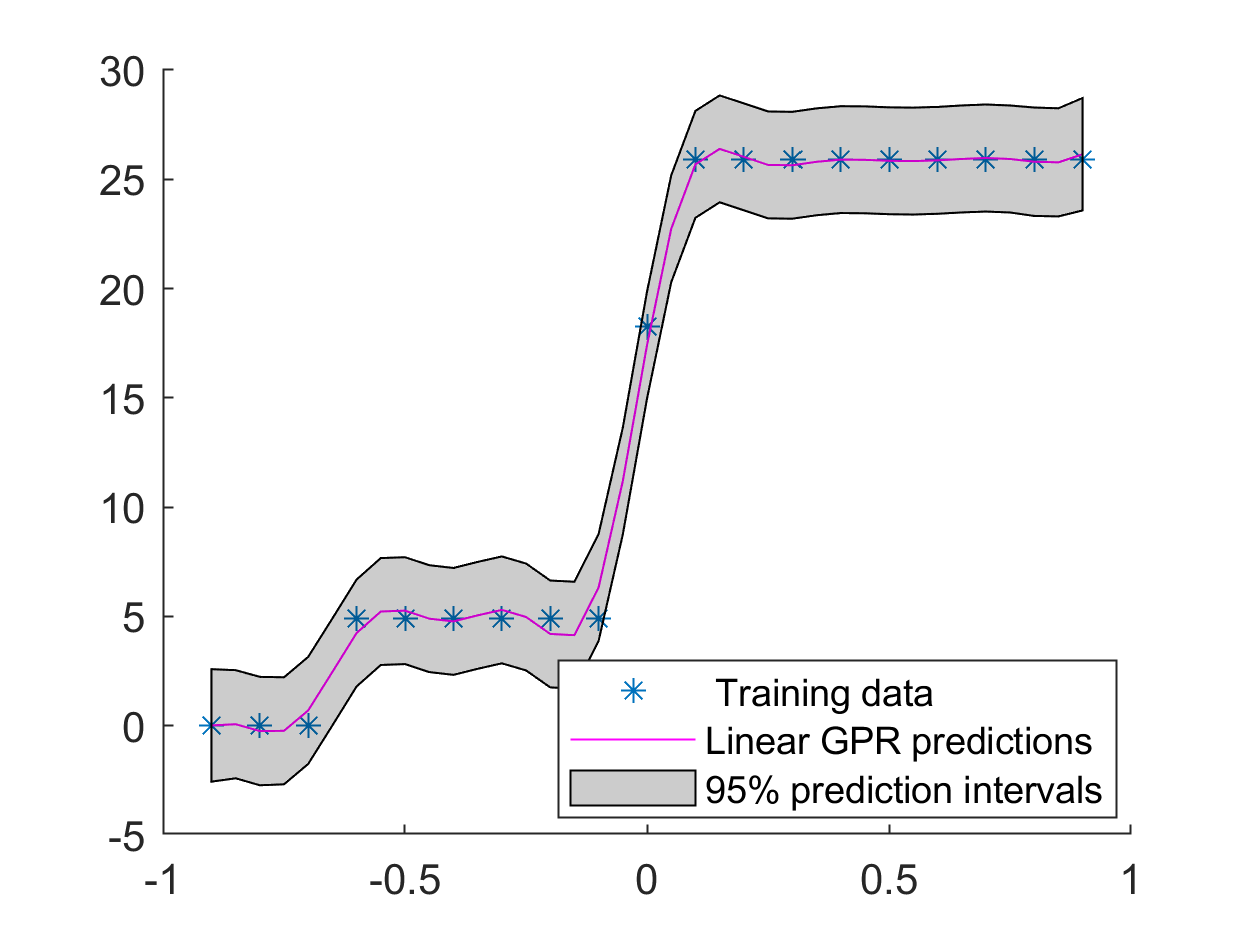}
\includegraphics[width=4.5cm]{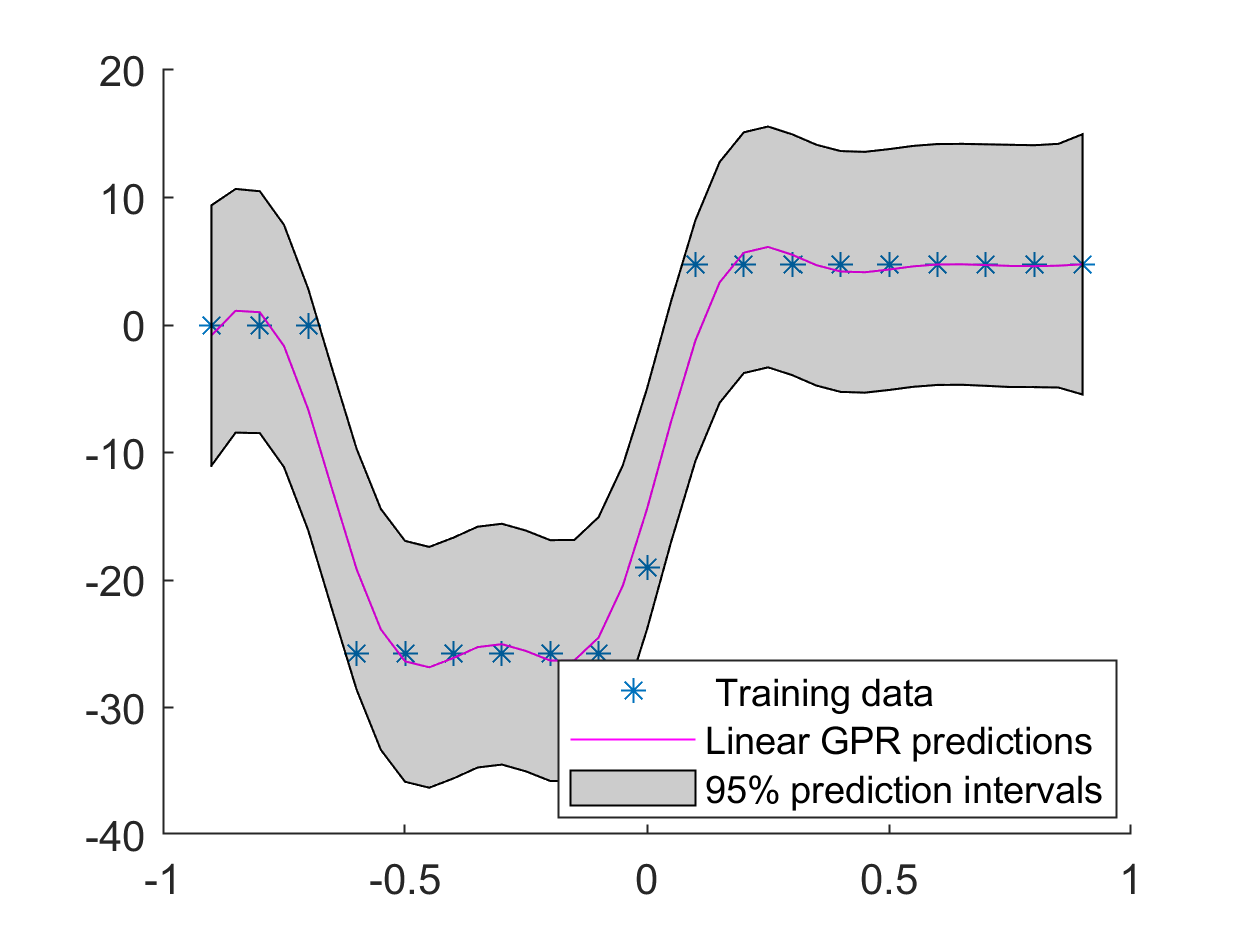}
\includegraphics[width=4.5cm]{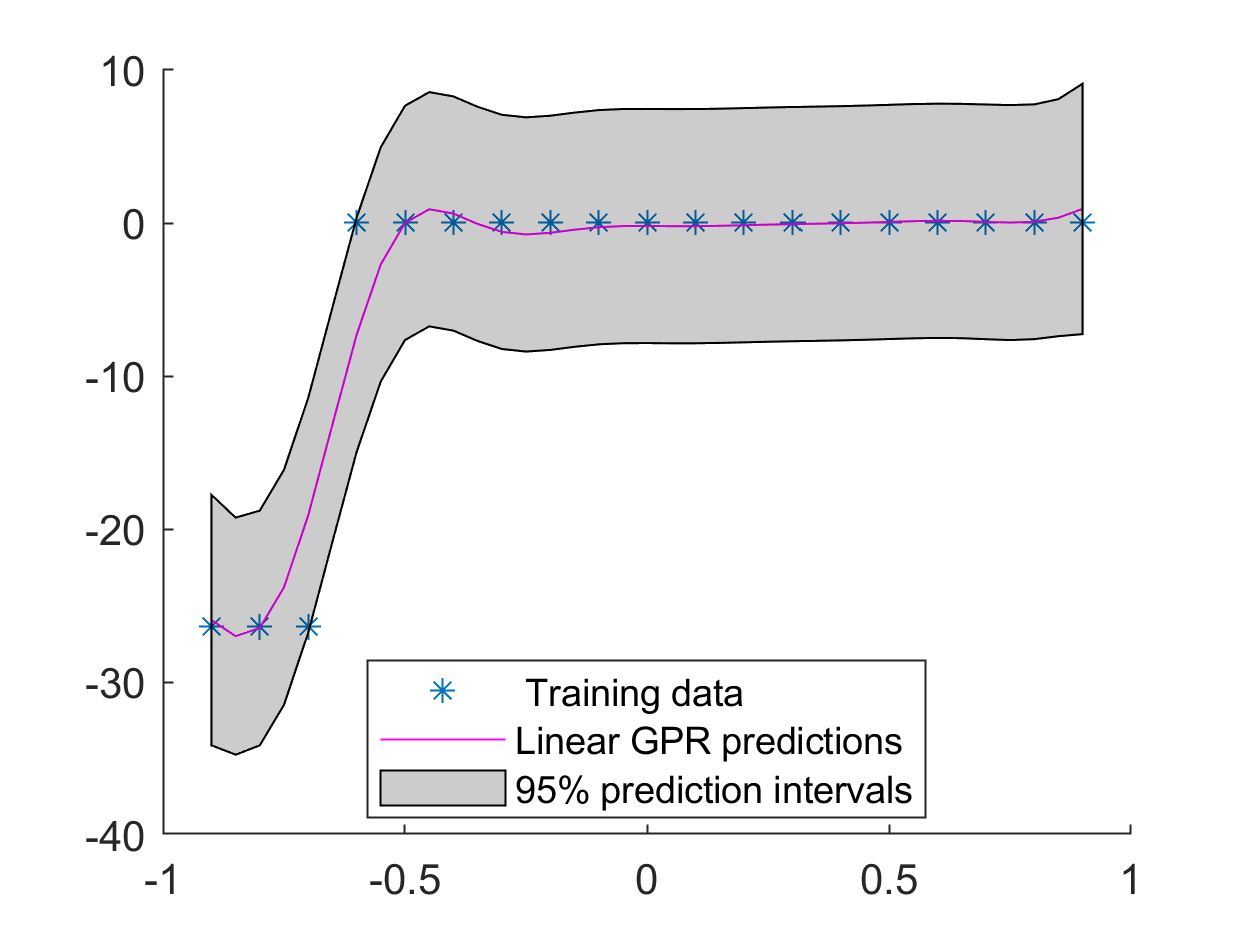}
\caption{Three GPR corresponding to the first three RB coefficients for the 3rd eigenvector of Problem~\eqref{mdl:crossing} with $h=0.05$. The training set is $\mu_{tr}=-0.9:0.1:0.9$.}
\label{figcrossing}
\end{figure}

In order to investigate how our model performs when a crossing occurs, we calculate the third eigenvalues and eigenvectors. The process of constructing the GPR model is the same but, instead of first eigenvalues and first RB eigenvectors, we used the third eigenvalues and third reduced eigenvectors as the training data at the sample points. The training points and the GPR for the three coefficients are presented in Figure~\ref{figcrossing} for the mesh size $h=0.05$. To draw the GPR we use the test points in the interval $[-0.9,0.9]$ with stepsize $0.05$. We show the confidence interval for all the three reduced coefficients. In this case the confidence interval is a little larger than in the previous tests, so the error is higher than for the first eigenvalues and eigenvectors.

\begin{figure}
\centering
\includegraphics[width=4.5cm]{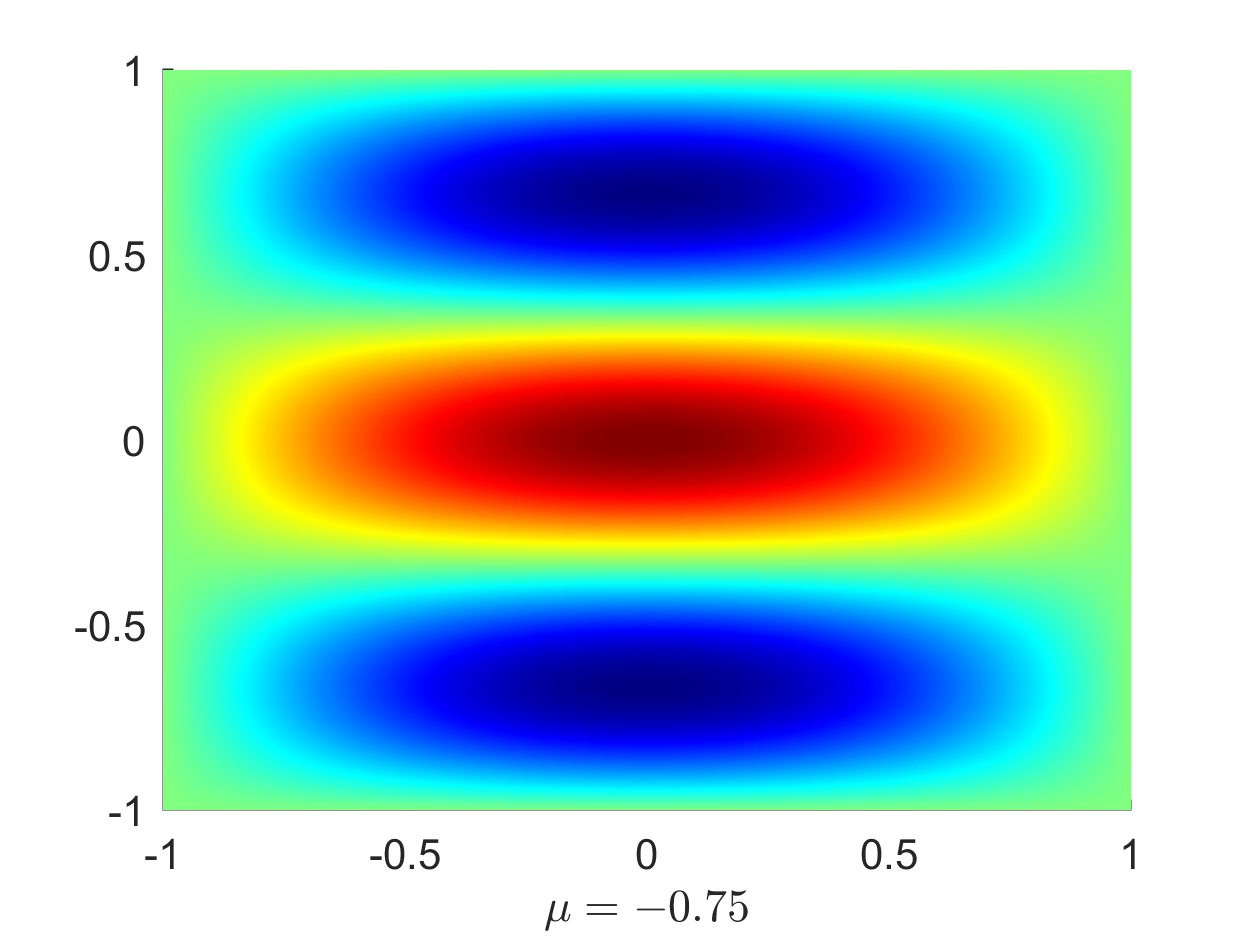}
\includegraphics[width=4.5cm]{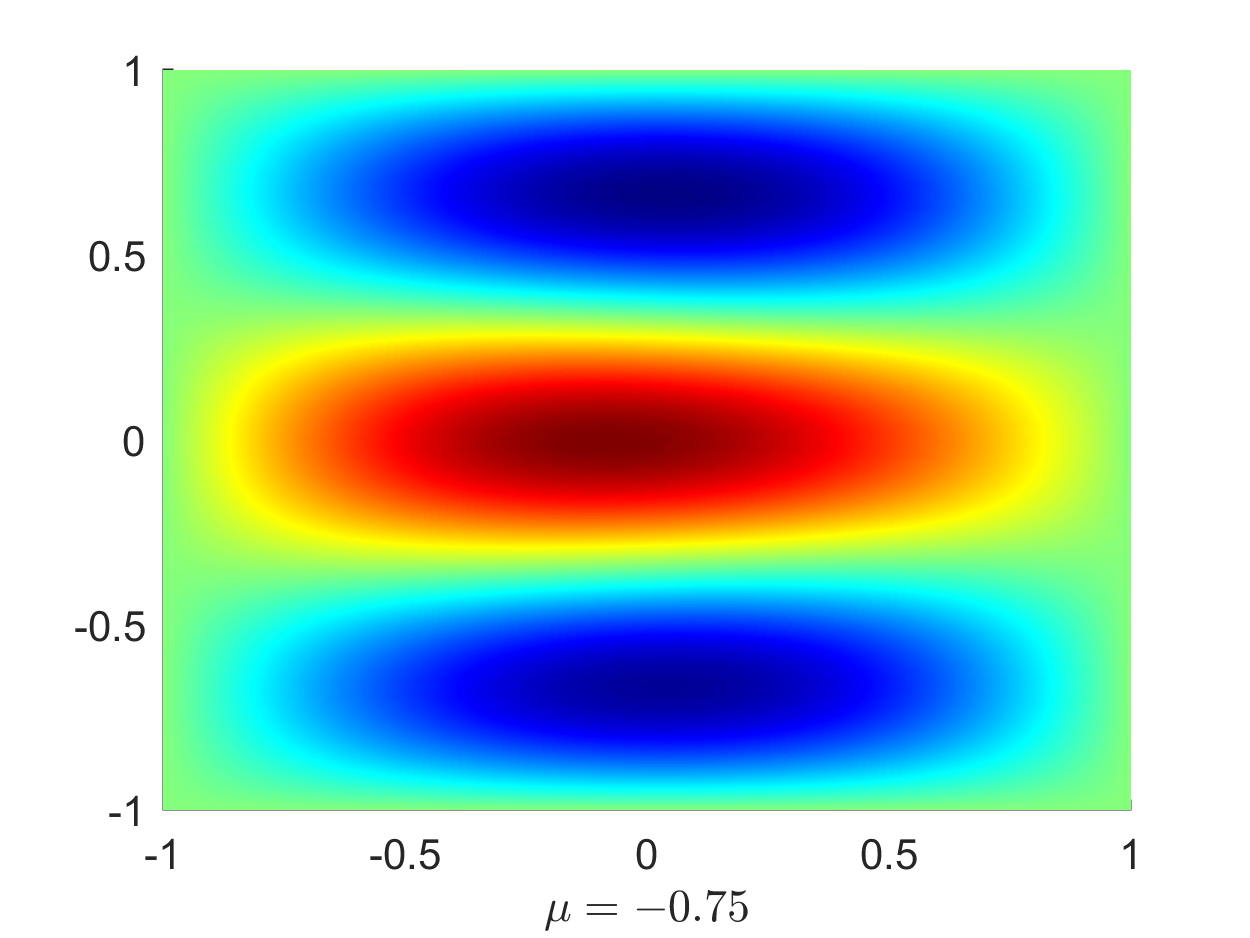}
\includegraphics[width=4.5cm]{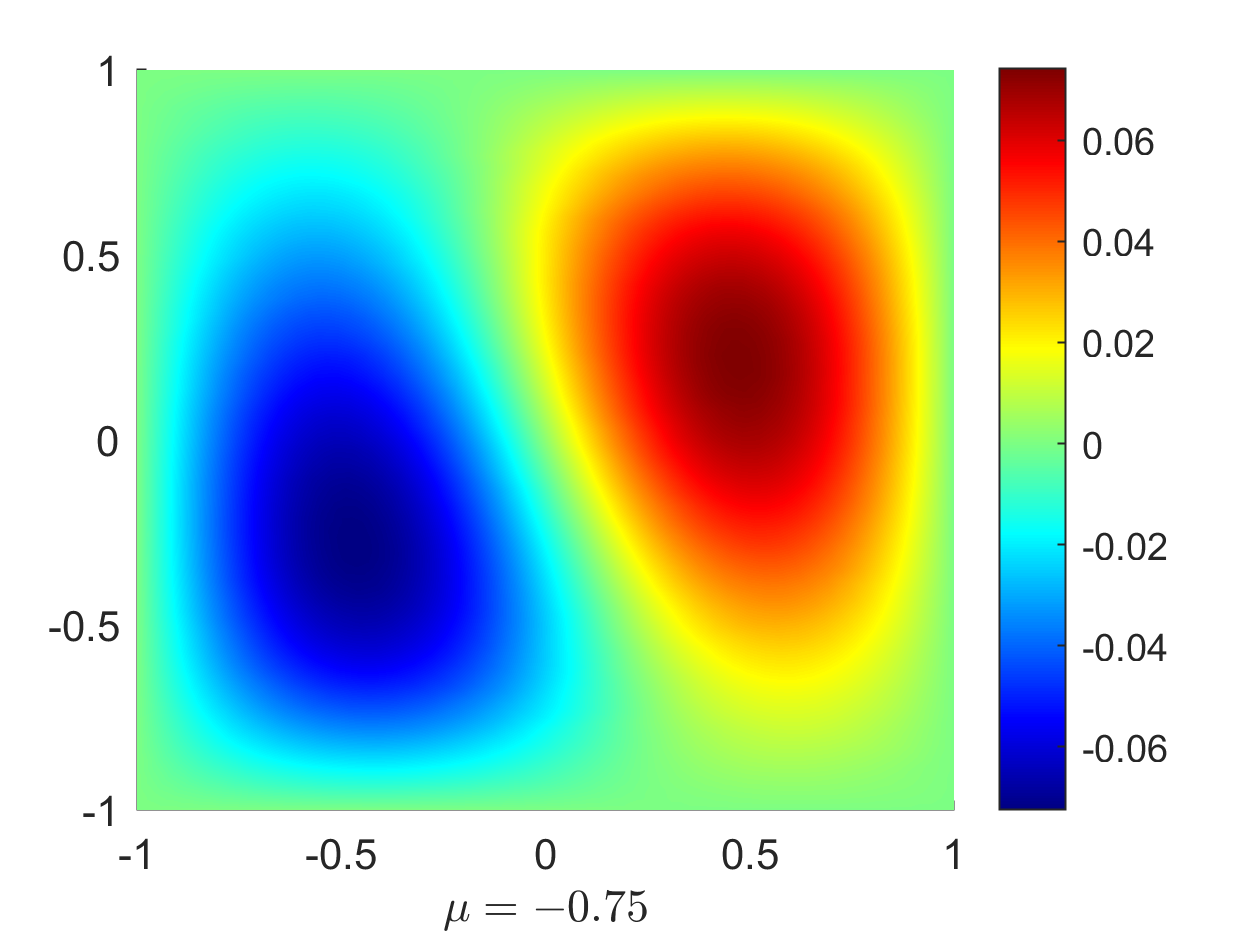}

\includegraphics[width=4.5cm]{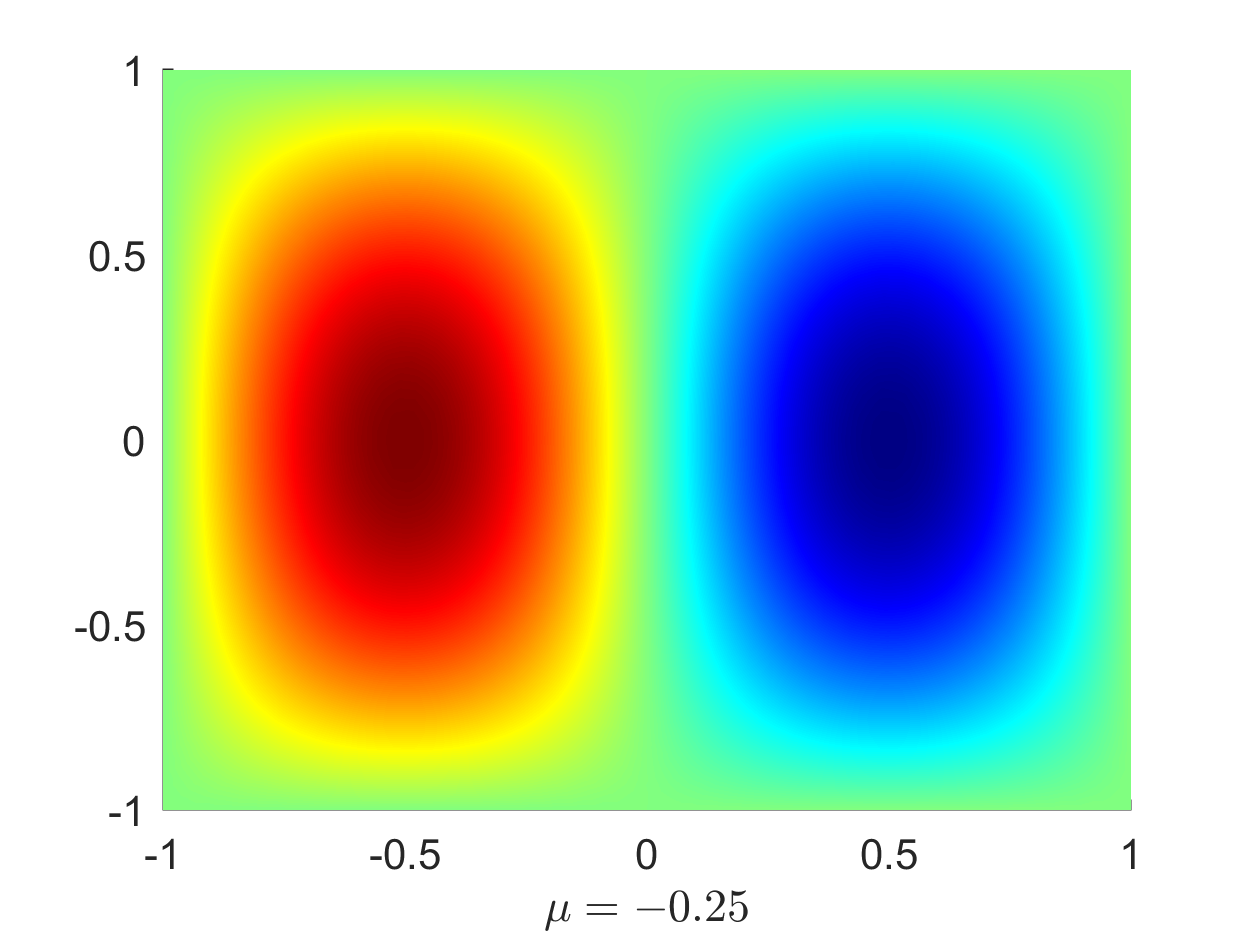}
\includegraphics[width=4.5cm]{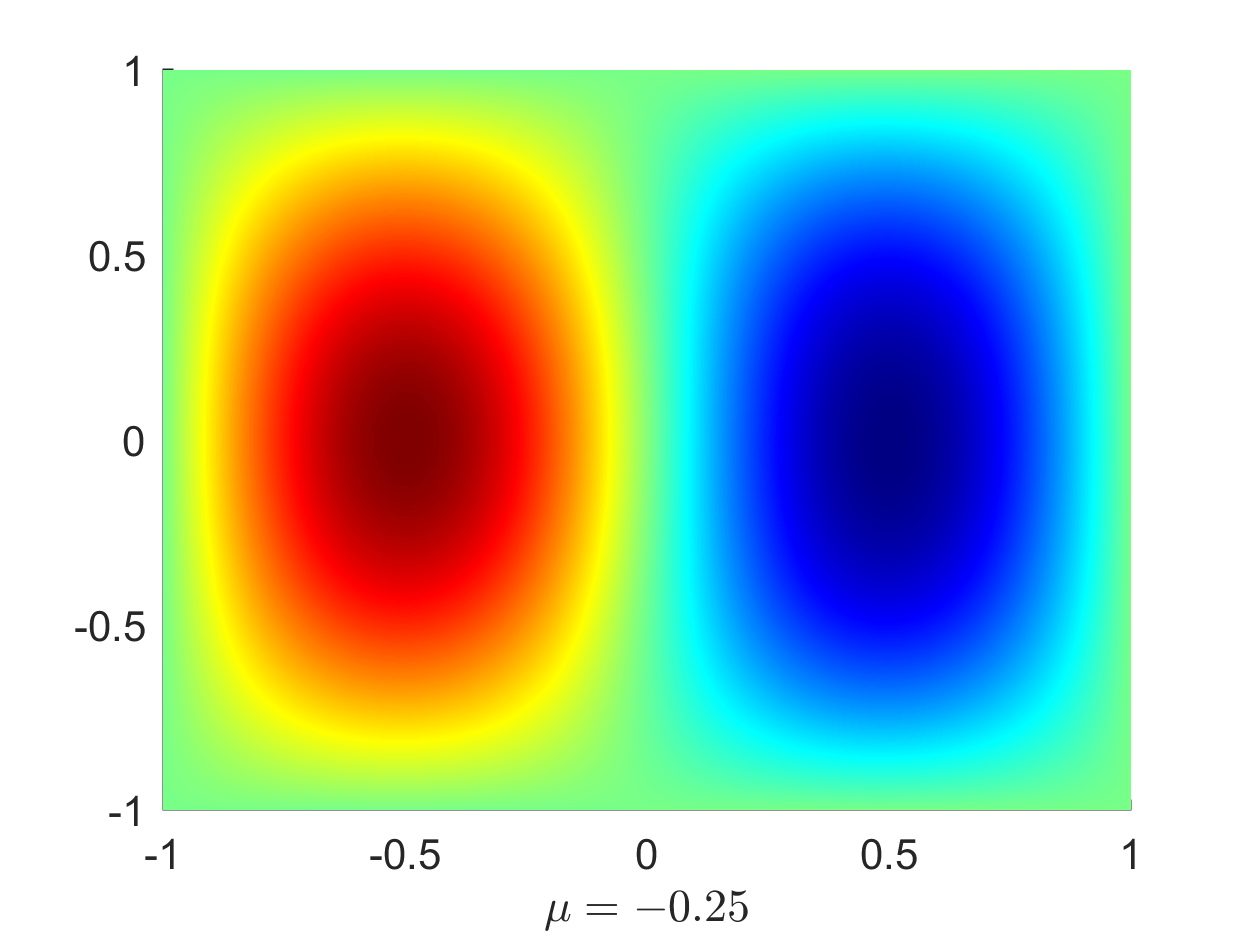}
\includegraphics[width=4.5cm]{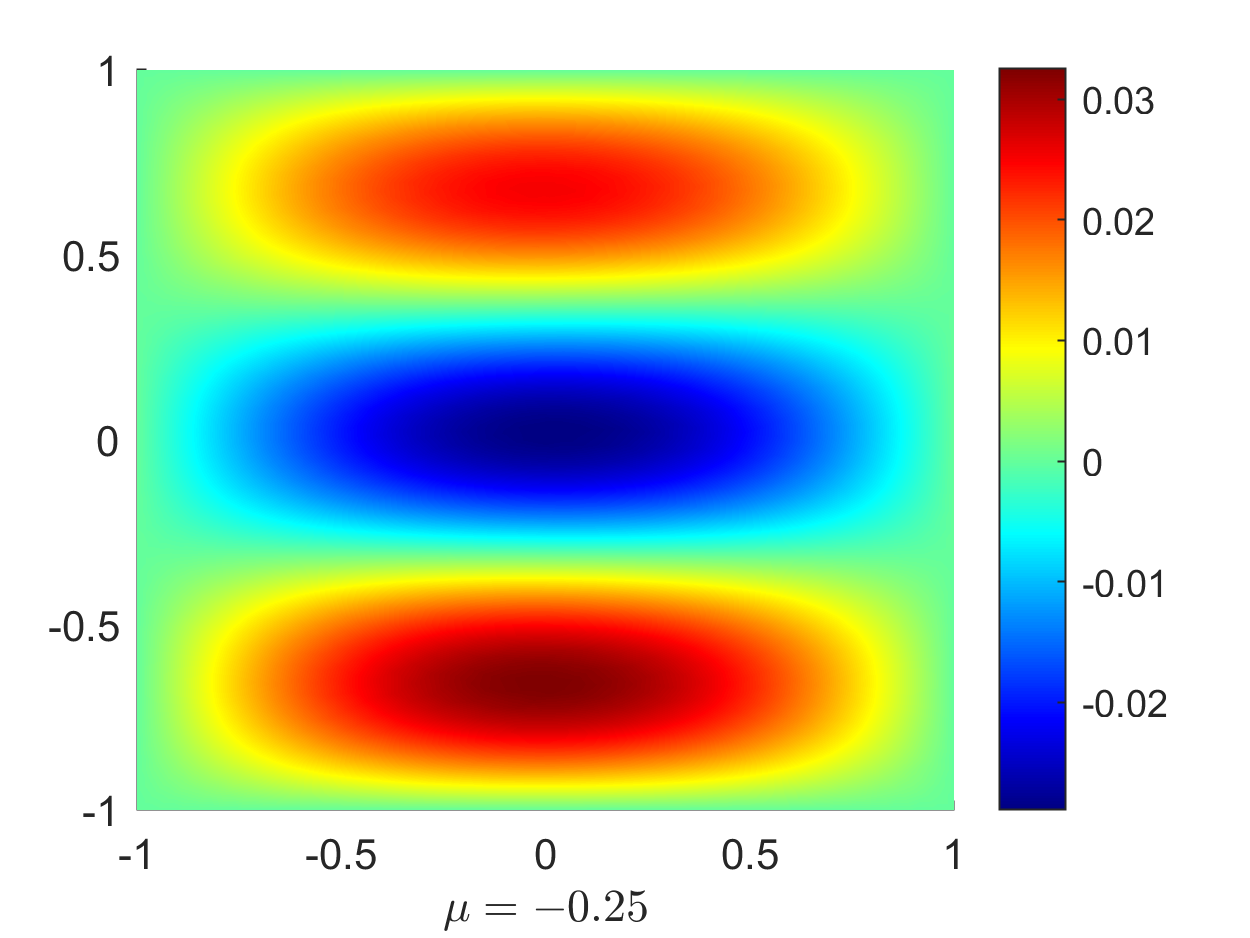}

\includegraphics[width=4.5cm]{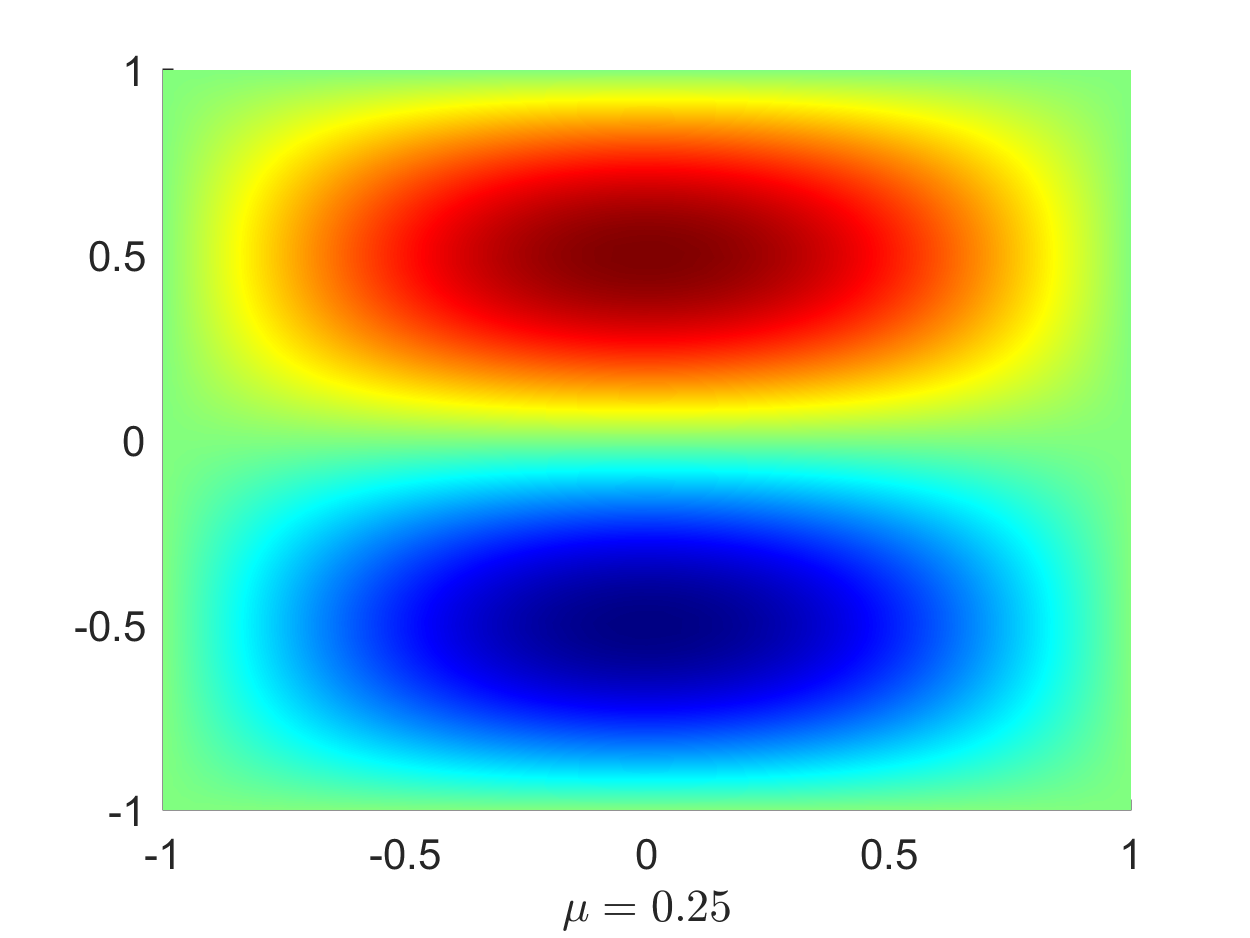}
\includegraphics[width=4.5cm]{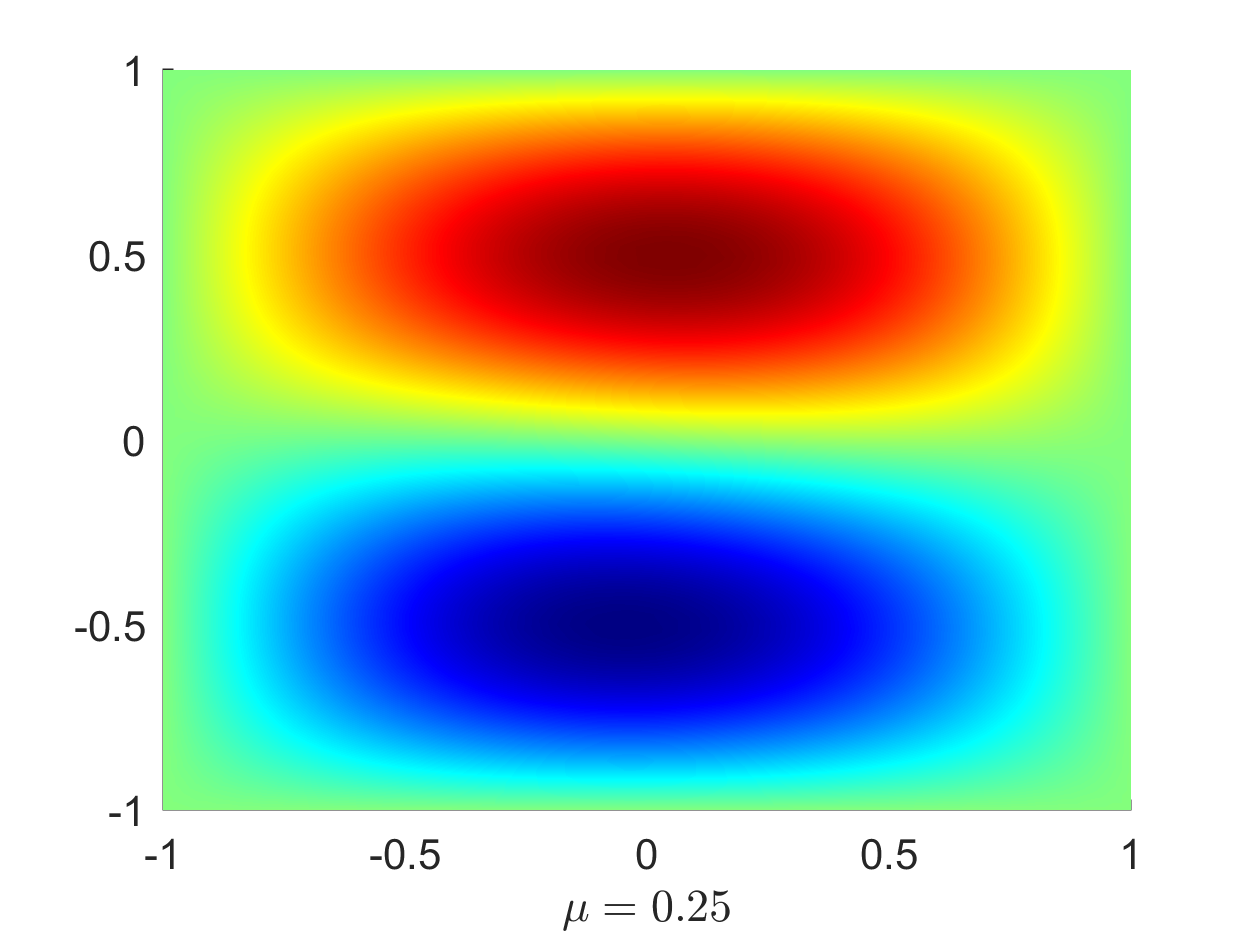}
\includegraphics[width=4.5cm]{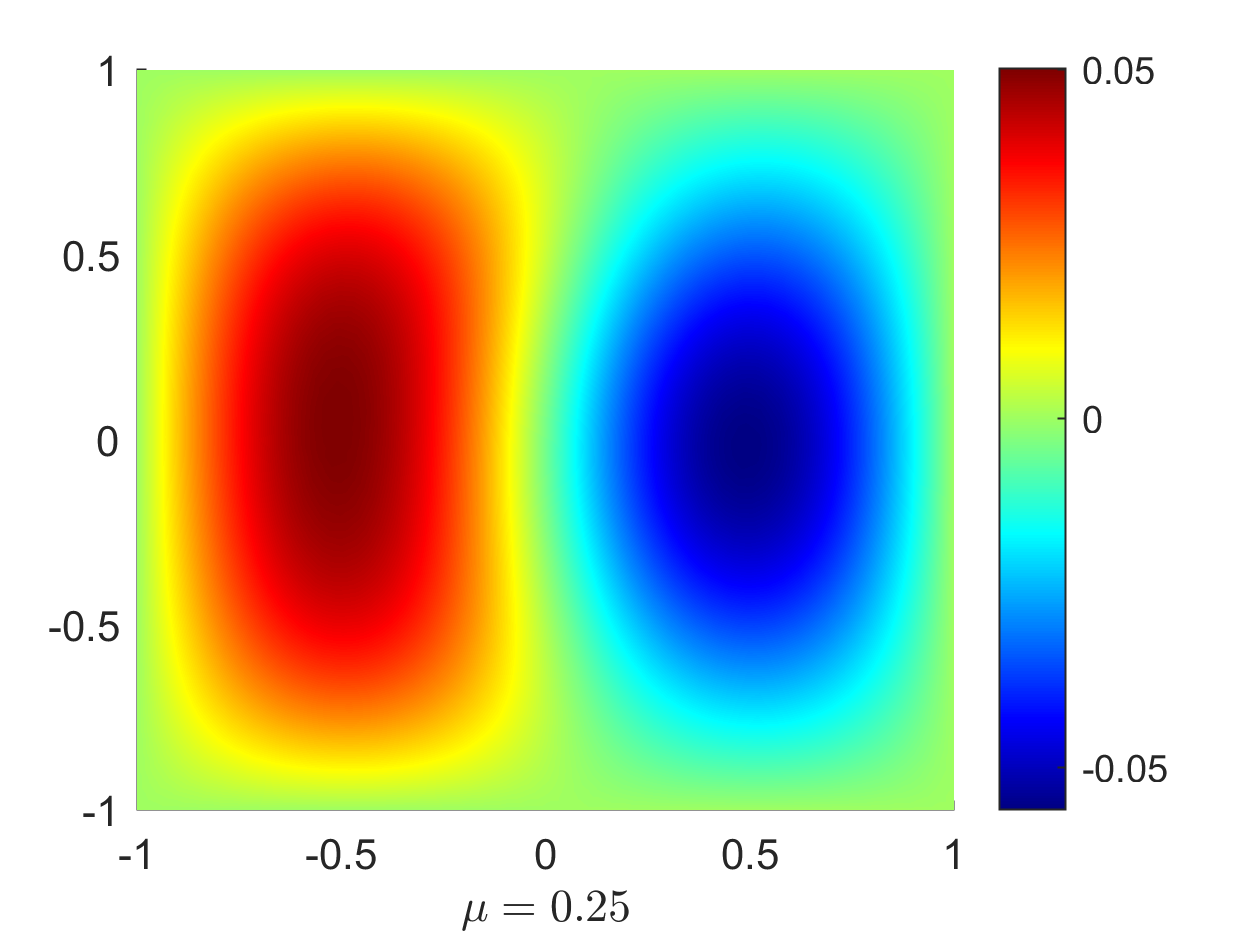}

\includegraphics[width=4.5cm]{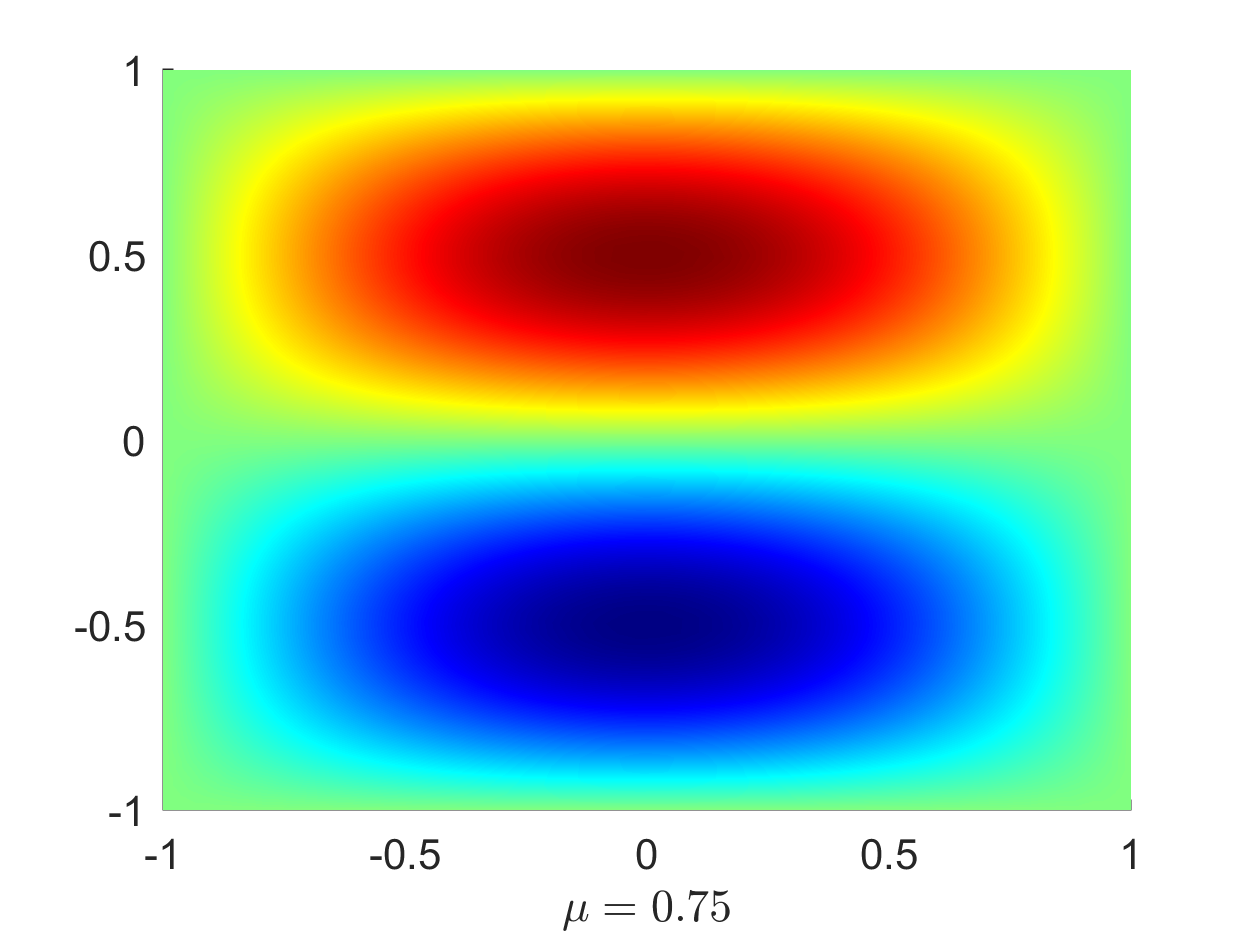}
\includegraphics[width=4.5cm]{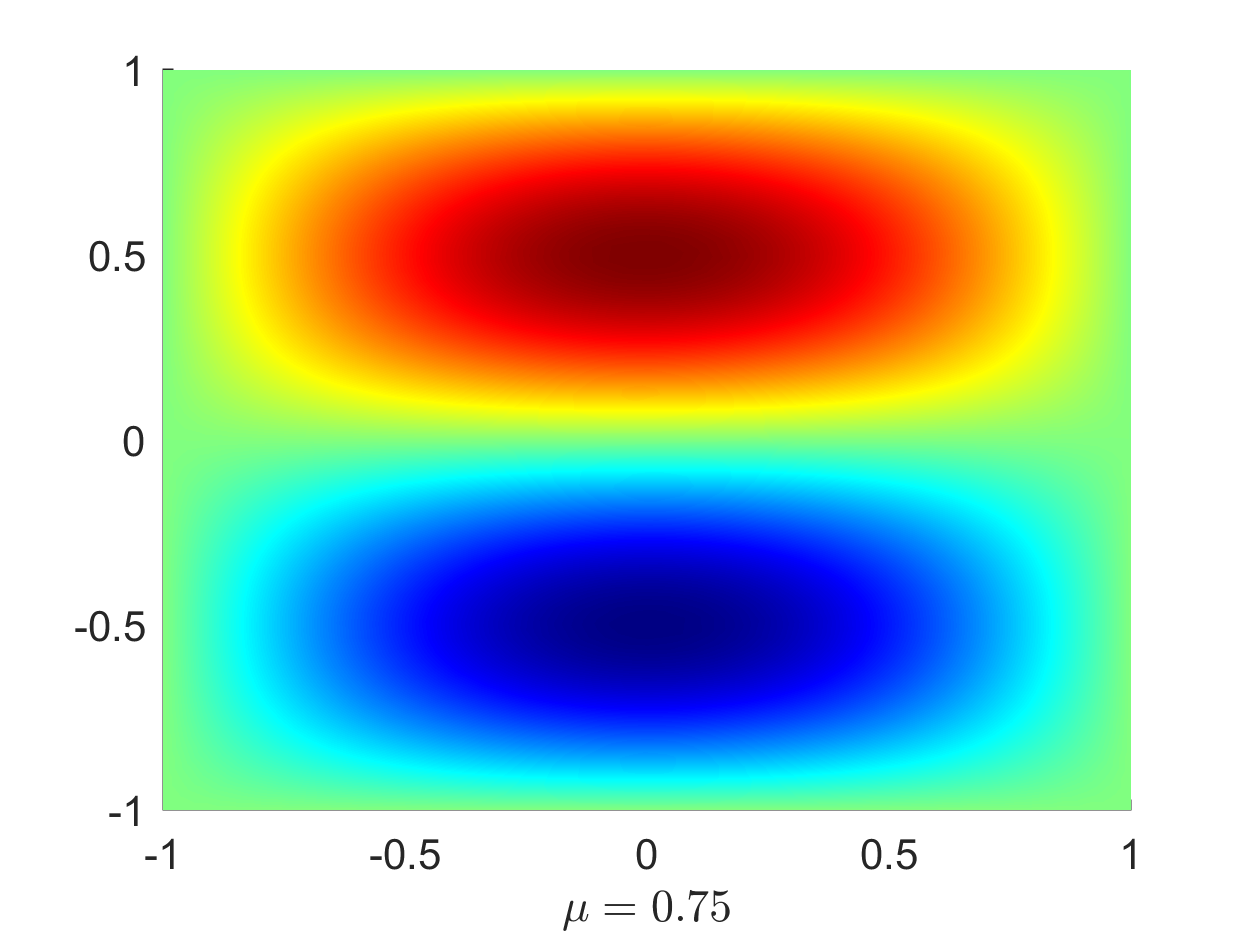}
\includegraphics[width=4.5cm]{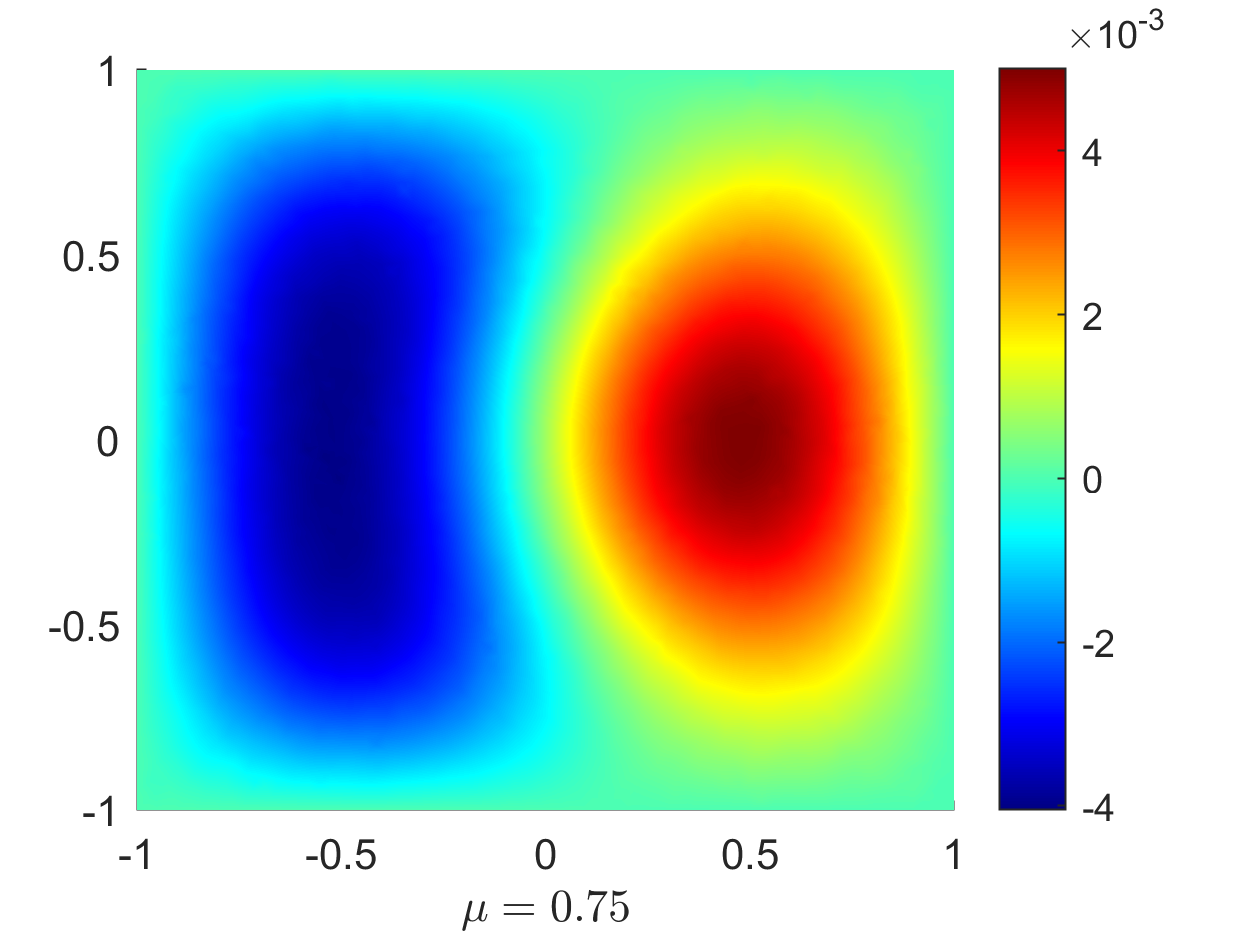}
\caption{3rd eigenvectors corresponding to mesh size $h=0.05$ using FEM (left) and DD model with 4 POD basis (middle) where the training set is $\mu_{tr}=-0.9:0.1:0.9$ and FEM. The difference is reported on the right column.}
\label{figcrossing3rdev}
\end{figure}

In Table~\ref{crossing:3rdev} we report the third eigenvalues obtained by our DD model at the four test points for the mesh size $h=0.1$, $0.05$, and $0.01$, respectively. For comparison we report the eigenvalues using FEM. We can see that the eigenvalues of our DD model are very close to the ones of the FEM method. In Figure~\ref{figcrossing3rdev} we display the third eigenvector at the four test parameters obtained by FEM and DD models. We display the difference between the two. The error is of order $10^{-2}$. Note that, although the eigenvector is independent of the parameter $\mu$, three different eigenvectors show up in the snapshot matrix. This is due to the crossings and it is the reason why the criterion expressed in~\eqref{criterion} gives three POD basis for calculating third eigenvector.

\begin{table}
\footnotesize
 	 	\centering
 	\begin{tabular}{|c|c|c|c|c|c|c|c|} 
 		\hline
 		h &  {\begin{tabular}[c]{@{}c@{}} Method \end{tabular}} &
		  {\begin{tabular}[c]{@{}c@{}} $\mu=-0.75$ \end{tabular}} &
 		 {\begin{tabular}[c]{@{}c@{}} $\mu=-0.25$ \end{tabular}} & 
 		 {\begin{tabular}[c]{@{}c@{}} $\mu=0.25$ \end{tabular}} &
 		{\begin{tabular}[c]{@{}c@{}}   $\mu=0.75$  \end{tabular}}\\
 \hline
 0.1&-FEM&8.14338931 &11.78888922&14.89196477&19.85303433\\
    &DD&8.21778394&11.81042527&14.91670069&19.85137938\\
  \hline
  0.05&FEM&8.05008647 &11.73700317 &14.82575077&19.76686171\\
      & DD&8.12327908&11.75458889&14.85046418&19.76559935\\
   \hline
    0.01&FEM&8.02024667 &11.72081569&14.80523805&19.74028492\\
        &DD&8.09313087&11.73711527 &14.82993019	&19.73915595\\         
  \hline
 	\end{tabular}
	\caption{3rd eigenvalues of \eqref{mdl:crossing} for different mesh using DD model with sample points $\mu_{tr}=-0.9:0.1:0.9$ and FEM model.}
	 	\label{crossing:3rdev}
\end{table}

\begin{figure}
\centering
\includegraphics[width=4.5cm]{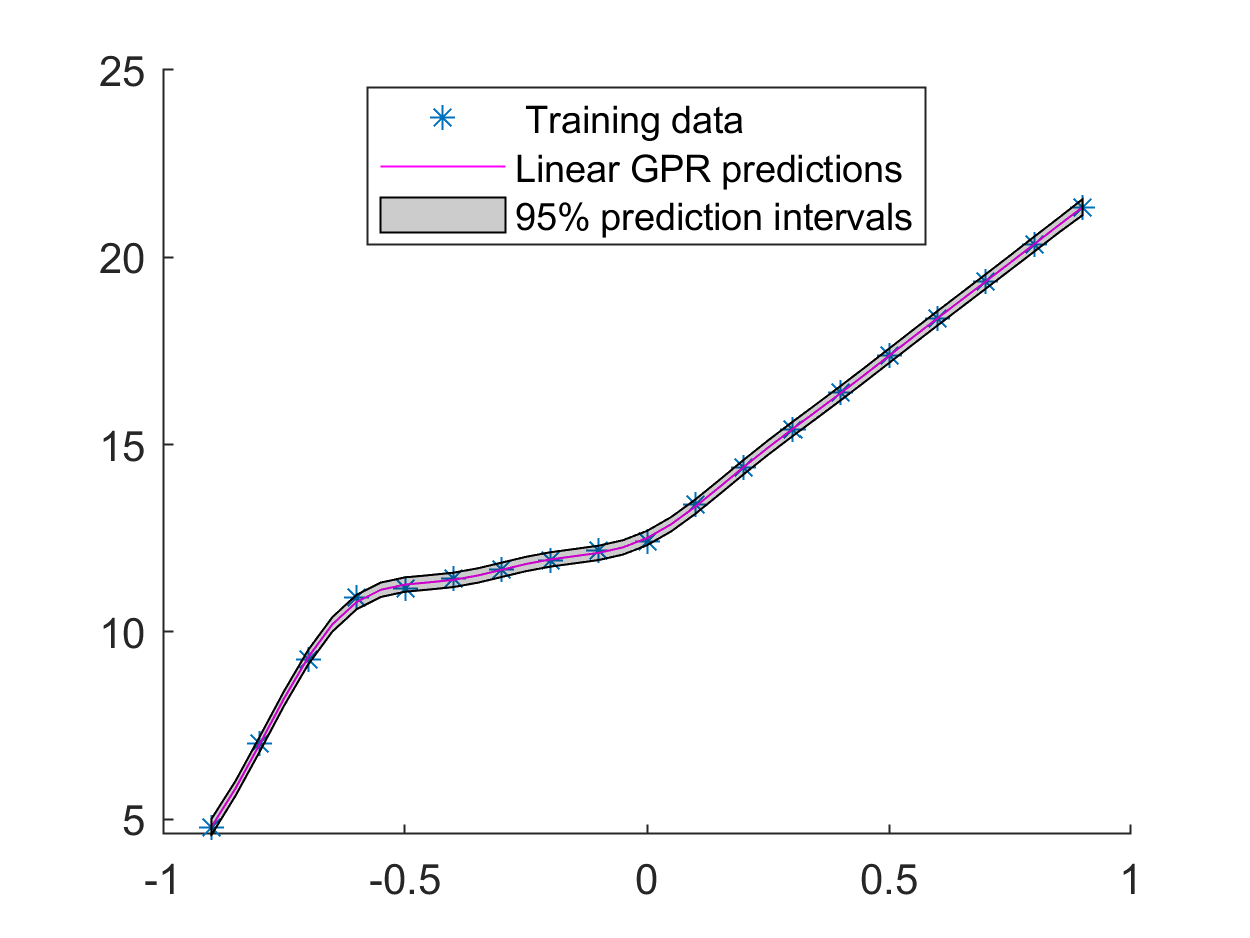}
\includegraphics[width=4.5cm]{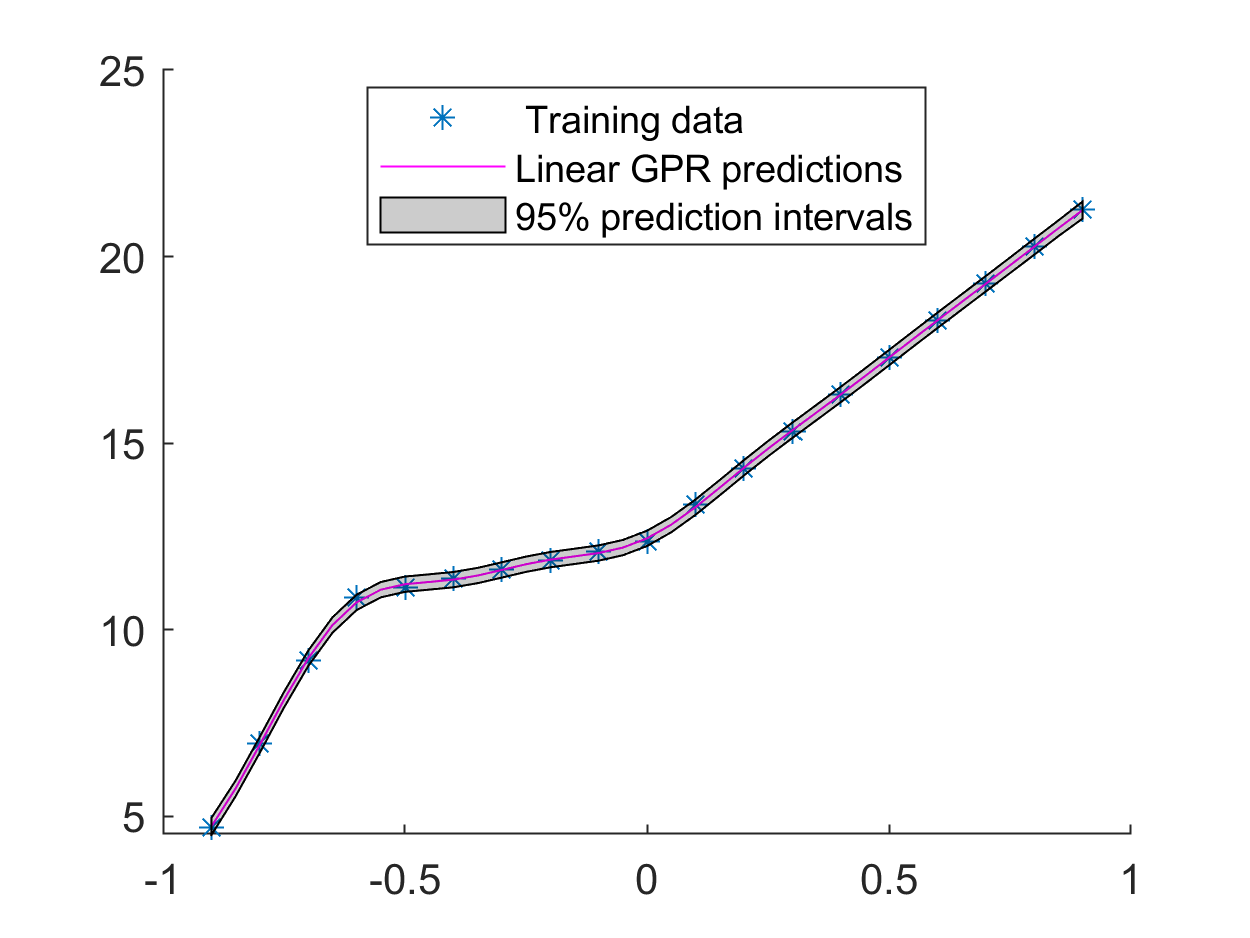}
\includegraphics[width=4.5cm]{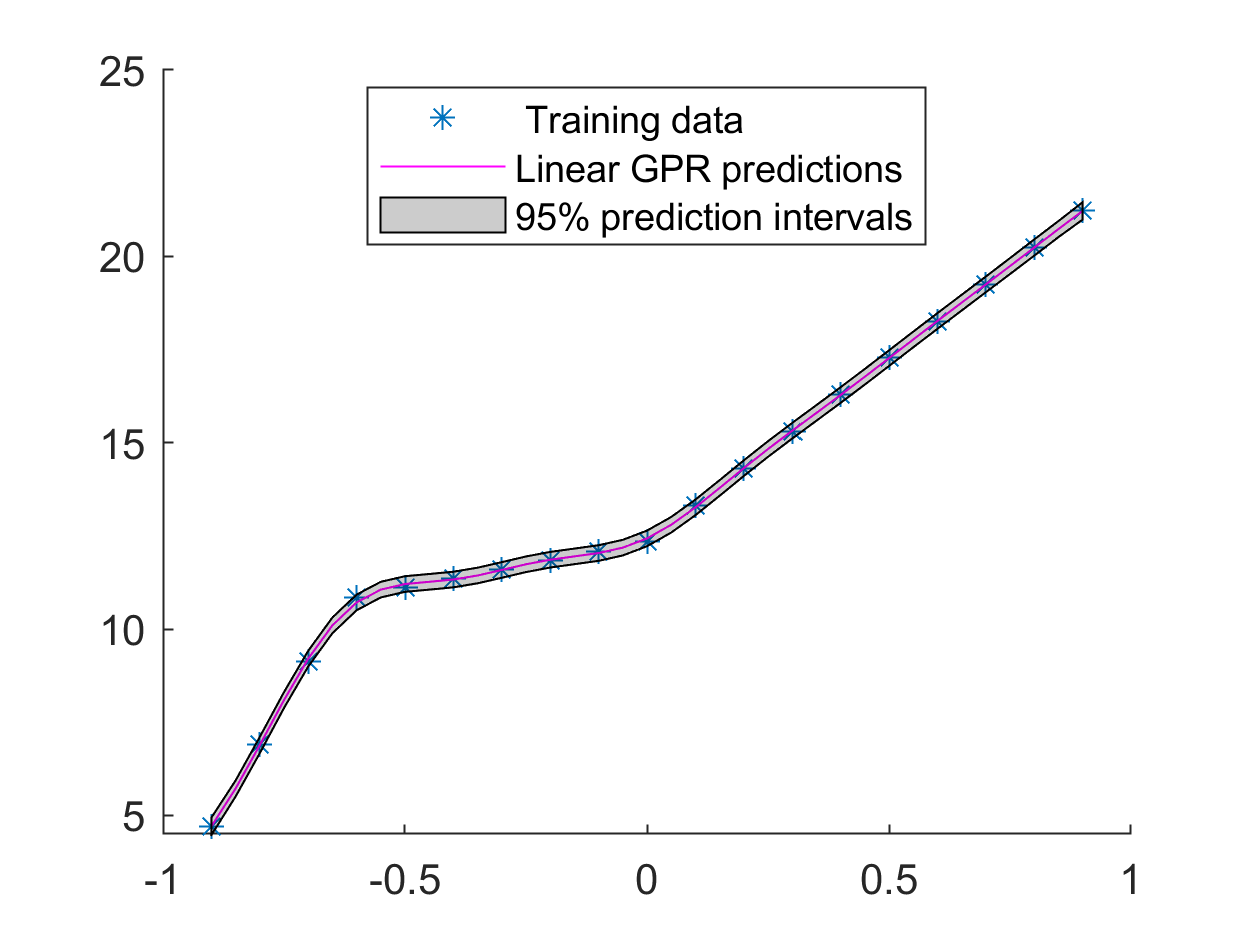}
\caption{The GPR corresponding to 3rd eigenvalues of \eqref{mdl:crossing} with mesh size $h=0.1,0.05$ and $0.01$ respectively. The training set is $\mu_{tr}=-0.9:0.1:0.9$.}
\label{fig:3rdcrossing}
\end{figure}

In Figure~\ref{fig:3rdcrossing} we plot the 3rd eigenvalues used for the training with star symbols and we plot in magenta color the mean of the GPR using the eigenvalues at the test points in the interval $[-0.9,0.9]$ with stepsize $0.05$, corresponding to mesh size $h=0.1$, $0.05$, and $0.01$. We show the $95$ percentage confidence interval for the 3rd eigenvalues. The confidence interval is very narrow, that means the eigenvalues we obtained using the DD model are matching well with the FEM ones.

\begin{figure}
\centering
\includegraphics[width=4.5cm]{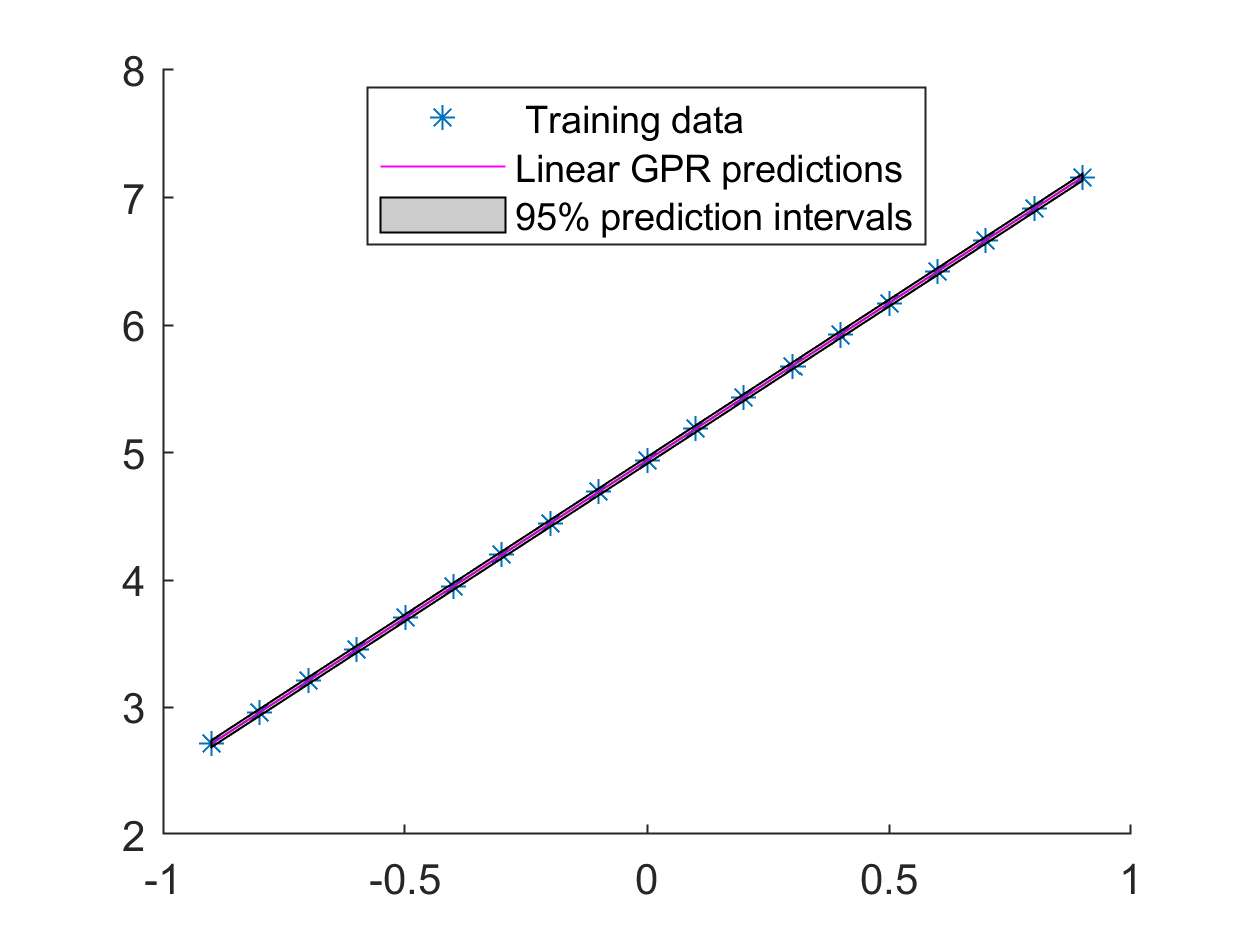}
\includegraphics[width=4.5cm]{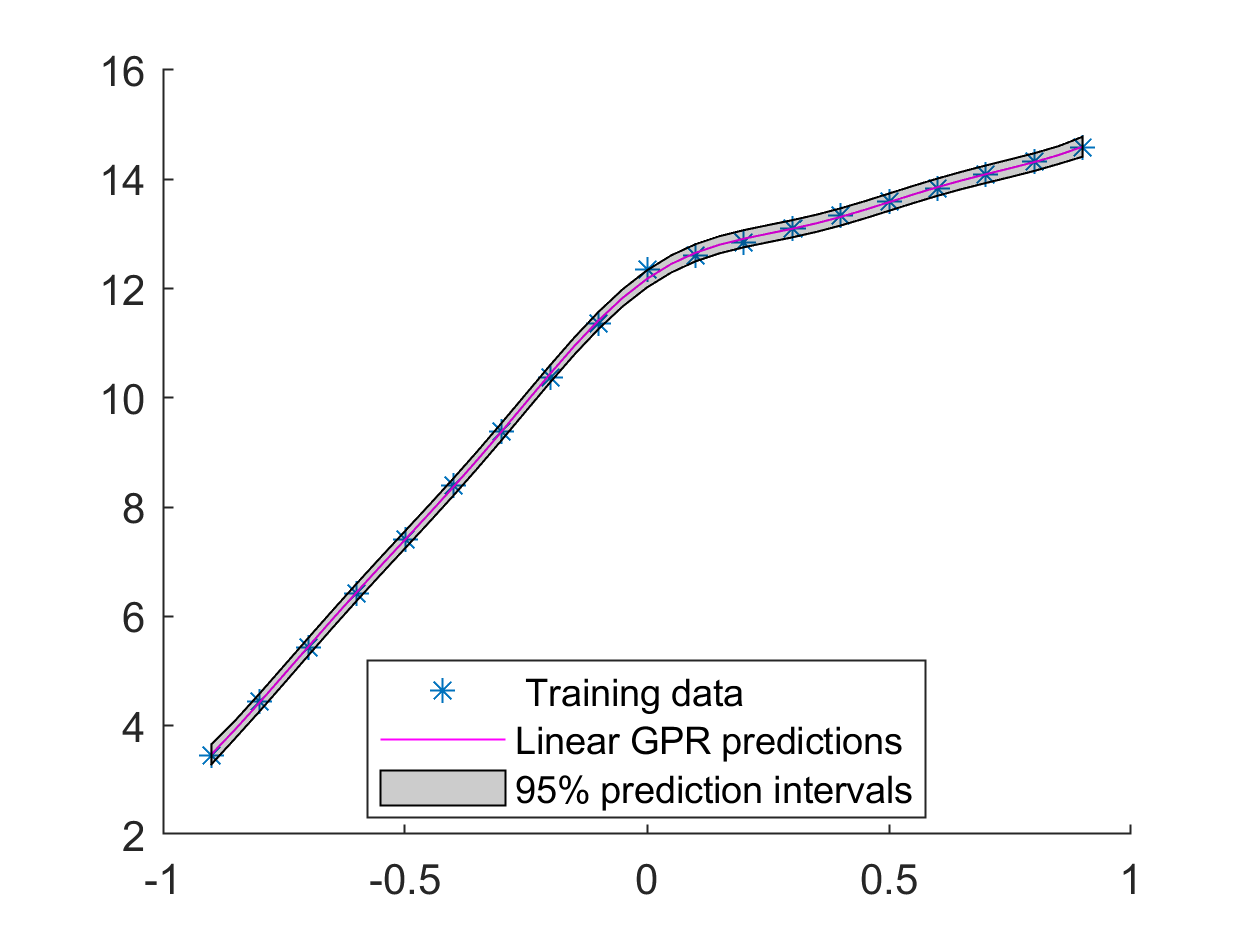}
\includegraphics[width=4.5cm]{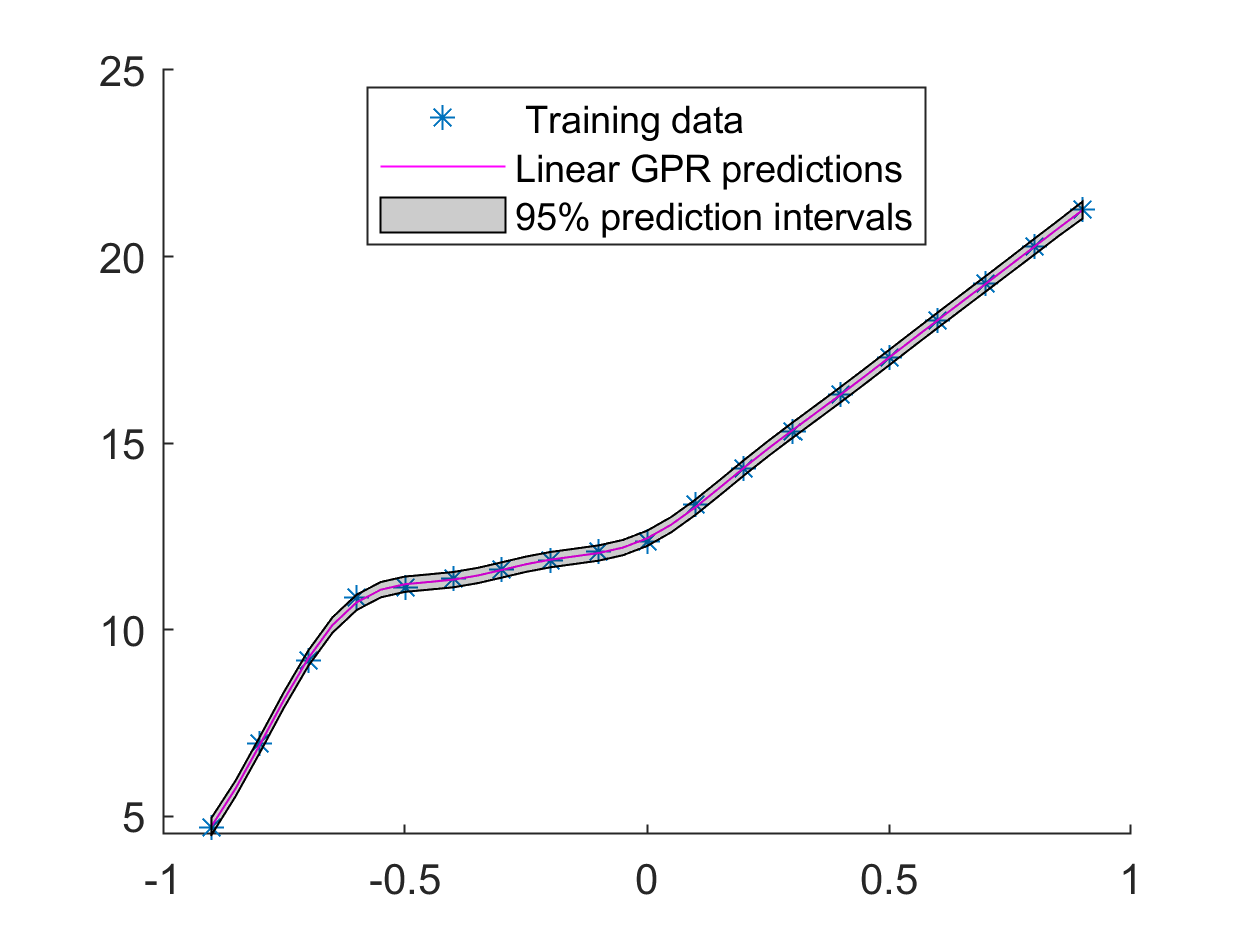}
\caption{The GPR corresponding to first three eigenvalues of \eqref{mdl:crossing} with mesh size $h=0.05$. The training set is $\mu_{tr}=-0.9:0.1:0.9$.}
\label{fig:ev_gpr}
\end{figure}

Finally, we have tested our DD model when multiple eigenfunctions are calculated simultaneously, by using the procedure described in Subsection~\ref{sec:simul}. We have computed first three eigenvalues and eigenvectors simultaneously for Problem~\eqref{mdl:crossing}. In Figure~\ref{fig:ev_gpr} we show the plot for the mean of the GPR and the training set, corresponding to the first three eigenvalues, with $95$ percentage of confidence interval. We observed that the training data are lying on the curve of the mean GPR for all the three cases. The confidence interval for the first eigenvalue is very narrow and for the other two it is a little larger than for the first one. All the three eigenvalues obtained by the DD model at the test points $\mu=-0.75$, $-0.25$, $0.25$, $0.75$, corresponding to mesh size $h=0,05$, are reported in Table~\ref{crossing:allev} and are compared with the FEM result. We can see that the first eigenvalues obtained by the DD model are more accurate than the second and third ones.

\begin{table}
\footnotesize
 	 	\centering
 	\begin{tabular}{|c|c|c|c|c|c|c|c|} 
 		\hline
 		 {\begin{tabular}[c]{@{}c@{}}   $\mu$  \end{tabular}}&
 		 {\begin{tabular}[c]{@{}c@{}} Method \end{tabular}} &
		  {\begin{tabular}[c]{@{}c@{}} 1st eigenvalue \end{tabular}} &
 		 {\begin{tabular}[c]{@{}c@{}}  2nd eigenvalue \end{tabular}} & 
 		 {\begin{tabular}[c]{@{}c@{}} 3rd eigenvalue  \end{tabular}} \\
 		 \hline
-0.75& FEM&3.08606437 &4.94333422 &8.05008647 \\
 &DD&3.08606913 &4.93514768 &8.12327908\\
 \hline
-0.25& FEM&4.32052203 &9.88459549 &11.73700317 \\
&DD& 4.32051579& 9.91536968& 11.75458889\\
 \hline
0.25& FEM &5.55496589 &12.97341757& 14.82575077\\
 &DD& 5.55496244 &13.00419699 &14.85046418\\
\hline
0.75& FEM &6.78940395 &14.20979040 &19.76686171\\
& DD& 6.78940910& 14.20160725 & 19.76559935\\
  \hline
 	\end{tabular}
	\caption{First three eigenvalues of \eqref{mdl:crossing} for different mesh $h=0.05$ using DD model with sample points $\mu_{tr}=-0.9:0.1:0.9$ and FEM model.}
\label{crossing:allev}
\end{table}

\begin{figure}
\centering
\subcaptionbox{GPR for 1st coefficient}{
\includegraphics[width=6.5cm]{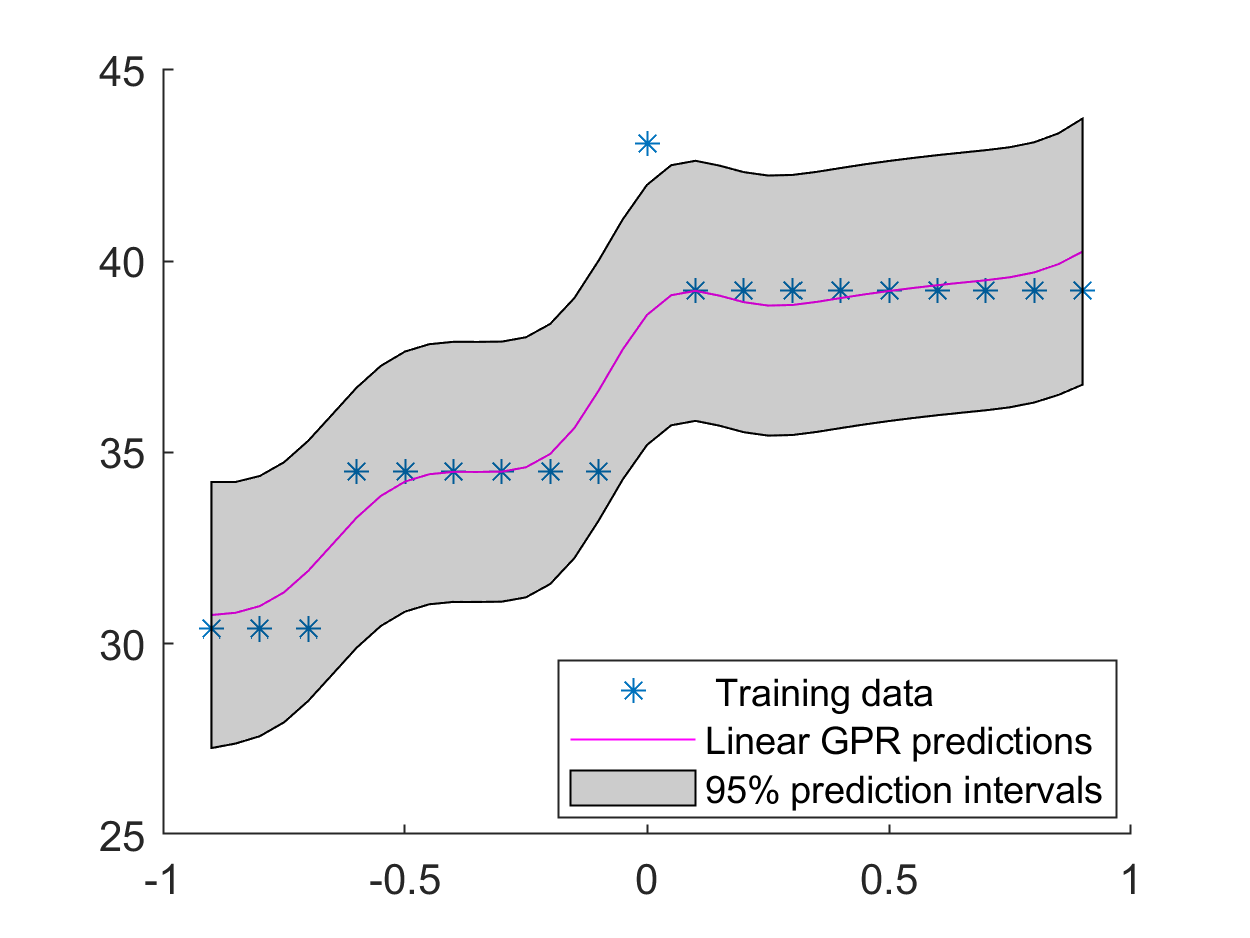}}
\subcaptionbox{GPR for 2nd coefficient}{
\includegraphics[width=6.5cm]{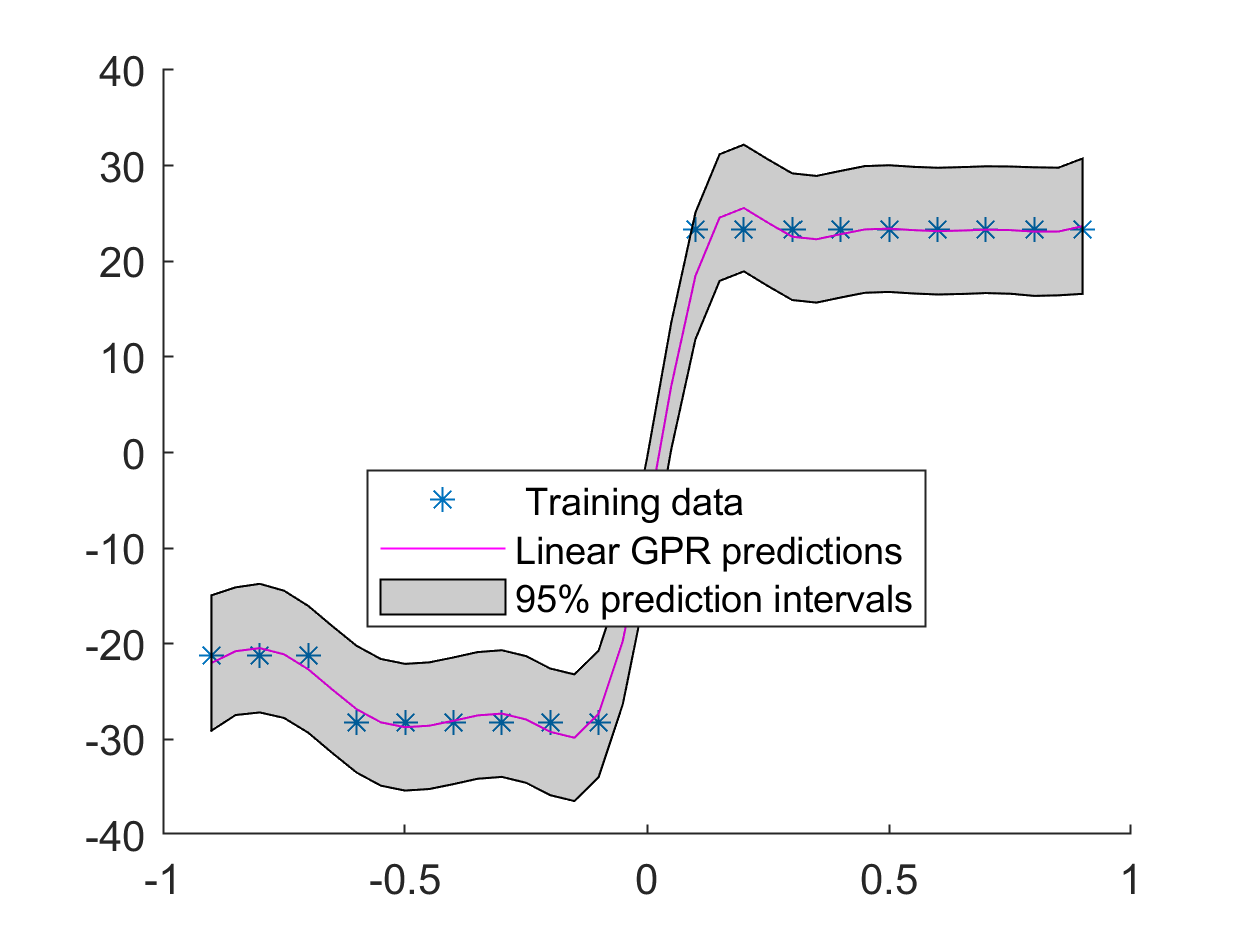}}

\subcaptionbox{GPR for 3rd coefficient}{
\includegraphics[width=6.5cm]{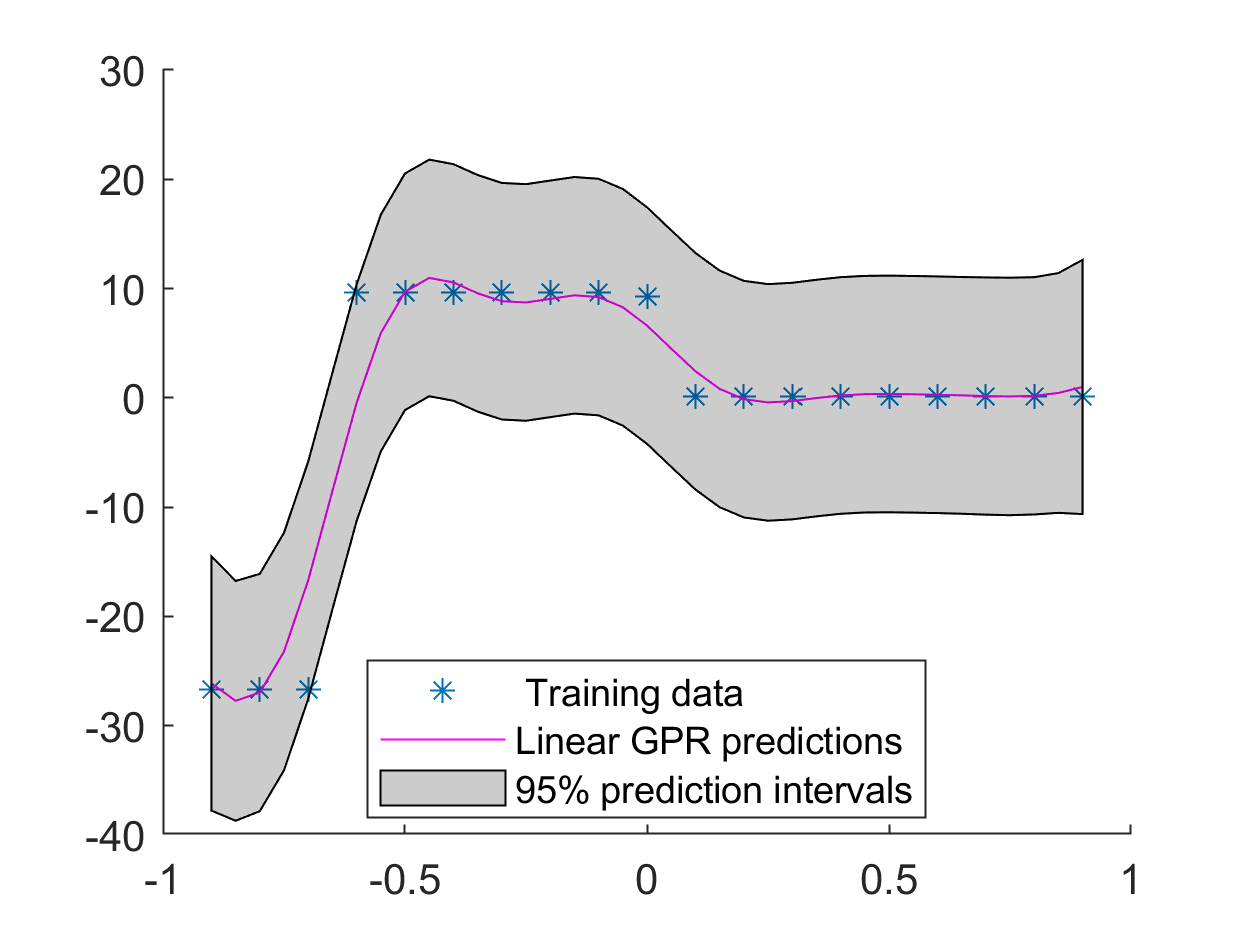}}
\subcaptionbox{GPR for 4th coefficient}{
\includegraphics[width=6.5cm]{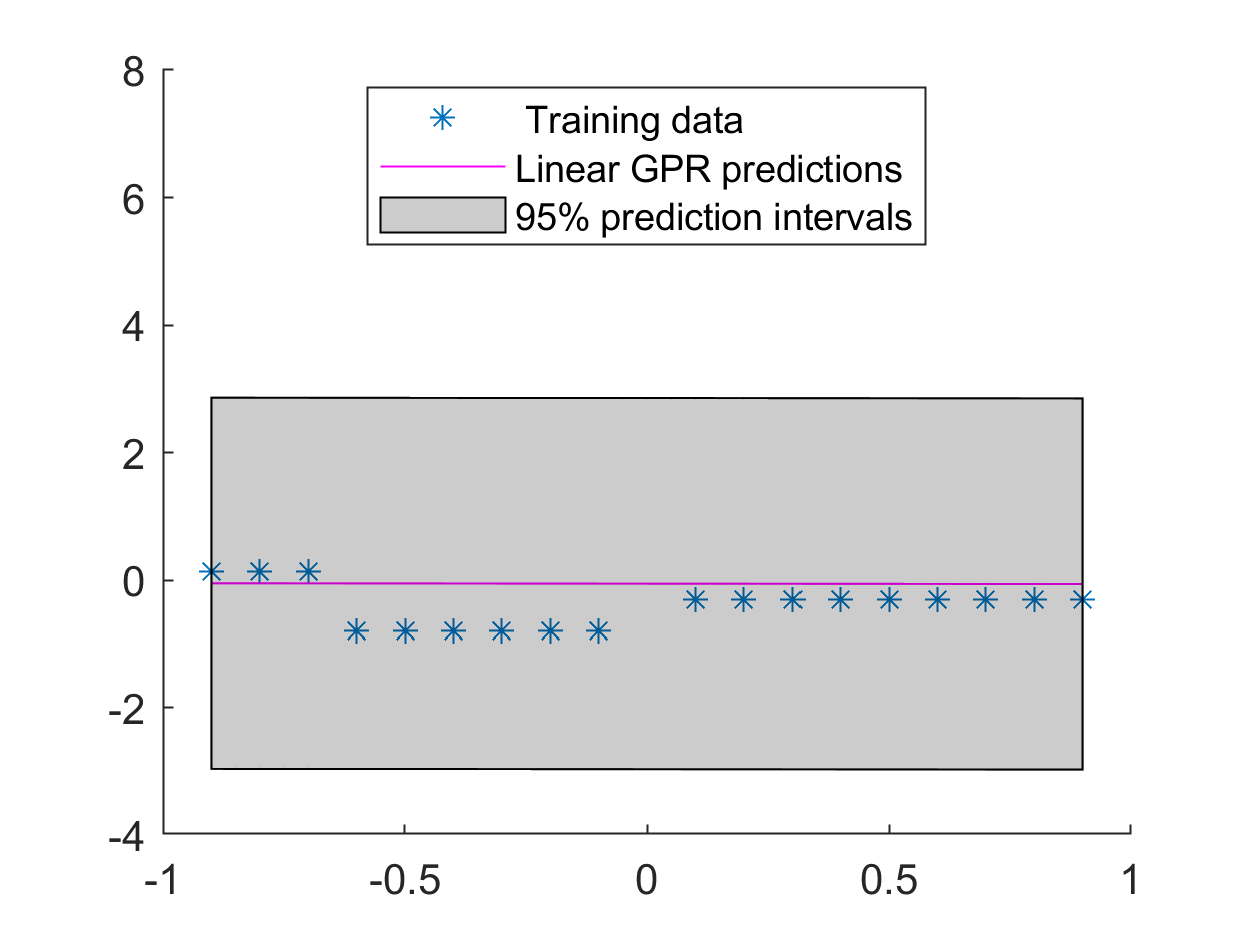}}
\caption{The GPR corresponding to the four projected coefficients \eqref{mdl:crossing} with mesh size $h=0.05$ when first three eigenvectors are considered simultaneously. The training set is $\mu_{tr}=-0.9:0.1:0.9$.}
\label{fig:si_ev_coff}
\end{figure}

In Figure~\ref{fig:si_ev_coff}, we presented the mean GPR corresponding to the four reduced coefficients and the training data for the simultaneous approximation case. We show the $95$ percentage of confidence interval for each GPR. The GPR for the first and third coefficients are a little large. Note that for the DD model the sign of the snapshots is very important. We have chosen the sign of the snapshots following the algorithm of~\cite{Pau07a}, reported in Figure~\ref{fg:alg}. In the algorithm, if $a$ and $b$ are column vectors then the notation $c=[a;b]$ means that concatenation of the two column vectors into one single column vector $c$. Analogously, the notation $a=[a;b]$ means that the vector $a$ is replaces by the concatenation of the two column vectors $a$ and $b$.

\begin{figure}
\begin{algorithmic}
\State $S_s = [u_{1,h}(\mu_1);u_{2,h}(\mu_1);u_{3,h}(\mu_1)]$
\For{$j=2:n_s$}
    \State $s_j=[]$
    \For{$i=1:3$} 
    \State $e_1=\|u_{i,h}(\mu_j)-u_{i,h}(\mu_{j-1}) \|$,
    \State $e_2=\|u_{i,h}(\mu_j)+u_{i,h}(\mu_{j-1}) \|$
    \If {$e_1\geq e_2$}
    \State $s_j=[s_j;-u_{i,h}(\mu_j)]$
   \Else
        \State $s_j=[s_j;u_{i,h}(\mu_j)]$
    \EndIf
    \EndFor
    \State $S_s= [S_s, s_j]$
\EndFor
\end{algorithmic}
\caption{Algorithm used for the choice of the sign of the snapshots}
\label{fg:alg}
\end{figure}

\begin{figure}
\centering
\includegraphics[width=4.5cm]{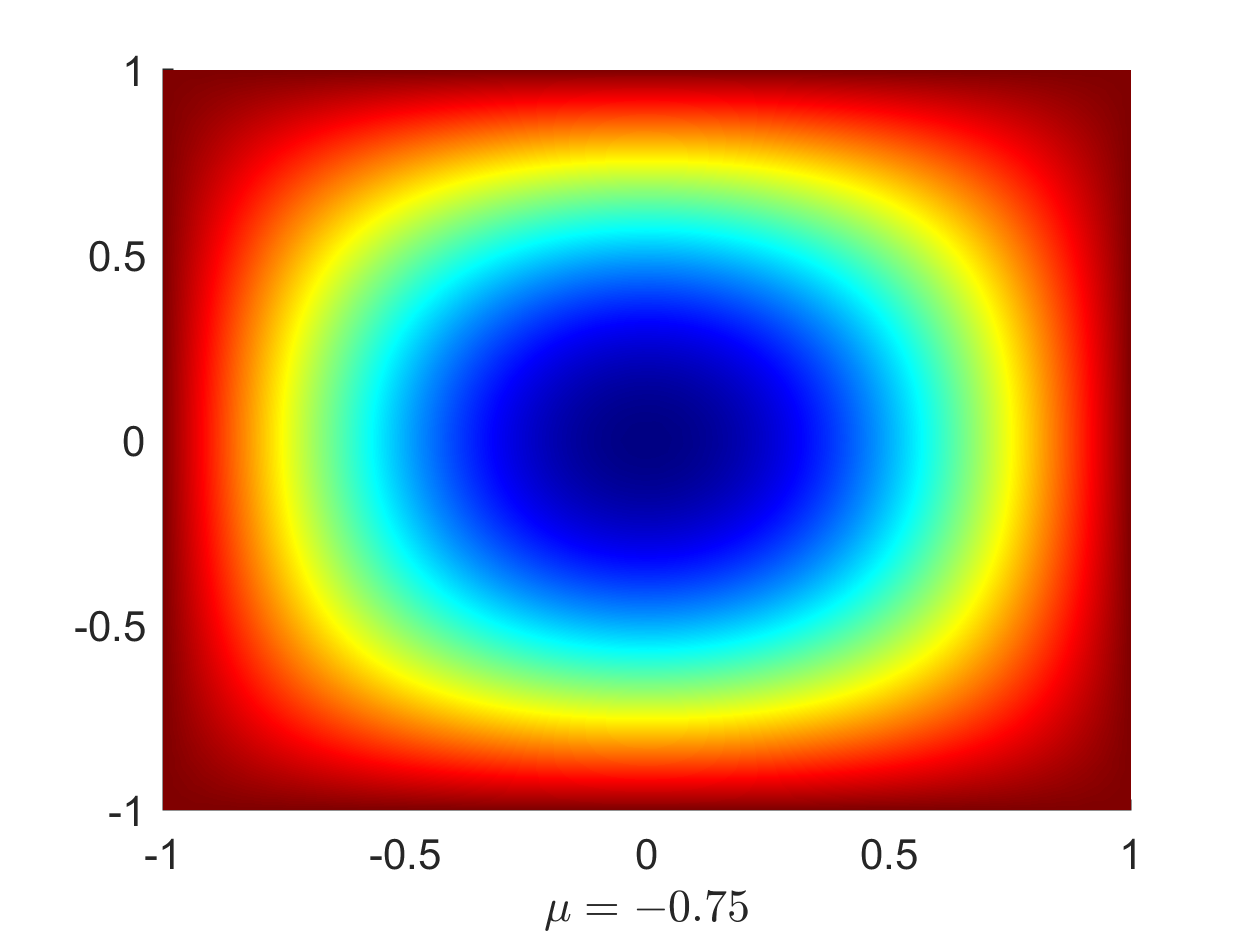}
\includegraphics[width=4.5cm]{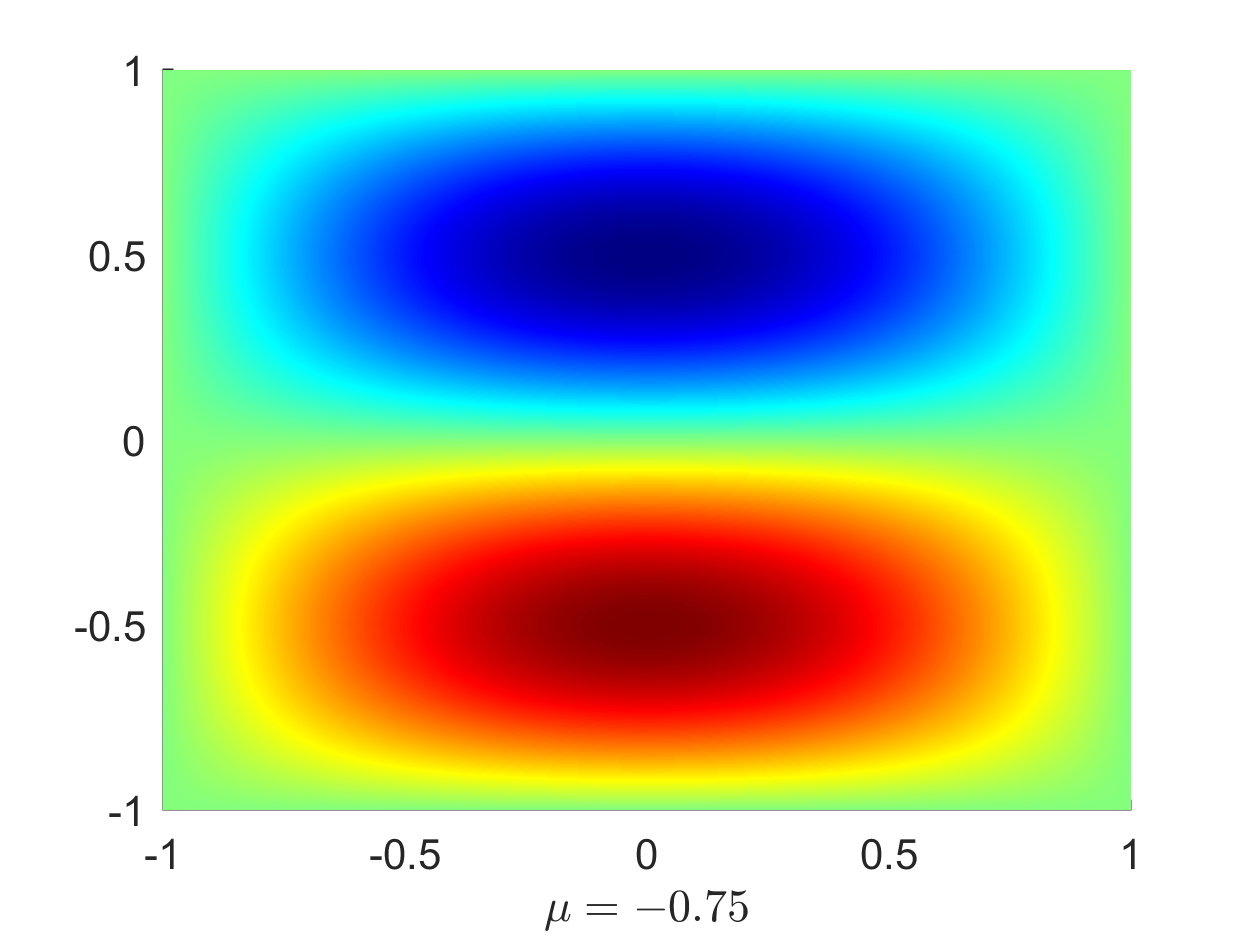}
\includegraphics[width=4.5cm]{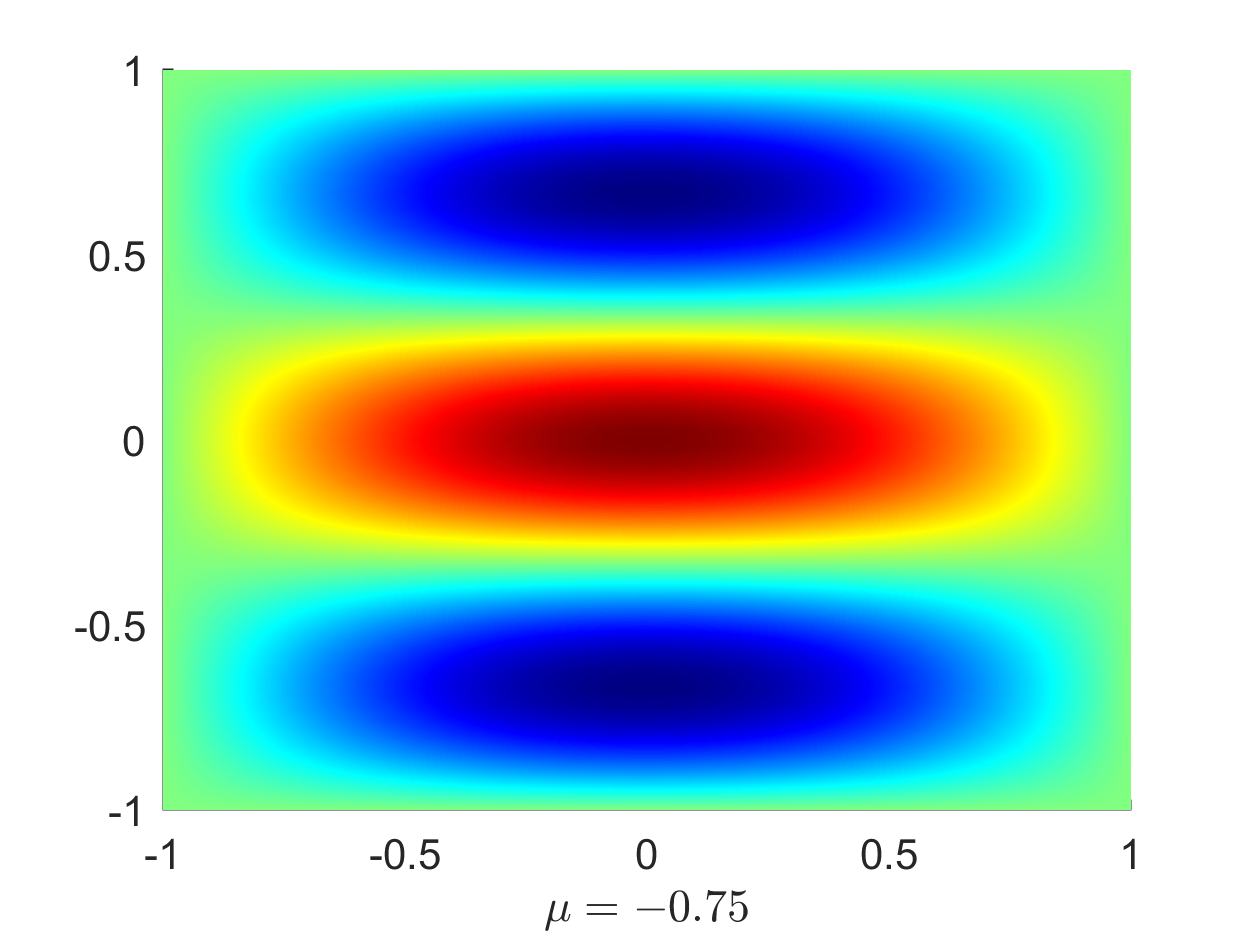}

\includegraphics[width=4.5cm]{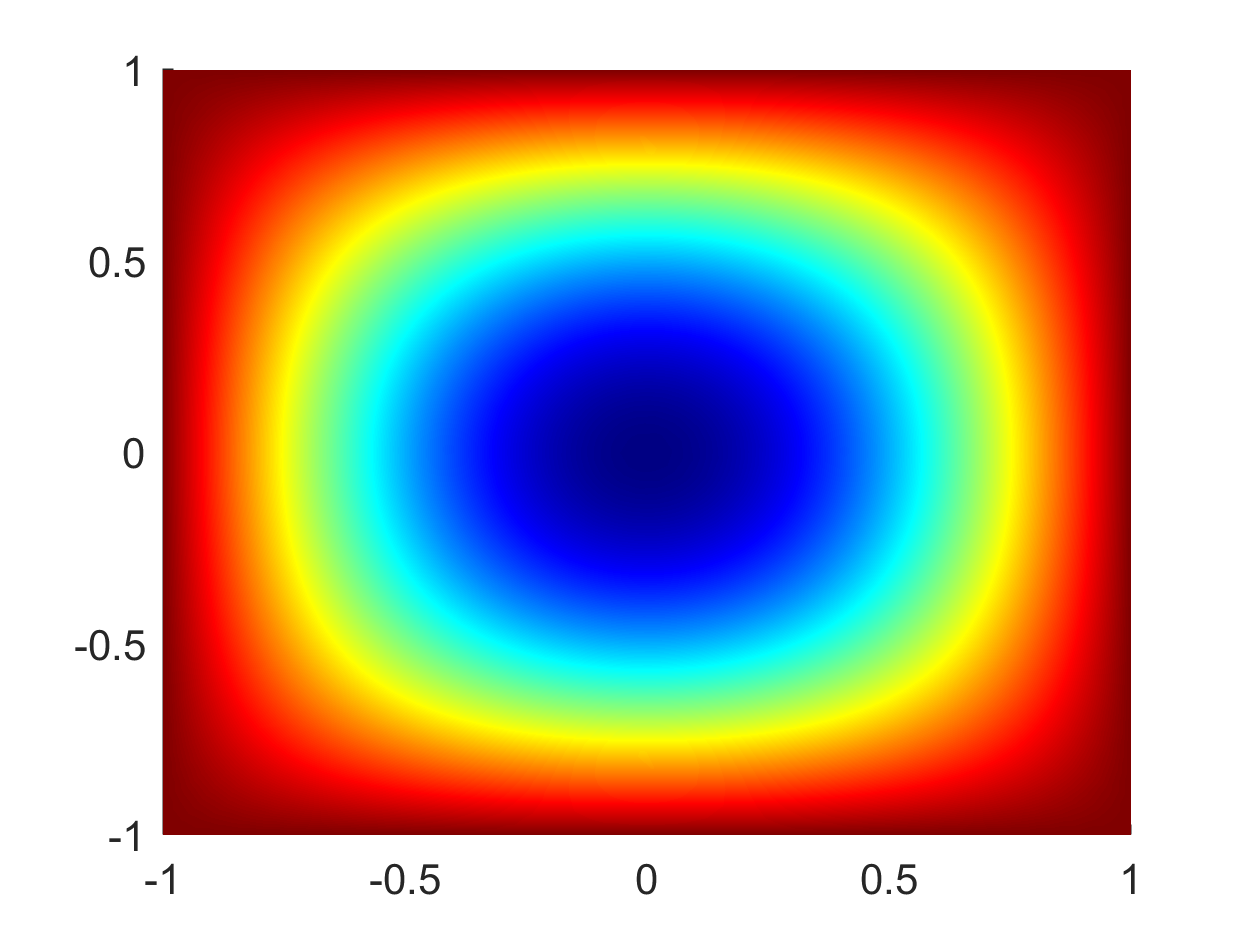}
\includegraphics[width=4.5cm]{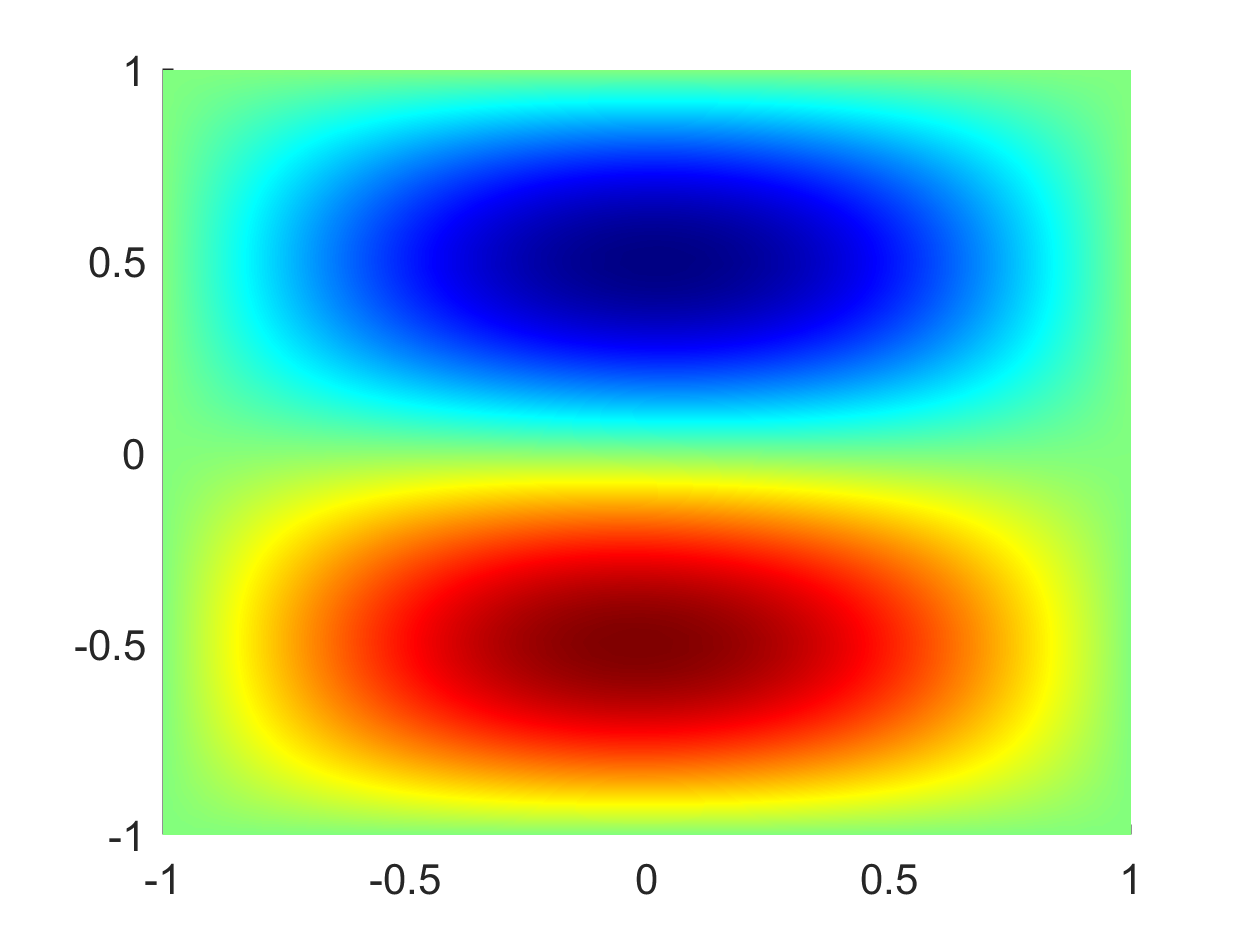}
\includegraphics[width=4.5cm]{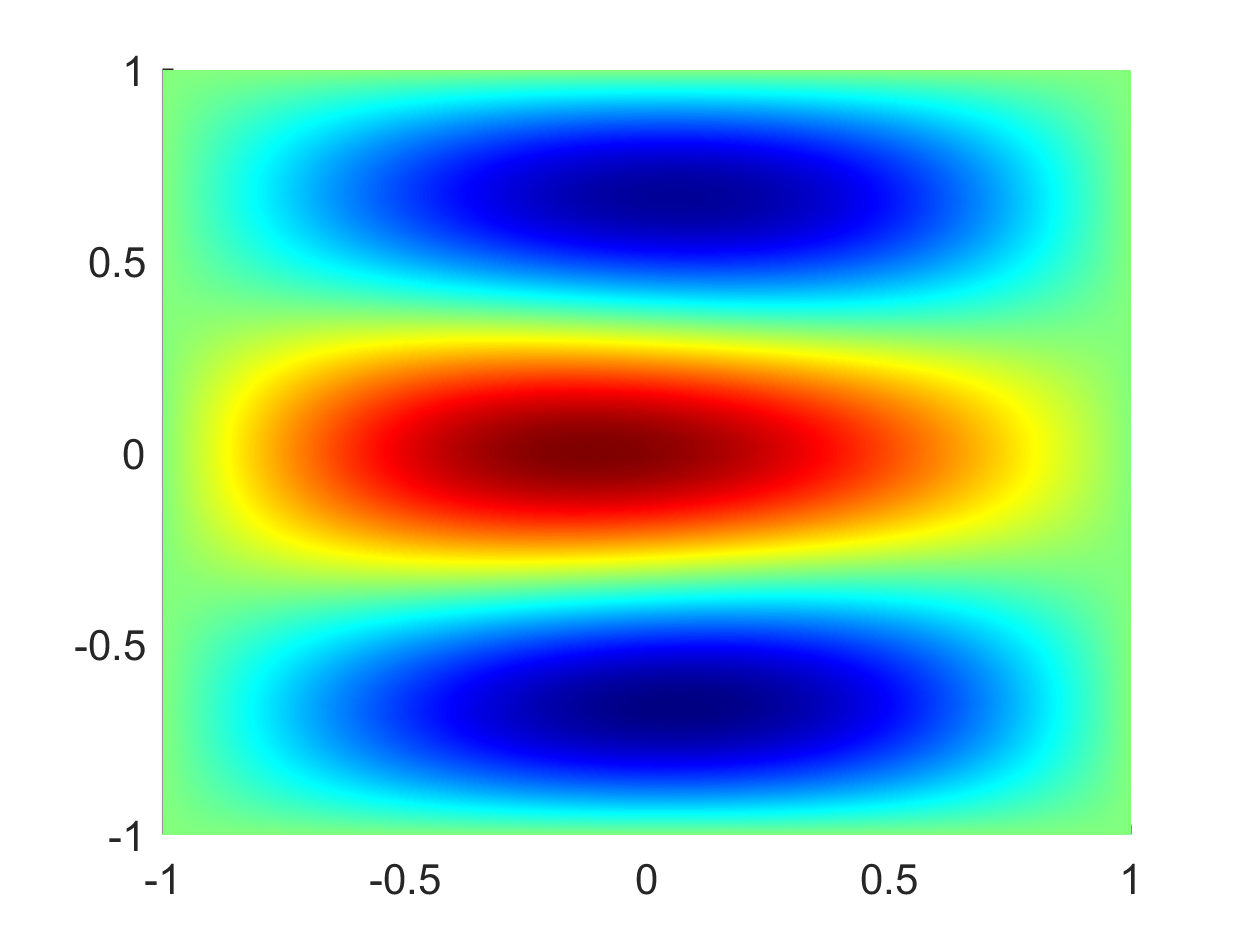}
\caption{Comparison of the first three eigenvectors of~\eqref{mdl:crossing} at $\mu=-0.75$ with $h=0.05$. FEM (top) and DD (bottom) solutions.}
\label{fig1:evct_simul}
\end{figure}

\begin{figure}
\centering
\includegraphics[width=4.5cm]{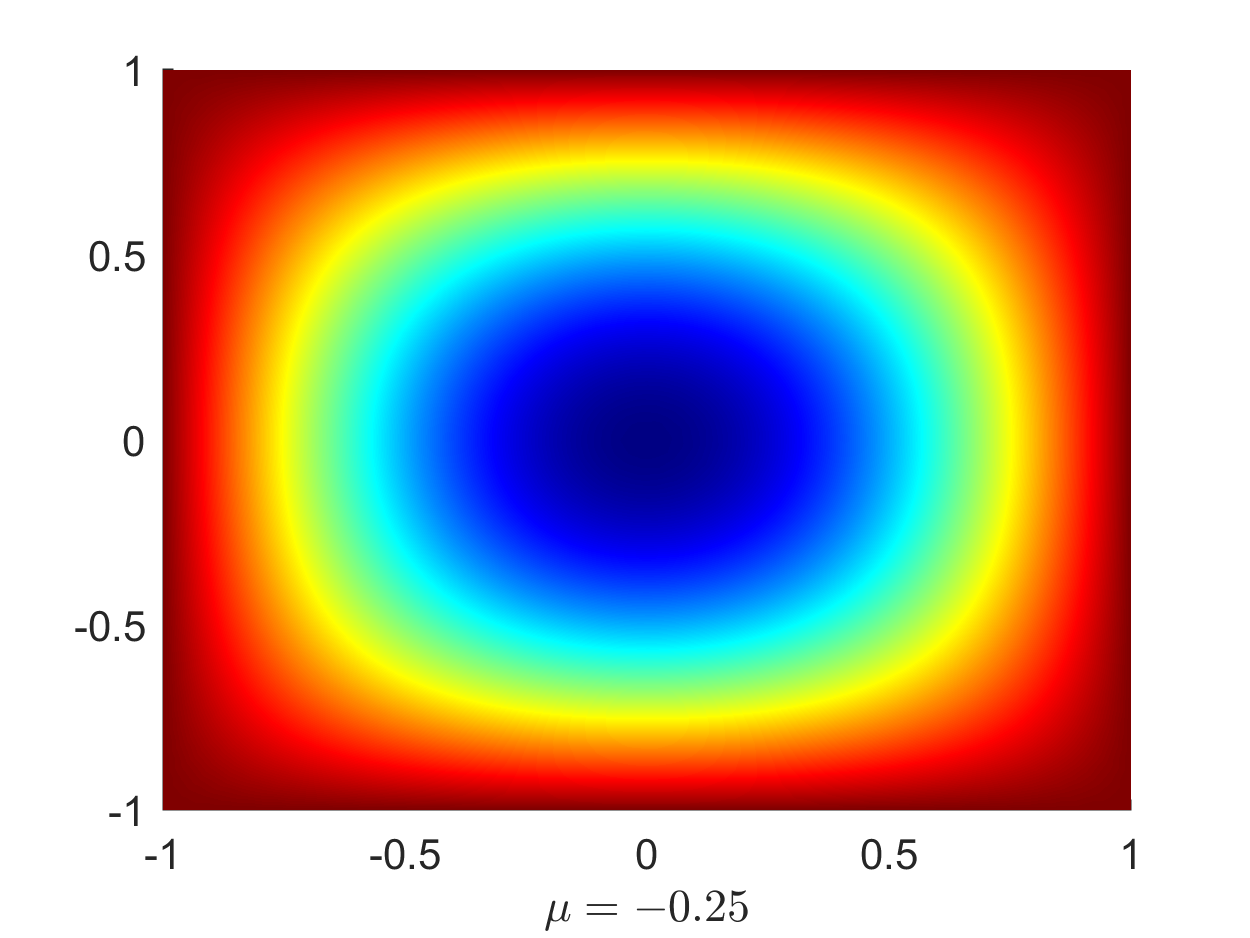}
\includegraphics[width=4.5cm]{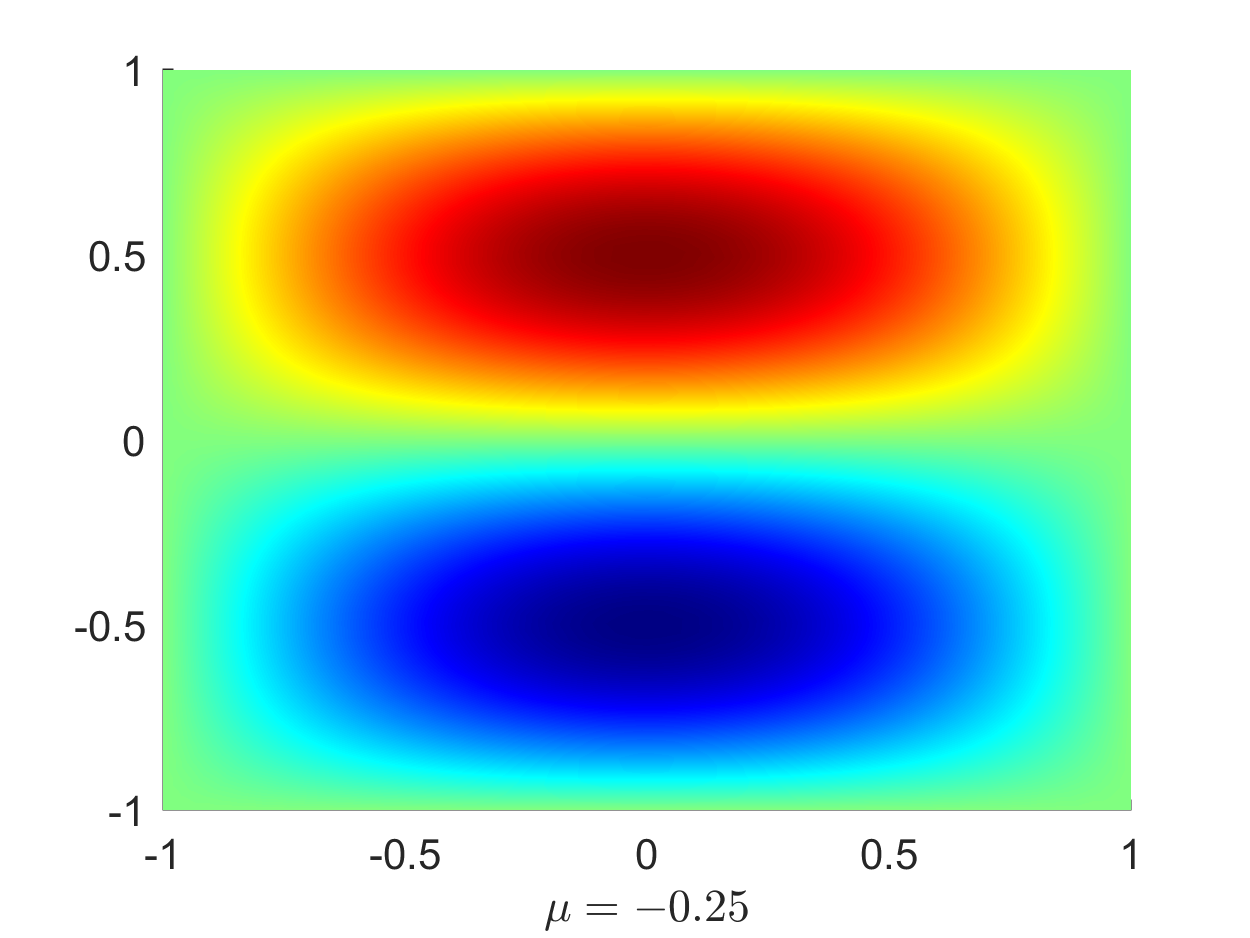}
\includegraphics[width=4.5cm]{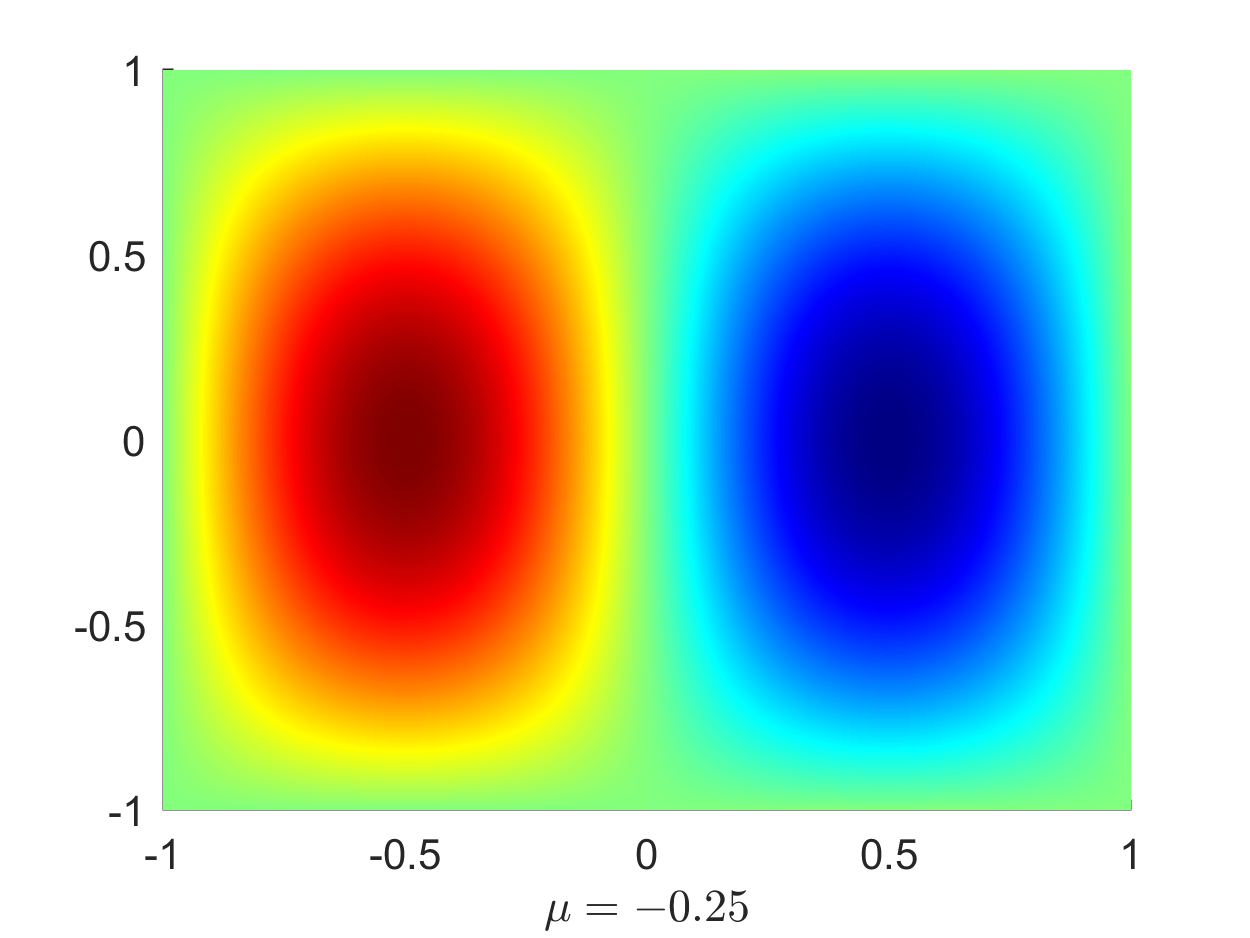}

\includegraphics[width=4.5cm]{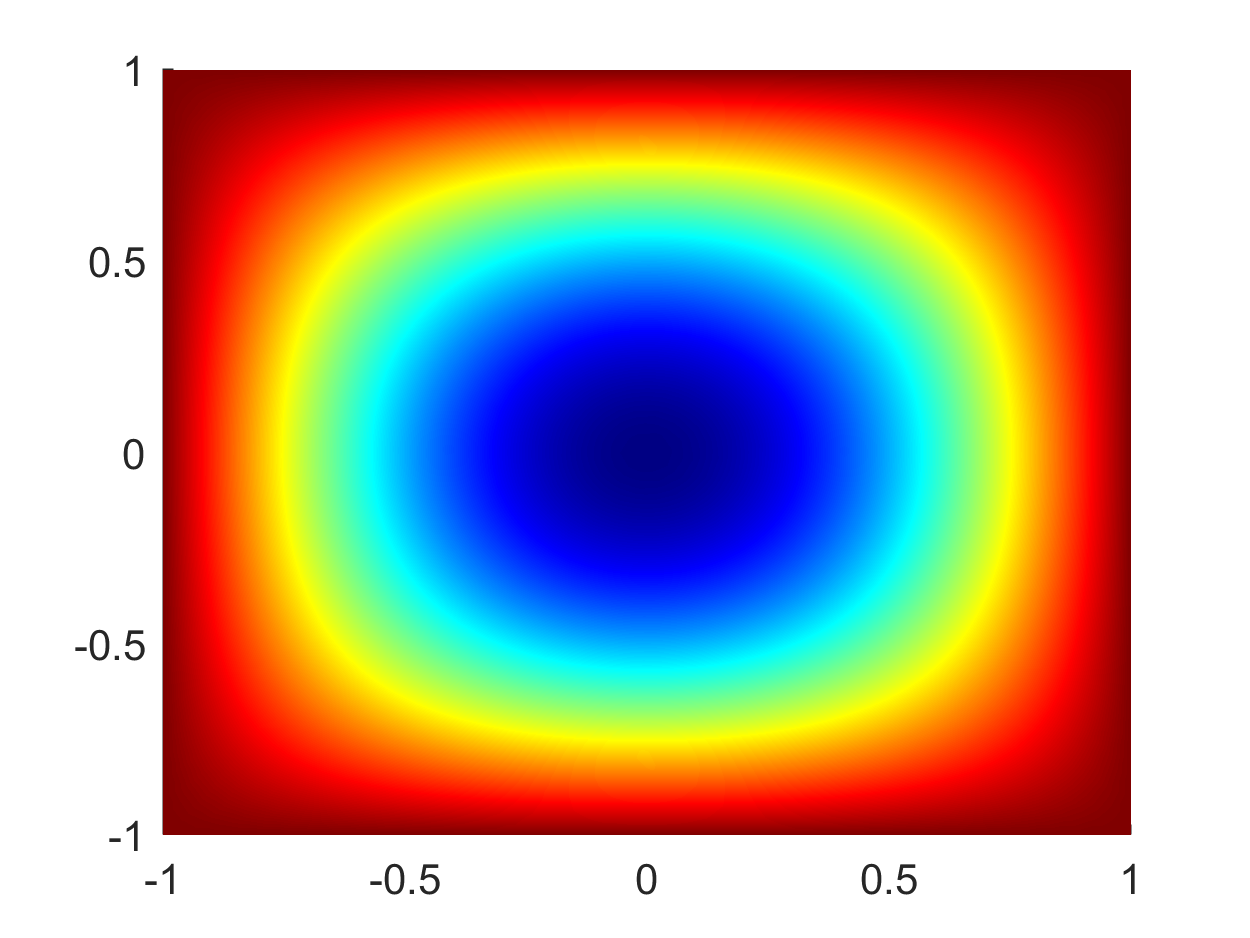}
\includegraphics[width=4.5cm]{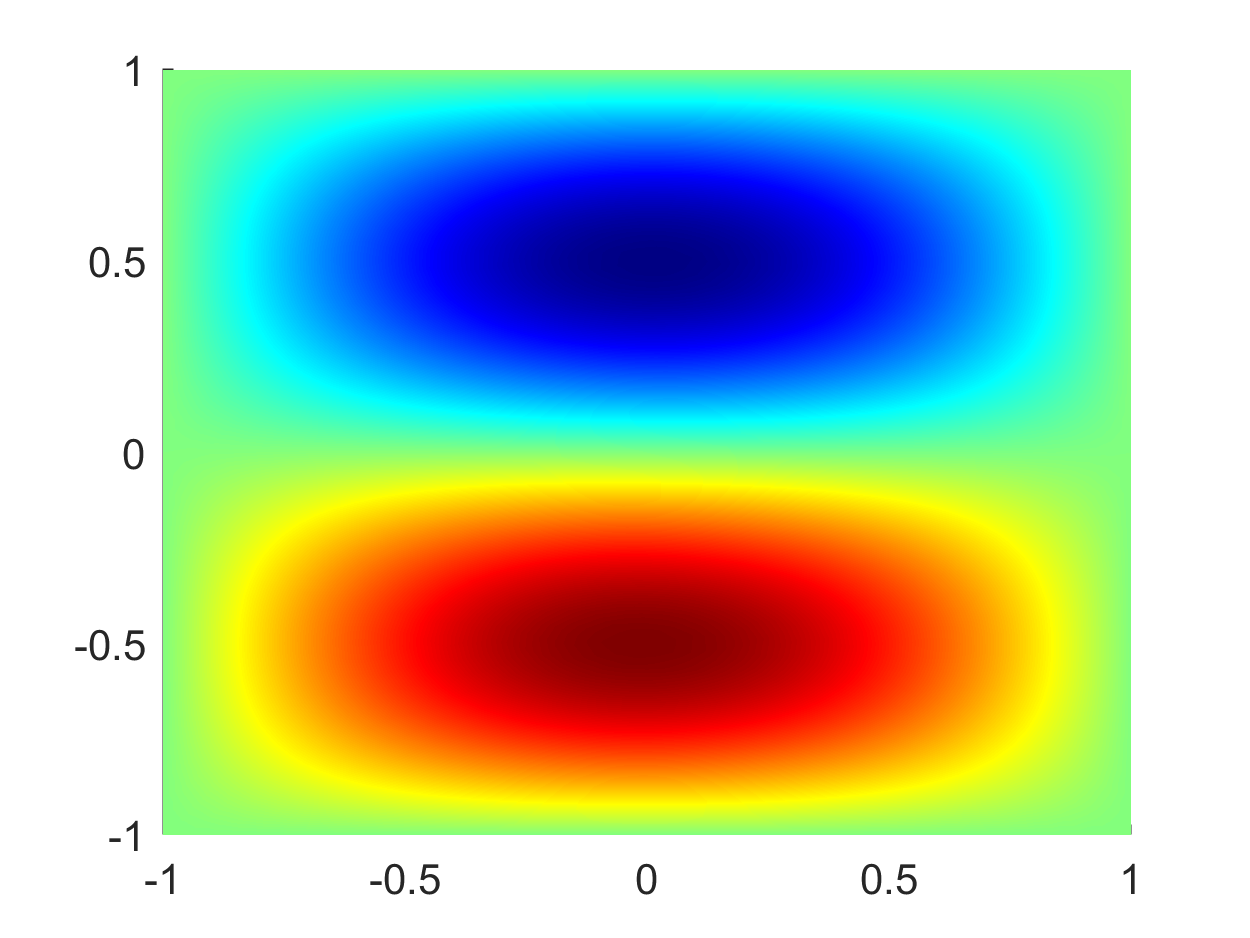}
\includegraphics[width=4.5cm]{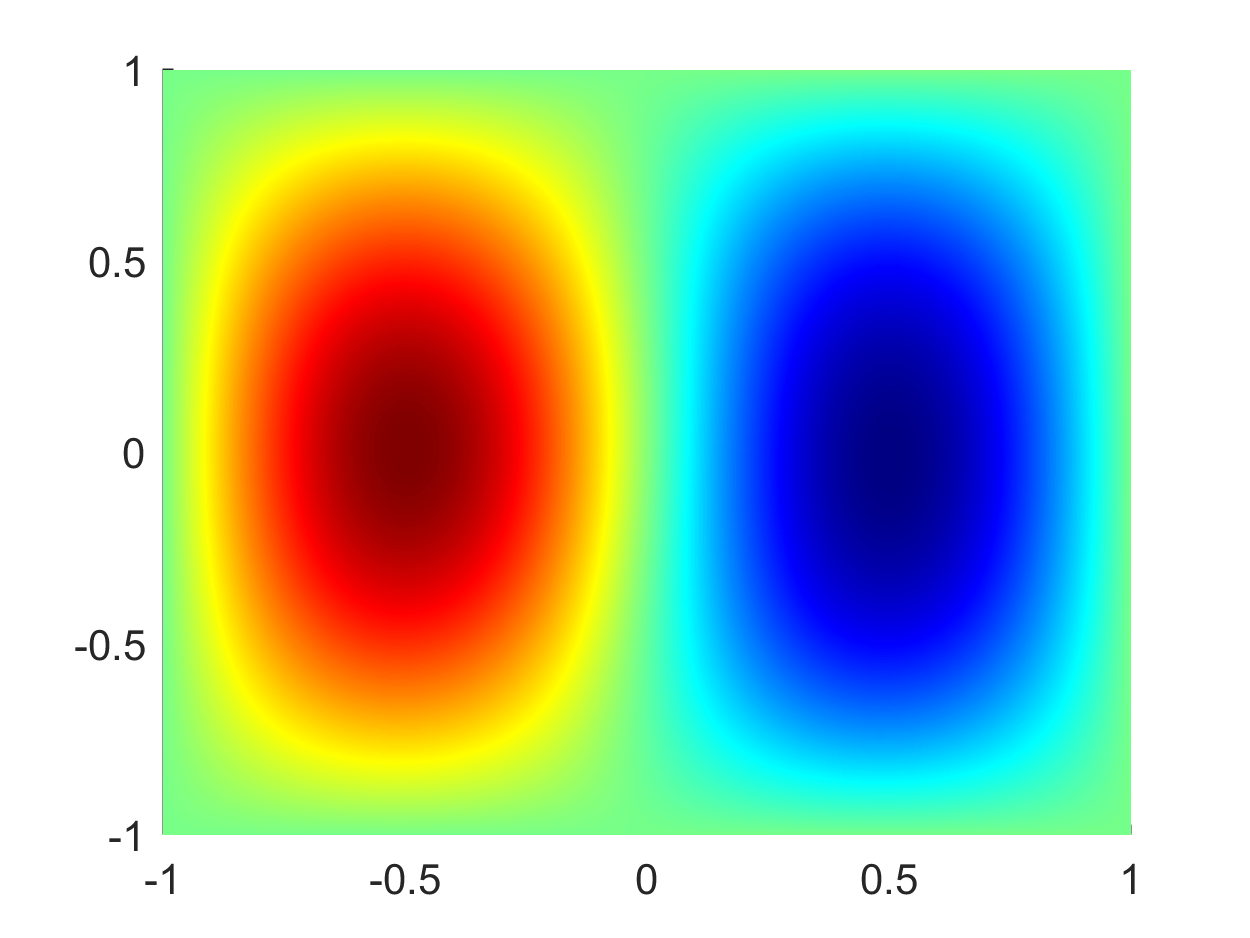}
\caption{Comparison of first three eigenvectors \eqref{mdl:crossing} at $\mu=-0.25$ with $h=0.05$. FEM (top) and DD (bottom) solutions.}
\label{fig2:evct_simul}
\end{figure}

\begin{figure}
\centering
\includegraphics[width=4.5cm]{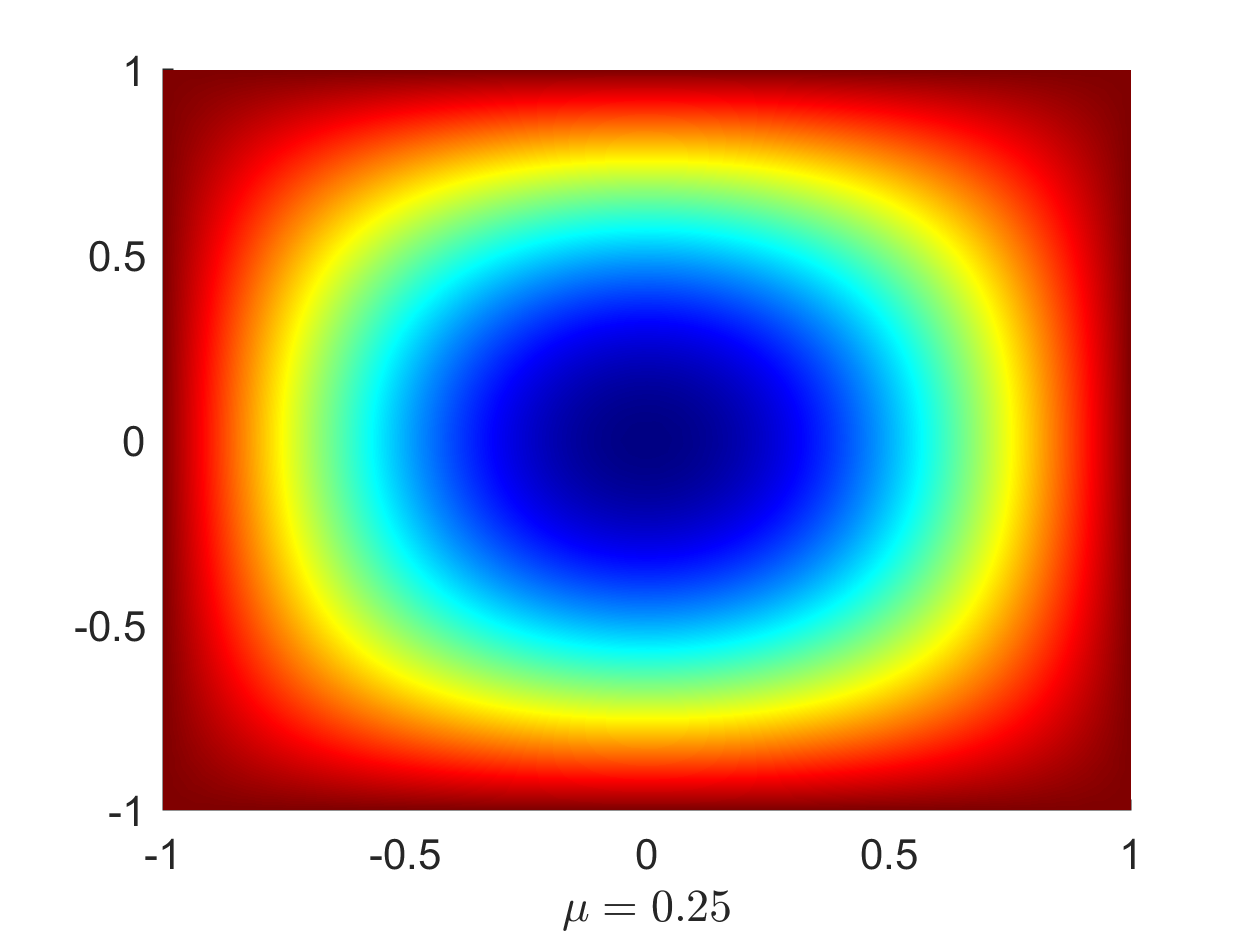}
\includegraphics[width=4.5cm]{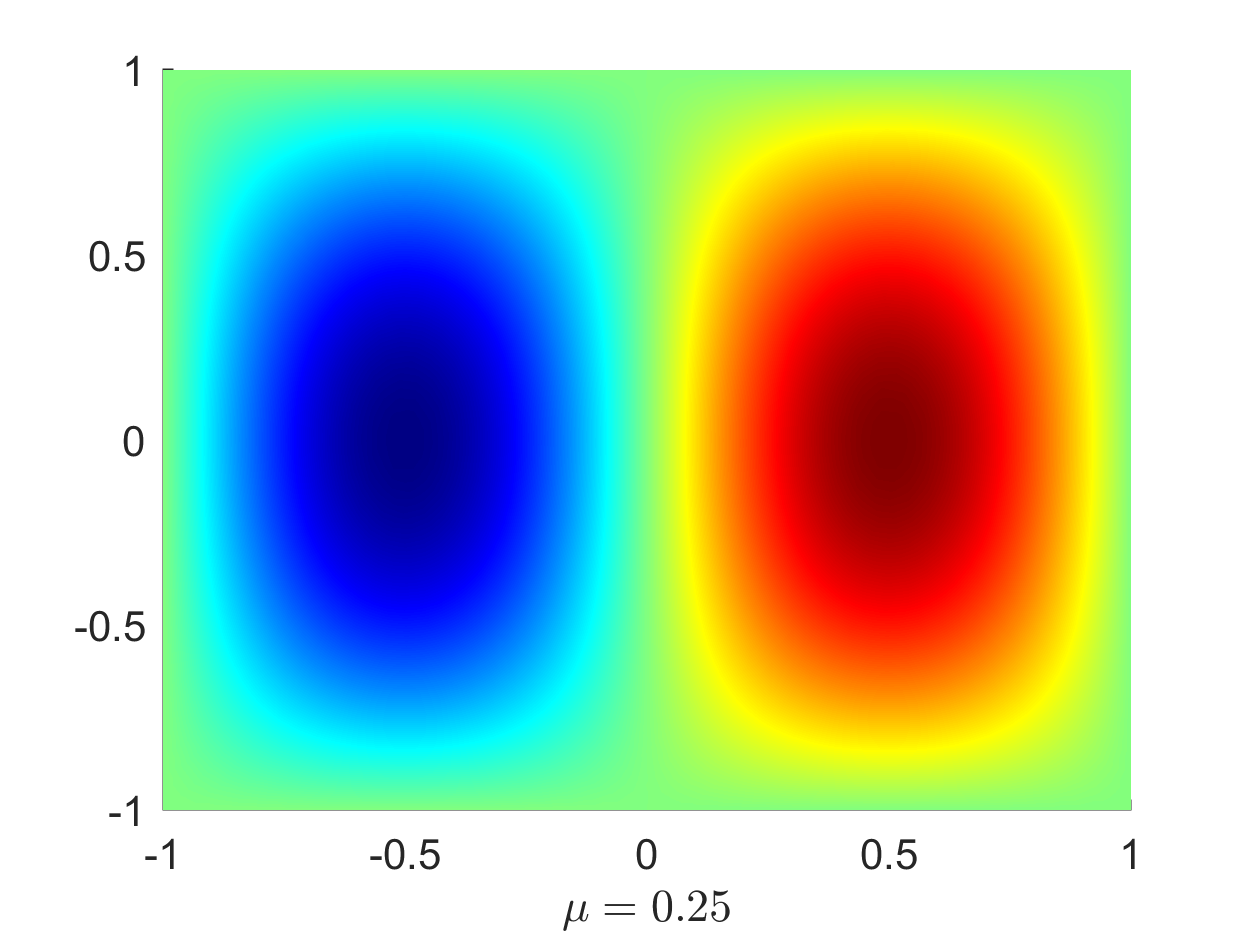}
\includegraphics[width=4.5cm]{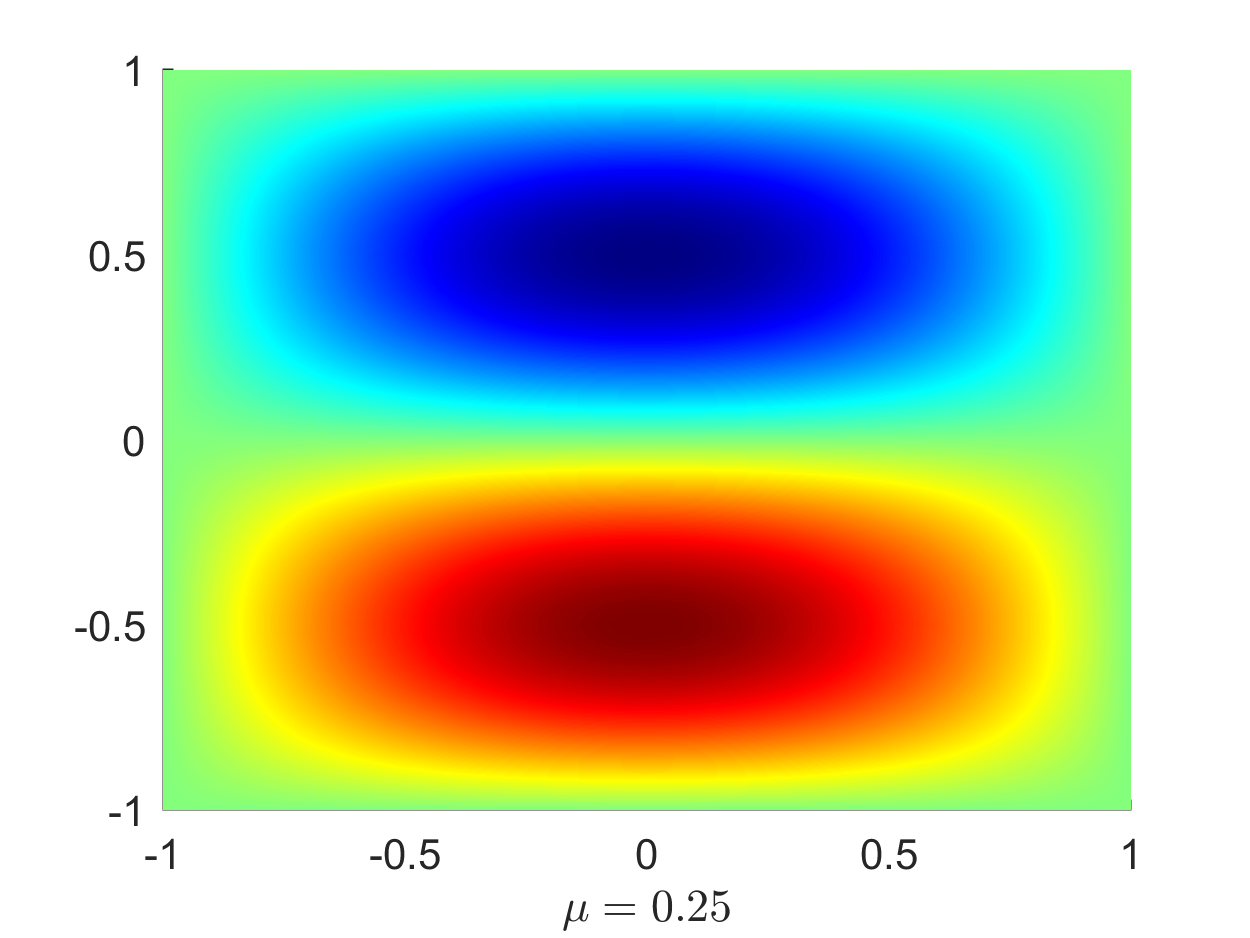}

\includegraphics[width=4.5cm]{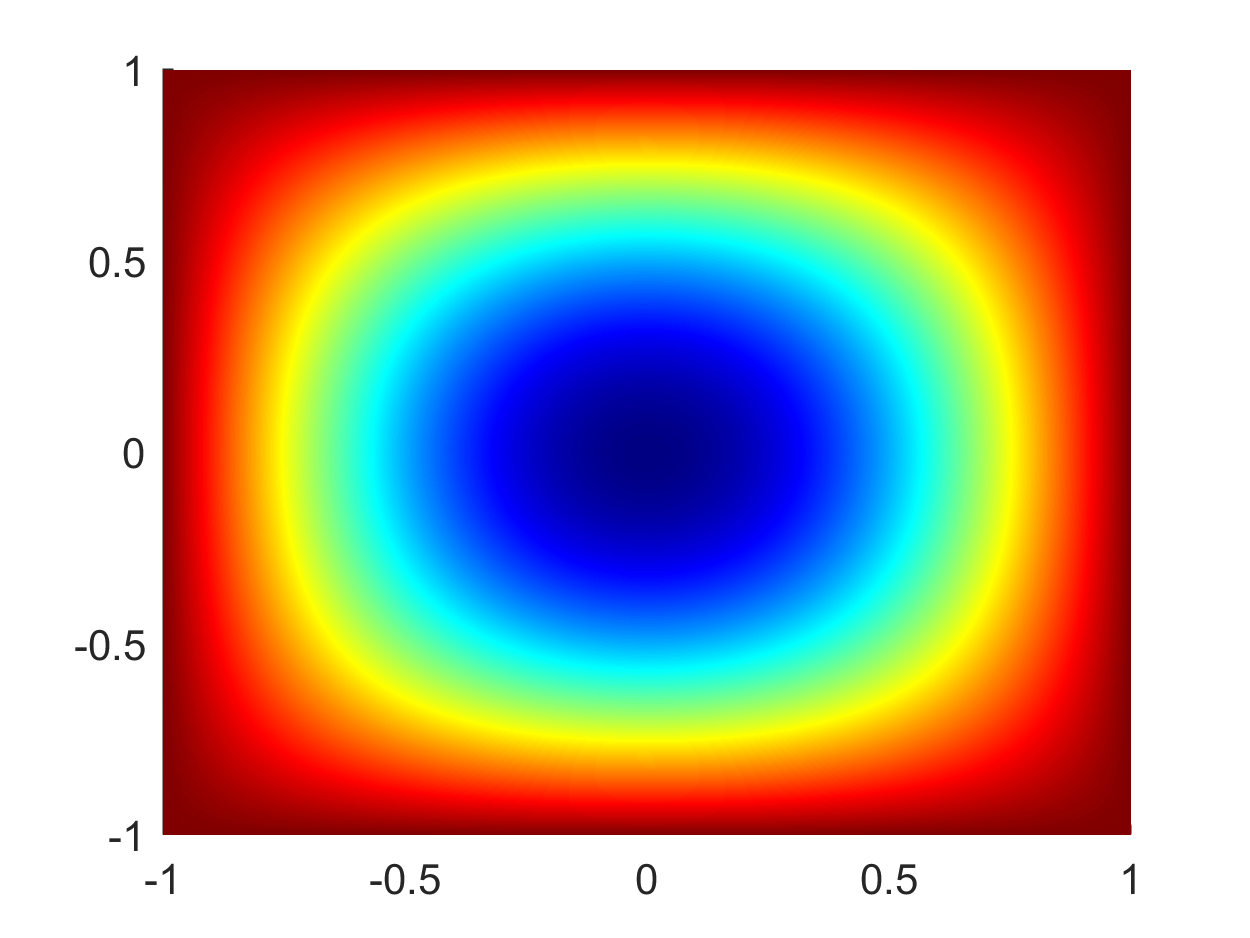}
\includegraphics[width=4.5cm]{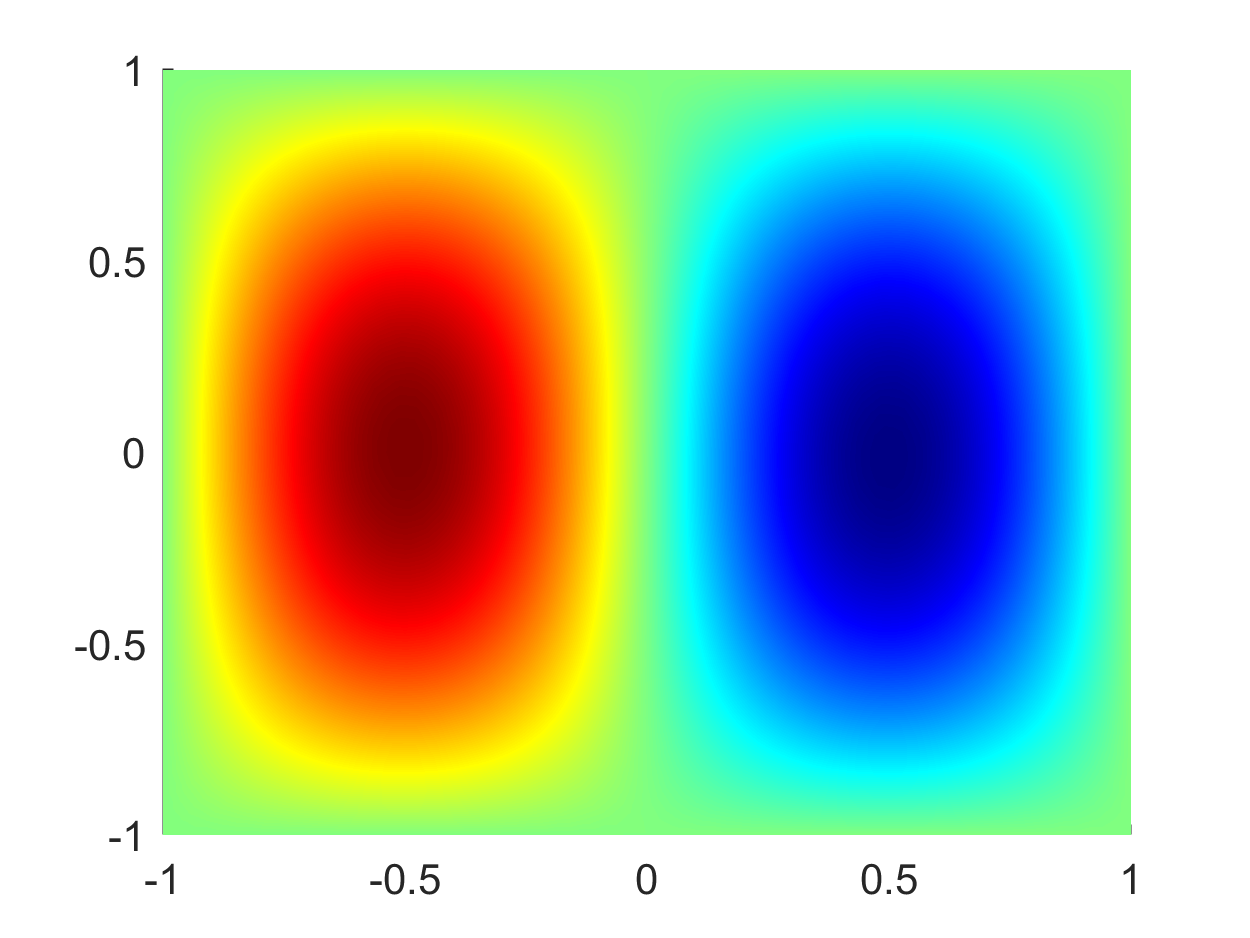}
\includegraphics[width=4.5cm]{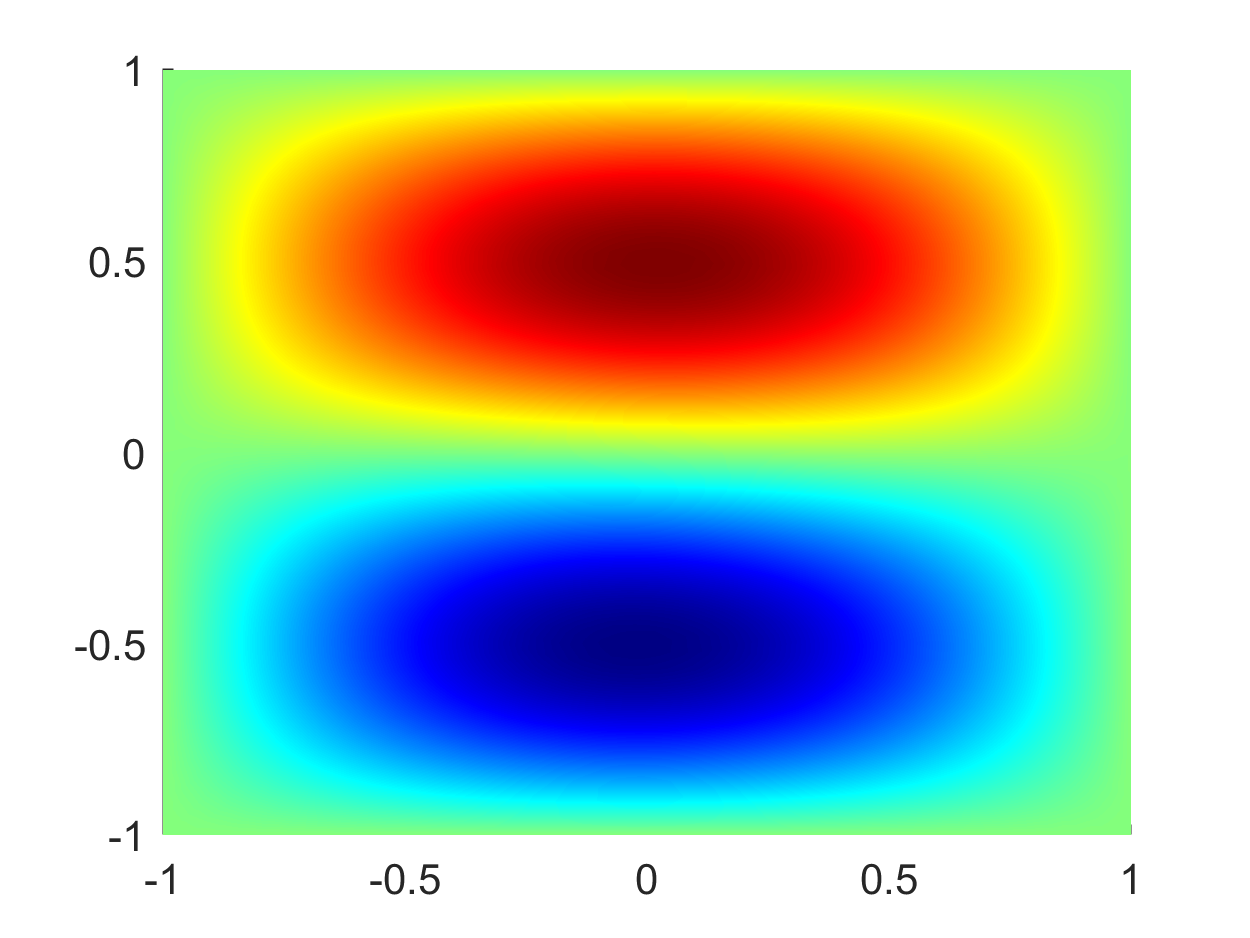}
\caption{Comparison of first three eigenvectors \eqref{mdl:crossing} at $\mu=0.25$ with $h=0.05$. FEM (top) and DD (bottom) solutions.}
\label{fig3:evct_simul}
\end{figure}

\begin{figure}
\centering
\includegraphics[width=4.5cm]{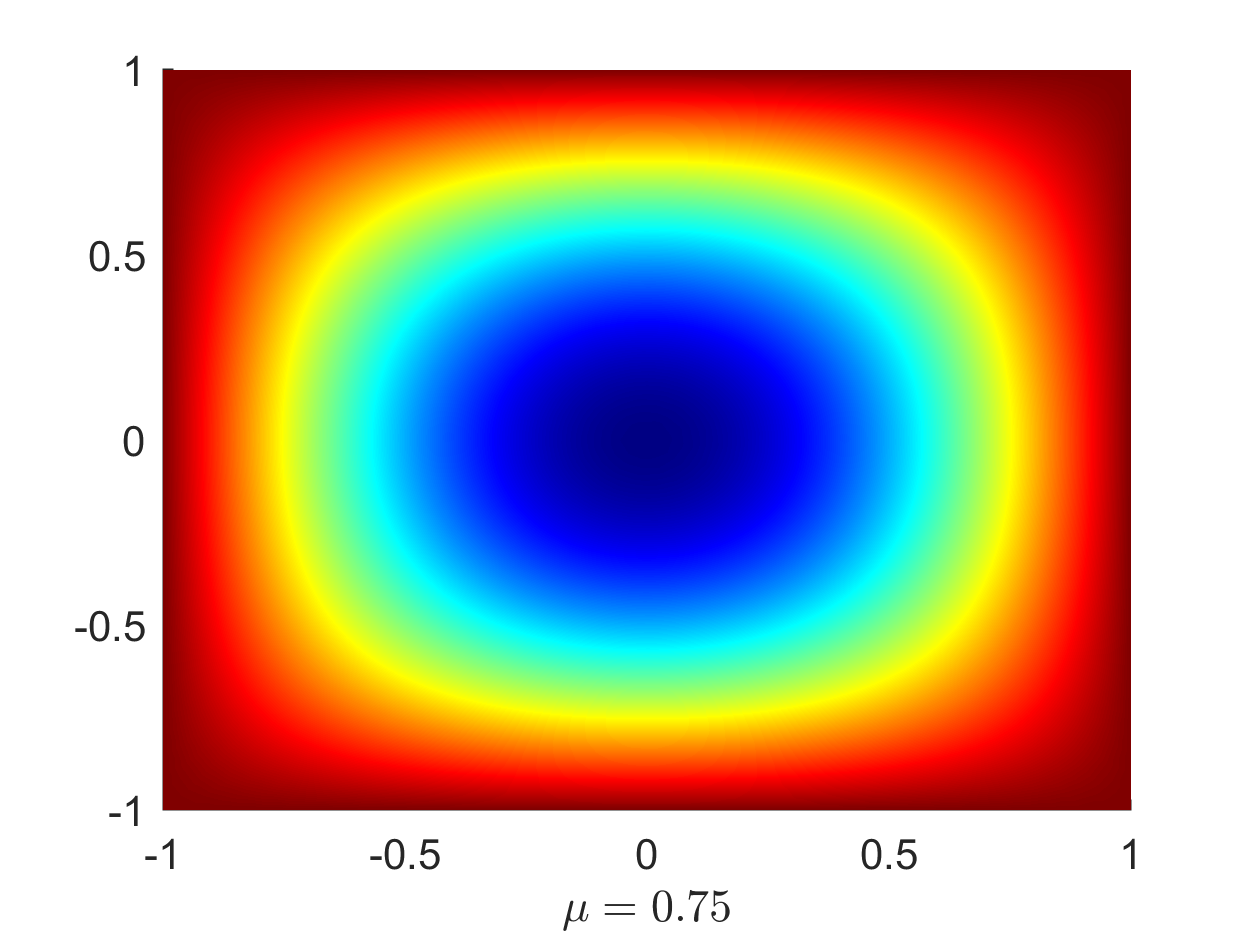}
\includegraphics[width=4.5cm]{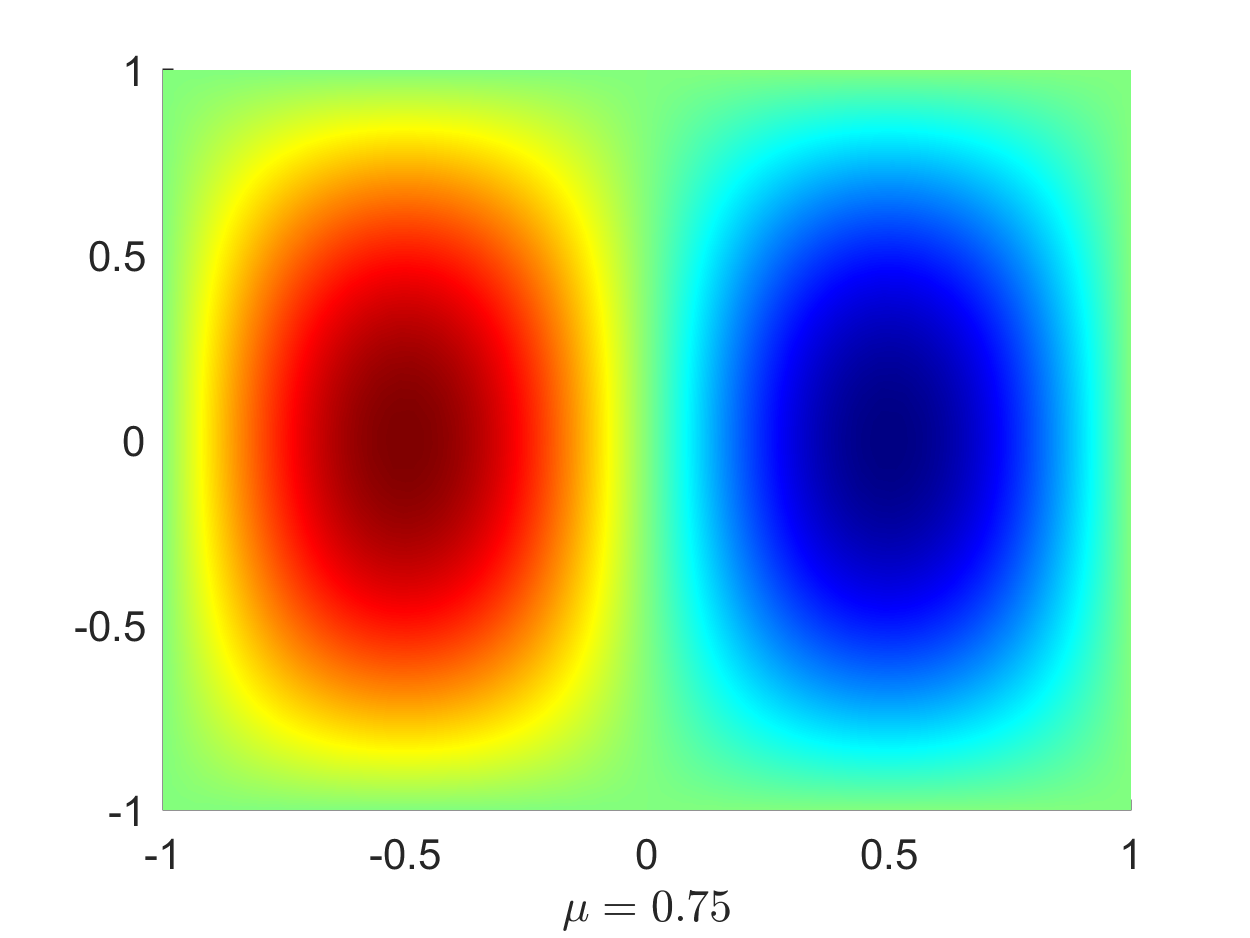}
\includegraphics[width=4.5cm]{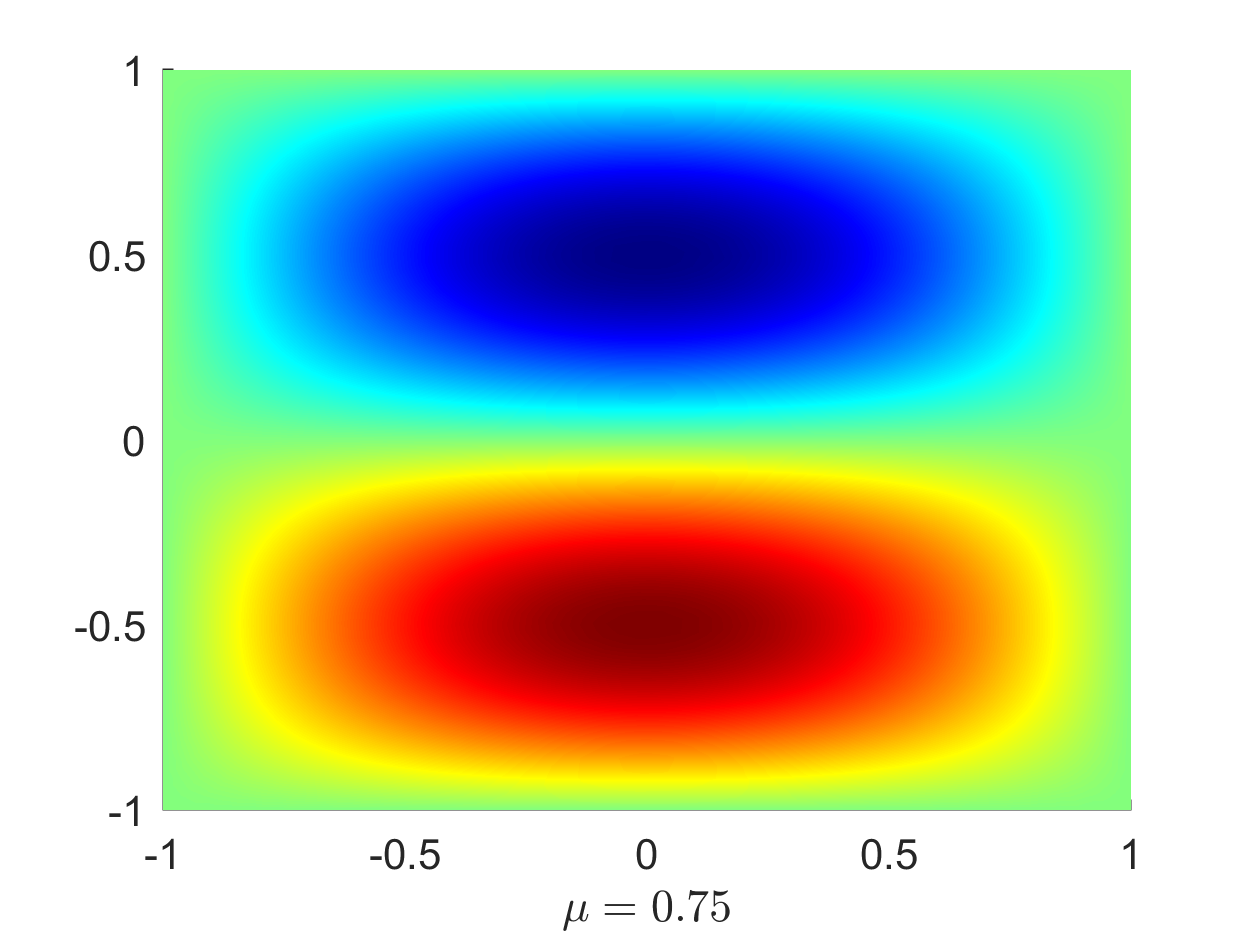}

\includegraphics[width=4.5cm]{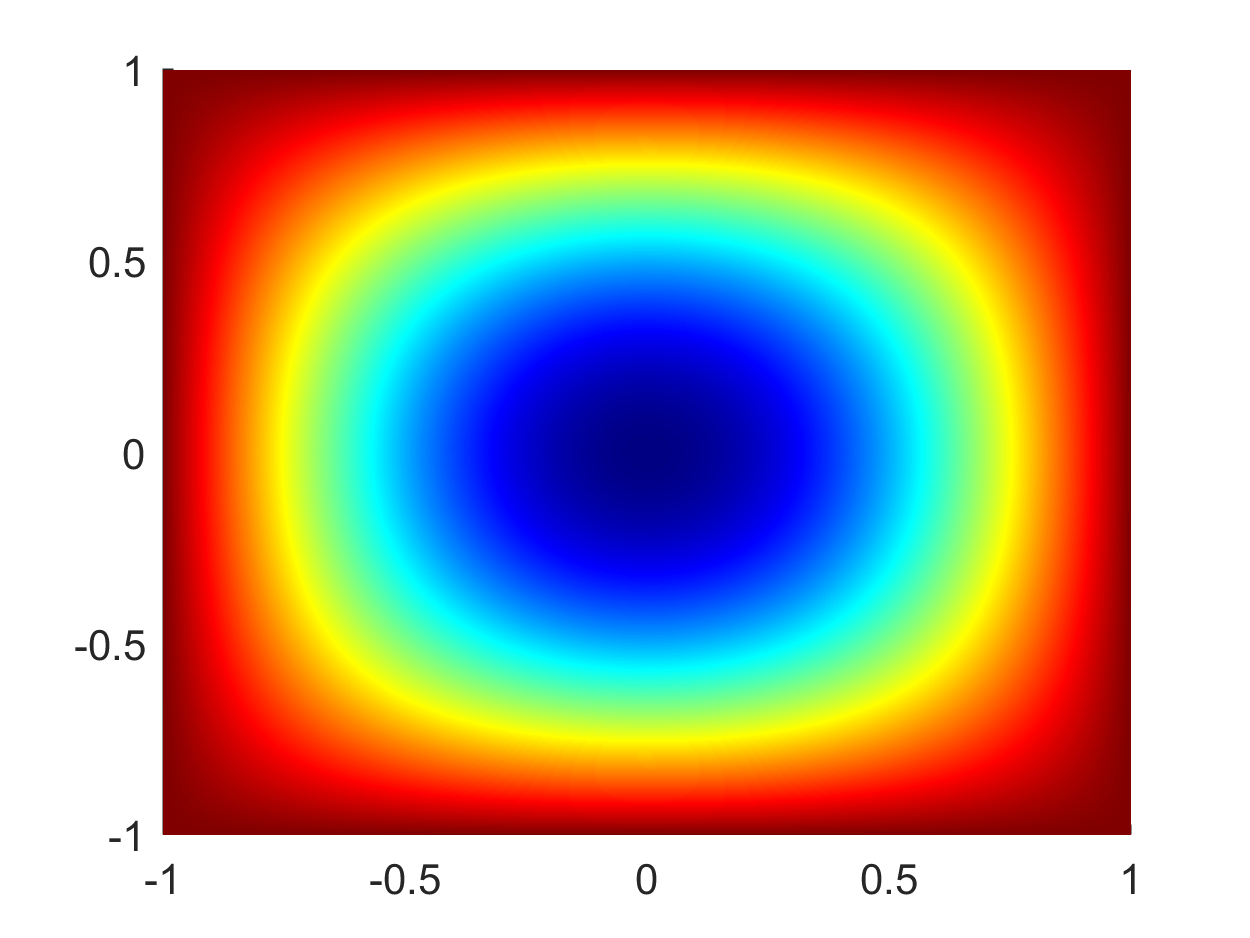}
\includegraphics[width=4.5cm]{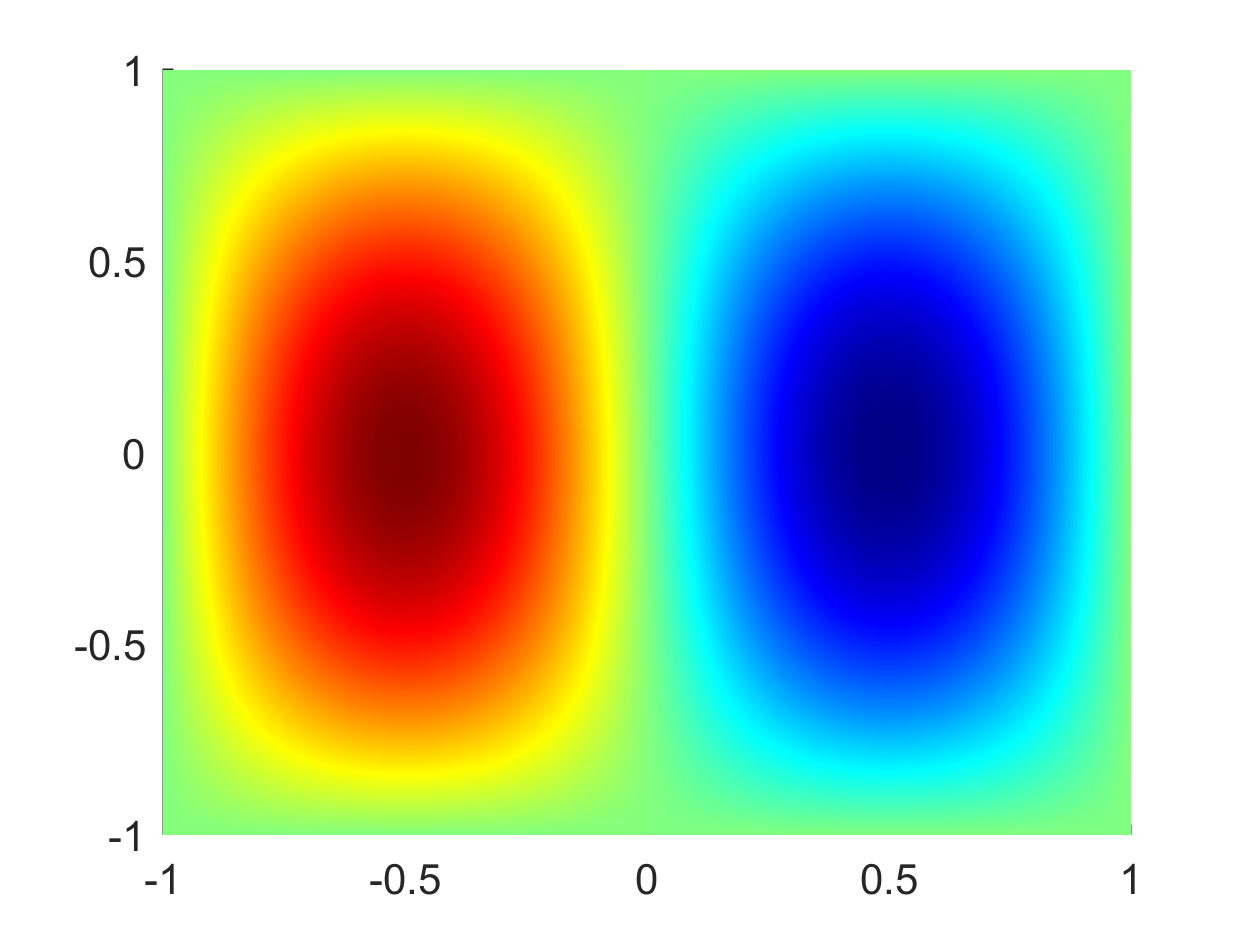}
\includegraphics[width=4.5cm]{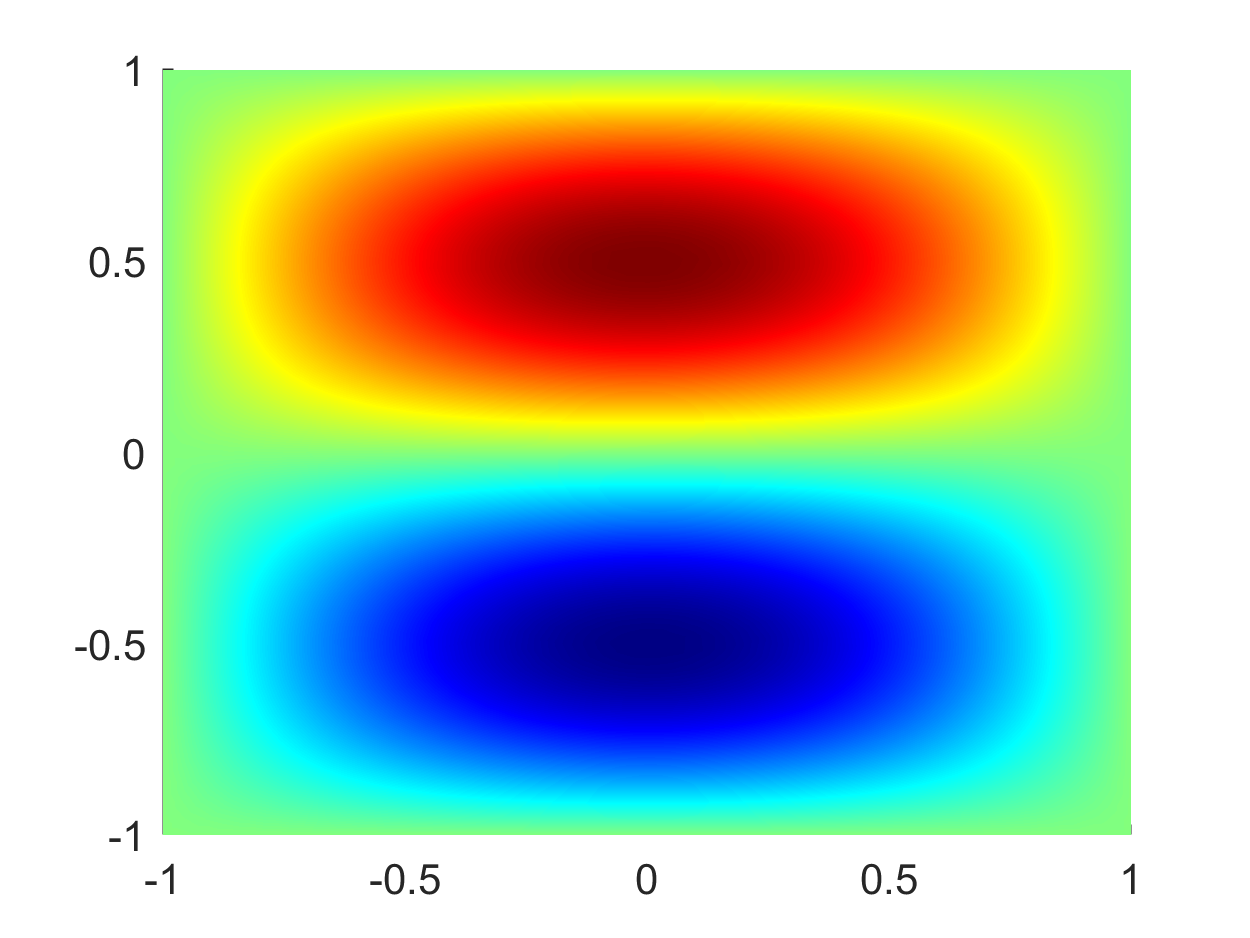}
\caption{Comparison of first three eigenvectors \eqref{mdl:crossing} at $\mu=0.75$ with $h=0.05$. FEM (top) and DD (bottom) solutions.}
\label{fig4:evct_simul}
\end{figure}

Finally, in Figures~\ref{fig1:evct_simul}, \ref{fig2:evct_simul}, \ref{fig3:evct_simul}, and~\ref{fig4:evct_simul} we show the first three eigenvectors obtained by the simultaneous DD model and compare them with that of FEM at the test points $\mu=-0.75$, $-0.25$, $0.25$, and $0.75$, respectively. We can see that the DD-based eigenvectors are matching with FEM-based eigenvectors at all the four points. Note that in this case we only had to train seven GPR: three corresponding to the three eigenvalues and four corresponding to the reduced coefficients.

\section*{Conclusions}

In this paper we introduced a new data driven reduced order model for the approximation of the eigensolution of parametric dependent partial differential equations. Our approach is based on the use of Gaussian process regressions which are trained by the eigenvalues and the coefficients of the eigenfunctions at given sample points of the parametric space.
We have described our scheme and we showed its effectiveness and robustness on several test cases.
Our DD algorithm can approximate successfully several model problems, including parametric PDEs that depend in a non affine way from the parameters. We also showed how the method can deal with the simultaneous approximation of eigenmodes and with crossings of eigenvalues.

\section*{Acknowledgments}

This research was supported by the Competitive Research Grants Program CRG2020 ``Synthetic data-driven model reduction methods for modal analysis'' awarded by the King Abdullah University of Science and Technology (KAUST).
Daniele Boffi is a member of the INdAM Research group GNCS and his research is partially supported by IMATI/CNR and by PRIN/MIUR.

%\bibliographystyle{elsarticle-harv}
%\bibliographystyle{elsarticle-num-names}
%\biboptions{square,sort,comma,numbers}
%\bibliographystyle{elsarticle-num}
\bibliographystyle{plain}
\bibliography{mybib}
\end{document}